\documentclass{amsart}
\usepackage{amsthm,amssymb,stmaryrd}
\usepackage{xparse,chngcntr}
\usepackage[table]{xcolor}
\usepackage{colortbl}
\usepackage[margin=1.1in]{geometry}
\usepackage{enumerate}
\usepackage{multicol}
\usepackage{amsfonts,mathrsfs,mathtools}
\usepackage{pgfplots}
\usepackage{pict2e,picture}
\usepackage[style=numeric, backend=biber]{biblatex}
\DeclareFieldFormat{postnote}{#1}
\DeclareFieldFormat{multipostnote}{#1}

\title{Borel combinatorics of $\mathbb{N}^{\ell}$-shift actions and hypersmooth equivalence relations}
\author{Scott Schneider}

\numberwithin{equation}{section}

\makeatletter
\newcommand{\pnrelbar}{%
  \linethickness{\dimen2}%
  \sbox\z@{$\m@th\prec$}%
  \dimen@=1.1\ht\z@
  \begin{picture}(\dimen@,.4ex)
  \roundcap
  \put(0,.2ex){\line(1,0){\dimen@}}
  \put(\dimexpr 0.5\dimen@-.2ex\relax,0){\line(1,1){.4ex}}
  \end{picture}%
}

\newcommand{\precneq}{\mathrel{\vcenter{\hbox{\text{\prec@neq}}}}}
\newcommand{\prec@neq}{%
  \dimen2=\f@size\dimexpr.04pt\relax
  \oalign{%
    \noalign{\kern\dimexpr.2ex-.5\dimen2\relax}
    $\m@th\prec$\cr
    \noalign{\kern-.5\dimen2}
    \hidewidth\pnrelbar\hidewidth\cr
  }%
}
\makeatother

\ExplSyntaxOn
\NewDocumentCommand{\newtheorembysubsection}{mom}
 {
  \IfNoValueTF{#2}
   {
    \newtheorem{#1}{#3}[subsection]
    \cs_set:cpn { the#1 } { \thesection. \deepesh_subsection: \arabic{#1} }
    \counterwithin*{#1}{section}
   }
   {
    \newtheorem{#1}[#2]{#3}
   }
 }
\cs_new:Nn \deepesh_subsection:
 {
  \int_compare:nT { \value{subsection}>0 } { \arabic{subsection}. }
 }
\ExplSyntaxOff

\definecolor{Gray}{gray}{0.9}
\newcolumntype{g}{>{\columncolor{Gray}}c}

\newtheorembysubsection{theorem}{Theorem}
\newtheorembysubsection{lemma}[theorem]{Lemma}
\newtheorembysubsection{fact}[theorem]{Fact}
\newtheorembysubsection{corollary}[theorem]{Corollary}
\newtheorembysubsection{proposition}[theorem]{Proposition}
\newtheorembysubsection{conjecture}[theorem]{Conjecture}
\newtheorembysubsection{claim}[theorem]{Claim}
\newtheorembysubsection{observation}[theorem]{Observation}
\theoremstyle{definition}
\newtheorembysubsection{definition}[theorem]{Definition}
\newtheorembysubsection{notation}[theorem]{Notation}
\newtheorembysubsection{remark}[theorem]{Remark}
\newtheorembysubsection{DefinitionNotation}[theorem]{Definition/Notation}
\newtheorembysubsection{example}[theorem]{Example}
\newtheorembysubsection{question}[theorem]{Question}

\newcommand{\smallnwarrow}{\text{\scalebox{.7}{\ensuremath{\nwarrow}}}}
\newcommand{\smallsearrow}{\text{\scalebox{.7}{\ensuremath{\searrow}}}}

\newcommand{\N}{\ensuremath{\mathbb{N}}}
\newcommand{\Z}{\ensuremath{\mathbb{Z}}}

\newcommand{\R}{\ensuremath{\mathcal{R}}}
\renewcommand{\S}{\ensuremath{\mathcal{S}}}

\renewcommand{\L}{\ensuremath{\mathcal{L}}}
\newcommand{\B}{\ensuremath{\mathcal{B}}}
\newcommand{\A}{\ensuremath{\mathcal{A}}}
\renewcommand{\H}{\ensuremath{\mathcal{H}}}

\newcommand{\npseJ}{\ensuremath{\neg\textup{Per}_{\q}(\JT)}}
\newcommand{\npseJm}{\ensuremath{\neg\textup{Per}_{\q}(J_{\m}^{\omega^2})}}

\newcommand{\Ol}{\ensuremath{\Omega^{\omega^\ell}}\!}
\newcommand{\Jl}{\ensuremath{J^{\omega^\ell}}\!}
\newcommand{\Tl}{\ensuremath{2^{\omega^\ell}}\!}
\newcommand{\OT}{\ensuremath{\Omega^{\omega^2}}\!}
\newcommand{\JT}{\ensuremath{J^{\omega^2}}\!}
\newcommand{\TT}{\ensuremath{2^{\omega^2}}\!}
\newcommand{\FOl}{\ensuremath{\textup{Fr}(\Omega^{\omega^\ell})}}

\newcommand{\FTl}{\ensuremath{\textup{Fr}(2^{\omega^\ell})}}

\newcommand{\FTT}{\ensuremath{\textup{Fr}(2^{\omega^2})}}

\newcommand{\POl}{\ensuremath{\textup{Per}(\Omega^{\omega^\ell})}}

\newcommand{\POT}{\ensuremath{\textup{Per}(\Omega^{\omega^2})}}

\newcommand{\PTT}{\ensuremath{\textup{Per}(2^{\omega^2})}}
\newcommand{\PmOl}{\ensuremath{\textup{Per}_{\m}(\Omega^{\omega^\ell})}}
\newcommand{\EOl}{\ensuremath{E_{\ell}^{\Omega}}}

\newcommand{\ETl}{\ensuremath{E_{\ell}^{2}}}
\newcommand{\EOT}{\ensuremath{E_{2}^{\Omega}}}
\newcommand{\EJT}{\ensuremath{E_{2}^{J}}}
\newcommand{\ETT}{\ensuremath{E_{2}^{2}}}
\newcommand{\FOln}{\ensuremath{F_{\ell,\d{n}}^{\Omega}}}

\renewcommand{\d}[1]{\ensuremath{\dot{#1}}}
\newcommand{\FeqOl}[1]{\ensuremath{F^{\Omega}_{\ell,{#1}}}}

\newcommand{\p}{\ensuremath{\bar{p}}}
\newcommand{\q}{\ensuremath{\bar{q}}}

\newcommand{\s}{\ensuremath{\bar{s}}}
\renewcommand{\t}{\ensuremath{\bar{t}}}
\newcommand{\m}{\ensuremath{\bar{m}}}
\renewcommand{\c}{\ensuremath{\bar{c}}}
\newcommand{\JTT}{\ensuremath{J^{\omega^2}}}

\newcommand{\EJmtwo}{\ensuremath{E_2^{J_{\m}}}}
\newcommand{\Ge}{\ensuremath{\G_{\bar{e}_1}}}
\newcommand{\Gee}{\ensuremath{\G_{\bar{e}_2}}}
\newcommand{\Gea}{\ensuremath{\vec{\G}_{\bar{e}_1}}}
\newcommand{\Geea}{\ensuremath{\vec{\G}_{\bar{e}_2}}}
\newcommand{\Gp}{\ensuremath{\G_{\p}}}
\newcommand{\Gm}{\ensuremath{\G_{\m}}}
\newcommand{\Gq}{\ensuremath{\G_{\q}}}
\newcommand{\Gse}{\ensuremath{\G_{\q}}}
\newcommand{\cse}{\ensuremath{c_{\q}}}
\newcommand{\Sse}{\ensuremath{S_{\q}}}
\newcommand{\Gsea}{\ensuremath{\vec{\G}_{\q}}}

\newcommand{\Wqnw}{\ensuremath{\Wq^{\smallnwarrow}}}

\newcommand{\tWqnw}{\ensuremath{\tilde{W}_{\q}^{\smallnwarrow}}}
\newcommand{\tWqse}{\ensuremath{\tilde{W}_{\q}^{\smallsearrow}}}

\newcommand{\Xse}{\ensuremath{X_{\q}}}
\newcommand{\Wse}{\ensuremath{W_{\q}}}
\newcommand{\Xp}{\ensuremath{X_{\p}}}
\newcommand{\Xq}{\ensuremath{X_{\q}}}
\newcommand{\Xm}{\ensuremath{X_{\m}}}
\newcommand{\Wp}{\ensuremath{W_{\p}}}
\newcommand{\Wq}{\ensuremath{W_{\q}}}
\newcommand{\Wm}{\ensuremath{W_{\m}}}

\newcommand{\Xe}{\ensuremath{X_{\bar{e}_1}}}
\newcommand{\We}{\ensuremath{W_{\bar{e}_1}}}
\newcommand{\Xee}{\ensuremath{X_{\bar{e}_2}}}

\newcommand{\dom}{\textup{dom}}
\newcommand{\ran}{\textup{ran}}
\newcommand{\res}{\restriction}
\newcommand{\st}{\ensuremath{\::\:}}
\newcommand{\acts}{\ensuremath{\curvearrowright}}

\newcommand{\ld}{\ensuremath{\ell}}

\newcommand{\M}{\ensuremath{\mathcal{M}}}
\newcommand{\G}{\ensuremath{\mathcal{G}}}
\newcommand{\GG}{\ensuremath{\mathscr{G}}}
\newcommand{\F}{\ensuremath{\mathcal{F}}}
\newcommand{\E}{\ensuremath{\mathcal{E}}}
\newcommand{\W}{\ensuremath{W}}

\newcommand{\emp}{\emptyset}
\renewcommand{\mod}{\textup{ mod }}

\allowdisplaybreaks

\begin{document}

\begin{abstract}
We conjecture that for each integer $\ell\geq 2$, the equivalence relation arising from the shift action of $\N^{\ell}$ on $2^{\omega^{\ell}}$ is hypersmooth. We prove this conjecture for $\ell=2$, reduce the general case to a problem in Borel graph colorings, and establish the general case in measure. We also show that under certain mild hypotheses on a hypersmooth equivalence relation $E$, every finite-index Borel super-equivalence relation of $E$ is hypersmooth.
\end{abstract}

\maketitle
\tableofcontents

\section{Introduction}\label{sec:intro}

The Borel equivalence relation $E$ is \emph{smooth} if $E$ Borel reduces to the equality relation on some standard Borel space, and \emph{hypersmooth} if $E$ is the union of an increasing sequence of smooth equivalence relations. In this paper, we show that the equivalence relation generated by the shift action of $\N^2$ on $\TT$ is hypersmooth. Some portions of our proof go through in higher dimensions, and we conjecture that the same conclusion holds there. Given a standard Borel space $\Omega$ and positive integer $\ell$, the shift action $\N^\ell\acts\Ol$ of $\N^{\ell}$ on $\Omega^{\omega^\ell}$ is defined by 
\begin{equation}
\big((a_1,\ldots,a_\ell)\cdot x\big)(k_1,\ldots,k_\ell) \ := \ x(k_1+a_1,\ldots,k_{\ell}+a_\ell).
\end{equation}
Let us write $\EOl$ for the Borel equivalence relation on $\Ol$ arising from this $\N^{\ell}$-action, so
\begin{equation}\label{eq:OmegaNotation}
x\mathrel{\EOl}y \quad :\Longleftrightarrow \quad (\exists\bar{a},\bar{b}\in\N^\ell)\,(\bar{a}\cdot x=\bar{b}\cdot y),
\end{equation}
where we write $\bar{a}=(a_1,\ldots,a_\ell)$ for elements of $\Z^{\ell}$. We conjecture the following:

\begin{conjecture}\label{con:HigherDimensions}
For every integer $\ell\geq 1$, the equivalence relation $\ETl$ is hypersmooth.
\end{conjecture}

That Conjecture \ref{con:HigherDimensions} holds for $\ell=1$ was shown in \cite[8.1]{DJK}, where in fact it was shown that $E^{\Omega}_1$ is hypersmooth for any $\Omega$.  
We will prove Conjecture \ref{con:HigherDimensions} in the case $\ell=2$ and provide a framework for understanding the general case, in particular reducing it to a problem in Borel graph colorings. This reduction, along with our work in the 2-dimensional case, will enable us to show that Conjecture \ref{con:HigherDimensions} holds on a full measure subset of $\Tl$ for all $\ell\geq 3$. We point out that Conjecture \ref{con:HigherDimensions} implies the more general statement that the equivalence relation generated by a finite family of commuting Borel endomorphisms of a standard Borel space is hypersmooth.

\begin{conjecture}\label{con:final} Let $X$ be a standard Borel space, let $\ell$ be a positive integer, and let $U_1,\ldots,U_\ell:X\to X$ be Borel functions on $X$ such that $U_i\circ U_j=U_j\circ U_i$ for all $1\leq i,j\leq\ell$. Let $E$ be the equivalence relation on $X$ defined by 
\[
x\mathrel{E}y \quad :\Longleftrightarrow \quad (\exists\bar{a},\bar{b}\in\N^{\ell})  
\big(U_1^{a_1}\circ \cdots\circ U_\ell^{a_\ell}\big)(x)=\big(U_1^{b_1}\circ\cdots\circ U_\ell^{b_\ell}\big)(y).
\]
Then $E$ is hypersmooth.
\end{conjecture}

\begin{proof}[Proof Sketch that \ref{con:HigherDimensions} implies \ref{con:final}]
A standard argument shows $E\sqsubseteq_B E^{2^\omega}_{\ell}$, and a pair of coding arguments show $E^{2^\omega}_{\ell}\sqsubseteq_B E^3_{\ell+1}$ and $E^3_{\ell+1}\sqsubseteq_B E^2_{\ell+1}$. (For a detailed proof, see Section \ref{sec:measure}.)
\end{proof}

Hypersmooth equivalence relations have received a great deal of attention since the 1990s, mostly devoted to the important subclass of hyperfinite equivalence relations (see Section \ref{sec:prelimBERs} for the basic definitions, facts, and references concerning these notions). Every hyperfinite equivalence relation is countable, and therefore arises as the orbit equivalence relation induced by a Borel action of a countable group by the celebrated Feldman-Moore Theorem \cite{FM}. In fact, as has long been known (\cite{Weiss84}, \cite{SS88}), the equivalence relation $E$ is hyperfinite if and only if $E$ arises as the orbit equivalence relation of a Borel $\Z$-action (equivalently, is generated by a single Borel automorphism), and finding ever larger classes of countable groups with the property that all their Borel actions yield hyperfinite orbit equivalence relations has been a significant theme of research in the subject, aimed at answering Weiss's question of whether the countable amenable groups form such a class (see \cite{AsymDim} for some of this project's history, and for its most recent and strongest result).  Analogously, every equivalence relation that arises from a Borel $\N$-action (equivalently, is generated by a single Borel endomorphism) is hypersmooth, and this paper initiates the parallel project for semigroups of Borel endomorphisms, with an ultimate goal, perhaps, of showing that Conjecture \ref{con:final} holds not just for finite commuting families of Borel endomorphisms but even for countable ones, mirroring the Gao-Jackson result for Borel actions of countable abelian groups \cite{GJ}.

Previous work on this project in the hyperfinite setting (see for instance \cite{JKL}, \cite{GJ}, \cite{Nilpotent}, and \cite{AsymDim}) has tended to follow the same general pattern of distinguishing the free and non-free portions of the action, making critical use of Borel marker sets, exploiting the local geometry of the acting group (often in terms of a Cayley graph), and finally appealing to the fact that a countable Borel equivalence relation that contains a hyperfinite sub-equivalence relation of finite index is itself hyperfinite.

Our arguments, too, will involve these elements.  We will show separately that the restriction of $\ETT$ to the periodic part $\textup{Per}(\TT):=\{x\in\TT\st (\exists\bar{a}\ne\bar{b}\in\N^2) \, \bar{a}\cdot x=\bar{b}\cdot x\}$ is hypersmooth, while focusing most of our attention on the free part $\textup{Fr}(\TT):=\TT\setminus\textup{Per}(\TT)$. Using the notion of Borel chromatic number from \cite{KST}, we will construct Borel marker sets in certain $\ETT$-invariant Borel subsets of $\FTT$ and use them to show $\ETT$ is hypersmooth on these sets in a manner reminiscent of the proof that Borel actions of finitely-generated groups of polynomial growth generate hyperfinite orbit equivalence relations (see \cite[11.1]{JKL}). We will make extensive use of the Euclidean planar geometry that can be associated to each $\ETT$-class $C\subseteq\FTT$, as described in Section \ref{subsec:geometry}.  And finally, since it is not currently known whether a Borel equivalence relation with a hypersmooth sub-equivalence relation of finite index is itself hypersmooth, we will need to establish this fact in the instances that concern us.

Despite these broad similarities, there are significant obstacles we will face that do not arise in the hyperfinite setting, and which provide some explanation for why such impressive progress has been made on Weiss's question while we are still stuck, forty years after Weiss showed Borel $\Z^n$-actions yield hyperfinite orbit equivalence relations, at the level of $\N^2$ in the context of semigroup actions. The primary source of these obstacles is the most obvious difference between the two contexts, namely that we will be dealing with a semigroup of endormorphisms with uncountable fibers rather than with a group of invertible automorphisms. The continuous maps $x\mapsto (1,0)\cdot x$ and $x\mapsto (0,1)\cdot x$ that generate the $\N^2$-action on $\TT$ are not injective or even countable-to-one; consequently, our equivalence classes are uncountable, the graphs we will attempt to color are locally uncountable, images of Borel sets under semigroup elements need not be Borel, saturations of Borel sets under our equivalence relations need not be Borel, and balls around Borel sets using the local geometry we will introduce need not be Borel. 

Indeed, by Theorems 1 and 4.1 of \cite{KL}, for any hypersmooth equivalence relation $E$ that, like $\ETT$, is not essentially hyperfinite, there is no Borel action of a Polish group that realizes $E$ as its orbit equivalence relation. Consequently, equivalence classes in $\FTT$ do not actually have the geometric structure of any Polish group, and the ``Euclidean planar geometry" on each $\ETT$-class $C\subseteq\FTT$ that we alluded to above and discuss in Section \ref{subsec:geometry} actually arises on a quotient of $C$ by a non-smooth sub-equivalence relation of $\ETT$ rather than on $C$ itself. We will say a bit more in Section \ref{sec:cascade} about how this impacts the construction of marker sets and the ``local Borel algorithms'' that feature prominently in some of the hyperfiniteness proofs, but the most notable consequence is that our constructions do not work, and cannot work, on the entirety of $\FTT$. 
Specifically, we will define a certain family of Borel marker sets, prove that the existence of such a family is equivalent to the existence of finite Borel colorings of certain graphs, and show that the existence of either in an invariant Borel set $X\subseteq\FTT$ implies hypersmoothness of $\ETT\res X$. Then we will also observe that no such family of Borel marker sets or finite Borel graph colorings can actually exist throughout all of $\FTT$.  What we will do, instead, is obtain a partition of $\FTT$ into a measure one ``tame'' set $X$ that admits the desired Borel marker sets and graph colorings, and a measure zero ``wild'' set $W$ that does not, but on which $\ETT$ is sufficiently simple that we can show it to be hypersmooth directly by reducing it to $E_1^{2^\omega}$.

For this reason, our claim to have reduced Conjecture \ref{con:HigherDimensions} to a problem in Borel graph colorings requires an asterisk, because while we do prove that hypersmoothness of $\ETl\res\FTl$ follows from the existence of finite Borel colorings of certain graphs on $\FTl$, we also point out that no such colorings exist. The challenge in extending our proof to dimensions $\ell\geq 3$, therefore, lies not just in coloring certain graphs, but in finding an appropriate partition of $\FTl$ into a colorable ``tame'' part and an uncolorable ``wild'' part that is nevertheless sufficiently structurally constrained that we can prove hypersmoothness of $\ETl$ on it by other means. For $\ell\geq 3$, we will be able to find a full measure colorable portion, but not one with a sufficiently nice complement to enable us to prove Conjecture \ref{con:HigherDimensions} outright.

\begin{remark}
Our coloring arguments rely on the base space $\Omega$ being finite, so our proof that $\ETT$ is hypersmooth goes through with any finite base space replacing $2=\{0,1\}$, but it does not work for an infinite base space $\Omega$. In particular, we do not yet know how to show $E^{2^\omega}_2$ is hypersmooth, which is unfortunate since hypersmoothness of $E^{2^\omega}_2$ would imply that any equivalence relation generated by a commuting pair of Borel endomorphisms is hypersmooth. In Example \ref{ex:stuck}, we describe a specific measure zero subset of $(2^\omega)^{\omega^2}$ on which our graphs have infinite Borel chromatic number but we do not yet know how to prove hypersmoothness of $E^{2^\omega}_2$ by other means.
\end{remark}

\begin{remark}
In recent, independent, and as yet unpublished work, Petr Naryshkin, Forte Shinko, Felix Weilacher, and Jing Yu have shown using the theory of Borel asymptotic dimension introduced in \cite{AsymDim} that the equivalence relation generated by a commuting family of countably many bounded-to-one Borel endomorphisms of a standard Borel space is hyperfinite, thus establishing a special case of the countable version of Conjecture \ref{con:final}. Their argument focuses on commuting pairs of bounded-to-one Borel endomorphisms, and extends to countable families via the machinery of Borel asymptotic dimension (see \cite{AsymDim}).  Neither our result not theirs implies the other, since theirs holds for countable families of Borel endormorphisms but assumes these endomorphisms to be bounded-to-one, so their equivalence relations are countable and they can use the Feldman-Moore Theorem and Borel asymptotic dimension. 
Notably, it is not yet clear whether or how the notion of Borel asymptotic dimension can be adapted to the setting of Borel equivalence relations with uncountable classes, although the equivalence relations $\ETl$ should serve as interesting test cases, and it seems plausible that ideas introduced in \cite{AsymDim} could play a role in attempting to prove the countably infinite versions of Conjectures \ref{con:HigherDimensions} and \ref{con:final}.
\end{remark}

The remainder of this paper is organized as follows. In Section \ref{sec:prelim}, we recall some basic definitions and facts we shall need concerning Borel equivalence relations and graph colorings, and then establish some definitions, results, and notations concerning $\N^{\ell}$-shift actions and the geometry of $\EOl$-classes. In Section \ref{sec:hs}, we recall various definitions and facts we will need concerning  hypersmooth equivalence relations, and obtain some conditions under which a finite-index super-equivalence relation of a hypersmooth equivalence relation is again hypersmooth.  In Section \ref{sec:cascade}, we discuss certain types of Borel marker sets in the context of $\N^{\ell}$-shift actions, and show how the existence of such marker sets can be used to prove hypersmoothness of $\EOl$ on invariant Borel subsets of $\Ol$. In Section \ref{sec:graphs}, we introduce various Borel graphs related to the $\N^{\ell}$-shift actions, and establish a connection between the existence of the Borel marker sets we use to prove hypersmoothness and finite Borel colorings of these graphs.  It is here that we state the main technical result of the paper (Theorem \ref{thm:main}), which gives a bifurcation of the free part of $\N^2\acts\TT$ into a ``tame'' portion where we can find good marker sets and a ``wild'' portion where our graphs have infinite Borel chromatic number but where we can show $\ETT$ to be hypersmooth by other means. Section \ref{sec:periodicity} is devoted to proving that the restriction of $\ETT$ to the periodic or ``non-free'' part of $\N^2\acts\TT$ is hypersmooth, after we first develop some general facts about periodicity in arbitrary dimensions.  Sections \ref{sec:ep}--\ref{sec:Gm} are then devoted to proving Theorem \ref{thm:main}, which involves considering in turn the various graphs we want to color. Finally, in Section \ref{sec:measure} we prove Conjecture \ref{con:HigherDimensions} in measure before concluding with a proof that Conjecture \ref{con:HigherDimensions} implies Conjecture \ref{con:final}. 

\bigskip

\noindent \textbf{Acknowledgements}. I would like to thank Brandon Seward for reading and providing helpful comments on an early draft of this paper.

\bigskip

\noindent \textbf{Use of AI}. AI was not used in any way to write any portion of or obtain any results in this paper.

\section{Preliminaries}\label{sec:prelim}

\subsection{Borel equivalence relations}\label{sec:prelimBERs}

Let $E$ and $F$ be equivalence relations on the sets $X$ and $Y$, respectively, and let $x\in X$ and $A\subseteq X$. We write $[x]_E:=\{y\in X\st x\mathrel{E}y\}$ for the $E$-equivalence class of $x$, and $[A]_E:=\{x\in X\st (\exists y\in A)\; x\mathrel{E}y\}$ for the \emph{$E$-saturation} of $A$. The set $A$ is \emph{$E$-invariant} if $[A]_E=A$, and is a \emph{complete section} for $E$ if $A$ meets every $E$-class.  We write $E\res A:=\{(x,y)\in A^2\st x\mathrel{E}y\}$ for the \emph{restriction} of $E$ to $A$, and $\Delta_X$ and $I_X$ for the equality and indiscrete equivalence relations on $X$, respectively, so $I_X:=X^2$. The product equivalence relation $E\times F$ on $X\times Y$ is defined by $(x,y)\mathrel{E\times F}(x',y') \ :\Leftrightarrow \ x\mathrel{E}x' \, \wedge \, y\mathrel{F}y'$.  If $X=Y$ and $E\subseteq F$, the quotient equivalence relation $F/E$ on $X/E$ is defined by $[x]_E\mathrel{F/E}[y]_E \ :\Leftrightarrow \ x\mathrel{F}y$.  A \emph{reduction} from $E$ to $F$ is a function $f:X\to Y$ such that for all $x,x'\in X$, we have $x\mathrel{E}x'$ if and only if $f(x)\mathrel{F}f(x')$. An injective reduction is called an \emph{embedding}, and a bijective reduction an \emph{isomorphism}. The reduction $f:X\to Y$ from $E$ to $F$ is called \emph{invariant} if its range $\ran(f)$ is $F$-invariant.

The above notions make sense for any equivalence relations, but all our work will be done in the Borel setting, meaning that our spaces are standard Borel and all our equivalence relations and reductions are Borel. A \emph{standard Borel space} is a measurable space $(X,\B)$ for which there exists a Polish (i.e., separable and completely metrizable) topology on $X$ that generates $\B$ as its $\sigma$-algebra of Borel sets. The equivalence relation $E$ on the standard Borel space $X$ is \emph{Borel} if $E$ is a Borel subset of $X\times X$.  If $E$ and $F$ are Borel equivalence relations on the standard Borel spaces $X$ and $Y$, respectively, then $E$ is \emph{Borel reducible} to $F$ if there exists a Borel reduction from $E$ to $F$, and the terms \emph{Borel embeddable}, \emph{invariantly Borel embeddable}, and \emph{Borel isomorphic} are defined analogously.  We will use the following standard notation:
\[
\begin{array}{lll}
E\leq_BF        & :\Longleftrightarrow & \mbox{$E$ is Borel reducible to $F$;} \\
E\sqsubseteq_BF & :\Longleftrightarrow & \mbox{$E$ is Borel embeddable in $F$;} \\
E\sqsubseteq_B^iF & :\Longleftrightarrow & \mbox{$E$ is invariantly Borel embeddable in $F$;} \\
E\cong_BF       & :\Longleftrightarrow & \mbox{$E$ is Borel isomorphic to $F$.} \\
\end{array}
\]

The Borel equivalence relation $E$ on the standard Borel space $X$ is \emph{finite} if each $E$-class is finite, \emph{countable} if each $E$-class is countable, and \emph{smooth} if $E\leq_B\Delta_Y$ for some (equivalently, any uncountable) standard Borel space $Y$. A \emph{transversal} for $E$ is a subset of $X$ that meets each $E$-class in exactly one point, and a \emph{selector} for $E$ is a reduction $f:X\to X$ from $E$ to $\Delta_X$ such that $f(x)\mathrel{E}x$ for all $x\in X$. We call $E$ \emph{selective} if $E$ admits a Borel selector, which is equivalent to admitting a Borel transversal. The equivalence relation $E$ is \emph{hyperfinite} if it can be expressed as the union of an increasing sequence of finite Borel equivalence relations, and \emph{hypersmooth} if it can be expressed as the union of an increasing sequence of smooth equivalence relations. Analogously, we call $E$ \emph{hyperselective} if it can be expressed as the union of an increasing sequence of selective equivalence relations.

Every finite Borel equivalence relation is selective, and every selective equivalence relation is smooth, so every hyperfinite equivalence relation is hyperselective and every hyperselective equivalence relation is hypersmooth; in fact, by \cite[5.1]{DJK}, the equivalence relation $E$ is hyperfinite if and only if $E$ is hypersmooth and countable.  Every smooth \emph{countable} Borel equivalence relation is selective, but not every smooth equivalence relation is selective, and not every hypersmooth equivalence relation is hyperselective. Our terms ``selective'' and ``hyperselective'' are not widely used, but we prefer them to the terms ``strongly smooth'' and ``strongly hypersmooth'' that were introduced for these notions in \cite{KL}, since there are other notable strengthenings of smoothness (see for instance \cite{Simultaneous}).  In fact, the distinction between smooth and selective is rarely emphasized, since most naturally occuring examples of smooth equivalence relations are selective.  Our proof that $\ETT$ is hypersmooth will actually show that it is hyperselective, but we emphasize the distinction primarily because in Section \ref{sec:hs} we pose a question and establish some results that are sensitive to it. 

Hyperfinite and hypersmooth equivalence relations were studied and classified in \cite{DJK} and \cite{KL}, respectively, and the former have received extensive subsequent attention (\cite{JKL}, \cite{KM}, \cite{BJ}, \cite{GJ}, \cite{Nilpotent}, \cite{CJMST2},  \cite{AsymDim}, \cite{GaoXiao}, \cite{FSV}). There exist universal instances of each, namely $E_0$ and $E_1$, the eventual equality relations on $2^\omega$ and on $(2^\omega)^\omega$, respectively. Formally these are defined as follows: 
\begin{align*}
\text{for all $x,y\in 2^\omega$,} \quad x\mathrel{E_0}y \ &:\Longleftrightarrow \ (\exists n)(\forall k)\, x(n+k)=y(n+k); \\
\text{for all $x,y\in (2^\omega)^\omega$,} \quad x\mathrel{E_1}y \ &:\Longleftrightarrow \ (\exists n)(\forall k)\, x(n+k)=y(n+k).
\end{align*}

\begin{fact}[see \cite{DJK} and \cite{KL}]
For every Borel equivalence relation $E$, we have the following:
\begin{enumerate}[(i)]
 \item[(i)] $E$ is hyperfinite \ $\Longleftrightarrow$ \ $E\leq_BE_0$;
 \item[(ii)] $E$ is hypersmooth \ $\Longleftrightarrow$ \ $E\leq_BE_1$ \ $\Longleftrightarrow$ \ $E\sqsubseteq_BE_1$. \hfill $\square$
\end{enumerate}
\end{fact}

Sequences are $E_0$ or $E_1$ equivalent if they agree up to finite initial segments of the same length; if we ask that they agree up to finite initial segments of possibly different lengths, we get the so-called ``tail'' analogues of $E_0$ and $E_1$ defined as follows:
\begin{align*}
\text{for all $x,y\in 2^\omega$,} \quad x\mathrel{E_{0,t}}y \ &:\Longleftrightarrow \ (\exists m,n)(\forall k)\, x(m+k)=y(n+k); \\
\text{for all $x,y\in (2^\omega)^\omega$,} \quad x\mathrel{E_{1,t}}y \ &:\Longleftrightarrow \ (\exists m,n)(\forall k)\, x(m+k)=y(n+k).
\end{align*}

The equivalence relation $E_{0,t}$ is usually just called $E_t$, while $E_{1,t}$ has no special name in the literature. $E_{0,t}$ is Borel bi-embeddable with $E_0$ (that is, each Borel embeds into the other), while $E_{1,t}\cong_BE_1$. Of course, in our notation \eqref{eq:OmegaNotation}, $E_{1,t}$ is just the equivalence relation $E_{1}^{2^\omega}$ arising from the shift action of $\N$ on $(2^\omega)^\omega$. That $E_{1,t}$ is hypersmooth was first proved in \cite{DJK}; for a proof that $E_{1,t}\cong_BE_1$, see the discussion following \cite[Proposition 3.4]{KL}. Since hypersmoothness of $E_{1,t}$ will be important in what follows, we record it here formally.

\begin{fact}[see \cite{DJK} and \cite{KL}]
$E_{1,t}$ is hypersmooth; in fact, $E_{1,t}\cong_BE_1$. \hfill $\square$
\end{fact}

Of course, $E_1$ (and consequently $E_{1,t}$) is hyperselective, since $E_1$ is the increasing union of the selective equivalence relations $F_n$ on $(2^\omega)^\omega$ defined by $x\mathrel{F_n}y \ :\Leftrightarrow \ (\forall k)\, x(n+k)=y(n+k)$, or, in our notation \eqref{eq:OmegaNotation}, $x\mathrel{F_n}y \ \Leftrightarrow \ n\cdot x=n\cdot y$. In Section \ref{sec:hs}, we will recall some well-known closure properties of smooth and hypersmooth equivalence relations, and prove some additional facts that we will need concerning equivalence relations that contain a hypersmooth sub-equivalence relation of finite index.

\subsection{Graphs, colorings, and Borel chromatic numbers}\label{subsec:graphs}

Our proof that $\ETT$ is hypersmooth will involve Borel colorings of graphs, which were introduced and studied in \cite{KST}.  Intuitively, a coloring of the graph $\G$ assigns to each vertex in $\G$ one of a given set of colors so that adjacent vertices receive different colors. The $\N^{\ell}$-shift action on the product space $\Ol$ naturally induces the directed graph $\GG_{\ell}$ on $\Ol$ defined by
\[
x\mathrel{\GG_{\ell}}y \quad :\Longleftrightarrow \quad (\exists i\leq\ell) \ \bar{e}_i\cdot x=y,
\]
where $\bar{e}_i$ is the $i$th standard basis vector in $\N^{\ell}$. We shall attempt to color restrictions of these graphs, and others related to them, to various subsets of $\Ol$. The graphs $\GG_{\ell}$ are not simple, since we can have $\bar{e}_i\cdot x=x$. Since non-simple graphs never admit colorings, discussions of graph colorings often assume all graphs under consideration to be simple, and we will follow this practice. In this subsection, we gather the general definitions and notation that we will use regarding graphs and colorings, and in Section \ref{sec:graphs} we introduce the specific graphs that we then attempt to color in Sections \ref{sec:ep}--\ref{sec:Gm}. 

It will be convenient for us to work with both graphs and directed graphs, and to view graphs as special cases of directed graphs. Formally, we regard a \emph{directed graph} on the set $X$ as an irreflexive binary relation on $X$, and a \emph{graph} on $X$ as a symmetric directed graph on $X$. Suppose $\G$ is a directed graph on the set $X$. Elements of $X$ are called \emph{vertices} of $\G$, elements of $\G$ are called \emph{edges} of $\G$, and elements $x$ and $y$ in $X$ are \emph{neighbors} in $\G$ if $x\mathrel{\G}y$ or $y\mathrel{\G}x$. A \emph{path} in $\G$ from $x_0$ to $x_r$ is a finite sequence $\langle x_0,\ldots,x_r\rangle$ in $X$ with $r\geq 1$ such that $x_i\mathrel{\G}x_{i+1}$ for each $0\leq i<r$ and $x_i\ne x_j$ for all $0\leq i<j\leq r$ except possibly $x_0=x_r$. For any $x\in X$, we write $\G(x) := \{y\in X\st x\mathrel{\G}y\}$ for the set of \emph{children} of $x$ and $\G^{-1}(x) := \{y\in X\st y\mathrel{\G}x\}$ for the set of \emph{parents} of $x$, and call $|\G(x)|$ the \emph{out-degree} of $x$, $|\G^{-1}(x)|$ the \emph{in-degree} of $x$, and $|\G(x)\cup\G^{-1}(x)|$ the \emph{degree} of $x$. For any subset $Y\subseteq X$, the \emph{restriction} of $\G$ to $Y$, written $\G\res Y$, is the graph $\G\cap Y^2$ on $Y$.

Every directed graph $\G$ on $X$ induces the graph $\overline{\G}$ on $X$, its \emph{symmetrization}, defined by 
\[
x\mathrel{\overline{\G}}y \quad :\Longleftrightarrow \quad x\mathrel{\G}y \ \vee \ y\mathrel{\G}x.
\]
Let $\G$ be a directed graph on $X$, with $\overline{\G}$ is associated graph. Two elements $x,y\in X$ are \emph{connected} in $\G$ if there is a path in $\overline{\G}$ from $x$ to $y$, and \emph{strongly connected} in $\G$ if there is a path in $\G$ from $x$ to $y$. Being connected in $\G$ is an equivalence relation on $X$, the equivalence classes of which are called the \emph{connected components} of $\G$. The directed graph $\G$ is called \emph{connected} if it has just one connected component, that is, if any two elements of $X$ are $\G$-connected.

A \emph{coloring} of the directed graph $\G$ on $X$ is a function $c:X\to Z$ into some set $Z$ of \emph{colors} such that for all $(x,y)\in\G$ we have $c(x)\ne c(y)$. Note that if $c$ is a coloring of the directed graph $\G$, then $c$ is also a coloring $\overline{\G}$, so in our discussions of graph colorings we need not distinguish between directed and undirected graphs. We often call a coloring taking values in a set of cardinality $\kappa$ a \emph{$\kappa$-coloring}, and one taking values in a finite set a \emph{finite coloring}. As a trivial example, the identity map on the set of vertices of any directed graph $\G$ is a coloring of $\G$. The smallest cardinal $\kappa$ for which 
$\G$ admits a $\kappa$-coloring is called the \emph{chromatic number} of $\G$, written $\chi(\G)$.

Suppose now that $X$ is a standard Borel space. The directed graph $\G$ on $X$ is \emph{Borel} if $\G$ is a Borel subset of $X\times X$. Every directed graph we consider will be a Borel directed graph on a standard Borel space. A coloring $c:X\to Z$ of the Borel directed graph $\G$ on the standard Borel space $X$ is called \emph{Borel} if $Z$ is a standard Borel space and the function $c$ is Borel. If $\G$ is a Borel directed graph on a standard Borel space, then $\G$ admits a Borel $2^{\aleph_0}$-coloring, and the \emph{Borel chromatic number} of $\G$, written $\chi_B(\G)$, is the smallest cardinal $\kappa$ for which there exists a Borel $\kappa$-coloring of $\G$. Similarly, if $X$ is a Polish space and $\G$ is a Borel graph on $X$, a \emph{continuous coloring} of $\G$ is a coloring $c:X\to Z$ of $\G$ where $c$ is a continuous function into some Polish space $Z$. The \emph{continuous chromatic number of $\G$}, written $\chi_C(\G)$, is the smallest cardinal $\kappa$ for which there exists a continuous $\kappa$-coloring of $\G$. Since every uncountable standard Borel space has cardinality $2^{\aleph_0}$, the only infinite values $\chi_B(\G)$ and $\chi_C(\G)$ could have are $\aleph_0$ and $2^{\aleph_0}$.

We conclude this subsection with a couple basic facts from \cite{KST} and an important and well-known example of a graph with finite chromatic number but infinite Borel chromatic number. For any standard Borel space $\Omega$ and integer $\ell\geq 1$, let $\vec{\G}_{\ell}(\Omega)$ be the directed graph on $\Omega^{\omega^\ell}$ defined by
\begin{equation}\label{eq:Gl1}
x\mathrel{\vec{\G}_{\ell}(\Omega)}y \quad :\Longleftrightarrow \quad x\ne y \ \wedge \ (\exists i\leq\ell) \, \bar{e}_i\cdot x=y.
\end{equation}
Note that the out-degree of any vertex of $\vec{\G}_{\ell}(\Omega)$ is at most $\ell$, and that for any Borel set $B\subseteq\Omega^{\omega^\ell}$ the set of $\vec{\G}_{\ell}(\Omega)$-parents of vertices in $B$ is Borel (though possibly uncountable). Therefore, Fact \ref{fact:countable} below follows immediately from \cite[4.5]{KST}, while Fact \ref{fact:KST5.2} is shown in \cite[5.2]{KST} and the discussion that follows.

\begin{fact}\label{fact:KST5.2}
For every integer $J\geq 2$, we have $\chi_B\big(\vec{\G}_{1}(J)\big)=3$. \hfill $\square$
\end{fact}

\begin{fact}\label{fact:countable}
For every standard Borel space $\Omega$ and integer $\ell\geq 1$, we have $\chi_B\big(\vec{\G}_{\ell}(\Omega)\big)\leq\aleph_0$. 
\hfill $\square$
\end{fact}

%

\begin{example}[See {\cite[3.2]{KST}}]\label{ex:BadAssGraph}
Let $[\omega]^\omega$ be the set of all strictly increasing sequences of natural numbers. Then $[\omega]^{\omega}$ is a closed subspace of $\omega^\omega$, hence is Polish in the subspace topology. Let $\sigma$ be the shift map on $[\omega]^{\omega}$, defined by $\sigma(x)(n)=x(n+1)$ for all $x\in [\omega]^{\omega}$ and $n\in\omega$. Let $\vec{\GG}$ be the Borel directed graph on $[\omega]^{\omega}$ generated by $\sigma$, so $x\mathrel{\vec{\GG}}y \ :\Leftrightarrow \ \sigma(x)=y$. Then $\chi_B(\vec{\GG})=\aleph_0$, even though $\chi(\vec{\GG})=2$.
\end{example}

%
%

\subsection{$\N^{\ell}$-actions on product spaces}

As is customary in logic, we write $\N=\omega$ for the set of natural numbers, include $0\in\N$, and identify each $n\in\N$ with the set $\{k\in\N\st k<n\}$. We tend to write $\N$ when we are thinking of the additive commutative monoid, and $\omega$ otherwise.  We write $\Z_+$ for the set of positive integers and $\Z_-$ for the set of negative integers. 
For $\ell\in\Z_+$ and $\bar{a}\in\Z^\ell$, let $|\bar{a}|:=\max\{|a_1|,\ldots,|a_\ell|\}$ and define the well-ordering $<$ on $\Z^\ell$ by
\begin{equation}\label{eq:WO}
\bar{a}<\bar{b} \quad :\Longleftrightarrow \quad |\bar{a}|<|\bar{b}| \mbox{ or $|\bar{a}|=|\bar{b}|\,\wedge\,\bar{a}<_{\mbox{\scriptsize{lex}}}\bar{b}$}\,.
\end{equation}
Frequently we will use the restriction of $<$ to subsets of $\Z^{\ell}$, especially $\N^{\ell}$. We also write
\begin{equation}\label{eq:WellOrder}
\N^{\ell}_{\leq n} \ := \ \{\bar{a}\in\N^\ell\st |\bar{a}|\leq n\} \qquad \text{and} \qquad \N^{\ell}_{<n} \ := \ \{\bar{a}\in\N^\ell\st |\bar{a}|< n\},
\end{equation}
define $\Z^{\ell}_{\leq n}$ and $\Z^{\ell}_{<n}$ analogously, and write $\Z_{\ne 0}=\Z\setminus\{0\}$.

The commutative monoids $\N^{\ell}$ act by shifts on spaces of the form $\Ol$, where $\ell\in\Z_+$ and $\Omega$ is a standard Borel space with $|\Omega|\geq 2$. For convenience we call any space of this form a \emph{product space}, with \emph{base space} $\Omega$ and \emph{dimension} $\ell$. The shift action $\N^{\ell}\curvearrowright\Ol$ is given by
\begin{equation}\label{eq:action}
\big((a_1,\ldots,a_\ell)\cdot x\big)(k_1,\ldots,k_\ell) \ := \ x\big(k_1+a_1,\ldots,k_\ell+a_\ell\big).
\end{equation}
The \emph{periodic} and \emph{free} parts of $\Ol$ under this action are the $\EOl$-invariant Borel sets
\begin{equation}
\POl \ := \ \big\{x\in\Ol\st (\exists \bar{a}\ne\bar{b}\in\N^{\ell})\,\bar{a}\cdot x=\bar{b}\cdot x\big\} \qquad \text{and} \qquad \FOl \ := \ \Ol\setminus \POl.
\end{equation}
Given such a product space $\Ol$, we write $\EOl$ in place of the more customary (but for us cumbersome) notation $E_{\N^{\ell}}^{\Omega^{\omega^\ell}}$ to denote the equivalence relation on $\Ol$ defined by
\begin{equation}\label{eq:GeneralShiftEquivRel}
x\mathrel{\EOl}y \quad :\Longleftrightarrow \quad (\exists\bar{a},\bar{b}\in\N^{\ell})(\bar{a}\cdot x=\bar{b}\cdot y).
\end{equation}
We also define the ``diagonal'' action of $\N$ on $\Ol$ by
\begin{equation}
(n,x) \ \mapsto \ (n,\ldots,n)\cdot x.
\end{equation}
It will be convenient to have notation for constant tuples in $\N^{\ell}$, so when the dimension $\ell$ is clear from context we will write $\d{n}:=(n,\ldots,n)\in\N^{\ell}$. This allows us to write $\d{n}\cdot x$ for the diagonal $\N$-action, and to write, e.g., $\d{n}+\bar{a}$ for 
$(n+a_1,\ldots,n+a_{\ell})$. 
Given $Y\subseteq\Ol$ and $\bar{a}\in\N^{\ell}$, we will use the notation
\begin{equation}
\bar{a}\cdot Y := \{\bar{a}\cdot y\st y\in Y\} \qquad \text{and} \qquad \bar{a}^{-1}[Y] := \{x\in\Ol\st \bar{a}\cdot x\in Y\}.
\end{equation}
In general, given any sets $X$ and $Y$ and function $f:X\to Y$ with $A\subseteq X$ and $B\subseteq Y$, we write $f[A]$ and $f^{-1}[B]$ for the forward and inverse images of $A$ and $B$ under $f$, respectively, and $Y^X$ for the set of all functions from $X$ to $Y$.

For each $\bar{c}\in\N^{\ell}$, let $F^{\Omega}_{\ell,\bar{c}}$ be the selective sub-equivalence relation of $\EOl$ defined on $\Ol$ by 
\begin{equation}
x\mathrel{F^{\Omega}_{\ell,\bar{c}}}y \quad :\Longleftrightarrow \quad \bar{c}\cdot x=\bar{c}\cdot y.
\end{equation}
As important special cases, for each $n\in\omega$ we have
\begin{equation}
x\mathrel{\FOln}y \quad \Longleftrightarrow \quad \d{n}\cdot x=\d{n}\cdot y.
\end{equation} 
Additionally write $F_{\ell}^{\Omega}:=\cup_n\FOln$, so $F_{\ell}^{\Omega}$ is a hyperselective sub-equivalence relation of $\EOl$.

For every product space $\Ol$ that we consider, we will have in mind a fixed Polish topology on $\Omega$ that is compatible with its Borel structure, and we view $\Ol$ as a Polish space in the product topology. For each $\bar{a},\bar{b},\bar{c}\in\N^{\ell}$, let us write
\begin{equation}\label{eq:prec}
\bar{a}\preceq\bar{b} \ :\Longleftrightarrow \ (\forall i) \, a_i\leq b_i, \quad 
\bar{a}\precneq\bar{b} \ :\Longleftrightarrow \ \bar{a}\preceq\bar{b} \, \wedge \, \bar{a}\ne\bar{b}, \quad\text{and}\quad 
\bar{a}\prec\bar{b} \ :\Longleftrightarrow \ (\forall i) \, a_i<b_i, 
\end{equation}
and define
\begin{equation}\label{eq:Rec}
\R_{\ell}(\c) \ := \ \{\bar{a}\in\N^{\ell}\st \bar{a}\prec\bar{c}\},
\end{equation}
omitting $\ell$ when its value is clear from context. Note that $|\R_{\ell}(\c)|=|c_1\cdots c_{\ell}|$. If $(V_i)_{i\in I}$ is a countable base for $\Omega$ and for each $k\in\Z_+$ and function $f:\R_{\ell}(\d{k})\to I$ we write
\begin{equation}\label{eq:BasicOpen}
U_f \ := \ \Big\{x\in\Omega^{\omega^\ell}\st \big(\forall\bar{a}\in\R_{\ell}(\d{k})\big)\, x(\bar{a})\in V_{f(\bar{a})}\Big\},
\end{equation}
then the collection $\big\{U_f\st k\in\Z_+ \text{ and } f\in I^{\R_{\ell}(\d{k})}\big\}$ is a base for the product topology on $\Omega^{\omega^{\ell}}$.  Note that in this topology $\POl$ is $F_\sigma$ and $\FOl$ is $G_\delta$ in $\Ol$.  We show now that with respect to this topology, the shift $\N^{\ell}$-action and diagonal $\N$-action on $\Ol$ are continuous and open, and closed when $\Omega$ is compact.

\begin{claim}\label{claim:NiceAction}
Let $\Omega$ be a Polish space, let $\ell\in\Z_+$ and $\c\in\N^{\ell}$, and let $g_{\bar{c}}:\Ol\to\Ol$ be the surjective function defined by  $g_{\bar{c}}(x)=\bar{c}\cdot x$.
\begin{enumerate}[(i)]
 \item The function $g_{\bar{c}}$ is continuous.
 \item The function $g_{\bar{c}}$ is open.
 \item If $\Omega$ is compact, then the function $g_{\bar{c}}$ is closed. In particular, if $\Omega$ is compact then for every clopen subset $Y\subseteq\Ol$, the set $\bar{c}\cdot Y$ is clopen in $\Ol$.
\end{enumerate}
\end{claim}

\begin{proof}
Fix $\c\in\N^{\ell}$, let $(V_i)_{i\in I}$ be a fixed countable base for the topology on $\Omega$ that includes $\Omega$, and fix $i_0\in I$ such that $V_{i_0}=\Omega$. 

For (i), let $V\subseteq\Ol$ be open, and fix $x\in\Ol$ such that $y:=\bar{c}\cdot x\in V$. Fix $k_y\in\Z_+$ and $f_y:\R(\d{k}_y)\to I$ such that $y\in U_{f_y}\subseteq V$. Let $k_x:=k_y+|\bar{c}|$, and define the function $f_x:\R(\d{k}_x)\to I$ by
\[
f_x(\bar{a}) \ := \ 
\begin{cases}
f_y(\bar{a}-\bar{c}) & \text{if $\bar{a}-\bar{c}\in\R(\d{k}_y)$;} \\
i_0 & \text{otherwise.}
\end{cases}
\]
Suppose $\bar{a}\in\R(\d{k}_x)$. If $\bar{a}-\bar{c}\not\in\R(\d{k}_y)$, then $x(\bar{a})\in\Omega=V_{f_x(\bar{a})}$, and if $\bar{a}-\bar{c}\in\R(\d{k}_y)$, then once again $x(\bar{a})=y(\bar{a}-\bar{c})\in V_{f_y(\bar{a}-\bar{c})}=V_{f_x(\bar{a})}$; thus $x\in U_{f_x}$. Now let $x'\in U_{f_x}$ be arbitrary. Then for all $\bar{a}\in\R(\d{k}_y)$ we have $\bar{c}+\bar{a}\in\R(\d{k}_x)$ and $(\bar{c}\cdot x')(\bar{a})=x'(\bar{c}+\bar{a})\in V_{f_x(\bar{c}+\bar{a})}=V_{f_y(\bar{a})}$, which shows $\bar{c}\cdot x'\in U_{f_y}$. Thus $x\in U_{f_x}\subseteq\bar{c}^{-1}[U_{f_y}]\subseteq\bar{c}^{-1}[V]$, and we conclude that $\bar{c}^{-1}[V]$ is open. This shows $g_{\c}$ is continuous.

For (ii), again let $V\subseteq\Ol$ be open, let $y\in\bar{c}\cdot V$, and fix $x\in V$ such that $\bar{c}\cdot x=y$. 
Fix $k\in\Z_+$ and $f_x:\R(\d{k})\to I$ such that $x\in U_{f_x}\subseteq V$. Define the function $f_y:\R(\d{k})\to I$ by
\[
f_y(\bar{a}) \ = \ 
\begin{cases}
f_x(\bar{c}+\bar{a}) & \text{if $\bar{c}+\bar{a}\in\R(\d{k})$;} \\
i_0 & \text{otherwise.}
\end{cases}
\]
For each $\bar{a}\in\R(\d{k})$, if $\c+\bar{a}\in\R(\d{k})$ then $y(\bar{a})=x(\bar{c}+\bar{a})\in V_{f_x(\bar{c}+\bar{a})}=V_{f_y(\bar{a})}$, and if $\c+\bar{a}\not\in\R(\d{k})$ then $V_{f_y(\bar{a})}=\Omega$ so $y(\bar{a})\in V_{f_y(\bar{a})}$ trivially; this shows $y\in U_{f_y}$. Now let $y'\in U_{f_y}$ be arbitrary, and define $x'\in\Ol$ by
\[
x'(\bar{i}) \ := \ 
\begin{cases}
y'(\bar{i}-\bar{c}) & \text{if $\bar{i}-\bar{c}\in\N^{\ell}$;} \\
x(\bar{i}) & \text{otherwise.}
\end{cases}
\]
Then $\c\cdot x'=y'$, and we claim $x'\in U_{f_x}$, showing $y'\in\c\cdot V$. To see this, let $\bar{i}\in\R(\d{k})$ be arbitrary. If $\bar{i}-\c\in\N^{\ell}$, then also $\bar{i}-\c\in\R(\d{k})$, so $x'(\bar{i})=y'(\bar{i}-\bar{c})\in V_{f_y(\bar{i}-\c)}=V_{f_x(\bar{i})}$. If $\bar{i}-\c\not\in\N^{\ell}$, then again $x'(\bar{i})=x(\bar{i})\in V_{f_x(\bar{i})}$. This shows $x'\in U_{f_x}$. Since $y'\in U_{f_y}$ was arbitrary, we have shown $U_{f_y}\subseteq\c\cdot V$, and we conclude that $\c\cdot V$ is open as desired.

Finally, for (iii), assume $\Omega$ is compact and $F\subseteq\Ol$ is closed. Then $\Ol$ is compact, so $F$ is compact, which by continuity of $g_{\bar{c}}$ implies $\c\cdot F$ is compact, hence closed.
\end{proof}

\begin{corollary}\label{cor:saturation} 
Let $\Ol$ be a product space, and let $Y\subseteq\Ol$ and $\c\in\N^{\ell}$. 
\begin{enumerate}[(i)]
 \item If $Y$ is open, then the $F^{\Omega}_{\ell,\c}$-saturation of $Y$ is also open.
 \item If $\Omega$ is compact and $Y$ is closed, then the $F^{\Omega}_{\ell,\c}$-saturation of $Y$ is also closed.
 \item If $\Omega$ is compact and $Y$ is clopen, then the $F^{\Omega}_{\ell,\c}$-saturation of $Y$ is also clopen.
\end{enumerate}
\end{corollary}

\begin{proof}
Since $[Y]_{F^{\Omega}_{\ell,\c}}=\c^{-1}[\c\cdot Y]$, (i) and (ii) follow from Claim \ref{claim:NiceAction}, and (iii) follows from (i) and (ii).
\end{proof}

For us, the base space $\Omega$ will always be either a finite discrete space which we take to be an integer $J\geq 2$, a countably infinite discrete space which we take to be $\omega$, or else Cantor space $2^\omega$ in the usual topology. 
In the first case, which will be most common, let 
\[
\S \ = \ \S_{\ell}^{J} \ := \ \bigcup_{k \mkern1mu \in \mkern1mu \Z_+}J^{\R_{\ell}(\d{k})},
\]
and for each $s\in \S$ let $|s|$ denote the unique integer $k$ such that $s\in J^{\R(\d{k})}$. For $k\in\Z_+$ and $s\in \S$ or $x\in\Jl$, we  write $s\res k$ and $x\res k$ as abbreviations for $s\res\R(\d{k})$ and $x\res\R(\d{k})$, respectively.  For each $s\in \S$, let $[s]\subseteq \Jl$ be the set consisting of all elements extending $s$, i.e.,
\[
[s] \ := \ \big\{x\in\Jl \st x\res |s|=s\big\}.
\]
The collection $\mathcal{B}_{\ell}^J:=\big\{[s]\st s\in \S_{\ell}^J\big\}$ is a countable clopen base for the product topology on $\Jl$, which is compact, zero-dimensional, and Polish with compatible metric $\rho_{\ell}$ given by
\begin{equation}\label{eq:rho}
\rho_{\ell}(x,y) \ := \ 2^{-k} \ \text{where $k$ is least such that $x\res k \ne y\res k$}.
\end{equation}
With this topology, $\Jl$ is just a homeomorphic copy of Cantor space.  Similar remarks apply if $\Omega=2^\omega$, since $(2^\omega)^{\omega^\ell}$ is homeomorphic to $2^{\omega^{\ell+1}}$.

Finally, for each $\ell\in\Z_+$, let $\mu_{\ell}$ be the $(\frac{1}{2},\frac{1}{2})$-probability measure on $\Tl$, satisfying
\[
\mu_{\ell}\Big(\big\{x\in\Tl\st x(\bar{i})=0\big\}\Big) \ = \ \mu_{\ell}\Big(\big\{x\in\Tl\st x(\bar{i})=1\big\}\Big) \ = \ \frac{1}{2}
\]
for all $\bar{i}\in\N^{\ell}$. Any comments we make about measures on $\Tl$ will always refer to these $\mu_{\ell}$.

\subsection{Geometry of $\EOl$-classes}\label{subsec:geometry}
Let $\Ol$ be a product space. We would like to have a standard way of picturing elements of $\FOl$ and the geometric structure of their $\EOl$-classes, and we describe this now with an emphasis on the $2$-dimensional case, both because that will be our focus in Sections \ref{sec:periodicity}--\ref{sec:Gm} and because the general picture can to some extent be understood by analogy with this case.

We view elements of $\OT$ as two-dimensional arrays of entries in the base space, arranged as in Figure \ref{fig:2d}. 
When $\Omega=2$, this becomes an infinite two-dimensional array of zeros and ones. Given $x\in\OT$ and $\bar{c}\in\N^2$, we can picture $\bar{c}\cdot x$ as the subarray obtained by advancing forward through this array by the vector $\bar{c}$.  Two elements $x$ and $y$ are $\EOT$-related if there exist $\bar{a},\bar{b}\in\N^2$ such that $\bar{a}\cdot x=\bar{b}\cdot y$; thus we can picture $\EOT$-equivalence as in Figure \ref{fig:equivalence}, where the shaded region represents the common outputs of $\bar{a}\cdot x$ and $\bar{b}\cdot y$.

\begin{figure}[h]
\[
\begin{array}{|cccc}
\vdots & & & \reflectbox{$\ddots$} \\
x(0,2) & x(1,2) & x(2,2) & \\
x(0,1) & x(1,1) & x(2,1) & \\
x(0,0) & x(1,0) & x (2,0) & \cdots \\
\hline
\end{array}
\]
\caption{An element $x\in \OT$.}
\label{fig:2d}
\end{figure}

\begin{figure}[h]
\begin{center}
\begin{tikzpicture}
\draw [thick] (0,1) -- (10,1);
\draw [thick] (0,1) -- (0,5);
\draw [thick] (1.5,0) -- (10,0);
\draw [thick] (1.5,0) -- (1.5,5);
\node at (-.15,.85) {$x$};
\node at (1.35,-.15) {$y$};
\draw [dashed] (0,1) -- (3,2);
\draw [dashed] (1.5,0) -- (3,2);
\fill [fill=gray!20] (3,2) rectangle (10,5); 
\node at (3,2) [circle,fill,inner sep=1.2pt]{};
\node at (0,1) [circle,fill,inner sep=1.2pt]{};
\node at (1.5,0) [circle,fill,inner sep=1.2pt]{};
\node at (3.9,1.75) {$\bar{a}\cdot x=\bar{b}\cdot y$};
\node at (.8,1.5) {$\bar{a}$};
\node at (2.2,.55) {$\bar{b}$};
\draw [dotted,thick] (0,1) -- (1.5,0);
\end{tikzpicture}
\end{center}
\caption{$\EOT$-equivalence.}
\label{fig:equivalence}
\end{figure}

Note that if $\bar{a}\cdot x=\bar{b}\cdot y$, then it is natural to position $x$ and $y$  as differing by the vector $\bar{a}-\bar{b}$, represented in Figure \ref{fig:equivalence} by the dotted segment. For this to be justified, one must check that if $\bar{a}\cdot x=\bar{b}\cdot y$ and $\bar{c}\cdot x=\bar{d}\cdot y$, then $\bar{a}-\bar{b}=\bar{c}-\bar{d}$.  This fails in $\POl$, but holds easily in $\FOl$.

\begin{claim}\label{claim:geometry}
Let $\Ol$ be a product space, and let $x,y\in\textup{Fr}(\Ol)$ and $\bar{a},\bar{b},\bar{c},\bar{d}\in\N^{\ell}$. If $\bar{a}\cdot x=\bar{b}\cdot y$ and $\bar{c}\cdot x=\bar{d}\cdot y$, then $\bar{a}-\bar{b}=\bar{c}-\bar{d}$.
\end{claim}

\begin{proof} Fix $x,y\in\textup{Fr}(\Ol)$, and assume $\bar{a}\cdot x=\bar{b}\cdot y$ and $\bar{c}\cdot x=\bar{d}\cdot y$. Then 
\[
(\bar{a}+\bar{d})\cdot x \ = \ (\bar{b}+\bar{d})\cdot y \ = \ (\bar{b}+\bar{c})\cdot x,
\]
so $\bar{a}+\bar{d}=\bar{b}+\bar{c}$ since $x\in\textup{Fr}(\Ol)$. \end{proof}

Claim \ref{claim:geometry} justifies the following definition: given $x,y\in\textup{Fr}(\Ol)$ such that $x\mathrel{\EOl}y$, write $x-y$ for the unique element $\bar{m}\in\Z^{\ell}$ such that $\bar{m}=\bar{b}-\bar{a}$ for all $\bar{a},\bar{b}\in\N^\ell$ such that $\bar{a}\cdot x=\bar{b}\cdot y$. We then have $x-y=\bar{0}$ provided there is $\bar{a}\in\N^{\ell}$ such that $\bar{a}\cdot x=\bar{a}\cdot y$, in which case
\[
\d{n}\cdot x \ = \ (\d{n}-\bar{a})\cdot (\bar{a}\cdot x) \ = \ (\d{n}-\bar{a})\cdot (\bar{a}\cdot y) \ = \ \d{n}\cdot y
\]
where $n=|\bar{a}|$. Thus $x-y=\bar{0}$ if and only if $x\mathrel{F_{\ell}^\Omega}y$. Furthermore, given $\EOl$-equivalent elements $x,y,x',y'\in\FOl$, if $x\mathrel{F^{\Omega}_{\ell}}x'$ and $y\mathrel{F^{\Omega}_{\ell}}y'$ then $x-y=x'-y'$. It follows that the collection of $F_{\ell}^{\Omega}$-classes within a given $\EOl$-class has the geometric structure of $\Z^{\ell}$, although we cannot in general pick an ``origin'' in a Borel fashion.

To make the observations above a bit more precise, fix for each base space $\Omega$ some distinguished element $x_{\Omega}\in\Omega$, and for each $x\in\FOl$ and $\m\in\Z^{\ell}$ define $\m\star x\in\Ol$ by
\begin{equation}\label{eq:star}
(\m\star x)(\bar{i}) \ = \ \begin{cases} x(\m+\bar{i}) & \text{if $\m+\bar{i}\in\N^{\ell}$;} \\ x_{\Omega} & \text{otherwise.} \end{cases}
\end{equation}
Then $\star$ induces a group action of $\Z^{\ell}$ on the quotient $\FOl/F^{\Omega}_{\ell}$, defined by $\m\cdot [x]_{F^{\Omega}_{\ell}}=[\m\star x]_{F^{\Omega}_{\ell}}$. This action is free, and for each $x\in\FOl$, the orbit of $[x]_{F^{\Omega}_{\ell}}$ under it is $\big\{[y]_{F^{\Omega}_{\ell}}\st x\mathrel{\EOl}y\big\}$; thus each $\EOl/F^{\Omega}_{\ell}$-class in $\FOl/F^{\Omega}_{\ell}$ looks like $\Z^{\ell}$. This geometric perspective will be quite useful, but we must be careful not to be misled by it, since the equivalence relation $F_{\ell}^{\Omega}$ is not smooth, $\FOl/F^{\Omega}_{\ell}$ is not a standard Borel space, and dividing out by $F_{\ell}^{\Omega}$ erases a great deal of structure, as we see below in Figure \ref{fig:EtClass} and the discussion around it. Note that $\star$ is not a $\Z^{\ell}$-action on $\FOl$, since for instance $-\m\star(\m\star x)$ need not equal $x$.

Our proof that $\ETT\res\FTT$ is hypersmooth will be based in part on the notion of a \emph{cascade} (see Section \ref{sec:cascade}), which will require that we have a metric providing us with a sense of how far apart two elements within an $\EOl$-class are from each other.  The $\Z^{\ell}$-geometry of $\EOl/F_{\ell}^{\Omega}$ suggests one potential notion of distance within an $\EOl$-class: namely, we could take the distance between $x$ and $y$ to be $|x-y|$.  But this does not take into account the magnitudes of minimal vectors $\bar{a}$ and $\bar{b}$ satisfying $\bar{a}\cdot x=\bar{b}\cdot y$, so it will not be of much use.  The metrics $\rho_{\ell}$ defined in Equation \eqref{eq:rho} are even less helpful, since two elements of $\FOl$ can be very close with respect to $\rho_{\ell}$ without even belonging to the same $\EOl$-class.  Instead, the primary notion of distance within an $\EOl$-class that we shall use is given in Definition \ref{df:d}.

\begin{definition}\label{df:d} 
Let $\Ol$ be a product space. For each $x,y\in\FOl$, define
\[
d_{\ell}(x,y) \ := \ \begin{cases} \min\Big\{\max\{|\bar{a}|,|\bar{b}|\}\st \bar{a},\bar{b}\in\N^\ell \,\wedge\,\bar{a}\cdot x=\bar{b}\cdot y\,\Big\} & \mbox{if $x\mathrel{\EOl}y$;} \\ \infty & \mbox{otherwise.} \end{cases}
\]
Also put $d_{\ell}(x,Y)=\min\{d_{\ell}(x,y)\st y\in Y\}$ for $x\in\Ol$ and $Y\subseteq\Ol$.
\end{definition}

The definition of $d_{\ell}$ does not depend on $\Omega$, so we omit it from the notation, and when $\ell$ is clear from context we will omit it as well. An \emph{extended metric} is a metric for which the value $+\infty$ is allowed. To any extended metric space $(X,d)$ we associate the finite-distance equivalence relation $E_d=\{(x,y)\in X^2\st d(x,y)<\infty\}$.

\begin{claim}\label{claim:metric} 
For every product space $\Ol$, the function $d_{\ell}$ is a Borel extended metric on $\FOl$ such that $E_{d_{\ell}}=\EOl\res\FOl$.
\end{claim}

\begin{proof} Fix a product space $\Ol$, and write $d=d_{\ell}$. To see that $d$ is Borel, note that $d(x,y)\leq n$ if and only if there exist $\bar{a},\bar{b}\in\N^{\ell}_{\leq n}$ such that $\bar{a}\cdot x=\bar{b}\cdot y$, where for each fixed $\bar{a},\bar{b}$ the set $$\big\{(x,y)\in \Ol\times\Ol\st \bar{a}\cdot x=\bar{b}\cdot y\big\}$$ is closed in $\Ol\times\Ol$ as a consequence of Claim \ref{claim:NiceAction}.

To see that the restriction of $d$ to each $\EOl$-class in $\FOl$ is a metric, we verify the triangle inequality, as the other axioms are trivial to check. Let $x,y,z\in\FOl$ and suppose
\[
 \bar{a}\cdot x \ = \ \bar{b}\cdot y, \quad \bar{c}\cdot x \ = \ \bar{d}\cdot z, \quad \mbox{and} \quad \bar{u}\cdot y \ = \ \bar{v}\cdot z,
\]
where
\[
 d(x,y) \ = \ \max\{|\bar{a}|,|\bar{b}|\}, \quad d(x,z) \ = \ \max\{|\bar{c}|,|\bar{d}|\}, \quad \mbox{and} \quad d(y,z) \ = \ \max\{|\bar{u}|,|\bar{v}|\}.
\]
Then
\[
 (\bar{c}+\bar{v})\cdot x \ = \ (\bar{d}+\bar{v})\cdot z \ = \ (\bar{d}+\bar{u})\cdot y,
\]
so
\[
 d(x,y) \ = \ \max\{|\bar{a}|,|\bar{b}|\} \ \leq \ \max\{|\bar{c}+\bar{v}|,|\bar{d}+\bar{u}|\}.
\]
Therefore, if $\max\{|\bar{c}+\bar{v}|,|\bar{d}+\bar{u}|\}=|\bar{c}+\bar{v}|$, then 
\[
 d(x,y) \ \leq \ |\bar{c}+\bar{v}| \ \leq \ |\bar{c}|+|\bar{v}| \ \leq \ \max\{|\bar{c}|,|\bar{d}|\}+\max\{|\bar{u}|,|\bar{v}|\} \ = \ d(x,z)+d(y,z),
\]
and if $\max\{|\bar{c}+\bar{v}|,|\bar{d}+\bar{u}|\}=|\bar{d}+\bar{u}|$, then again
\[
 d(x,y) \ \leq \ |\bar{d}+\bar{u}| \ \leq \ |\bar{d}|+|\bar{u}| \ \leq \ \max\{|\bar{c}|,|\bar{d}|\}+\max\{|\bar{u}|,|\bar{v}|\} \ = \ d(x,z)+d(y,z).
\] 
Finally, the fact that $E_{d_{\ell}}=\EOl\res\FOl$ is immediate from the definition of $d_{\ell}$.
\end{proof}

Throughout our arguments, all references to the distance between two $\EOl$-related points $x,y\in\FOl$ will be with respect to $d_{\ell}$. We will never use $\rho_{\ell}$ when considering the geometry of $\EOl$-classes, even though $\rho_{\ell}$ generates the topology of $\Ol$ while $d_{\ell}$ does not. The metric $d_{\ell}$ is more closely related to the graphs we define in Section \ref{sec:graphs}, although it is not exactly the graph metric arising from any of these graphs, either. %
In any metric space $(X,d)$, we write
\[
B_r(x_0) \ := \ \{x\in X\st d(x,x_0)<r\} \qquad \text{and} \qquad B_r[x_0] \ := \ \{x\in X\st d(x,x_0)\leq r\}
\]
for the $d$-open and $d$-closed balls of radius $r$ about $x_0$ in $X$; similarly, given $Y\subseteq X$, we write $B_r(Y):=\{x\in X\st d(x,Y)<r\}$ and $B_r[Y]:=\{x\in X\st d(x,Y)\leq r\}$, where $d(x,Y):=\inf\{d(x,y)\st y\in Y\}$. 
Since the metrics $d_{\ell}$ of Definition \ref{df:d} take values in $\N$, we have $d_{\ell}(x,Y)=\min\{d_{\ell}(x,y)\st y\in Y\}$, and every $d_{\ell}$-open ball in $\FOl$ is $d_{\ell}$-closed and vice versa, with $B_r[x]=B_{r+1}(x)$. 
Any mention of balls in the free part of a product space $\Ol$ will refer to the metric $d_{\ell}$ of Definition \ref{df:d}, so we must be careful with terminology: since $d$ does not generate the topology of $\Ol$, ``$d$-open" balls in $\Ol$ need not be topologically open, and ``$d$-closed'' balls need not be topologically closed. We have the following:

\begin{claim}\label{claim:dSat} 
Let $\Ol$ be a product space, and let $Y\subseteq\FOl$ and $r\in\omega$. 
\begin{enumerate}[(i)]
 \item If $Y$ is open, then $B_r(Y)$ is also open.
 \item If $\Omega$ is compact and $Y$ is closed, then $B_{r}(Y)$ is also closed.
\end{enumerate}
\end{claim}

\begin{proof} This follows from Claim \ref{claim:NiceAction}, since $d(x,Y)<r$ if and only if there exist $\bar{a},\bar{b}\in\N^\ell_{<r}$ such that $\bar{a}\cdot x\in\bar{b}\cdot Y$.
\end{proof}

In order to picture the metric $d$ and balls $B_{r}[x]$, we first consider the 1-dimensional case. When $\Omega=2$ and $\ell=1$, $E^2_1$ is just the tail equivalence relation $E_{0,t}$ on Cantor space, and an $E^2_1$-class in $\textup{Fr}(2^\omega)$ looks like an infinite binary tree with no root as in Figure \ref{fig:EtClass}, where we display a local portion of the directed graph $\vec{\G}_1$ on $\textup{Fr}(2^\omega)$ defined by $x\mathrel{\vec{\G}_1}y \ :\Leftrightarrow \ 1\cdot x=y$. Each node corresponds to some point $x$ in $\textup{Fr}(2^\omega)$, and is connected to its lone $\vec{\G}_1$-child $1\cdot x$ on the right and its two $\vec{\G}_1$-parents $0^\frown x$ and $1^\frown x$ on the left. We have labeled points $x_0$ and $y_0$ in $\textup{Fr}(2^\omega)$, and shaded their $F^2_{1,3}$-classes. Note that $3\cdot x_0=5\cdot y_0$ and there exist no $a<3$ and $b<5$ for which $a\cdot x_0=b\cdot y_0$, so $d_1(x_0,y_0)=5$. We have also shaded the closed $d_1$-ball of radius $2$ about $x_0$. Note that $B_{2}[x_0]=B_{3}(x_0)$.

\begin{figure}[h]
\begin{center}
\begin{tikzpicture}
 \fill [fill=gray!28] (-1.5,-4) rectangle (-0.5,4);
 \fill [fill=gray!28] (-5.5,-4.2) rectangle (-4.5,-2); 
 \foreach \x in {-3.75,-3.5,...,3.75,4} \node at (-5,\x-.125) [circle,fill,inner sep=1pt]{};
 \foreach \x in {-3.75,-3.25,...,3.25,3.75} \node at (-3,\x) [circle,fill,inner sep=1.1pt]{};
 \foreach \x in {-3.5,-2.5,...,2.5,3.5} \node at (-1,\x) [circle,fill,inner sep=1.2pt]{};
 \foreach \x in {-3,-1,1,3} \node at (1,\x) [circle,fill,inner sep=1.3pt]{};
 \foreach \x in {-2,2} \node at (3,\x) [circle,fill,inner sep=1.4pt]{};
 \node at (5,0) [circle,fill,inner sep=1.5pt]{};
 \foreach \x in {-3.75,-3.25,...,3.25,3.75} \draw (-5,\x-.125) -- (-3,\x);
 \foreach \x in {-3.75,-3.25,...,3.25,3.75} \draw (-5,\x+.125) -- (-3,\x);
 \foreach \x in {-3.75,-2.75,...,2.25,3.25} \draw (-3,\x) -- (-1,\x+.25);
 \foreach \x in {-3.75,-2.75,...,2.25,3.25} \draw (-3,\x+.5) -- (-1,\x+.25);
 \foreach \x in {-3,-1,...,1,3} \draw (-1,\x+.5) -- (1,\x);
 \foreach \x in {-3,-1,...,1,3} \draw (-1,\x-.5) -- (1,\x);
 \foreach \x in {-2,2} \draw (1,\x+1) -- (3,\x);
 \foreach \x in {-2,2} \draw (1,\x-1) -- (3,\x);
 \draw (3,2) -- (5,0);
 \draw (3,-2) -- (5,0);
 \draw [dotted,thick] (5,0) -- (6.5,0);
 \node at (-1,0.75) {$x_0$};
 \node at (-5.18,-4) {$y_0$};
 \node at (1,2) {$B_2[x_0]$};
 \node at (-5.8,0) {\ldots};
\draw plot [smooth cycle,tension=0.5] coordinates {(-5.1,0) (-5.1,0.94) (-4.2,1.1) (-3.8,1.9) (-2,2.1) (-1.1,3.75) (1,3.4) (3.28,2.05) (1.8,0.9) (-2,0.1)};%
\fill [gray!35, opacity=0.3] plot [smooth cycle,tension=0.5] coordinates {(-5.1,0) (-5.1,0.94) (-4.2,1.1) (-3.8,1.9) (-2,2.1) (-1.1,3.75) (1,3.4) (3.28,2.05) (1.8,0.9) (-2,0.1)};
 \filldraw[fill=white] (-1,0.5) circle (1.75pt);
 \filldraw[fill=white] (-5,-3.875) circle (1.35pt);
\end{tikzpicture}
\end{center}
\caption{Local geometry within an $E^2_1$-class.}
\label{fig:EtClass}
\end{figure}
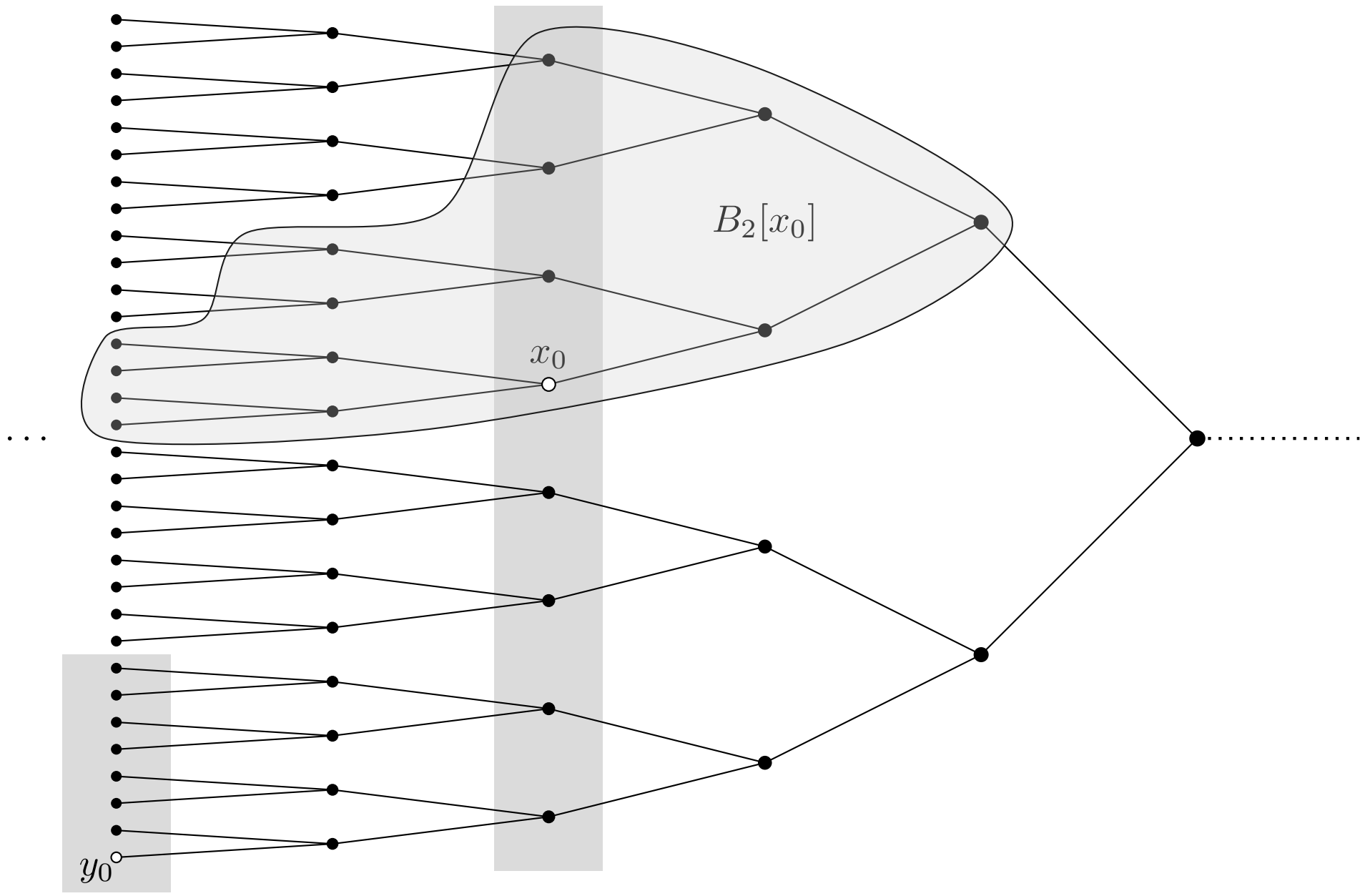

One may picture an entire $E^2_1$-class by imagining Figure \ref{fig:EtClass} continuing indefinitely in both directions; then the $F^2_1$-classes are (infinite) vertical strips, which are ordered within an $E^2_1$-class like $\Z$. From that perspective, $B_{2}[x_0]$ is contained in the closed interval of radius 2 centered at $[x_0]_{F^2_1}$ in the $\Z$-ordering of $F^2_1$-classes within $[x_0]_{E^2_1}$, as shown in Figure \ref{fig:1dimProj}.

\begin{figure}[h]
\begin{center}
\begin{tikzpicture}
 \fill [fill=gray!25] (-5,-0.15) rectangle (3,0.15);
 \foreach \x in {-5,-3,...,5} \node at (-\x,0) [circle,fill,inner sep=1.5pt]{};
 \draw (-6,0) -- (6,0);
 \draw [dotted] (-7,0) -- (-6,0);
 \draw [dotted] (6,0) -- (7,0);
 \node at (-0.8,-0.6) {$[x_0]_{F^2_1}$};
 \node at (-4.8,-0.6) {$[y_0]_{F^2_1}$};
 \node at (-5,0) {$\Big[$};
 \node at (3,0) {$\Big]$};
\end{tikzpicture}
\end{center}
\caption{The projection of $B_{2}[x_0]$ in $\textup{Fr}(2^\omega)/F^2_1$.}
\label{fig:1dimProj}
\end{figure}
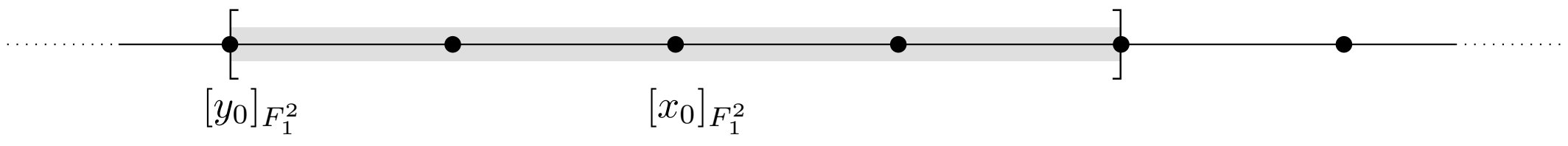

\begin{figure}[h]
\begin{center}
\begin{tikzpicture}
 \fill [fill=gray!25] (-1.6,-1.6) rectangle (1.6,1.6);
 \foreach \x in {-2,-1.6,...,2} \foreach \y in {-2,-1.6,...,2} \node at (\x,\y) [circle,fill,inner sep=0.6pt]{};
 \draw [dashed] (0,0) -- (1.6,1.2);
 \draw [dashed] (0.8,-0.4) -- (1.6,1.2);
 \node at (-0.1,-0.18) {$x_0$};
 \node at (0.7,-0.58) {$y_0$};
 \node at (1.22,0.18) {$\c$};
  \filldraw[fill=white] (0,0) circle (1.4pt);
  \filldraw[fill=white] (0.8,-0.4) circle (1.4pt);
  \draw (-1.6,-1.6) -- (-1.6,1.6);
  \draw (-1.6,1.6) -- (1.6,1.6);
  \draw (-1.6,-1.6) -- (1.6,-1.6);
  \draw (1.6,-1.6) -- (1.6,1.6);
\end{tikzpicture}
\end{center}
\caption{The projection of $B_{4}[x_0]$ in $\textup{Fr}(2^{\omega^2})/F^2_2$.}
\label{fig:2dimProj}
\end{figure}
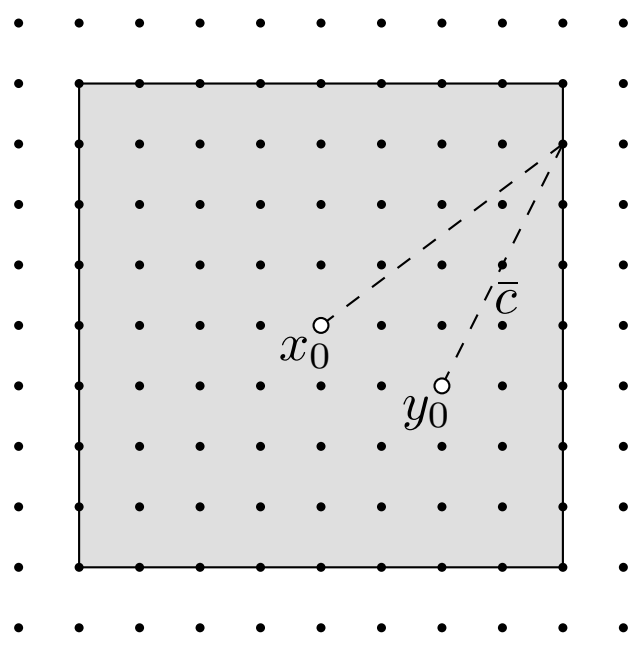
Similarly, in the 2-dimensional case, the closed $d_2$-ball of radius $r$ about a point $x_0\in \textup{Fr}(2^{\omega^2})$ is contained in the ``square'' of radius $r$ about $x_0$ in the quotient $\textup{Fr}(2^{\omega^2})/F^2_2$, which looks like $\Z^2$ (see Figure \ref{fig:2dimProj}). 
However, Figure \ref{fig:2dimProj} collapses a great deal of structure, in the same way that Figure \ref{fig:1dimProj} collapses the structure appearing in Figure \ref{fig:EtClass}. In order to picture $B_{4}[x_0]$ using Figure \ref{fig:2dimProj}, one must imagine that each node within the shaded square represents many points in $\textup{Fr}(2^{\omega^2})$; specifically, a given node $y_0$ within the shaded square represents an entire $F^2_{2,\c}$-class, where $\c\in\N^2_{\leq r}$ is maximal such that $\c\cdot y_0$ still belongs to the shaded square. Note that whenever the dimension $\ell$ is at least $2$, and even if the base space is finite, every nontrivial $F^{\Omega}_{\ell,\c}$-class is uncountable, so in particular every nontrivial $d_{\ell}$-ball is uncountable. Nevertheless, for any base space $\Omega$ and dimension $\ell$, the $d_{\ell}$-ball $B_{r}[x]$ is covered by the finitely many $F_{\ell,\d{r}}^{\Omega}$-classes $[\m\star x]_{F_{\ell,\d{r}}^{\Omega}}$ where $\m\in\Z^\ell_{\leq r}$. 

We remark that for any base space $\Omega$ and Borel set $Y\subseteq\Ol$, neither $\c\cdot Y$ nor $B_{r}[Y]$ need be Borel. Furthermore, when $\Omega$ is not compact, the maps $g_{\c}$ of Claim \ref{claim:NiceAction} need not be closed, so neither $\c\cdot Y$ nor $B_r[Y]$ need be Borel even when $Y$ is closed. However, see Lemma \ref{lem:ForwardBorel} for some conditions on the Borel set $Y\subseteq\Ol$ that guarantee $\bar{c}\cdot Y$ and $B_r[Y]$ are Borel.

\section{Hypersmooth equivalence relations}\label{sec:hs}

We continue here the discussion begun in Section \ref{sec:prelimBERs}, with an emphasis on some of the closure properties of subclasses of hypersmooth equivalence relations. Recall that the classes of smooth and hypersmooth equivalence relations are closed under products, Borel reducibility, and restrictions to Borel subsets, while the class of hyperselective equivalence relations is closed under restrictions to \emph{invariant} Borel subsets, but not under Borel reductions or restrictions to arbitrary Borel subsets. We record the facts we shall need formally below; all are well-known, or are obvious from the definitions. For more context and details, see \cite{DJK}, \cite{JKL}, \cite{KL}, and \cite{Gao}.

\begin{fact}
\label{fact:BasicClosure}
Let $E$ and $F$ be Borel equivalence relations on the standard Borel spaces $X$ and $Y$, respectively.
\begin{enumerate}[(i)]
 \item If $E$ is finite, then $E$ is smooth.
 \item If $F$ is smooth and $E\leq_BF$, then $E$ is smooth.
 \item If $E$ is smooth and $B\subseteq X$ is Borel, then $E\res B$ is smooth.
 \item If $F$ is hypersmooth and $E\leq_BF$, then $E$ is hypersmooth.
 \item If $E$ is hypersmooth and $B\subseteq X$ is Borel, then $E\res B$ is hypersmooth.
 \item If $E$ is hyperselective and $B\subseteq X$ is an $E$-invariant Borel subset of $X$, then $E\res B$ is hyperselective.
 \item If $E$ and $F$ are hypersmooth, then $E\times F$ is hypersmooth. \hfill $\square$
\end{enumerate}
\end{fact}

The classes of hypersmooth and hyperselective equivalence relations are also closed under countable disjoint unions, and we separately state and prove a slight variant of this fact which will be used repeatedly in what follows.

\begin{lemma}\label{lem:piecewiseHS}
Let $E$ be a Borel equivalence relation on the standard Borel space $X$, and let $(X_n)_{n\in\omega}$ be a sequence of $E$-invariant Borel subsets of $X$ such that $\cup_nX_n=X$. 
\begin{enumerate}[(i)]
 \item If $E\res X_n$ is hypersmooth for each $n$, then $E$ is hypersmooth.
 \item If $E\res X_n$ is hyperselective for each $n$, then $E$ is hyperselective.
\end{enumerate}
\end{lemma}

\begin{proof}
By replacing each $X_n$ with $X_n\setminus \cup_{k<n}X_k$, we may assume the sets $X_n$ are pairwise disjoint. 

To prove (i), fix for each $n\in\omega$ a Borel reduction $f_n:X_n\to (2^\omega)^\omega$ from $E\res X_n$ to $E_1$, and define the Borel function $g_n:X_n\to (2^\omega)^\omega$ by 
\[
g_n(x)(m) \ := \ \begin{cases} f_n(x)(m/2) & \text{if $m$ is even;} \\ {1^n}{}^\frown 0^\infty & \text{if $m$ is odd.} \end{cases}
\]
It is easy to check that $g:=\sqcup_ng_n$ is a Borel reduction from $E$ to $E_1$. 

For (ii), suppose that for each $n$,  $\big(F_k^{(n)}\big)_{k\in\omega}$ is an increasing sequence of selective equivalence relations on $X_n$ such that $\cup_kF^{(n)}_k=E\res X_n$.  For each $k$, let $F_k:=\bigcup_{n\in\omega}F_k^{(n)}$. Then $(F_k)_{k\in\omega}$ is an increasing sequence of selective equivalence relations on $X$ with $\cup_kF_k=E$.
\end{proof}

As the class of hyperfinite equivalence relations is precisely the class of countable hypersmooth equivalence relations, it too is closed, within the universe of countable Borel equivalence relations, under Borel reducibility, restrictions to Borel subsets, products, and countable disjoint unions. But hypefinite equivalence relations enjoy another closure property that is not known to extend to the hypersmooth setting, and which will occupy us for the rest of this section. Given equivalence relations $E$ and $F$ on $X$, we call $E$ a \emph{sub-equivalence relation} of $F$, and $F$ a \emph{super-equivalence relation} of $E$, if $E\subseteq F$, that is, if $x\mathrel{E}y$ implies $x\mathrel{F}y$ for all $x,y\in X$. If $E$ is a sub-equivalence relation of $F$, we say that $E$ has \emph{finite index} in $F$, and write $[F:E]<\infty$, if every $F$-class contains finitely many $E$-classes. It was shown in \cite{JKL} that if the countable Borel equivalence relation $F$ has a hyperfinite sub-equivalence relation of finite index, then $F$ itself is hyperfinite. However, we do not know whether this generalizes to hypersmooth or hyperselective equivalence relations.

\begin{question}\label{ques:HSFI}
Let $E\subseteq F$ be Borel equivalence relations on the standard Borel space $X$ such that $E$ is hypersmooth and $[F:E]<\infty$. Is $F$ also hypersmooth?
\end{question}

\begin{question}\label{ques:HSelFI}
Let $E\subseteq F$ be Borel equivalence relations on the standard Borel space $X$ such that $E$ is hyperselective and $[F:E]<\infty$. Is $F$ also hyperselective?
\end{question}

We will return briefly to Question \ref{ques:HSelFI} in its full generality at the end of this section, but our main goal will be to provide positive answers to these questions in certain special cases. To state and prove our results, we need a couple more definitions from \cite{KL}. Recall that the subset $B$ of the Polish space $X$ is $K_{\sigma}$ if $B$ can be expressed as the union of a countable family of compact subsets of $X$.

\begin{definition}[\cite{KL}]
\begin{enumerate}[(i)]
 \item The Borel equivalence relation $E$ on the Borel set $B$ in the Polish space $Y$ \emph{has $K_\sigma$ classes} if every equivalence class of $E$ is $K_{\sigma}$ in $Y$.
 \item The Borel equivalence relation $E$ is \emph{uniformly continuous} if $E\cong_BE\times I_{2^\omega}$.
\end{enumerate}
\end{definition}

By \cite[3.4]{KL}, uniform continuity is equivalent to containing a smooth Borel sub-equivalence relation with uniformly continuum-size equivalence classes; that is, $E$ on $X$ is uniformly continuous if and only if there is a smooth sub-equivalence relation $F\subseteq E$ and a Borel function $f:X\times 2^\omega\to X$ such that for all $(x,y)\in F$ and $\alpha\ne\beta$ in $2^\omega$, we have $f(x,\alpha)\mathrel{F}x$, $f(x,\alpha)=f(y,\alpha)$, and $f(x,\alpha)\ne f(x,\beta)$. Note that $E_1$ is uniformly continuous with witness $x\mathrel{F}y \ \Leftrightarrow \ 1\cdot x=1\cdot y$, so $E_{1,t}$ is also uniformly continuous since $E_{1,t}\cong_BE_1$. Note also that every Borel super-equivalence relation of a uniformly continuous equivalence relation is itself uniformly continuous, with the same witness $F$. Let us make a couple more elementary observations.

\begin{observation}\label{ob:UniformlyContinuous}
Let $\Ol$ be a product space, and $X$ an $\EOl$-invariant Borel subset of $\FOl$.
\begin{enumerate}
 \item[(i)] If $\Omega$ is compact, then $\EOl\res X$ has $K_{\sigma}$ classes.
 \item[(ii)] If $\ell\geq 2$ or $\Omega$ is uncountable, then $\EOl\res X$ is uniformly continuous.
\end{enumerate}
\end{observation}

\begin{proof}
Assume the hypotheses, and write $E:=\EOl\res X$. For (i), assume $\Omega$ is compact, so $\Ol$ is also compact, and let $C=[x_0]_E$ be an $E$-class. For each $\bar{a},\bar{b}\in\N^{\ell}$, the set $\bar{b}^{-1}[\{\bar{a}\cdot x_0\}]$ is closed in $\Ol$ by Claim \ref{claim:NiceAction}, and is therefore compact, so
\[
C \ = \ \bigcup_{\bar{a},\bar{b} \mkern1mu \in \mkern1mu \N^{\ell}}\bar{b}^{-1}\big[\{\bar{a}\cdot x_0\}\big]
\]
is $K_{\sigma}$ in $\Ol$. For (ii), if either $\ell\geq 2$ or $\ell=1$ but $\Omega$ is uncountable, then $F_{\ell,\d{1}}^{\Omega}\res X$ is a smooth sub-equivalence relation of $E$ that has uniformly continuum-sized classes. To see this when $\ell\geq 2$, fix distinct elements $x_0,x_1\in\Omega$ and take $f:X\times 2^\omega\to X$ to be given by 
\[
f(x,\alpha)(\bar{i}) \ := \ 
\begin{cases}
 x_{\alpha(i_1)} & \text{if $i_j=0$ for all $1<j\leq\ell$;} \\
 x(\bar{i}) & \text{if $\d{1}\preceq\bar{i}$;} \\
 x_0 & \text{otherwise.}
\end{cases}
\]
If $\ell=1$ but $\Omega=2^\omega$, we can put $f(x,\alpha)(0):=\alpha$ and $f(x,\alpha)(n):=x(n)$ for all $n\geq 1$.
\end{proof}

\begin{observation}\label{ob:UniformlyContinuousAndKSigma}
Every uniformly continuous hypersmooth equivalence relation with $K_{\sigma}$ classes is hyperselective.
\end{observation}

\begin{proof}
Let $E$ be a uniformly continuous hypersmooth equivalence relation on the Borel set $B$ in the Polish space $Y$ such that every $E$-class is $K_{\sigma}$ in $Y$. By \cite[Lemma 3.2]{KL}, we have $E\times I_{2^\omega}\sqsubseteq_B^i E_1$. Since $E$ is uniformly continuous, this means $E\sqsubseteq_B^i E_1$. Thus $E$ is Borel isomorphic to the restriction of $E_1$ to an $E_1$-invariant Borel subset of $(2^\omega)^\omega$, so $E$ is hyperselective by Fact \ref{fact:BasicClosure}(vi).
\end{proof}

Now we are ready to state and prove the main results of this section.  We begin by describing conditions that allow us to mimic, in the hypersmooth setting, the standard proof that finite-index super-equivalence relations of hyperfinite equivalence relations are hyperfinite.


\begin{lemma}\label{lem:FIoverHS}
Let $E\subseteq F$ be Borel equivalence relations on the standard Borel space $X$ such that $E$ is hypersmooth and $[F:E]<\infty$. Suppose $B\subseteq X$ is a Borel complete section for $F$ such that for each $k\in\Z_+$, the set
\[
X_k \ := \ \big\{x\in X\st [x]_F\cap B \textup{ meets exactly $k$ $E$-classes}\big\}
\]
is Borel and there exist Borel functions $f^{(k)}_0,\ldots,f^{(k)}_{k-1}:X_k\to B$ with the property that for all $x\in X_k$ and distinct $i,j\in k$, we have $f_i^{(k)}(x)\mathrel{F}x$ and $\neg\mkern1mu f_i^{(k)}(x)\mathrel{E}f_j^{(k)}(x)$. Then $F$ is hypersmooth.
\end{lemma}

\begin{proof}
Assume the hypotheses and notation in the statement of the lemma. Note that the sets $X_k$ are $F$-invariant.  By Lemma \ref{lem:piecewiseHS}, it will suffice to show that for each $k\in\Z_+$ the restriction $F^{(k)}$ of $F$ to $X_k$ is hypersmooth. So fix $k\in\Z_+$, and let $E^{(k)}$ be the restriction of $E$ to $X_k$, so $E^{(k)}$ is hypersmooth. Let $(D^{(k)}_n)_{n\in\omega}$ be an increasing sequence of smooth equivalence relations on $X_k$ such that $E^{(k)}=\cup_nD_n^{(k)}$. For each $n\in\omega$, let $g_n^{(k)}:X_k\to 2^{\omega}$ be a Borel reduction from $D_n^{(k)}$ to $\Delta(2^\omega)$, and let $h^{(k)}_n:X_k\to (2^\omega)^k$ be the Borel map defined by
\[
h^{(k)}_n(x) \ := \ \Big(\big(g_n^{(k)}\circ f_0^{(k)}\big)(x), \ldots, \big(g_n^{(k)}\circ f_{k-1}^{(k)}\big)(x)\Big).
\]
Let $\Delta_k$ be the finite (hence smooth) Borel equivalence relation on $(2^{\omega})^k$ defined by
\[
(x_0,\ldots,x_{k-1})\mathrel{\Delta_k}(y_0,\ldots,y_{k-1}) \quad :\Longleftrightarrow \quad (\exists\pi\in \textup{Sym}_k) \, (x_{\pi(0)},\ldots,x_{\pi(k-1)})=(y_0,\ldots,y_{k-1}),
\]
where $\textup{Sym}_k$ is the group of permutations of $k$. If we now let $F_n^{(k)}$ be the Borel equivalence relation on $X_k$ defined by
\[
x\mathrel{F_n^{(k)}}y \quad :\Longleftrightarrow \quad h^{(k)}_n(x)\mathrel{\Delta_k}h^{(k)}_n(y),
\]
so $h_n^{(k)}$ Borel reduces $F_n^{(k)}$ to $\Delta_k$,
then each $F_n^{(k)}$ is smooth and $F^{(k)}$ is the increasing union of the $F_n^{(k)}$, showing $F^{(k)}$ is hypersmooth as desired.
\end{proof}

\begin{theorem}\label{thm:E1fi}
Let $F$ be a Borel super-equivalence relation of $E_1$ on $(2^\omega)^\omega$ such that $[F:E_1]<\infty$. Then for each $k\in\omega$, the set
\(
X_k \ := \ \{x\in (2^\omega)^\omega \st [x]_F \text{ contains exactly $k$ $E_1$-classes}\}
\)
is Borel, and $F$ is hyperselective.
\end{theorem}

\begin{proof}
We prove by induction on $k\in\omega$ that each $X_k$ is Borel, while inductively constructing Borel functions $f_k$ satisfying the following:
\begin{enumerate}
 \item[(i)] $\dom(f_k)=(2^\omega)^\omega\setminus \cup_{i\leq k}X_i$ and $f_k(x)\in [x]_F$ for all $x\in\dom(f_k)$; and 
 \item[(ii)] for all $x\in X_{k+1}$, the points $f_0(x),\ldots,f_{k}(x)$ are pairwise $F$-equivalent and $E_1$-inequivalent.
\end{enumerate}
For the induction base $k=0$, we have $X_0=\emp$, so $X_0$ is Borel, and we set $f_0:=\textup{id}_{(2^\omega)^\omega}$, noting that (i) and (ii) then hold for $k=0$. Now fix $k\in\omega$, assume $X_i$ is Borel for all $i\leq k$, and suppose we have constructed Borel functions $f_0,\ldots,f_{k}$ satisfying (i) and (ii). Let 
\[
P \ := \ \big\{(x,y)\in \dom(f_k)\times (2^\omega)^\omega \st x\mathrel{F}y \, \wedge \, (\forall 0\leq i\leq k) \, \neg\, f_i(x)\mathrel{E_1}y \big\}.
\]
Then $P$ is a Borel subset of $\dom(f_k)\times (2^\omega)^\omega$, and each nonempty vertical section $\{y\in (2^\omega)^\omega\st (x,y)\in P\}$ of $P$ is a union of finitely many $E_1$-classes and is therefore $K_\sigma$ in $(2^\omega)^\omega$. It follows by the Arsenin-Kunugui Uniformization Thereom (see \cite[35.46]{Kechris} or \cite[5.12.1]{Srivastava}) that $P$ admits a Borel uniformization and has Borel projection to its first coordinate. As $X_{k+1}$ is the complement of this projection in $\dom(f_k)$, $X_{k+1}$ is also Borel. Furthermore, we can take $f_{k+1}$ to be the Borel uniformization of $P$ given by the Arsenin-Kunugui Theorem. It is straightforward to check that $f_{k+1}$ satisfies (i) and (ii), completing the inductive construction.

Now, for each $k\in\Z_+$ and $j\in k$, let $f_j^{(k)}:=f_{j}\res X_k$. Then the set $B=(2^\omega)^\omega$ and the functions $\big\{f_j^{(k)}\st k\in\Z_+,\, j\in k\big\}$ satisfy the hypotheses of Lemma \ref{lem:FIoverHS}, so $F$ is hypersmooth. Since $F$ has $K_{\sigma}$ classes and is uniformly continuous, it follows by Observation \ref{ob:UniformlyContinuousAndKSigma} that $F$ is hyperselective.
\end{proof}

Using Theorem \ref{thm:E1fi}, we can give a positive answer to Question \ref{ques:HSelFI} for uniformly continuous equivalence relations.

\begin{corollary}\label{cor:fi1}
Let $E\subseteq F$ be Borel equivalence relations on the standard Borel space $X$ such that $[F:E]<\infty$. If $E$ is hyperselective and uniformly continuous, then $F$ is hyperselective.
\end{corollary}

\begin{proof}
Assume the hypotheses. Since $E$ is hypersmooth, so is $E\times I_{2^\omega}$, so $E\times I_{2^{\omega}}\sqsubseteq_B E_1$. In fact, since $E$ is hyperselective, $E\times I_{2^\omega}\sqsubseteq_B^i E_1$ by \cite[Lemma 3.3]{KL}, and therefore $E\sqsubseteq_B^iE_1$ since $E$ is uniformly continuous. Fix an invariant Borel embedding $\phi:X\to (2^\omega)^\omega$ from $E$ to $E_1$, and let $Y=\ran(\phi)$, so $Y$ is Borel and $E_1$-invariant. Let $F^\phi$ be the Borel equivalence relation on $Y$ defined by
\[
x\mathrel{F^{\phi}}y \quad :\Longleftrightarrow \quad \phi^{-1}(x)\mathrel{F}\phi^{-1}(y).
\]
Then $F^{\phi}\cup E_1$ is a Borel equivalence relation on $(2^\omega)^\omega$ that contains $E_1$ as a finite index sub-equivalence relation, so $F^{\phi}\cup E_1$ is hyperselective by Theorem \ref{thm:E1fi}. Therefore $F^{\phi}$ is hyperselective by Fact \ref{fact:BasicClosure}, which means $F$ itself is hyperselective since $\phi$ is a Borel isomorphism from $F$ to $F^{\phi}$.
\end{proof}

\begin{corollary}\label{cor:fi2}
Let $\Ol$ be an $\ell$-dimensional product space, let $X$ be an $\EOl$-invariant Borel subset of $\FOl$, and write $F=\EOl\res X$. Suppose $F$ has a hypersmooth sub-equivalence relation of finite index. Then $F$ is hypersmooth, and if $\Omega$ is compact with $\ell\geq 2$, then $F$ is hyperselective.
\end{corollary}

\begin{proof} Assume the hypotheses, and let $E$ be a hypersmooth finite-index sub-equivalence relation of $F$. Call an $E$-class $C$ \emph{dense} if for all $x\in [C]_F$ there is $\bar{a}\in\N^{\ell}$ such that $\bar{a}\cdot x\in C$, and let $B$ be the union of all dense $E$-classes. Since $[F:E]<\infty$, every $F$-class contains at least one dense $E$-class, so for every $x\in X$ there is $\bar{a}\in\N^{\ell}$ such that $\bar{a}\cdot x\in B$; in particular, $B$ is a complete section for $F$. The set $B$ is also Borel, since $B=\{x\in X\st (\forall\bar{a}\in\N^{\ell})(\exists\bar{b}\in\N^{\ell})\, (\bar{a}+\bar{b})\cdot x\mathrel{E}x\}$.

Now, inductively define Borel functions $f_k$ for $k\in\omega$ as follows. Given $x\in X$ and $k\in\omega$, let $f_k(x)=\bar{a}\cdot x$ where $\bar{a}\in\N^{\ell}$ is least such that $\bar{a}\cdot x\in B$ and for all $i\in k$ we have $\neg\,f_i(x)\mathrel{E}\bar{a}\cdot x$, provided such $\bar{a}$ exists. Note that $\dom(f_0)=X$. For each $k\in\Z_+$, let $X_k$ be the set of all $x\in X$ such that $k$ is least for which $f_k$ is undefined, so $X_k$ is Borel and $F$-invariant. For each $k\in\Z_+$ and $j\in k$, put $f^{(k)}_j:=f_j\res X_k$. Then the set $B$, the sets $X_k$, and the functions $\big\{f_j^{(k)}\st k\in\Z_+,\, j\in k\big\}$ satisfy the hypotheses of Lemma \ref{lem:FIoverHS}, so $F$ is hypersmooth. If $\Omega$ is compact and $\ell\geq 2$, then $F$ is uniformly continuous with $K_{\sigma}$ classes by \ref{ob:UniformlyContinuous}, so $F$ is hyperselective by \ref{ob:UniformlyContinuousAndKSigma}.
\end{proof}


We close this section with a conjecture and an observation that will not be needed in what follows.

\begin{conjecture}\label{con:OTHERPROJECT}
For every Borel equivalence relation $E$, if $E\times I_{2^\omega}$ is hyperselective then $E$ is hyperselective.
\end{conjecture}

\begin{observation}
Conjecture \ref{con:OTHERPROJECT} implies a positive answer to Question \ref{ques:HSelFI}.
\end{observation}

\begin{proof}
Let $E\subseteq F$ be Borel equivalence relations with $E$ hyperselective and $[F:E]<\infty$.  Then $E\times I_{2^\omega}$ is a finite-index sub-equivalence relation of $F\times I_{2^\omega}$ that is both hyperselective and uniformly continuous, so $F\times I_{2^\omega}$ is hyperselective by Corollary \ref{cor:fi1}, and $F$ itself is hyperselective by Conjecture \ref{con:OTHERPROJECT}.
\end{proof}


\section{Marker sets and hypersmoothness}\label{sec:cascade}

Most proofs of hypersmoothness (in particular, hyperfiniteness) involve the use of \emph{marker sets}, the precise definition of which varies by context but which ultimately are just subsets of the underlying space that enjoy nice structural or geometric properties within equivalence classes. For our purposes, we will be especially interested in Borel marker sets in $\FOl$ having any of the following properties:

\begin{definition}\label{df:MarkerSets}
Let $\Ol$ be an $\ell$-dimensional product space, and $X$ an $\EOl$-invariant Borel subset of $\FOl$. Let $p,q\in\Z_+$. The subset $Y\subseteq\FOl$ is:
\begin{itemize}
 \item \emph{$p$-sparse} if for all $y,y'\in Y$, $d_{\ell}(y,y')<p$ implies $y\mathrel{F^{\Omega}_{\ell,\d{p}}}y'$;
 \item \emph{dense in $X$} if for all $x\in X$ there is $\bar{a}\in\N^{\ell}$ for which $\bar{a}\cdot x\in Y$; 
 \item \emph{$q$-dense in $X$} if for all $x\in X$ there is $\bar{a}\in\N^{\ell}_{<q}$ for which $\bar{a}\cdot x\in Y$;
 \item \emph{$(p,q)$-regular in $X$} if $Y$ is $p$-sparse and $q$-dense in $X$.
\end{itemize}
We say $Y$ is \emph{dense} if $Y$ is dense in its $\EOl$-saturation, and similarly define \emph{$q$-dense} and \emph{$(p,q)$-regular}.
\end{definition}

Note that if $Y$ is $(p,q)$-regular, then necessarily $p\leq q$. Given an $\EOl$-invariant subset $X\subseteq\FOl$, any transversal for $\EOl\res X$ naturally yields $(p,p)$-regular subsets $Y_p$ in $X$, for $p\geq 1$, as follows. Suppose $T\subseteq X$ is a transversal for $\EOl\res X$ with corresponding selector $\tau$, let $p\in\Z_+$, and put $x\in Y_p$ if and only if $p$ divides each component of $x-\tau(x)$. Then $Y_p$ is $(p,p)$-regular in $X$. Of course, $Y_p$ need not be Borel if $X$ or $T$ is not Borel, and we expect that for all $\ell\geq 1$ and $p\geq 2$ there is no $(p,p)$-regular Borel set in any $\ETl$-invariant Borel set $X\subseteq\Tl$ of positive $\mu_{\ell}$-measure, though we do not need this fact and will not attempt to prove it. In general, we have found constructing $(p,q)$-regular Borel marker sets in nontrivial subsets of $\FTl$ to be quite challenging, and in the next section we note that when $\ell\geq 2$ it can only be done in proper subsets of $\FTl$ (see Corollary \ref{cor:NoNiceMarkerSets}).

The following lemma shows that $p$-sparse marker sets and $q$-dense marker sets can always be made as ``$F$-invariant'' as we like:

\begin{lemma}\label{lem:RegularInvar}
Let $\Ol$ be an $\ell$-dimensional product space, and $X$ an $\EOl$-invariant Borel subset of $\FOl$. Let $p,q\in\Z_+$, and suppose $Y\subseteq X$ is Borel. 
\begin{enumerate}[(i)]
 \item If $Y$ is $p$-sparse, then for every $r\in\omega$ the set $\d{r}^{-1}[Y]$ is $F^{\Omega}_{\ell,\d{r}}$-invariant and $p$-sparse.
 \item If $Y$ is $q$-dense in $X$, then for every $r\in\omega$ the set $\d{r}^{-1}[Y]$ is $F^{\Omega}_{\ell,\d{r}}$-invariant and $q$-dense in $X$.
 \item If $Y$ is $(p,q)$-regular in $X$, then for every $r\in\omega$ the set $\d{r}^{-1}[Y]$ is $F^{\Omega}_{\ell,\d{r}}$-invariant and $(p,q)$-regular in $X$.
\end{enumerate}
\end{lemma}

\begin{proof}
Let $p$, $q$, and $Y$ be given, and fix arbitrary $r\in\omega$ and write $Z:=\d{r}^{-1}[Y]$. That $Z$ is Borel and $F^{\Omega}_{\ell,\d{r}}$-invariant is clear. For (i), assume $Y$ is $p$-sparse. Fix $x,y\in Z$ such that $d(x,y)<p$, along with $\bar{a},\bar{b}\in\N^{\ell}_{<p}$ such that $\bar{a}\cdot x=\bar{b}\cdot y$. Then $\bar{a}\cdot (\d{r}\cdot x)=\bar{b}\cdot (\d{r}\cdot y)$, so $d(\d{r}\cdot x, \d{r}\cdot y)<p$, which implies $\d{r}\cdot x\mathrel{F^{\Omega}_{\ell,\d{p}}}\d{r}\cdot y$ since $Y$ is $p$-sparse. This means $x\mathrel{F^{\Omega}_{\ell,\d{r}+\d{p}}}y$, and therefore $x\mathrel{F^{\Omega}_{\ell,\d{p}}}y$ since $d(x,y)<p$. For (ii), assume $Y$ is $q$-dense in $X$, and let $x\in X$ be arbitrary. Fix $\bar{a}\in\N_{<q}$ such that $\bar{a}\cdot (\d{r}\cdot x)\in Y$. Then $\d{r}\cdot (\bar{a}\cdot x)\in Y$, so $\bar{a}\cdot x\in Z$, showing $Z$ is $q$-dense in $X$. Finally, (iii) is immediate from (i) and (ii).
\end{proof}

Our goal now will be to show how the existence of regular Borel marker sets in $\EOl$-invariant Borel subsets $X\subseteq\FOl$ implies hypersmoothness of $\EOl\res X$. For this we use a general fact, stated below as Lemma \ref{lem:cascade}, that is based on the proof from \cite{JKL} showing that Borel actions of finitely generated groups of polynomial growth yield hyperfinite orbit equivalence relations (see \cite[Lemma 1.19]{JKL}). Our statement is particularly inspired by the presentation of this result in \cite[Ch 11]{KM} using the notion of a \emph{cascade}. 
Recall that to any extended metric space $(X,d)$ we associate the equivalence relation $E_d=\{(x,y)\in X^2\st d(x,y)<\infty\}$. Furthermore, we say that the function $f:X\to X$ on the (extended) metric space $(X,d)$ has \emph{bounded displacement} if there is some $M\in\mathbb{R}$ such that $d\big(x,f(x)\big)<M$ for all $x\in X$. 

\begin{lemma}[{See \cite[Ch 11]{KM}} for comparison]\label{lem:cascade}
Let $d:X\times X\to [0,\infty]$ be a Borel extended metric on the standard Borel space $X$, and let $(g_n)_{n\in\omega}$ be a sequence of Borel functions from $X$ to $X$ with bounded displacement. Suppose there is a fixed constant $C\in\Z_+$ such that for every $n\in\omega$ and $d$-open ball $B_n\subseteq X$ of radius $n$, the image of $B_n$ under $g_n$ has cardinality at most $C$. Then $E_d$ has a hypersmooth sub-equivalence relation of index at most $C$.
\end{lemma}

\begin{proof}
Assume the hypotheses, and in particular let $C$ be as stated. For each $n\in\omega$, fix $m_n>0$ such that $d\big(x,g_n(x)\big)<m_n$ for all $x\in X$. Let $(n_k)_{k\in\omega}$ be an increasing sequence of positive integers such that
\[
n_{k+1} \ > \ 2(m_{n_0}+\cdots +m_{n_k})+k
\]
for all $k\in\omega$. For each $k\in\omega$, let $f_k:=g_{n_k}\circ\cdots\circ g_{n_0}$, and define the equivalence relation $F_k$ on $X$ by
\[
x\mathrel{F_k}y \quad :\Longleftrightarrow \quad f_k(x)=f_k(y).
\]
Then $(F_k)_{k\in\omega}$ is an increasing sequence of smooth sub-equivalence relations of $E_d$, and we claim that the hypersmooth sub-equivalence relation $F:=\cup_kF_k$ of $E_d$ has index at most $C$ in $E_d$. To see this, we first show by induction on $k$ that for all $k\in\omega$,
\begin{equation*}\tag{$\ast$}
d\big(x,\,f_k(x)\big) \ < \ m_{n_0}+\cdots+m_{n_k} \quad \text{for all $x\in X$}.
\end{equation*}
For the base case $k=0$ we have
\[
d\big(x,\,f_0(x)\big) \ = \ d\big(x,g_{n_0}(x)\big) \ < \ m_{n_0} \quad \text{for all $x\in X$}.
\]
For the inductive step, fix $k\geq 0$ and assume $(\ast)$ holds. Then for all $x\in X$ we have
\begin{align*}
d\big(x,\,f_{k+1}(x)\big) \ &\leq \ d\big(x,\,f_{k}(x)\big)+d\big(f_k(x),\,f_{k+1}(x)\big) \\
&< \ m_{n_0}+\cdots+m_{n_k} + d\big(f_{k}(x),\,g_{n_{k+1}}(f_k(x))\big) \ < \ m_{n_0}+\cdots+m_{n_k}+m_{n_{k+1}},
\end{align*}
completing the induction.

Suppose now for contradiction that $x_0,\ldots,x_C$ are pairwise $F$-inequivalent points in the same $E_d$-class in $X$. Fix $k$ such that $d(x_0,x_i)\leq k$ for all $0\leq i\leq C$. Then by $(\ast)$ we have that for each $0\leq i\leq C$,
\begin{align*}
d\big(f_k(x_0),\,f_k(x_i)\big) \ &\leq \ d\big(f_k(x_0), \, x_0\big) + d(x_0,x_i) + d\big(x_i,\,f_k(x_i)\big) \\
&< \ (m_{n_0}+\cdots+m_{n_k})+k+(m_{n_0}+\cdots+m_{n_k}) \ < \ n_{k+1}. 
\end{align*}
So $\{f_{k}(x_i)\st 0\leq i\leq C\}$ is contained in the open ball of radius $n_{k+1}$ centered at $f_k(x_0)$, which implies that its image under $g_{n_{k+1}}$ has  cardinality at most $C$. Thus $\{f_{k+1}(x_i)\st 0\leq i\leq C\}$ has cardinality at most $C$, contradicting the fact that the $x_i$ are pairwise $F_{k+1}$-inequivalent.
\end{proof}

Using Lemma \ref{lem:cascade}, we can now prove one of the key results to which we will ultimately appeal in our proof that $\ETT$ is hypersmooth. 

%

\begin{theorem}\label{thm:RegularHS}
Let $\Ol$ be an $\ell$-dimensional product space, and $X$ an $\EOl$-invariant Borel subset of $\FOl$. Suppose there is a fixed positive integer $c$ such that for every $n\in\Z_+$ there is a Borel subset $Z_n\subseteq X$ that is $(n,cn)$-regular in $X$. Then $\EOl\res X$ is hyperselective.
\end{theorem}

\begin{proof}
Assume the hypotheses, and in particular let the constant $c$ and sets $Z_n$ be as stated. For each $n\in\Z_+$ and $x\in X$, let $\tau_n(x):=\bar{a}\cdot x$ where $\bar{a}\in \N^{\ell}_{<cn}$ is least such that $\bar{a}\cdot x\in Z_n$, and let $g_n(x):=2\d{n}\cdot\tau_{2n}(x)$. Note that for each $n\in\Z_+$ and $x\in X$ we have $d\big(x,\,g_n(x)\big)<2cn+2n$, so each $g_n$ has bounded displacement. Let $n\in\Z_+$ and $x\in X$ be arbitrary, and write $B:=B_n(x)$; we claim that $\tau_{2n}[B]$ is contained in the union of the open balls of radius $n$ about the points $n\bar{u}\cdot x$, where $\bar{u}\in\N^{\ell}_{\leq 2c}$.  To see this, fix $y\in B$, along with $\bar{a},\bar{b}\in\N^{\ell}_{<n}$ such that $\bar{a}\cdot x=\bar{b}\cdot y$ and $\c\in\N^{\ell}_{<2cn}$ such that $\tau_{2n}(y)=\bar{c}\cdot y$.  For each integer $1\leq i\leq \ell$, let $0\leq u_i\leq 2c$ and $0\leq r_i<n$ be the unique integers satisfying $a_i+c_i=u_in+r_i$, so
\[
\bar{a}+\bar{c} \ = \ (u_1n+r_1,\ldots,u_{\ell}n+r_{\ell}) \ = \ n\bar{u}+\bar{r}.
\]
Then 
\[
\bar{b}\cdot (\bar{c}\cdot y) \ = \ (\bar{a}+\bar{c})\cdot x \ = \ \bar{r}\cdot(n\bar{u}\cdot x),
\]
so
\[
d\Big(\tau_{2n}(y), \, n\bar{u}\cdot x\Big) \ < \ n.
\]
This shows $\tau_{2n}[B]$ is indeed covered by the $(2c+1)^{\ell}$ open balls of radius $n$ centered at the points $n\bar{u}\cdot x$, where $\bar{u}\in\N^{\ell}_{\leq 2c}$. Since $\tau_{2n}$ takes values in the $2n$-sparse set $Z_{2n}$, each of these $(2c+1)^\ell$ balls can meet $Z_{2n}$ in at most one $F^{\Omega}_{\ell,2\d{n}}$-class. 
But this implies that $g_n[B]$ has size at most $(2c+1)^\ell$. We can now apply Lemma \ref{lem:cascade} using the sequence of functions $(g_n)$ and the constant $C:=(2c+1)^\ell$ to obtain a hypersmooth subequivalence relation of $\EOl\res X$ of index at most $C$; then Corollary \ref{cor:fi2} implies that $\EOl\res X$ is hyperselective.
\end{proof}

We conclude this section with some additional comments about marker sets in $\FTl$ that will not be needed later, but may help motivate Definition \ref{df:MarkerSets}; specifically, we discuss what happens when one attempts to show $\ETl\res\FTl$ is hypersmooth by following as far as possible the proof from \cite{JKL} that Borel actions of finitely generated groups of polynomial growth have hyperfinite orbit equivalence relations. All our comments here refer exclusively to the free parts of any actions considered.  In Theorem 1.16 of \cite{JKL}, Chapter 11 of \cite{KM}, and Section 2 of \cite{GJ}, we find the following general pattern of argument (see also \cite[Theorem 7.1]{AsymDim}).  First, there is a locally finite Borel graph $\G$ that generates the equivalence relation $E$ under consideration, along with a metric $\rho$ on $E$-classes that is closely related to $\G$. For each distance $n$, the graph 
\[
x\mathrel{\G_n}y \quad :\Longleftrightarrow \quad x\ne y \ \wedge \ \rho(x,y)<n
\]
admits a Borel $\aleph_0$-coloring, which may be used in combination with a greedy algorithm to produce ``Borel kernels" for the graphs $\G_n$, that is, Borel sets $S_n$ satisfying: 
\begin{enumerate}[(i)]
 \item if $x,y\in S_n$ are distinct, then $\rho(x,y)\geq n$; and
 \item for every point $x$, 
 $\rho(x,S_n)<n$.
\end{enumerate}
Then for every marker distance $n$, there is a Borel function $g_n$ taking each $x$ to some nearest element of $S_n$, yielding a cascade of functions $g_n$ and marker sets $S_n$ which can be used to define a hypersmooth sub-equivalence relation $\mathcal{F}$ of $E$ as in \eqref{eq:cascade} below. If for each $n$ the marker set $S_n$ is sufficiently sparsely distributed, a rapidly increasing sequence of marker distances can be chosen as in the proof of Lemma \ref{lem:cascade} so that $\mathcal{F}$ has finite index in $E$. 

For the most part, this entire approach can be carried out for $\N^{\ell}\curvearrowright 2^{\omega^{\ell}}$, even though we have no group action or locally finite graph available, but it fails to show $E_{\ell}^2$ is hypersmooth because the resulting marker sets are not ``sufficiently sparsely distributed." Specifically, we can show the following, although we omit a proof since these results will not be needed in what follows.

\begin{proposition}\label{prop:extra}
Let $\ell$ and $n$ be arbitrary positive integers. 
\begin{enumerate}[(i)]
 \item There exists an $n$-sparse 
 relatively clopen set $S_n\subseteq\FTl$ such that $d_{\ell}(x,S_n)<n$ for every $x\in\FTl$.
 \item With $S_n$ as in (i), there is a continuous function $g_n:\FTl\to S_n$ such that $d_{\ell}\big(g_n(x),x\big)<n$ for every $x\in\FTl$.
\end{enumerate}
\end{proposition}

From these functions $g_n$ we obtain, for any rapidly growing sequence $\alpha$ of positive integers, a hypersmooth sub-equivalence relation $\mathcal{F}_{\alpha}:=\bigcup_n\mathcal{F}_{\alpha,n}$ of $\ETl\res\FTl$, where
\begin{equation}\label{eq:cascade}
x\mathrel{\mathcal{F}_{\alpha,n}}y \quad :\Longleftrightarrow \quad (g_{\alpha(n)}\circ\cdots\circ g_{\alpha(0)})(x) = (g_{\alpha(n)}\circ\cdots\circ g_{\alpha(0)})(y).
\end{equation} 
But no matter how $\alpha$ is chosen, we cannot expect $\mathcal{F}_{\alpha}$ to have finite index in $\ETl\res\FTl$ since there is no fixed constant $C$ bounding the number of $F^{2}_{\ell,\d{n}}$-classes that intersect $g_n[B_n(x)]$. To see why, let us consider again the geometry of the $\ETl$-classes. With Figure \ref{fig:2dimProj} in mind, recall that after being projected onto $\FTl/F^2_{\ell}$, $d_{\ell}$-balls in $\FTl$ look like $\ell$-cubes in $\Z^{\ell}$. Given a $d_{\ell}$-ball $B_r(x)$ of radius $r$, $g_n[B_r(x)]$ will be contained in $B_{r+n}(x)$, and we must ask how many $F^2_{\ell,\d{n}}$-classes can meet $S_n\cap B_{r+n}(x)$. Counting $F^2_{\ell,\d{n}}$-classes that meet $S_n\cap B_{r+n}(x)$ in \emph{front} of $x$ yields the same answer as if we were dealing with orbits of $\Z^{\ell}\curvearrowright 2^{\Z^{\ell}}$: it is not too hard to show that for each $c\in\Z_+$, $\N^{\ell}_{<cn}\cdot x$ can meet $S_n$ in at most $c^{\ell}$ distinct $F^2_{\ell,\d{n}}$-classes. But due to the geometry of $d_{\ell}$-balls that we described in Section \ref{subsec:geometry} and illustrated using Figures \ref{fig:EtClass}, \ref{fig:1dimProj}, and \ref{fig:2dimProj}, for each $r\geq n$ the set $S_n\cap \{y\in B_r(x)\st x-y\in\N^{\ell}\}$ can meet all $r^{\ell}-1$ different $F^2_{\ell,\d{r}}$-classes that lie strictly \emph{behind} $x$, and we run into similar trouble if we consider $F^2_{\ell,\d{r}}$-classes containing points $y\in S_n\cap B_r(x)$ for which $x-y\in\Z^{\ell}\setminus\pm\N^{\ell}$. Thus the number of $F^2_{\ell,\d{r}}$-classes meeting $S_n\cap g_n[B_r(x)]$ grows without bound as $n\to\infty$ no matter how $r$ is chosen in relation to $n$, and we can have no hope of applying anything like Lemma \ref{lem:cascade} to these functions $g_n$.

This discussion explains how applying a cascade argument to $\N^\ell\curvearrowright 2^{\omega^\ell}$ must differ from applying one to $\Z^\ell\curvearrowright 2^{\Z^{\ell}}$: in the latter case, the functions $g_n$ can safely move points in any direction, but in the former case they must move points forward. Thus at stage $n$, for some appropriate distance $r_n$, our marker sets need to meet not just all balls $B_{r_n}(x)$ but all regions of the form $\N^{\ell}_{<r_n}\cdot x$, and this is what motivates our definitions of \emph{dense} and \emph{$(p,q)$-regular} marker sets. Indeed, the sets $S_n$ described in Proposition \ref{prop:extra} are not dense in $\FTl$. One might hope to carry out some sort of local Borel algorithm to transform them into regular marker sets or even just dense sparse ones, but we have not been able to make progress in this direction due to the peculiar geometry of $\ETl$-classes. Indeed, we will see in Corollary \ref{cor:NoNiceMarkerSets} that for all $\ell,p\geq 2$ and $q\geq 1$, there does not exist a Borel subset $Y\subseteq\FTl$ that is both $p$-sparse and $q$-dense in all of $\FTl$. We \emph{will} be able to construct regular Borel marker sets in large invariant subsets of $\FTT$, but our arguments do not use Proposition \ref{prop:extra} or any sort of local Borel algorithm that adjusts marker points or boundaries of marker regions geometrically in the spirit of \cite{GJ} or \cite{Nilpotent}. Instead, we will employ a combinatorial argument using Borel graph colorings, which is the main topic of the next section.


\section{Borel graph colorings and regular marker sets}\label{sec:graphs}

Let $\Ol$ be an $\ell$-dimensional product space, and $X$ an $\EOl$-invariant Borel subset of $\FOl$. In the previous section, we reduced the problem of showing $\EOl\res X$ is hypersmooth to that of constructing regular Borel marker sets in $X$. In this section, we establish a close connection between regular Borel marker sets in $X$ and finite Borel colorings of certain graphs on $X$ arising from the $\N^{\ell}$-shift action. Roughly speaking, we show how to produce regular Borel marker sets from Borel graph colorings, and regular clopen marker sets from continuous graph colorings. We begin by defining various graphs on $\EOl$-classes and introducing notation for them, continuing the discussion begun in Section \ref{subsec:graphs}.

As we defined in Equation \eqref{eq:Gl1}, for each $\ell\in\Z_+$ let $\vec{\G}_{\ell}$ be the directed graph on $\Ol$ given by
\[
x\mathrel{\vec{\G}_{\ell}}y \quad :\Longleftrightarrow \quad x\ne y \ \wedge \ (\exists i\leq\ell) \, \bar{e}_i\cdot x=y,
\]
and let $\G_{\ell}$ be its symmetrization.  Additionally define, for each $\ell\in\Z_+$ and nonzero $\m\in\N^{\ell}$, the directed graph $\vec{\G}_{\m}$ on $\Ol$ by
\begin{equation}\label{eq:Gl}
x\mathrel{\vec{\G}_{\m}}y \quad :\Longleftrightarrow \quad x\ne y \ \wedge \ \m\cdot x=y,
\end{equation}
and let $\G_{\m}$ be its symmetrization. Note that 
\[
\vec{\G}_{\ell}=\bigcup_{1 \leq i\leq\ell}\vec{\G}_{\bar{e}_i} \qquad \text{and} \qquad {\G}_{\ell}=\bigcup_{1\leq i\leq\ell}{\G}_{\bar{e}_i}.
\]
Now we generalize \eqref{eq:Gl} as follows. For each $m\in\Z$, define 
\[
m^+ \ := \ \begin{cases} m & \text{if $m\geq 0$} \\ 0 & \text{if $m<0$} \end{cases} \qquad \text{and} \qquad 
m^- \ := \ \begin{cases} 0 & \text{if $m\geq 0$} \\ |m| & \text{if $m<0$,} \end{cases} 
\]
and for each $\m=(m_1,\ldots,m_{\ell})\in\Z^\ell$, let $\m^+$ and $\m^-$ be the $\ell$-tuples in $\N^\ell$ defined by
\begin{equation}\label{eq:mplusminus}
\m^+ \ := \ (m_1^+,\ldots,m_{\ell}^+) \qquad \text{and} \qquad \m^- \ := \ (m_1^-,\ldots,m_\ell^-).
\end{equation}
Then, given arbitrary nonzero $\m\in\Z^{\ell}$, define the directed graph $\vec{\G}_{\m}$ on $\Ol$ by
\begin{equation}\label{eq:mq}
x\mathrel{\vec{\G}_{\m}}y \quad :\Longleftrightarrow \quad x\ne y \ \wedge \ \m^+\cdot x=\m^-\cdot y,
\end{equation}
and let $\G_{\m}$ be its symmetrization.  Note that this agrees with our previous definition of $\vec{\G}_{\m}$ and $\G_{\m}$ when $\m\in\N^{\ell}$, since in that case $\m^+=\m$ and $\m^-=\bar{0}$. Note also that if $x\mathrel{\vec{\G}_{\m}}y$ then $y-x=\m$ and $d(x,y)=|\m|$. Of the graphs $\Gm$ for $\m\in\Z^2\setminus\{\bar{0}\}$, two will be particularly important for us in Sections \ref{sec:ep}--\ref{sec:Gm}, namely those where $\m=(1,\pm 1)$. For this reason, in Sections \ref{sec:ep}--\ref{sec:Gm} we will use special notation for these two vectors, writing $\p:=(1,1)$ and $\q:=(1,-1)$.

Continuing with our definitions, given $\ell,n\in\Z_+$, define for each nonzero $\m\in\Z^{\ell}$ the (undirected) graph $\G_{\m,n}$ on $\Ol$ by
\[
x\mathrel{\G_{\m,n}}y \quad :\Longleftrightarrow \quad x\ne y \, \wedge \, \big(\exists \bar{a},\bar{b}\in\N^{\ell}_{<n}\big) \, \Big(\bar{a}\cdot x=\bar{b}\cdot y \, \wedge \, \bar{b}-\bar{a}=\pm\m\Big),
\]
and let 
\begin{equation}\label{eq:Gln}
\mathrel{\G_{\ell,n}} \ := \ \bigcup_{\m \mkern1mu \in \mkern1mu \Z^{\ell}_{<n}\setminus\{\bar{0}\}}\G_{\m,n} \ = \ \Big\{(x,y)\in \Ol\times\Ol\st \neg\,  x\mathrel{F^{\Omega}_{\ell}}y \, \wedge \, d_{\ell}(x,y)<n\Big\}.
\end{equation}
Note that $\G_{\m,n}=\G_{-\m,n}$ for each $\m\in\Z^{\ell}\setminus\{\bar{0}\}$. In the case $\ell=2$, write
\[
\G_{2,n,+} \ := \ \bigcup_{\m \mkern1mu \in \mkern1mu \N^{2}_{<n}\setminus\{\bar{0}\}}\G_{\m,n} \qquad \text{and} \qquad 
\G_{2,n,-} \ := \ \bigcup_{\m \mkern1mu \in \mkern1mu (\Z_+\times\Z_-)_{<n}}\G_{\m,n},
\]
so $\G_{2,n}=\G_{2,n,+}\cup\G_{2,n,-}$.

All the graphs above are defined the same way over any base space, so we omit the base space from the notation when convenient. When we do wish to specify the base space, we will include it in parentheses after the graph name, so, for instance, $\vec{\G}_{2}(3)$ is the directed graph on $3^{\omega^2}$ defined by
\[
x\mathrel{\vec{\G}_{2}(3)}y \quad \Longleftrightarrow \quad x\ne y \ \wedge \ \Big(\bar{e}_1\cdot x=y \ \vee \ \bar{e}_2\cdot x=y\Big),
\]
and $\vec{\G}_{2}(2^\omega)$ is the same graph defined on $(2^\omega)^{\omega^2}$.  When we wish to restrict any of the graphs defined above to the free part of the base space, we will write, e.g., $\G_{\m}^{\textup{Fr}}(2):=\G_{\m}(2)\res\FTl$ for this restriction.

Now we turn our attention to establishing a connection between regular Borel marker sets and finite Borel graph colorings. We start by pointing out some nice properties enjoyed by $\FeqOl{\c}$-invariant sets, for $\c\in\N^{\ell}$, in the next three lemmas.

\begin{lemma}\label{lem:ForwardBorel}
Let $r\in\Z_+$ and $\bar{a},\c\in\N^{\ell}$ with $\d{r}\preceq\bar{c}$ and $\bar{a}\preceq\c$. Let $X$ be an $\EOl$-invariant Borel subset of the product space $\Ol$, and suppose $Y$ is an $F_{\ell,\c}^{\Omega}$-invariant subset of $X$.
\begin{enumerate}[(i)]
 \item If $Y$ is Borel, then $\bar{a}\cdot Y$ and $B_r[Y]$ are Borel.
 \item If $Y$ is clopen in $X$, then $\bar{a}\cdot Y$ and $B_r[Y]$ are clopen in $X$.
\end{enumerate}
\end{lemma}

\begin{proof}
Assume the hypotheses, and in particular assume $Y\subseteq X$ is $\FeqOl{\c}$-invariant, which means $Y$ is $F^{\Omega}_{\ell,\bar{a}}$-invariant since $\bar{a}\preceq\c$. Let $x_{\Omega}$ be some fixed element of $\Omega$, and define the continuous function $h_{\bar{a}}:\Ol\to\Ol$ by
\[
h_{\bar{a}}(x)(\bar{i}) \ := \ \begin{cases}
    x(\bar{i}-\bar{a}) & \text{if $\bar{i}-\bar{a}\in\omega^{\ell}$;} \\
    x_{\Omega} & \text{otherwise.}
\end{cases}
\]
Note that for all $x\in\Ol$ we have $\bar{a}\cdot h_{\bar{a}}(x)=x$, and also $h_{\bar{a}}(\bar{a}\cdot x)\mathrel{F^{\Omega}_{\ell,\bar{a}}}x$. Note further that $h_{\bar{a}}^{-1}[X]=X$, since $X$ is $\EOl$-invariant.  We claim now that $\bar{a}\cdot Y=h_{\bar{a}}^{-1}[Y]$. Indeed, if $x\in h_{\bar{a}}^{-1}[Y]$, so $h_{\bar{a}}(x)\in Y$, then $x=\bar{a}\cdot h_{\bar{a}}(x)\in\bar{a}\cdot Y$; conversely, if $x=\bar{a}\cdot y$ where $y\in Y$, then $h_{\bar{a}}(x)=h_{\bar{a}}(\bar{a}\cdot y)\mathrel{F^{\Omega}_{\ell,\bar{a}}}y$, so $h_{\bar{a}}(x)\in Y$ since $Y$ is $F^{\Omega}_{\ell,\bar{a}}$-invariant. Thus $\bar{a}\cdot Y=h_{\bar{a}}^{-1}[Y]$. As $h_{\bar{a}}$ is continuous, this means $\bar{a}\cdot Y$ is Borel if $Y$ is Borel, and clopen in $h^{-1}_{\bar{a}}[X]=X$ if $Y$ is clopen in $X$.

Furthermore, since
\[
B_r[Y] \ = \ \bigcup_{\bar{u},\bar{v}\,\in\,\N^{\ell}_{\leq r}}\bar{v}^{-1}\big[\bar{u}\cdot Y\big]
\]
and $\d{r}\preceq\c$, it follows from what we just showed together with Claim \ref{claim:NiceAction} that $B_r[Y]$ is Borel if $Y$ is Borel, and clopen in $X$ if $Y$ is clopen in $X$.
\end{proof}

\begin{lemma}\label{lem:ForwardInvar}
Let $r\in\Z_+$ and $\bar{a},\c\in\N^{\ell}$, and let $Y$ be an $\FeqOl{\c}$-invariant subset of the product space $\Ol$.
\begin{enumerate}[(i)]
 \item The set $\bar{a}^{-1}[Y]$ is $F^{\Omega}_{\ell,\bar{a}+\c}$-invariant.
 \item If $\bar{a}\preceq\c$, then the set $\bar{a}\cdot Y$ is $F^{\Omega}_{\ell,\c-\bar{a}}$-invariant.
 \item If $\d{r}\preceq\c$, then the set $B_r[Y]$ is $F^{\Omega}_{\ell,\bar{c}-\d{r}}$-invariant.
\end{enumerate}
\end{lemma}

\begin{proof}
Assume the hypotheses. For (i), suppose $x\in\bar{a}^{-1}[Y]$ and $x\mathrel{F^{\Omega}_{\ell,\bar{a}+\c}}x'$, and write $y=\bar{a}\cdot x\in Y$ and $y'=\bar{a}\cdot x'$. Then
\[
\c\cdot y \ = \ (\bar{a}+\c)\cdot x \ = \ (\bar{a}+\c)\cdot x' \ = \ \c\cdot y',
\]
so $y'\in Y$ since $Y$ is $F^{\Omega}_{\ell,\c}$-invariant, which means $x'\in\bar{a}^{-1}[Y]$. This shows $\bar{a}^{-1}[Y]$ is $F^{\Omega}_{\ell,\bar{a}+\c}$-invariant.

For (ii), assume $\bar{a}\preceq\c$. Fix $z\in\bar{a}\cdot Y$, along with $y\in Y$ and $z'\in\Ol$ such that $z=\bar{a}\cdot y$ and $z\mathrel{F^{\Omega}_{\ell,\c-\bar{a}}}z'$. Define $y'\in\Ol$ by 
\[
y'(\bar{i}) \ := \ \begin{cases}
z'(\bar{i}-\bar{a}) & \text{if $\bar{i}-\bar{a}\in\omega^{\ell}$;} \\ 
y(\bar{i}) & \text{otherwise.}
\end{cases}
\]
Then $\bar{a}\cdot y'=z'$ and 
\[
\c\cdot y \ = \ (\c-\bar{a})\cdot (\bar{a}\cdot y) \ = \ (\c-\bar{a})\cdot z \ = \ (\c-\bar{a})\cdot z' \ = \ (\c-\bar{a})\cdot (\bar{a}\cdot y') \ = \ \c\cdot y',
\]
so $y'\in Y$ since $Y$ is $\FeqOl{\c}$-invariant. Thus $z'\in\bar{a}\cdot Y$, showing $\bar{a}\cdot Y$ is indeed $F^{\Omega}_{\ell,\c-\bar{a}}$-invariant.

Finally, for (iii), assume $\d{r}\preceq\c$. By (i) and (ii), for all $\bar{a},\bar{b}\in\N^{\ell}_{\leq r}$ the set $\bar{b}^{-1}\big[\bar{a}\cdot Y\big]$ is $\FeqOl{\c-\bar{a}+\bar{b}}$-invariant, hence $\FeqOl{\c-\d{r}}$-invariant since $\c-\d{r}\preceq \c-\bar{a}+\bar{b}$. Since 
\[
B_r[Y] \ = \ \bigcup_{\bar{a},\bar{b}\,\in\,\N^{\ell}_{\leq r}}\bar{b}^{-1}\big[\bar{a}\cdot Y\big] 
\]
and unions of invariant sets are invariant, the claim follows.
\end{proof}

\begin{lemma}\label{lem:ForwardSparse}
Let $p\in\Z_+$ and $\bar{a},\bar{c}\in\N^{\ell}$ with $\bar{a}+\d{p}\preceq\bar{c}$, and let $Y$ be an $\FeqOl{\bar{c}}$-invariant subset of the product space $\Ol$. If $Y$ if $p$-sparse, then $\bar{a}\cdot Y$ is $p$-sparse.
\end{lemma}

\begin{proof}
Let $p$, $\bar{a}$, and $\bar{c}$ be as stated, and suppose $Y\subseteq\Ol$ is $\FeqOl{\bar{c}}$-invariant and $p$-sparse. Let $z_0,z_1\in\bar{a}\cdot Y$ with $d(z_0,z_1)<p$.  Fix $\bar{u},\bar{v}\in\N^{\ell}_{<p}$ such that $\bar{u}\cdot z_0=\bar{v}\cdot z_1$, along with $y_0,y_1\in Y$ such that $\bar{a}\cdot y_0=z_0$ and $\bar{a}\cdot y_1=z_1$. Fix some $x_{\Omega}\in\Omega$, and define $y_1'\in\Ol$ by
\[
y_1'(\bar{i}) \ = \ \begin{cases} 
y_0(\bar{i}+\bar{u}-\bar{v}) & \text{if $\bar{i}+\bar{u}-\bar{v}\in\omega^{\ell}$;} \\ x_{\Omega} & \text{otherwise.}
\end{cases}
\]
Then $\bar{u}\cdot y_0=\bar{v}\cdot y_1'$, so $d(y_0,y_1')<p$. Furthermore, for all $\bar{i}\in\omega^{\ell}$ we have
\[
y_1'(\dot{p}+\bar{a}+\bar{i}) \ = \ y_0(\dot{p}+\bar{a}+\bar{i}+\bar{u}-\bar{v}) \ = \ (\bar{u}\cdot z_0)(\dot{p}-\bar{v}+\bar{i}) \ = \ (\bar{v}\cdot z_1)(\dot{p}-\bar{v}+\bar{i}) \ = \ y_1(\dot{p}+\bar{a}+\bar{i}),
\]
which shows both $\bar{a}\cdot y_1\mathrel{F^{\Omega}_{\ell,\d{p}}}\bar{a}\cdot y_1'$ and $y_1\mathrel{F^{\Omega}_{\ell,\bar{c}}}y_1'$, since $\dot{p}+\bar{a}\preceq\bar{c}$. This means $y_1'\in Y$ since $Y$ is $F^{\Omega}_{\ell,\bar{c}}$-invariant, which implies $y_0\mathrel{F^{\Omega}_{\ell,\d{p}}}y_1'$ since $Y$ is $p$-sparse. Thus $\bar{a}\cdot y_0\mathrel{F^{\Omega}_{\ell,\d{p}}}\bar{a}\cdot y_1'$, and therefore $z_0\mathrel{F^{\Omega}_{\ell,\d{p}}}z_1$, showing that $\bar{a}\cdot Y$ is indeed $p$-sparse as claimed.
\end{proof}

Now we come to one of the key lemmas in this paper, which states that $(p,q)$-regular marker sets can be converted into $(p,3p)$-regular marker sets. The ``$3p$" here is not necessarily optimal, but it will suffice for our purposes.

\begin{lemma}\label{lem:SparseBoundedlyDense}
Let $\Ol$ be an $\ell$-dimensional product space, and $X$ an $\EOl$-invariant Borel subset of $\FOl$. Let $p,q\in\Z_+$, and suppose there exists a Borel subset $Y\subseteq X$ that is $(p,q)$-regular in $X$. Then there exists a Borel subset $Z\subseteq X$ that is $(p,3p)$-regular in $X$. Furthermore, if $Y$ is clopen in $X$ then $Z$ can be taken to be clopen in $X$ as well.
\end{lemma}

\begin{proof}
Assume the hypotheses, and let $r:=p+p(2p+q)^\ell$. By replacing $Y$ with $\d{r}^{-1}[Y]$ if necessary, we may assume by Lemma \ref{lem:RegularInvar} that $Y$ is $F^{\Omega}_{\ell,\d{r}}$-invariant. Now, let
\[
L \ := \ \{\c\in\R(2\d{p}+\d{q})\st (\forall i) \ p\mid c_i\}.
\]
As a subset of $\Z^{\ell}$, $L$ inherits the well-order that was defined in Equation \eqref{eq:WO}. Using this well-order, list the elements of $L$ as $(\c_n)_{n<|L|}$, noting that $|L|\leq (2p+q)^\ell<r$. We construct by induction on $L$, for each $\c_n\in L$, an $\FeqOl{\d{r}-n\d{p}}$-invariant Borel subset $Z_n$ of $X$. As we construct these sets we write $Z_{<n}:=\cup_{m<n}Z_m$ and $Z_{\leq n}:=\cup_{m\leq n}Z_m$. For the inductive definition, fixing $\c_n\in L$ and supposing $Z_m$ has been defined for all $m<n$, define
\[
Z_{n} \ := \ \c_n\cdot Y \setminus B_{p}\big(Z_{<n}\big),
\]
so in particular $Z_0=Y$ since $\c_0=\bar{0}$ and $Z_{<0}=\emp$. We claim that the set $Z:=\bigcup_{\c_n\in L}Z_{n}$ has the desired properties.

First we show by induction that $Z_n$ is $\FeqOl{\d{r}-n\d{p}}$-invariant for each $\c_n\in L$. For the base case $\c_0=\bar{0}$, $Z_0=Y$ is $\FeqOl{\d{r}}$-invariant by assumption. Now fix $n>0$ and assume $Z_m$ is $\FeqOl{\d{r}-m\d{p}}$-invariant for all $m<n$. Since $|\c_n|\leq np$, the set $\c_n\cdot Y$ is $\FeqOl{\d{r}-n\d{p}}$-invariant by Lemma \ref{lem:ForwardInvar}(ii). And since $B_p(Z_m)$ is $\FeqOl{\d{r}-(m+1)\d{p}}$-invariant for each $m<n$ by Lemma \ref{lem:ForwardInvar}(iii), the set $B_p(Z_{<n})=\cup_{m<n}B_p(Z_m)$ is $\FeqOl{\d{r}-n\d{p}}$-invariant.  It follows that $Z_n$ is $\FeqOl{\d{r}-n\d{p}}$-invariant, completing the inductive argument. Therefore, since $r-|L|p\geq p$, every $Z_n$, along with $Z$ itself, is $\FeqOl{\d{p}}$-invariant.

Since $Y$ is $\FeqOl{\d{r}}$-invariant and $|\c_n|<2p+q<r$ for each $n$, each set $\c_n\cdot Y$ is Borel by Lemma \ref{lem:ForwardBorel}(i), and clopen in $X$ if $Y$ is clopen in $X$ by Lemma \ref{lem:ForwardBorel}(ii). Using this, along with the fact just proved that every $Z_n$ is $\FeqOl{\d{p}}$-invariant, an easy induction on $L$ using Lemma \ref{lem:ForwardBorel} again shows that each $Z_n$, and hence $Z$ itself, is Borel, and in fact clopen in $X$ if $Y$ is clopen in $X$. 

For $p$-sparsity, we argue by induction on $L$ that for all $\bar{c}_n\in L$, the set $Z_{\leq n}$ is $p$-sparse.  Since $Z=Z_{<|L|}$, this shows $Z$ is $p$-sparse as well. The base case $\bar{c}_0=\bar{0}$ holds by assumption, since $Z_0=Y$. For the inductive step, fix $n>0$ and assume $Z_{<n}$ is $p$-sparse. Note that $Z_{\leq n}=Z_{<n}\cup Z_{n}$. Since $Z_{n}\subseteq\c_n\cdot Y$ and $\c_n\cdot Y$ is $p$-sparse by Lemma \ref{lem:ForwardSparse}, $Z_{n}$ itself is $p$-sparse. So in order to show $Z_{\leq n}$ is $p$-sparse, it will suffice to show that for all $x\in Z_{n}$ and $z\in Z_{<n}$, if $d(x,z)<p$ then $x\mathrel{F^{\Omega}_{\ell,\d{p}}}z$. But in fact there do not exist $x\in Z_{n}$ and $z\in Z_{<n}$ such that $d(x,z)<p$, by definition of $Z_{n}$. So $Z_{\leq n}$ is indeed $p$-sparse, completing the induction and showing $Z$ is $p$-sparse as well. 

To complete the proof, fix $x\in X$; we must find $\bar{u}\in\N^{\ell}_{<3p}$ such that $\bar{u}\cdot x\in Z$. Fix some $x_0\in X$ such that $\dot{q}\cdot x_0=x$, along with $\bar{a}\in\N^{\ell}_{<q}$ such that $y:=\bar{a}\cdot x_0\in Y$. Let $\bar{b}:=\dot{q}-\bar{a}$, so 
$\bar{b}\cdot y=x$ and $|\bar{b}|\leq q$. Let $\bar{d}\in\R(\dot{p})$ be the unique tuple such that $p$ divides each component of $\bar{c}:=\bar{b}+\d{p}+\bar{d}$. Then $|\bar{c}|<2p+q$, so $\bar{c}\in L$, say $\c=\c_n$. Now we consider two cases. If $\c\cdot y\in Z_{n}$, then since $\c\cdot y=(\d{p}+\bar{d})\cdot x$, we can take $\bar{u}=\d{p}+\bar{d}\in\N^{\ell}_{<2p}$. 
Otherwise, we must have $\c\cdot y\in B_p(Z_{<n})$, so fix $z\in Z_{<n}$ and $\bar{s},\bar{t}\in\N^{\ell}_{<p}$ such that $\bar{s}\cdot (\c\cdot y)=\bar{t}\cdot z$. Then
\[
(\bar{b}+\bar{d}+2\d{p}-\bar{t}+\bar{s})\cdot y \ = \ (\d{p}-\bar{t}+\bar{s})\cdot (\bar{c}\cdot y) \ = \ (\d{p}-\bar{t})\cdot (\bar{t}\cdot z) \ = \ \d{p}\cdot z.
\]
Since $Z_{<n}$ is $\FeqOl{\d{p}}$-invariant, this means
\[
(\bar{d}+\d{p}-\bar{t}+\bar{s})\cdot x \ = \ (\bar{b}+\bar{d}+\d{p}-\bar{t}+\bar{s})\cdot y \ \in \ Z,
\]
so we can take $\bar{u}=\bar{d}+\d{p}-\bar{t}+\bar{s}\in\N^{\ell}_{<3p}$. This completes the proof.
\end{proof}

\begin{theorem}\label{thm:connection}
Let $\Ol$ be an $\ell$-dimensional product space, and $X$ an $\EOl$-invariant Borel subset of $\FOl$. Let $p,q,K\in\Z_+$. 
\begin{enumerate}[(i)]
 \item If there exists a Borel set $Y\subseteq X$ that is $(p,q)$-regular in $X$, then $\chi_B\big(\G_{\ell,p}\res X\big)\leq |\R(\d{q})|$.
 \item If there exists a relatively clopen set $Y\subseteq X$ that is $(p,q)$-regular in $X$, then $\chi_C\big(\G_{\ell,p}\res X\big)\leq |\R(\d{q})|$.
 \item If $\chi_B\big(\G_{\ell,p}\res X\big)\leq K$, then there exists a $(p,2Kp)$-regular Borel set $Y\subseteq X$.
 \item If $\chi_C\big(\G_{\ell,p}\res X\big)\leq K$, then there exists a $(p,2Kp)$-regular clopen set $Y\subseteq X$.
\end{enumerate}
\end{theorem}

\begin{proof}
For (i), fix $p,q\in\Z_+$ and suppose $Y$ is a Borel subset of $X$ that is $(p,q)$-regular in $X$. For each $x\in X$, let $\gamma(x)$ be the least $\bar{c}\in\R(\d{q})$ such that $\bar{c}\cdot x\in Y$. We claim that the Borel function $\gamma:X\to\R(\d{q})$ is a coloring of $\G_{\ell,p}\res X$. To see this, fix $x,y\in X$ such that $x\mathrel{\G_{\ell,p}}y$, along with $\bar{a},\bar{b}\in\N^{\ell}_{<p}$ such that $\bar{a}\cdot x=\bar{b}\cdot y$. Suppose for contradiction that $\gamma(x)=\gamma(y)$, and write $\bar{c}\in\R(\d{q})$ for their common value. Then $\bar{a}\cdot (\bar{c}\cdot x)=\bar{b}\cdot (\bar{c}\cdot y)$, so $d\big(\bar{c}\cdot x,\, \bar{c}\cdot y\big)<p$. Since $\bar{c}\cdot x$ and $\bar{c}\cdot y$ belong to $Y$ and $Y$ is $p$-sparse, this implies $\bar{c}\cdot x\mathrel{F^{\Omega}_{\ell,\d{p}}}\bar{c}\cdot y$. Thus $x\mathrel{F^{\Omega}_{\ell}}y$, contradicting $x\mathrel{\G_{\ell,p}}y$. So $\gamma$ is indeed a coloring of $\G_{\ell,p}\res X$, which shows $\chi_B\big(\G_{\ell,p}\res X\big)\leq |\R(\d{q})|$. Note that if $Y$ is clopen in $X$ then $\gamma$ is continuous, proving (ii) as well.

To prove (iii), we employ a greedy algorithm using an $\FeqOl{\d{r}}$-invariant coloring for a sufficiently large $r$, putting $x\in Y$ if and only if there is no point in $B_p(x)$ with a smaller color. For the details, suppose $c_0:X\to K$ is a Borel coloring of $\G_{\ell,p}\res X$, let $r:=Kp$, and define $c:X\to K$ by $c(x)=c_0(\d{r}\cdot x)$, so $c$ is an $\FeqOl{\d{r}}$-invariant Borel $K$-coloring of $\G_{\ell,p}\res X$. Now, writing $Y_{<k}:=\cup_{i<k}Y_i$ and $Y_{\leq k}:=\cup_{i\leq k}Y_i$, inductively define the sequence $(Y_k)_{k\in K}$ of Borel subsets of $X$ by
\[
Y_k \ := \ c^{-1}[\{k\}]\setminus B_p(Y_{<k}),
\]
so $Y_0:=c^{-1}\left[\{0\}\right]$. An easy induction using Lemmas \ref{lem:ForwardBorel} and \ref{lem:ForwardInvar}(iii) shows that for each $k<K$, $Y_k$ is $F^{\Omega}_{\ell,{\d{r}-k\d{p}}}$-invariant and Borel, and is clopen in $X$ if $c_0$ is continuous. We claim that the $F^{\Omega}_{\ell,\d{p}}$-invariant Borel set $Y:=\bigcup_{k\in K}Y_k$ is $(p,2Kp)$-regular in $X$.

To prove this, we first show by induction on $k\in K$ that each set $Y_{\leq k}$ is $p$-sparse. That the set $Y_{\leq 0}=c^{-1}[\{0\}]$ is $p$-sparse follows immediately from the fact that $c$ colors $\G_{\ell,p}\res X$. For the induction, fix $0<k<K$ and assume $Y_{<k}$ is $p$-sparse. Since $c$ colors $\G_{\ell,p}\res X$ and $Y_{k}\subseteq c^{-1}[\{k\}]$, the set $Y_{k}$ is also $p$-sparse, so to show $Y_{\leq k}$ is $p$-sparse it will suffice to show that for all $x\in Y_{k}$ and $y\in Y_{<k}$, if $d(x,y)<p$ then $x\mathrel{F^{\Omega}_{\ell,\d{p}}}y$. But there do not exist $x\in Y_k$ and $y\in Y_{<k}$ such that $d(x,y)<p$, by the definition of $Y_k$. Thus each $Y_{\leq k}$ is indeed $p$-sparse, and it follows that $Y$ itself is $p$-sparse.

To complete the proof that $Y$ is $(p,2Kp)$-regular in $X$, let $x\in X$ be arbitrary; we must find $\bar{c}\in\N^{\ell}_{<2Kp}$ such that $\bar{c}\cdot x\in Y$. We inductively construct an infinite sequence $(\bar{c}_k)_{k\in\omega}$ in $\N^{\ell}_{<2Kp}$, yielding points $x_k:=\bar{c}_k\cdot x\in X$, having the following properies for all $k\in\omega$:
\begin{enumerate}
 \item[(i)] $|\bar{c}_k-K\d{p}|\leq kp$; and
 \item[(ii)] if $k\geq 1$ and $x_{k-1}\not\in Y$, then $\bar{c}_{k}\ne \bar{c}_{k-1}$ and $c(x_k)<c(x_{k-1})$.
\end{enumerate}
Given such a sequence, there must be some $k<K$ such that $x_k\in Y$, since otherwise $\big(c(x_k)\big)_{k\leq K}$ is a strictly decreasing sequence in $K$ of length greater than $K$. Letting $k_0<K$ be least such that $x_{k_0}\in Y$, we then have $|\bar{c}_{k_0}-K\d{p}|\leq k_0p<Kp$, so $\bar{c}_{k_0}\in\N^{\ell}_{<2Kp}$ and $\bar{c}_{k_0}\cdot x\in Y$, as desired. To construct the sequence, we start with $\bar{c}_0:=K\d{p}$, so (i) and (ii) hold for $k=0$. Then, continuing inductively, fix $k>0$ and suppose $\bar{c}_0,\ldots,\bar{c}_{k-1}\in\N^{\ell}_{<2Kp}$ have beed defined and satisfy (i) and (ii). 
If $x_{k-1}\in Y$, then we let $\bar{c}_{k}:=\bar{c}_{k-1}$, so $x_k=x_{k-1}$ and (i) and (ii) continue to hold. Otherwise, write $j:=c(x_{k-1})$ and note that since $x_{k-1}\not\in Y_j$, we must have $k\geq 2$ and $x_{k-1}\in B_p(Y_{<j})$. So fix $y\in Y_{<j}$ such that $d(x_{k-1},y)<p$, and put $\bar{c}_k:=\bar{c}_{k-1}+(y-x_{k-1})$. Note that $y-x_{k-1}\ne\bar{0}$ since $c(x_{k-1})=j$ and $Y_{<j}\subseteq c^{-1}[\{0,\ldots,j-1\}]$ is $F^{\Omega}_{\ell,\d{p}}$-invariant. Thus $\bar{c}_k\ne\bar{c}_{k-1}$. Furthermore, $|\bar{c}_{k}-K\d{p}| \leq |\bar{c}_{k-1}-K\d{p}|+p<kp$, and
\[
c(x_k) \ = \ c(y) \ < \ j \ = \ c(x_{k-1}),
\]
so again (i) and (ii) continue to hold. This completes the inductive construction, and the proof of (iii).

Finally, note that in our proof of (iii), if $c_0$ is continuous then so is $c$ and the sets $Y_k$ and $Y$ itself are clopen in $X$, so (iv) follows immediately from our proof of (iii).
\end{proof}

\begin{corollary}\label{cor:main}
Let $\Ol$ be an $\ell$-dimensional product space, and $X$ an $\EOl$-invariant Borel subset of $\FOl$. If $\chi_B\big(\G_{\ell,n}\res X\big)$ is finite for every $n\in\Z_+$, then $\EOl\res X$ is hyperselective.
\end{corollary}

\begin{proof}
Assume the hypotheses, and in particular assume $\chi_B\big(\G_{\ell,n}\res X\big)$ is finite for each $n$. Then for every $n\in\Z_+$, there exists a Borel subset $Y_n\subseteq X$ that is $(n,2Kn)$-regular in $X$ by Thereom \ref{thm:connection}, so in fact there exists a Borel subset $Z_n\subseteq X$ that is $(n,3n)$-regular in $X$ by Lemma \ref{lem:SparseBoundedlyDense}. It now follows from Theorem \ref{thm:RegularHS} (with $c=3$) that $\EOl\res X$ is hyperselective.
\end{proof}

Given an $\EOl$-invariant Borel set $X\subseteq\FOl$, we have reduced the problem of showing $\EOl\res X$ is hypersmooth to that of obtaining finite Borel colorings of the graphs $\G_{\ell,n}(\Omega)\res X$. In particular, if $\chi_B\big(\G^{\textup{Fr}}_{\ell,n}(\Omega)\big)$ is finite for each $n\in\Z_+$, then $\EOl\res\FOl$ is hypersmooth. Indeed, $\chi_B\big(\G^{\textup{Fr}}_{1,n}(2)\big)$ \emph{is} finite for each $n$, and $E_{1}^{2}$ (that is, the tail equivalence relation $E_t$) is hypersmooth. Unfortunately, it turns out that $\chi_B\big(\G^{\textup{Fr}}_{\ell,n}(2)\big)=\aleph_0$ for all $\ell\geq 2$ and $n>1$. Indeed, the argument below shows that when $\ell=2$, the graph $\G_{\bar{e}_1}(2)$ has infinite Borel chromatic number, and with a bit of care it can be adapted to show the same for $\G^{\textup{Fr}}_{\bar{e}_1}(2)$ and all the graphs $\G_{\m}$ that we consider in Sections \ref{sec:ep}--\ref{sec:Gm}.

\begin{observation}\label{claim:NoColoring}
For every base space $\Omega$ and integer $\ell\geq 2$, we have $\chi_B\big(\vec{\G}_{\ell}(\Omega)\big)=\aleph_0$.
\end{observation}

\begin{proof}
Let $\Omega$ and $\ell\geq 2$ be given. 
We already know $\chi_B\big(\vec{\G}_{\ell}(\Omega)\big)\leq\aleph_0$ by Fact \ref{fact:countable}. Therefore, since $\vec{\G}_{\bar{e}_1}(2)$ embeds into $\vec{\G}_{\ell}(\Omega)$, it will suffice to show $\chi_B\big(\vec{\G}_{\bar{e}_1}(2)\big)\geq\aleph_0$. 
As in Example \ref{ex:BadAssGraph}, let $[\omega]^\omega$ be the Polish space of all strictly increasing sequences of natural numbers, let $\sigma$ be the shift map on $[\omega]^\omega$, and let $\vec{\GG}$ be the directed graph generated by $\sigma$, so $\chi_B(\vec{\GG})=\aleph_0$. Define the Borel map $f:[\omega]^\omega\to 2^{\omega^2}$ by 
\[
f(x)(m,n) \ = \ \begin{cases} 1 & \text{if $x(m)=n$;} \\ 0 & \text{otherwise,} \end{cases}
\]
so $f(x)$ is the characteristic function of the graph of the function $x:\omega\to\omega$. 
Note that $\bar{e}_1\cdot f(x)=f\big(\sigma(x)\big)$ for all $x\in [\omega]^{\omega}$. Therefore, if $c:2^{\omega^2}\to Z$ is a Borel coloring of $\vec{\G}_{\bar{e}_1}(2)$, then for all $x\in [\omega]^\omega$ we have $c\big(f(x)\big)\ne c\big(f(\sigma(x))\big)$, so $c\circ f$ is a Borel coloring of $\vec{\GG}$. It follows that $Z$ is infinite, showing $\chi_B\big(\vec{\G}_{\bar{e}_1}(2)\big)\geq\aleph_0$. 
\end{proof}

Since $\G_{\ell}^{\textup{Fr}}$ embeds into $\G_{\ell,p}^{\textup{Fr}}$ whenever $p\geq 2$, it follows that $\aleph_0\leq\chi_B\big(\G_{\ell}^{\textup{Fr}}\big)\leq\chi_B\big(\G_{\ell,p}^{\textup{Fr}}\big)$ for all $\ell\geq 2$, so we immediately obtain the following from Theorem \ref{thm:connection}(i):

\begin{corollary}\label{cor:NoNiceMarkerSets}
Let $\Ol$ be an arbitrary product space with $\ell\geq 2$, and let $p,q\in\Z_{\geq 2}$ be arbitrary. There does not exist a Borel set $Y\subseteq\FOl$ that is $(p,q)$-regular in $\FOl$. \hfill $\square$
\end{corollary}


Despite Observation \ref{claim:NoColoring} and Corollary \ref{cor:NoNiceMarkerSets}, finite Borel graph colorings will play a crucial role in our proof that $\ETT$ is hypersmooth. 
Indeed, we will see that even though $\chi_B\big(\G_2(2))\big)=\aleph_0$, the graphs $\G_{2,n}(2)$ have finite Borel chromatic numbers on suitably ``large'' subsets of $\FTT$. Specifically, as the main technical result of this paper we will prove the following:

\begin{theorem}\label{thm:main}
For every $n\in\Z_+$, there is a partition of $\FTT$ into a pair of $\ETT$-invariant Borel sets $\FTT=X_n\sqcup W_n$ such that $\ETT\res W_n$ is hypersmooth and $\G_{2,n}\res X_n$ admits a Borel $K_n$-coloring, where $K_n=3^{n^2}\cdot 4^{n^4}$.
\end{theorem}

\begin{corollary}\label{cor:final}
$\ETT\res\FTT$ is hyperselective.
\end{corollary}

\begin{proof}[Proof of \ref{cor:final} from \ref{thm:main}]
For each $n\in\Z_+$, fix sets $X_n$ and $W_n$ as given in Theorem \ref{thm:main}, and let $X:=\bigcap_nX_n$ and $W:=\bigcup_nW_n$, so $X$ and $W$ are $\ETT$-invariant Borel sets such that $\FTT=X\sqcup W$. Then for each $n\in\Z_+$ we have $\chi_B\big(\G_{2,n}\res X\big)\leq 3^{n^2}\cdot 4^{n^4}$ since $X\subseteq X_n$, so $\ETT\res X$ is hyperselective by Corollary \ref{cor:main}. Furthermore, $\ETT\res W$ is hypersmooth by Lemma \ref{lem:piecewiseHS}, hence hyperselective by Corollary \ref{cor:fi2}, so $\ETT\res\FTT$ itself is hyperselective by another application of Lemma \ref{lem:piecewiseHS}.
\end{proof}

Before presenting a rather lengthy and technical proof of Theorem \ref{thm:main} in Sections \ref{sec:ep}--\ref{sec:Gm}, we turn our attention in Section \ref{sec:periodicity} to the periodic part of $\ETT$, and show that $\ETT\res\PTT$ is also hyperselective.

\section{Periodicity and symmetry}\label{sec:periodicity}

Most recent proofs that Borel actions of a given class of countable groups generate hyperfinite orbit equivalence relations focus on the free part of the action as the main case, while either omitting the ``non-free'' part or treating it as a (sometimes quite tedious) secondary case (see for instance \cite{GJ}, \cite{Nilpotent}, and \cite{AsymDim}).  In this section, we will address the non-free, or periodic, part of the $\N^{2}$-action on $\OT$, both because we will be able to show without too much difficulty that $\ETT\res\PTT$ is hypersmooth, but also because we will use some facts about periodicity in our proof of Theorem \ref{thm:main}. We begin with some general definitions and facts in the higher dimensional setting.

\begin{definition}
Let $\Ol$ be a product space, let $x\in\Ol$, and let $\m\in\Z^\ell$ and $H\leq\Z^\ell$.
\begin{enumerate}[(i)]
 \item $x$ is \emph{$\m$-periodic} if there exist $\bar{a},\bar{b}\in\N^{\ell}$ such that $\bar{a}\cdot x=\bar{b}\cdot x$ and $\bar{a}-\bar{b}=\m$.
 \item $x$ is \emph{$H$-periodic} if $x$ is $\m$-periodic for all $\m\in H$.
 \item $x$ is \emph{maximally $H$-periodic} if $x$ is $H$-periodic and for any subgroup $G\leq\Z^{\ell}$ such that $x$ is $G$-periodic, we have $G\leq H$.
 \item $x$ is \emph{periodic} if $x$ is $\m$-periodic for some nonzero $\m\in\Z^{\ell}$.
\end{enumerate}
\end{definition}

We have already been writing $\POl$ for the set of all periodic elements of $\Ol$, and $\FOl$ for the set of non-periodic elements. For each $\m\in\Z^{\ell}$, let us write
\[
\PmOl \ := \ \big\{x\in\Ol\st \text{$x$ is $\m$-periodic}\big\} \qquad \text{and} \qquad \neg\PmOl \ := \ \Ol\setminus\PmOl.
\]
Note that every $x\in\Ol$ is $\bar{0}$-periodic, so $\textup{Per}_{\bar{0}}(\Ol)=\Ol$.

\begin{claim}\label{claim:PerInvar}
 For every product space $\Ol$ and $\m\in\Z^{\ell}$, the set $\PmOl$ is an $\EOl$-invariant Borel subset of $\Ol$.
\end{claim}

\begin{proof}
Fix $\m\in\Z^{\ell}$. For each $\bar{a},\bar{b}\in\N^{\ell}$, the set of all $x\in\Ol$ such that $\bar{a}\cdot x=\bar{b}\cdot x$ is closed in $\Ol$ by continuity of the $\N^{\ell}$-action, so the set
\[
\PmOl \ = \ \bigcup_{\{\bar{a},\bar{b}\,\in\,\N^{\ell}\st \bar{a}-\bar{b}\,=\,\m\}} \big\{x\in\Ol\st \bar{a}\cdot x=\bar{b}\cdot x\big\}
\]
is $F_{\sigma}$. For $\EOl$-invariance, fix $x\in\PmOl$ and $\bar{a},\bar{b}\in\N^{\ell}$ such that $\bar{a}\cdot x=\bar{b}\cdot x$ and $\bar{a}-\bar{b}=\m$. Let $y\in\Ol$, and suppose $\bar{c},\bar{d}\in\N^{\ell}$ are such that $\bar{c}\cdot x=\bar{d}\cdot y$. Then $(\bar{a}+\bar{d})-(\bar{b}+\bar{d})=\bar{m}$ and
\[
(\bar{a}+\bar{d})\cdot y \ = \ \bar{a}\cdot (\bar{c}\cdot x) \ = \ \bar{c}\cdot(\bar{b}\cdot x) \ = \ \cdot y) \ = \ (\bar{b}+\bar{d})\cdot y,
\]
so $y$ is $\m$-periodic as well.
\end{proof}

\begin{definition}
Let $\Ol$ be a product space, let $x\in\Ol$, and let $\m\in\Z^\ell$ and $H\leq\Z^\ell$.
\begin{enumerate}[(i)]
 \item $x$ is \emph{strongly $\m$-periodic} if $x(\bar{i})=x(\bar{j})$ for all $\bar{i},\bar{j}\in\N^{\ell}$ such that $\bar{i}-\bar{j}=\m$.
 \item $x$ is \emph{strongly $H$-periodic} if $x$ is strongly $\m$-periodic for all $\m\in H$.
 \item $x$ is \emph{strongly-maximally $H$-periodic} if $x$ is both strongly $H$-periodic and maximally $H$-periodic.
\end{enumerate}
\end{definition}

The following lemma establishes the basic properties of the various types of periodicity and their relationships to each other. Most of these facts are what one would expect, although strong $\langle\s,\t\rangle$-periodicity need not follow from strong $\s$-periodicity and strong $\t$-periodicity. Parts (vi) and (vii) of the lemma address this subtlety.

\begin{lemma}\label{lem:per}
Let $\Ol$ be a product space, let $x\in\Ol$, and let $\m\in\Z^\ell$ and $H\leq\Z^\ell$.
\begin{enumerate}[(i)]
 \item If $x$ is strongly $\bar{m}$-periodic, then $x$ is $\bar{m}$-periodic.
 \item If $x$ is $\bar{m}$-periodic, then so is $\bar{c}\cdot x$ for any $\bar{c}\in\mathbb N^\ell$; the same holds for strong $\bar{m}$-periodicity, $H$-periodicity, strong $H$-periodicity, maximal $H$-periodicity, and strongly-maximal $H$-periodicity.
 \item If $x$ is $\bar{m}$-periodic, then there is $\bar{c}\in\mathbb N^\ell$ such that $\bar{c}\cdot x$ is strongly $\bar{m}$-periodic.
 \item The set of $\bar{m}\in\mathbb Z^\ell$ for which $x$ is $\bar{m}$-periodic is a subgroup of $\mathbb Z^\ell$.
 \item If $x$ is strongly $\bar{m}$-periodic, then $x$ is strongly $(c\bar{m})$-periodic for all $c\in\mathbb Z$. 
 \item If there is some $\bar{m}\in\Z_+^\ell$ such that $x$ is strongly $\bar{m}$-periodic, then the set of all $\bar{m}\in\mathbb Z^\ell$ for which $x$ is strongly $\bar{m}$-periodic is a subgroup of $\mathbb Z^\ell$.
 \item Suppose $H=\langle\bar{m}_1,\ldots,\bar{m}_n\rangle$. If $x$ is strongly $\bar{m}_j$-periodic for each $1\leq j\leq n$, then there is $\bar{c}\in\mathbb N^\ell$ with $|\bar{c}|\leq n\cdot\max\{|\bar{m}_1|,\ldots,|\bar{m}_n|\}$ such that $\bar{c}\cdot x$ is strongly $H$-periodic.
 \item If $x$ is $H$-periodic, then there is $\bar{c}\in\mathbb N^\ell$ such that $\bar{c}\cdot x$ is strongly $H$-periodic.
 \item If $x$ is maximally $H$-periodic, then there is $\bar{c}\in\N^{\ell}$ such that $\bar{c}\cdot x$ is strongly-maximally $H$-periodic.
\end{enumerate}
\end{lemma}

\begin{proof} For (i), if $x(\bar{i})=x(\bar{j})$ for all $\bar{i},\bar{j}\in\N^\ell$ such that $\bar{i}-\bar{j}=\bar{m}$, then for any particular $\bar{a},\bar{b}\in\N^\ell$ such that $\bar{a}-\bar{b}=\bar{m}$, we have $\bar{a}\cdot x=\bar{b}\cdot x$. 

For (ii), if $\bar{a}\cdot x=\bar{b}\cdot x$ with $\bar{a}-\bar{b}=\m$, then also $\bar{a}\cdot (\bar{c}\cdot x)=\bar{b}\cdot (\bar{c}\cdot x)$ with $(\bar{a}+\c)-(\bar{b}+\c)=\m$, so $\bar{c}\cdot x$ is $\m$-periodic if $x$ is. And if $x$ is strongly $\m$-periodic with $\bar{i},\bar{j}\in\N^{\ell}$ arbitrary such that $\bar{i}-\bar{j}=\m$, then $(\bar{c}\cdot x)(\bar{i})=x(\bar{c}+\bar{i})=x(\bar{c}+\bar{j})=(\bar{c}\cdot x)(\bar{j})$ since $(\bar{c}+\bar{i})-(\bar{c}+\bar{j})=\bar{m}$. The remaining claims follow easily.

For (iii), suppose $x$ is $\bar{m}$-periodic, and fix $\bar{a},\bar{b}\in\N^\ell$ such that $\bar{a}\cdot x=\bar{b}\cdot x$ and $\bar{a}-\bar{b}=\bar{m}$. Let $c_i=\max(a_i,b_i)$ for each $1\leq i\leq\ell$. We show that $\bar{c}\cdot x$ is strongly $\bar{m}$-periodic. Let $\bar{i},\bar{j}\in\N^\ell$ be such that $\bar{i}-\bar{j}=\bar{m}$. Then
\[
(\bar{c}\cdot x)(\bar{i}) \ = \ (\bar{c}\cdot x)(\bar{j}+\bar{m}) \ = \ (\bar{c}\cdot x)(\bar{j}+\bar{a}-\bar{b}) \ = \ (\bar{a}\cdot x)(\bar{j}+\bar{c}-\bar{b}) \ = \ (\bar{b}\cdot x)(\bar{j}+\bar{c}-\bar{b}) \ = \ (\bar{c}\cdot x)(\bar{j}).
\]

For (iv), it is immediate from the definition that if $x$ is $\bar{m}$-periodic then $x$ is also $(-\bar{m})$-periodic. Suppose now for $\s,\t\in\mathbb Z^\ell$  that $x$ is both $\s$-periodic and $\t$-periodic, and fix $\bar{a},\bar{b},\bar{c},\bar{d}\in\mathbb N^\ell$ such that $\bar{a}\cdot x=\bar{b}\cdot x$ and $\bar{c}\cdot x=\bar{d}\cdot x$ with $\s=\bar{a}-\bar{b}$ and $\t=\bar{c}-\bar{d}$. Then
\[
(\bar{a}+\bar{c})\cdot x \ = \ \bar{c}\cdot (\bar{a}\cdot x) \ = \ \bar{c}\cdot (\bar{b}\cdot x) \ = \ \bar{b}\cdot (\bar{c}\cdot x) \ = \ \bar{b}\cdot (\bar{d}\cdot x) \ = \ (\bar{b}+\bar{d})\cdot x,
\]
where $\s+\t=(\bar{a}+\bar{c})-(\bar{b}+\bar{d})$, so $x$ is $(\s+\t)$-periodic.

For (v), fix $\bar{m}\in\mathbb Z^\ell$ and suppose $x$ is strongly $\bar{m}$-periodic. It is immediate from the definition that $x$ is strongly $(-c\m)$-periodic if $x$ is strongly $c\m$-periodic, and $x$ is strongly $\bar{0}$-periodic trivially, so it will suffice to show $x$ is strongly $c\m$-periodic for all $c\in\Z_+$. We show this by induction on $c$, the base case $c=1$ being our assumption. Fix $c\geq 1$ and assume $x$ is $c\m$-periodic, and let $\bar{i},\bar{j}\in\N^{\ell}$ be arbitrary such that $\bar{i}-\bar{j}=(c+1)\m$, so $(\bar{i}-\m)-\bar{j}=c\m$. Note that $\bar{i}-\m\in\N^{\ell}$ by convexity of $\N^{\ell}$. Therefore $x(\bar{i}) = x(\bar{i}-\m)$ since $x$ is strongly $\m$-periodic, and $x(\bar{i}-\m)=x(\bar{j})$ since $x$ is strongly $c\m$-periodic, so $x(\bar{i})=x(\bar{j})$ to complete the induction and the proof.

For (vi), fix $\bar{m}\in \Z_+^\ell$ such that $x$ is strongly $\bar{m}$-periodic, and let $\s,\t\in\Z^\ell$ be arbitrary such that $x$ is strongly $\s$-periodic and strongly $\t$-periodic. By (v), it will suffice to show $x$ is strongly $(\s+\t)$-periodic. Let $\bar{i},\bar{j}\in\N^\ell$ satisfy $\bar{i}-\bar{j}=\s+\t$, and fix $c\in\N$ such that $\bar{j}+\s+c\bar{m}\in\N^\ell$. Note that $\bar{i}+c\m\in\N^{\ell}$ and $\bar{j}+c\m\in\N^{\ell}$ since $\m\in\Z_+^\ell$, and also that
\(
(\bar{i}+c\bar{m})-(\bar{j}+\s+c\bar{m}) = \t.
\)
Thus
\[
x(\bar{i}+c\bar{m}) \ = \ x(\bar{j}+\s+c\bar{m}) \ = \ x(\bar{j}+c\bar{m})
\]
since $x$ is strongly $\t$-periodic and strongly $\s$-periodic, and since $x$ is strongly $(c\bar{m})$-periodic by (v) it follows that
\[
x(\bar{i}) \ = \ x(\bar{i}+c\bar{m}) \ = \ x(\bar{j}+c\bar{m}) \ = \ x(\bar{j}),
\]
completing the proof.

Proving (vii) takes a bit more work than the others. For any $n\in\Z_+$, write $\|\cdot\|$ for the Euclidean norm and $\|\cdot\|_{\infty}$ for the sup norm on $\mathbb{R}^{n}$, along with $d_E$ and $d_{\infty}$, respectively, for the corresponding metrics. Call any subset of $\mathbb{R}^{n}$ of the form $\{\bar{a}+\bar{i}\st \|\bar{i}\|_{\infty}\leq 1\}$ where $\bar{a}\in\Z^{n}$ a \emph{unit cell}. Let $\Gamma_n$ be the grid graph on $\Z^{n}$, given by $\bar{a}\mathrel{\Gamma_n}\bar{b} \ :\Leftrightarrow \ d_E(\bar{a},\bar{b})\leq 1$.  Our proof is based on the fact that every line segment in $\mathbb{R}^{n}$ can be approximated by a path in $\Gamma_n$. Such approximations 
have been well studied in the computer science literature (see for instance \cite{DigitalStraightness}), and there are many ways of constructing them; for us, the following simplistic argument will suffice. Suppose $\bar{a}$ and $\bar{b}$ are two fixed points in $\Z^{n}$. Let $L$ be the Euclidean line segment between them, and let $P$ be the intersection of $\Z^{n}$ with the union of all unit cells that meet $L$, so $P:=\{\bar{c}\in\Z^{n}\st d_{\infty}(\bar{c},L)\leq 1\}$. Then $\bar{a},\bar{b}\in P$ and $\Gamma_n\res P$ is connected, so there is a path $(\bar{v}_0,\ldots,\bar{v}_s)$ from $\bar{a}$ to $\bar{b}$ in $\Gamma_n$ such that $d_{\infty}(\bar{v}_i,L)\leq 1$ for each $0\leq i\leq s$.

Now, suppose $H=\langle\m_1,\ldots,\m_n\rangle$, where for each $1\leq k\leq n$, $x$ is strongly $\m_k$-periodic and thus also strongly $(-\m_k)$-periodic by (v). Let $S:=\{\pm\bar{m}_1,\ldots,\pm\bar{m}_n\}$, and put $c:=n\cdot\max\{|\m|\st \m\in S\}$. We claim that $\d{c}\cdot x$ is strongly $H$-periodic. To see this, let $\bar{h}\in H$ be arbitrary and fix $\bar{i},\bar{j}\in\N^{\ell}$ such that $\bar{j}-\bar{i}=\bar{h}$; we must show $(\d{c}\cdot x)(\bar{i})=(\d{c}\cdot x)(\bar{j})$. 
View the shifted lattice $\bar{i}+H$ in $\Z^{\ell}$ as an undirected graph where two vertices are connected provided they differ by one of the generators in $S$. 
Fix $\bar{c}\in\Z^n$ such that $\bar{h}=\sum_{i=1}^nc_i\bar{m}_i$, and let $T:\mathbb{R}^n\to\mathbb{R}^{\ell}$ be the affine linear map defined by 
\[
T(\bar{x}) \ := \ \bar{i}+\sum_{i=1}^n x_i\bar{m}_i,
\]
so $T(\bar{0})=\bar{i}$ and $T(\bar{c})=\bar{j}$. A straightforward calculation shows that for all $\bar{x},\bar{y}\in\mathbb{R}^n$ we have
\[
d_{\infty}\big(T(\bar{x}),T(\bar{y})\big) \ \leq \ c\cdot d_{\infty}(\bar{x},\bar{y}).
\]
Let $L$ be the Euclidean line segment from $\bar{0}$ to $\bar{c}$ in $\mathbb{R}^n$, and let $(\bar{v}_0,\ldots,\bar{v}_s)$ be a path from $\bar{0}$ to $\bar{c}$ in $\Gamma_n$ such that $d_{\infty}(\bar{v}_i,L)\leq 1$ for each $0\leq i\leq s$. For each $0\leq i\leq s$, let $\bar{w}_i:=T(\bar{v}_i)$. Then $L':=T[L]$ is the Euclidean line segment from $\bar{i}$ to $\bar{j}$ in $\mathbb{R}^{\ell}$, and $(\bar{w}_0,\ldots,\bar{w}_s)$ is a path in $\bar{i}+H$ from $\bar{i}$ to $\bar{j}$ such that $d_{\infty}(\bar{w}_i,L')\leq c$ for each $0\leq i\leq s$. Since $L'\subseteq\mathbb{R}_{\geq 0}^{\ell}$, this implies that $\d{c}+\bar{w}_i\in\N^{\ell}$ for each $0\leq i\leq s$. Therefore, since $x$ is strongly $\m$-periodic for each $\m\in S$, we have
\[
(\d{c}\cdot x)(\bar{w}_i) \ = \ x(\d{c}+\bar{w}_i) \ = \ x(\d{c}+\bar{w}_{i+1}) \ = \ (\d{c}\cdot x)(\bar{w}_{i+1})
\]
for each $0\leq i<n$, so $(\d{c}\cdot x)(\bar{i})=(\d{c}\cdot x)(\bar{j})$ by induction, completing the proof of (vii).

For (viii), suppose $x$ is $H$-periodic, and fix generators $\m_1,\ldots,\m_r$ for $H$. By (ii) and (iii), there is $\bar{a}\in\mathbb N^\ell$ such that $\bar{a}\cdot x$ is strongly $\bar{m}_k$-periodic for each $1\leq k\leq r$, and then by (vii) there is $\bar{b}\in\mathbb N^\ell$ such that $\bar{b}\cdot (\bar{a}\cdot x)$ is strongly $H$-periodic.

Finally, for (ix), if $x$ is maximally $H$-periodic, then by (viii) there is $\c\in\N^{\ell}$ such that $\c\cdot x$ is strongly $H$-periodic, but also maximally $H$-periodic by (ii).
\end{proof}

For each $x\in\Ol$, we write $H_x$ for the unique subgroup $H$ of $\Z^{\ell}$ for which $x$ is maximally $H$-periodic, where the existence and uniqueness of such a subgroup follows immediately from the definitions and the fact that $\Z^{\ell}$ is Noetherian. Let us also write, for each $H\leq\Z^{\ell}$, 
\begin{equation}
\textup{Per}_{H}(\Ol) \ := \ \Big\{x\in\POl\st H_x=H\Big\}.
\end{equation}
Note that $x\in\FOl$ if and only if $H_x$ is trivial.

\begin{lemma}\label{lem:HBorel}
For every product space $\Ol$ and subgroup $H\leq\Z^{\ell}$, the set $\textup{Per}_{H}(\Ol)$ is Borel and $\EOl$-invariant.
\end{lemma}

\begin{proof}
Since
\[
\textup{Per}_H(\Ol) \ = \ \left(\bigcap_{\m\mkern1mu\in H}\PmOl\right) \setminus \left(\bigcup_{\m\mkern1mu\in\mkern2mu\Z^{\ell}\setminus H}\PmOl\right),
\]
$\textup{Per}_H(\Ol)$ is Borel and $\EOl$-invariant by Claim \ref{claim:PerInvar}.
\end{proof}

Before showing $\ETT\res\PTT$ is hypersmooth, 
let us examine two additional notions of symmetry that will be relevant in Sections \ref{sec:ep}--\ref{sec:Gm}. First, we consider an analogue of periodicity for ``rectangular patches," that is, for maps of the form $f:\R_{\ell}(\c)\to\Omega$ where $\Ol$ is a product space and $\c\in\N^{\ell}$. For each $\ell\in\Z_+$, write
\begin{equation}
\R_{\ell} \ := \ \bigcup_{\bar{c}\mkern1mu\in\mkern1mu\N^{\ell}}\R_{\ell}(\bar{c}),
\end{equation}
where $\R_{\ell}(\c):=\{\bar{a}\in\N^{\ell}\st \bar{a}\prec\c\}$ as in \eqref{eq:Rec}. Given $s,t\in\R_{\ell}$ and $\bar{c}\in\N^{\ell}$, let $\bar{c}\cdot s$ be the function defined by $(\bar{c}\cdot s)(\bar{i})=s(\bar{c}+\bar{i})$ for all $\bar{i}\in\omega^\ld$ such that $\bar{c}+\bar{i}\in\dom(s)$, and write $s\simeq t$ to mean $s(\bar{i})=t(\bar{i})$ for all $\bar{i}\in\dom(s)\cap\dom(t)$. Since the analogy between periodicity in $\Ol$ and in $\Omega^{\R_{\ell}}$ can be misleading, we use different terminology: given $s\in\Omega^{\R_{\ell}}$ and $\m\in\Z^{\ell}$, we call $s$ \emph{$\m$-symmetric} if $s(\bar{i})=s(\bar{j})$ for all $\bar{i},\bar{j}\in\dom(s)$ such that $\bar{i}-\bar{j}=\m$. This notion will be used in Section \ref{sec:SubSec3}.

Second, in the two-dimensional case we will frequently exploit the symmetry between the two coordinates. For each $\bar{a}=(a_1,a_2)\in\N^2$ and base space $\Omega$, write $\bar{a}^{\delta}:=(a_2,a_1)$ and define the coordinate-swapping involution $\delta_{\Omega}:\OT\to\OT$ by
\begin{equation}\label{eq:delta}
\delta_{\Omega}(x)(\bar{i}) \ := \ x(\bar{i}^\delta).
\end{equation}
We will always omit the base space from the notation here, and just write $\delta$. Basic properties of this map include the following.

\begin{lemma}\label{lem:delta}
Let $\OT$ be a two-dimensional product space, and write $\q:=(1,-1)$.
 \begin{enumerate}[(i)]
  \item The map $\delta:\OT\to\OT$ is a homeomorphism.
  \item For all $\bar{a}\in\N^2$ and $x\in\OT$, we have $\delta(\bar{a}\cdot x)=\bar{a}^{\delta}\cdot\delta(x)$. In particular, $\delta(\bar{e}_1\cdot x)=\bar{e}_2\cdot\delta(x)$ and $\delta(\bar{e}_2\cdot x)=\bar{e}_1\cdot\delta(x)$.
  \item $\delta$ is an isomorphism from $\EOT$ to itself.
  \item For all $x,y\in\OT$, we have $x\mathrel{\Gea}y$ if and only if $\delta(x)\mathrel{\Geea}\delta(y)$, and $x\mathrel{\Gsea}y$ if and only if $\delta(y)\mathrel{\Gsea}\delta(x)$.
 \end{enumerate}
\end{lemma}

\begin{proof}
Let $(V_i)_{i\in I}$ be a fixed countable base for the topology on $\Omega$, let $k\in\Z_+$ and $f:\R(\d{k})\to I$ be arbitrary, and consider the basic open set $U_f\subseteq\OT$ as defined in Equation \eqref{eq:BasicOpen}. Define $g:\R(\d{k})\to I$ by $g(\bar{a})=f(\bar{a}^\delta)$. 
Then for all $x\in\OT$, we have
\begin{align*}
\delta(x)\in U_f \quad &\Longleftrightarrow \quad \big(\forall\bar{a}\in\R(\d{k})\big) \, \delta(x)(\bar{a})\in V_{f(\bar{a})} \\
&\Longleftrightarrow \quad \big(\forall\bar{a}\in\R(\d{k})\big) \, x(\bar{a}^\delta)\in V_{g(\bar{a}^\delta)} \\
&\Longleftrightarrow \quad \big(\forall\bar{a}\in\R(\d{k})\big) \, x(\bar{a})\in V_{g(\bar{a})} \quad \Longleftrightarrow \quad x\in U_g.
\end{align*}
Thus $\delta^{-1}[U_f]=U_g$, showing $\delta$ is continuous. Since $\delta^{-1}=\delta$, it follows that $\delta$ is a homeomorphism.

For (ii), let $\bar{a}=(a_1,a_2)\in\N^2$ and $x\in\OT$. Then for all $i,j\in\omega$, we have
\[
\delta(\bar{a}\cdot x)(i,j) \ = \ (\bar{a}\cdot x)(j,i) \ = \ x(a_1+j, \, a_2+i) \ = \ \delta(x)(a_2+i,\, a_1+j) \ = \ \big(\bar{a}^\delta\cdot\delta(x)\big)(i,j).
\]
Finally, (iii) and (iv) follow easily from (ii).
\end{proof}

Now we are ready to conclude this section with a proof that when $\ell=2$, the restriction of $\EOl$ to the periodic part $\POl$ is hypersmooth. We work with an arbitrary fixed base space $\Omega$, since the proof does not depend on it. Let $\H$ be the countable set of all nontrivial subgroups of $\Z^2$. Since
\[
\POT \ = \ \bigsqcup_{H\in\H}\textup{Per}_H(\OT),
\]
by Lemmas \ref{lem:HBorel} and \ref{lem:piecewiseHS} we may work separately on the $\EOT$-invariant Borel sets $\textup{Per}_H(\OT)$. We will consider several cases, depending on $H$. The three main cases are: (i) $H$ has rank $2$; (ii) $H=\langle\m\rangle$ where $\m\in\Z_+^2$; and (iii) $H=\langle\m\rangle$ where $\m\in\Z_+\times -\N$. The few remaining cases are either trivial or can easily be reduced to one of these three.

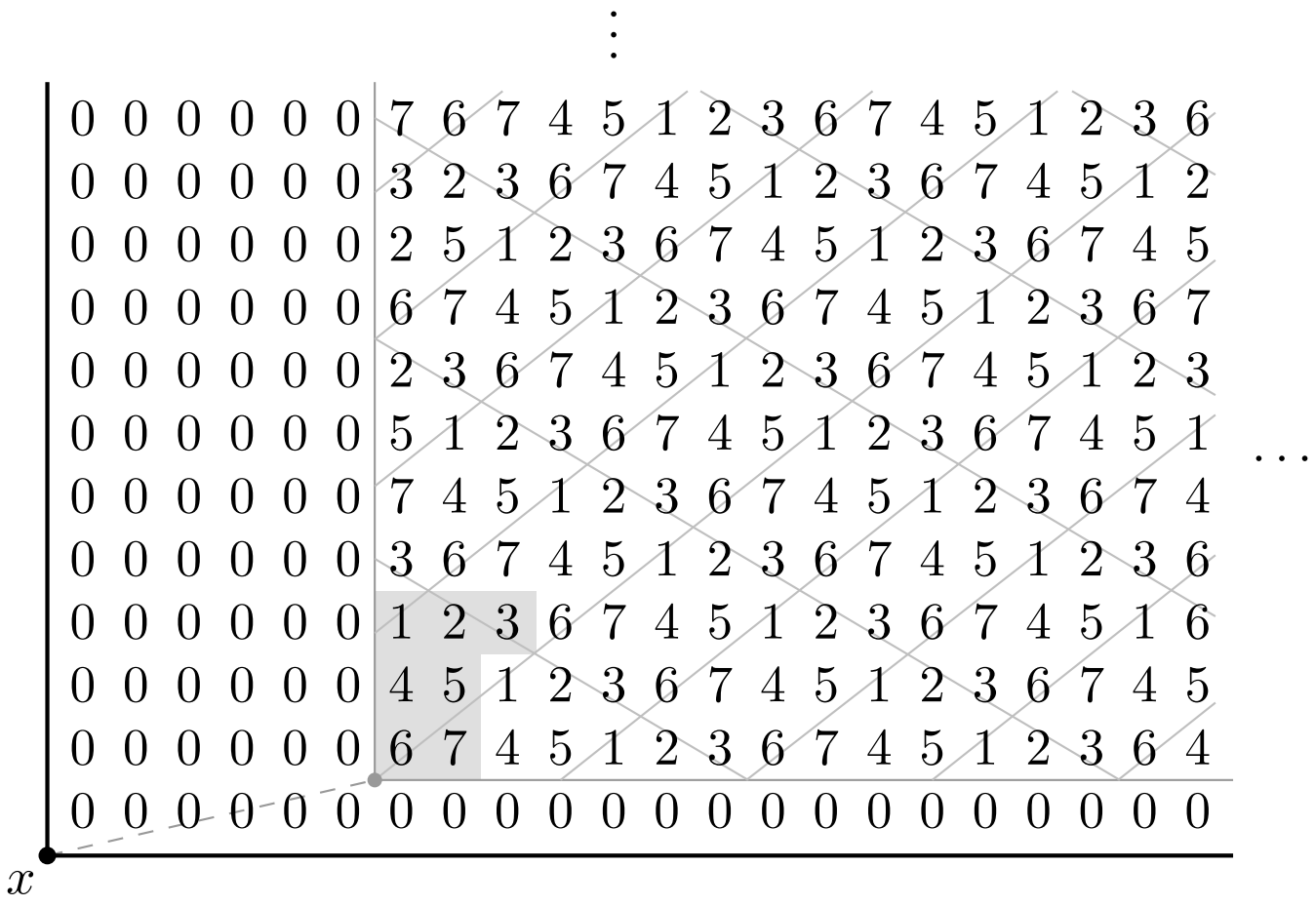
\begin{figure}[h]
\begin{center}
\begin{tikzpicture}[x=1.2cm, y=1.22cm]
 \draw [dashed,gray!80] (0,.7) -- (1.85,1.12);
 \draw [thick] (0,.7) -- (6.7,.7);
 \draw [thick] (0,.7) -- (0,5);
  \fill [fill=gray!25] (1.85,1.12) rectangle (2.45,2.17); 
  \fill [fill=gray!25] (2.44,1.82) rectangle (2.76,2.17); 
  \node at (0,.7) [circle,fill,inner sep=1.2pt]{};
  \draw [gray!80] (1.85,1.12) -- (6.7,1.12);
  \draw [gray!80] (1.85,1.12) -- (1.85,5);
    \draw [gray!50] (1.85,1.12) -- (6.6,4.83);
    \draw [gray!50] (2.9,1.12) -- (6.6,4.01);
    \draw [gray!50] (3.95,1.12) -- (6.6,3.15);
    \draw [gray!50] (5,1.12) -- (6.6,2.37);
    \draw [gray!50] (6.05,1.12) -- (6.6,1.55);
     \foreach \y in {7/3,14/3,7,28/3} \draw [gray!50] (1.85,1.12+.35*\y) -- (6.754-.448*\y,4.95);
    \draw [gray!50] (1.85,2.35) -- (3.96,1.12);
    \draw [gray!50] (1.85,3.575) -- (6.06,1.12);
    \draw [gray!50] (1.85,4.8) -- (6.6,2.03);
    \draw [gray!50] (3.69,4.95) -- (6.6,3.26);
     \draw [gray!50] (5.79,4.95) -- (6.6,4.485);
  \node at (1.85,1.12) [circle,fill=gray!80,inner sep=1pt]{};
 \foreach \x in {0,.3,...,6.6} \foreach \y in {.7} \node at (.2+\x,.25+\y) {$0$};
 \foreach \x in {0,.3,.6,.9,1.2,1.5} \foreach \y in {1.05,1.4,...,4.55} \node at (.2+\x,.25+\y) {$0$};
  \foreach \x in {.6,.9} \node at (.2+1.2+\x,.25+1.4+3.5*\x) {$6$};
  \foreach \x in {0,.3,.6,.9} \node at (.2+1.8+\x,.25+1.05+3.5*\x) {$6$};
  \foreach \x in {.3,.6,.9} \node at (.2+1.8+.6+\x,.25+.7+3.5*\x) {$6$};
  \foreach \x in {.3,.6,.9,1.2} \node at (.2+1.8+1.2+\x,.25+.35+3.5*\x) {$6$};
  \foreach \x in {.3,.6,.9,1.2} \node at (.2+1.8+1.8+\x,.25+3.5*\x) {$6$};
  \foreach \x in {.6,.9,1.2} \node at (.2+1.8+2.4+\x,.25-.35+3.5*\x) {$6$};
  \foreach \x in {.6,.9,1.2,1.5} \node at (.2+1.8+3+\x,.25-.7+3.5*\x) {$6$};
  \foreach \x in {.6,.9} \node at (.2+2.4+3+\x,.25-1.05+3.5*\x) {$6$};
  \node at (2,4.1) {$2$};
  \foreach \x in {0,.3,.6,.9} \node at (.2+1.8+\x,.25+1.4+3.5*\x) {$4$};
  \foreach \x in {0,.3,.6,.9} \node at (.2+1.8+.6+\x,.25+1.05+3.5*\x) {$4$};
  \foreach \x in {.3,.6,.9} \node at (.2+1.8+1.2+\x,.25+.7+3.5*\x) {$4$};
  \foreach \x in {.3,.6,.9,1.2} \node at (.2+1.8+1.8+\x,.6+3.5*\x) {$4$};
  \foreach \x in {.3,.6,.9,1.2} \node at (.2+1.8+2.4+\x,.25+3.5*\x) {$4$};
  \foreach \x in {.6,.9,1.2} \node at (.2+1.8+3+\x,.25-.35+3.5*\x) {$4$};
  \foreach \x in {.6,.9} \node at (.2+2.4+3+\x,.25-1.05+.35+3.5*\x) {$4$};
  \node at (2,4.45) {$3$};
  \foreach \x in {0,.3,.6} \node at (.2+1.8+\x,.25+.35+1.4+3.5*\x) {$1$};
  \foreach \x in {0,.3,.6,.9} \node at (.2+1.8+.6+\x,.25+.35+1.05+3.5*\x) {$1$};
  \foreach \x in {0,.3,.6,.9} \node at (.2+1.8+1.2+\x,.25+.35+.7+3.5*\x) {$1$};
  \foreach \x in {.3,.6,.9} \node at (.2+1.8+1.8+\x,.6+.35+3.5*\x) {$1$};
  \foreach \x in {.3,.6,.9,1.2} \node at (.2+1.8+2.4+\x,.25+.35+3.5*\x) {$1$};
  \foreach \x in {.3,.6,.9,1.2} \node at (.2+1.8+3+\x,.25-.35+.35+3.5*\x) {$1$};
  \foreach \x in {.6,.9} \node at (.2+2.4+3+\x,.25-1.05+.7+3.5*\x) {$1$};
  \foreach \x in {.3,.6,.9} \node at (.5+1.2+\x,.25+1.4+3.5*\x) {$7$};
  \foreach \x in {0,.3,.6,.9} \node at (.5+1.8+\x,.25+1.05+3.5*\x) {$7$};
  \foreach \x in {.3,.6,.9} \node at (.5+1.8+.6+\x,.25+.7+3.5*\x) {$7$};
  \foreach \x in {.3,.6,.9,1.2} \node at (.5+1.8+1.2+\x,.25+.35+3.5*\x) {$7$};
  \foreach \x in {.3,.6,.9,1.2} \node at (.5+1.8+1.8+\x,.25+3.5*\x) {$7$};
  \foreach \x in {.6,.9,1.2} \node at (.5+1.8+2.4+\x,.25-.35+3.5*\x) {$7$};
  \foreach \x in {.6,.9,1.2} \node at (.5+1.8+3+\x,.25-.7+3.5*\x) {$7$};
  \node at (6.5,1.3) {$4$};
  \foreach \x in {.3,.6} \node at (.5+1.2+\x,.25+.35+1.4+3.5*\x) {$5$};
  \foreach \x in {0,.3,.6,.9} \node at (.5+1.8+\x,.25+1.4+3.5*\x) {$5$};
  \foreach \x in {0,.3,.6,.9} \node at (.5+1.8+.6+\x,.25+1.05+3.5*\x) {$5$};
  \foreach \x in {.3,.6,.9} \node at (.5+1.8+1.2+\x,.25+.7+3.5*\x) {$5$};
  \foreach \x in {.3,.6,.9,1.2} \node at (.5+1.8+1.8+\x,.6+3.5*\x) {$5$};
  \foreach \x in {.3,.6,.9,1.2} \node at (.5+1.8+2.4+\x,.25+3.5*\x) {$5$};
  \foreach \x in {.6,.9,1.2} \node at (.5+1.8+3+\x,.25-.35+3.5*\x) {$5$};
  \node at (6.5,1.65) {$5$};
  \foreach \x in {.3,.6} \node at (.5+1.2+\x,.25+.7+1.4+3.5*\x) {$2$};
  \foreach \x in {0,.3,.6} \node at (.5+1.8+\x,.25+.35+1.4+3.5*\x) {$2$};
  \foreach \x in {0,.3,.6,.9} \node at (.5+1.8+.6+\x,.25+.35+1.05+3.5*\x) {$2$};
  \foreach \x in {0,.3,.6,.9} \node at (.5+1.8+1.2+\x,.25+.35+.7+3.5*\x) {$2$};
  \foreach \x in {.3,.6,.9} \node at (.5+1.8+1.8+\x,.6+.35+3.5*\x) {$2$};
  \foreach \x in {.3,.6,.9,1.2} \node at (.5+1.8+2.4+\x,.25+.35+3.5*\x) {$2$};
  \foreach \x in {.3,.6,.9,1.2} \node at (.5+1.8+3+\x,.25-.35+.35+3.5*\x) {$2$};
  \node at (6.5,2) {$6$};
  \node at (2,4.8) {$7$};
  \foreach \x in {.3,.6,.9} \node at (.5+1.2+\x,.25-.35+1.4+3.5*\x) {$3$};
  \foreach \x in {0,.3,.6} \node at (.8+1.8+\x,.25+.35+1.4+3.5*\x) {$3$};
  \foreach \x in {0,.3,.6,.9} \node at (.8+1.8+.6+\x,.25+.35+1.05+3.5*\x) {$3$};
  \foreach \x in {0,.3,.6,.9} \node at (.8+1.8+1.2+\x,.25+.35+.7+3.5*\x) {$3$};
  \foreach \x in {.3,.6,.9} \node at (.8+1.8+1.8+\x,.6+.35+3.5*\x) {$3$};
  \foreach \x in {.3,.6,.9,1.2} \node at (.8+1.8+2.4+\x,.25+.35+3.5*\x) {$3$};
  \foreach \x in {.3,.6,.9} \node at (.8+1.8+3+\x,.25-.35+.35+3.5*\x) {$3$};
 \node at (-.15,.55) {$x$};
 \node at (7,2.9) {$\cdots$};
 \node at (3.2,5.35) {$\vdots$};
\end{tikzpicture}
\end{center}
\caption{A sample $x\in\textup{Per}_H(\OT)$, where $H=\big\langle (3,2), \, (1,3)\big\rangle$ and $\c^x=(6,1)$.}
\label{fig:rank2}
\end{figure}

The basic ideas for the proofs are fairly straightforward in all three cases, although the details are a bit tedious to verify in cases (ii) and (iii). We begin with case (i), where rank$(H)=2$, so $\Z^2/H$ is finite. Suppose $H\leq\Z^2$ has rank 2, and let $x\in\textup{Per}_H(\OT)$. Intuitively, the $\EOT$-class of $x$ is completely determined by the values of $x$ at points $\c^x+\bar{i}$ where $\c^x\cdot x$ is strongly $H$-periodic and $\bar{i}\in\N^{\ell}$ ranges over a set of coset representatives for $H$ in $\Z^2$, considered up to a shift modulo $H$. A visual example is given in Figure \ref{fig:rank2}, where different numbers represent different elements of $\Omega$. In this figure, we have taken $H=\big\langle (3,2), \, (1,3)\big\rangle$ and $\c^x=(6,1)$, and we have shaded a set of coset representatives for $H$ in $\Z^2$, namely $\R\big((2,3)\big)\cup\{(2,2)\}$. With this picture as a reference, we now turn to the formal proof.

\begin{proposition}\label{prop:Case1}
For any rank $2$ subgroup $H\leq\Z^2$, $\EOT\res\textup{Per}_H(\OT)$ is smooth.
\end{proposition}

\begin{proof}
Suppose $H\leq\Z^2$ has rank $2$, so $\Z^2/H$ is finite. Let $Q\subseteq\N^2$ be a fixed transversal containing $\bar{0}$, and $q$ the corresponding selector, for the coset equivalence relation on $\Z^2$ induced by $H$, so $Q$ is finite. Let $F$ be the finite (hence smooth) Borel equivalence relation on $\Omega^{Q}$ defined by
\[
\alpha\mathrel{F}\beta \quad :\Longleftrightarrow \quad (\exists\s\in Q)(\forall\bar{i}\in Q) \ \alpha(\bar{i})=\beta\big(q(\s+\bar{i})\big).
\]
Using Lemma \ref{lem:per}(viii), for each $x\in \textup{Per}_H(\OT)$ let $\bar{c}^x\in\N^2$ be least such that $\bar{c}^x\cdot x$ is strongly $H$-periodic, so $x\mapsto\bar{c}^x$ is Borel, and define $f(x)\in\Omega^{Q}$ by
\[
f(x)(\bar{i}) \, := \, x(\c^x+\bar{i}) \qquad \text{for each $\bar{i}\in Q$,}
\]
so $f:\textup{Per}_H(\OT)\to\Omega^{Q}$ is Borel. We claim that $f$ Borel reduces $\EOT\res\textup{Per}_H(\OT)$ to $F$.

To see this, fix $x,y\in\textup{Per}_H(\OT)$, and assume first that $x\mathrel{\EOT}y$. Fix $\bar{a},\bar{b}\in\N^2$ such that $\bar{a}\cdot x=\bar{b}\cdot y$. Let $\t$ be the unique element of $Q$ such that $q(\bar{c}^y+\bar{a}+\t)=\bar{0}$, and let $\s:=q(\bar{c}^x+\bar{b}+\t)$. 
Using strong $H$-periodicity of $\c^x\cdot x$ and $\c^y\cdot y$, for all $\bar{i}\in Q$ we have
\begin{align*}
f(x)(\bar{i}) \ &= \ x(\bar{c}^x+\bar{i}) \ = \ x(\bar{c}^x+\bar{i}+\bar{c}^y+\bar{a}+\t) \ = \ (\bar{a}\cdot x)(\bar{c}^x+\bar{i}+\bar{c}^y+\t) \ = \ (\bar{b}\cdot y)(\bar{c}^x+\bar{i}+\bar{c}^y+\t) \\
&= \ y(\bar{c}^x+\bar{i}+\bar{c}^y+\bar{b}+\t) \ = \ y(\bar{c}^y+\s+\bar{i}) \ = \ y\big(\bar{c}^y+q(\s+\bar{i})\big) \ = \ f(y)\big(q(\s+\bar{i})\big),
\end{align*}
which shows $f(x)\mathrel{F}f(y)$. 

Conversely, suppose $f(x)\mathrel{F}f(y)$, and fix $\s\in Q$ such that $f(x)(\bar{i})=f(y)\big(q(\s+\bar{i})\big)$ for all $\bar{i}\in Q$. Let $\bar{i}\in\N^2$ be arbitrary, and write $\bar{i}'=q(\bar{i})$. Then
\[
x(\bar{c}^x+\bar{i}) \ = \ x(\c^x+\bar{i}') \ = \  f(x)(\bar{i}') \ = \ f(y)\big(q(\s+\bar{i}')\big) \ = \ y\big(\bar{c}^y+q(\s+\bar{i}')\big) \ = \ y(\bar{c}^y+\s+\bar{i}),
\]
so $\bar{c}^x\cdot x=(\bar{c}^y+\s)\cdot y$, showing $x\mathrel{\EOT}y$. We have shown that $f$ Borel reduces $\EOT\res\textup{Per}_H(\OT)$ to the smooth equivalence relation $F$, completing the proof in this case.
\end{proof}

Now we turn our attention to cases (ii) and (iii), which are similar to each other. In case (ii), where $H=\langle\m\rangle$ with $\m\in\Z_+^2$, the $\EOT$-class of $x\in\textup{Per}_H(\OT)$ is determined by a bi-infinite vertical strip of values of width $m_1$, considered up to a horizontal shift.  A visual example is given in Figure \ref{fig:Case2}, where again different symbols represent different elements of $\Omega$. Here $\m=(3,1)$, and $\c^x\cdot x$ is strongly $H$-periodic where $\c^x=(6,3)$. In the picture, one should imagine completing the bi-infinite vertical strip below the $abc$ row by making use of the $(3,1)$-periodicity, so the row directly below $abc$ is $\alpha\beta\gamma$, and the row below that is $\delta\epsilon\zeta$, etc. Completed in this way, this bi-infinite vertical strip together with the vector $\m$ contains all the information necessary to determine $[x]_{\EOT}$. In fact it contains too much information, since we should not know which column in the strip comes ``first,'' so we will need to mod out by a horizontal shift.  This must be done carefully: by ``horizontal shift'' we do not quite mean a cyclic permutation of the columns, since there is a vertical step to consider as well.  For instance, if in Figure \ref{fig:Case2} we shift the shaded strip two places to the right, the $def$ row becomes $fab$ rather than $fde$, since the $de$ entries have now moved up one row. This motivates the definition of the maps $\sigma_{r,\m}$, which should be read with Figure \ref{fig:Case2} in mind.

\begin{figure}[h]
\begin{center}
\begin{tikzpicture}[x=1.2cm, y=1.22cm]
     \draw[->,gray!120] (4.55,1.15) -- (2.79,.42);
     \fill [fill=gray!25] (1.85,-.3) rectangle (2.75,3.6); 
     \draw(1.85,.1) rectangle (2.75,.4);
     \draw(4.55,1.15) rectangle (5.45,1.45);
 \draw [dashed,gray!80] (0,0) -- (1.85,1.12);
 \draw [thick] (0,0) -- (6.7,0);
 \draw [thick] (0,0) -- (0,3.6);
  \node at (0,0) [circle,fill,inner sep=1.2pt]{};
  \draw [gray!80] (1.85,1.12) -- (6.7,1.12);
  \draw [gray!80] (1.85,1.12) -- (1.85,3.6);
    \node at (1.85,1.12) [circle,fill=gray!80,inner sep=1pt]{};
 \foreach \x in {0,.3,...,6.6} \foreach \y in {0,.35,.7} \node at (.2+\x,.25+\y) {$0$};
 \foreach \x in {0,.3,.6,.9,1.2,1.5} \foreach \y in {1.05,1.4,...,3.15} \node at (.2+\x,.25+\y) {$0$};
  \foreach \x in {0,.3,.6,.9,1.2,1.5} \node[anchor=south] at (.2+1.8+3*\x,.25+.858+21*\x/18) {$a$};
  \foreach \x in {0,.3,.6,.9,1.2} \node[anchor=south] at (.5+1.8+3*\x,.25+.858+21*\x/18) {$b$};
  \foreach \x in {0,.3,.6,.9,1.2} \node[anchor=south] at (.8+1.8+3*\x,.25+.858+21*\x/18) {$c$};
  \foreach \x in {0,.3,.6,.9,1.2,1.5} \node[anchor=south] at (.2+1.8+3*\x,.25+1.208+21*\x/18) {$d$};
  \foreach \x in {0,.3,.6,.9,1.2} \node[anchor=south] at (.5+1.8+3*\x,.25+1.208+21*\x/18) {$e$};
  \foreach \x in {0,.3,.6,.9,1.2} \node[anchor=south] at (.8+1.8+3*\x,.25+1.15+21*\x/18) {$f$};
  \foreach \x in {0,.3,.6,.9,1.2} \node[anchor=south] at (.2+1.8+3*\x,.25+1.51+21*\x/18) {$g$};
  \foreach \x in {0,.3,.6,.9,1.2} \node[anchor=south] at (.5+1.8+3*\x,.25+1.558+21*\x/18) {$h$};
  \foreach \x in {0,.3,.6,.9,1.2} \node[anchor=south] at (.8+1.8+3*\x,.25+1.558+21*\x/18) {$i$};
  \foreach \x in {0,.3,.6,.9} \node[anchor=south] at (.2+1.8+3*\x,.25+1.86+21*\x/18) {$j$};
  \foreach \x in {0,.3,.6,.9} \node[anchor=south] at (.5+1.8+3*\x,.25+1.908+21*\x/18) {$k$};
  \foreach \x in {0,.3,.6,.9} \node[anchor=south] at (.8+1.8+3*\x,.25+1.908+21*\x/18) {$l$};
  \foreach \x in {0,.3,.6} \node[anchor=south] at (.2+1.8+3*\x,.25+2.258+21*\x/18) {$m$};
  \foreach \x in {0,.3,.6} \node[anchor=south] at (.5+1.8+3*\x,.25+2.258+21*\x/18) {$n$};
  \foreach \x in {0,.3,.6} \node[anchor=south] at (.8+1.8+3*\x,.25+2.258+21*\x/18) {$o$};
  \foreach \x in {0,.3} \node[anchor=south] at (.2+1.8+3*\x,.25+2.56+21*\x/18) {$p$};
  \foreach \x in {0,.3} \node[anchor=south] at (.5+1.8+3*\x,.25+2.56+21*\x/18) {$q$};
  \foreach \x in {0,.3} \node[anchor=south] at (.8+1.8+3*\x,.25+2.608+21*\x/18) {$r$};
  \foreach \x in {0} \node[anchor=south] at (.2+1.8+3*\x,.25+2.958+21*\x/18) {$s$};
  \foreach \x in {0} \node[anchor=south] at (.5+1.8+3*\x,.25+2.958+21*\x/18) {$t$};
  \foreach \x in {0} \node[anchor=south] at (.8+1.8+3*\x,.25+2.958+21*\x/18) {$u$};
  \foreach \x in {0,.3,.6,.9,1.2} \node[anchor=south] at (1.1+1.8+3*\x,.25+.858+21*\x/18) {$\alpha$};
  \foreach \x in {0,.3,.6,.9} \node[anchor=south] at (1.4+1.8+3*\x,.25+.8+21*\x/18) {$\beta$};
  \foreach \x in {0,.3,.6,.9} \node[anchor=south] at (1.7+1.8+3*\x,.25+.8+21*\x/18) {$\gamma$};
  \foreach \x in {0,.3,.6,.9} \node[anchor=south] at (2+1.8+3*\x,.25+.855+21*\x/18) {$\delta$};
  \foreach \x in {0,.3,.6} \node[anchor=south] at (2.3+1.8+3*\x,.25+.858+21*\x/18) {$\epsilon$};
  \foreach \x in {0,.3,.6} \node[anchor=south] at (2.6+1.8+3*\x,.25+.8+21*\x/18) {$\zeta$};
  \foreach \x in {0,.3,.6} \node[anchor=south] at (2.9+1.8+3*\x,.25+.82+21*\x/18) {$\eta$};
  \foreach \x in {0,.3} \node[anchor=south] at (3.2+1.8+3*\x,.25+.855+21*\x/18) {$\theta$};
  \foreach \x in {0,.3} \node[anchor=south] at (3.5+1.8+3*\x,.25+.858+21*\x/18) {$\iota$};
  \foreach \x in {0,.3} \node[anchor=south] at (3.8+1.8+3*\x,.25+.858+21*\x/18) {$\kappa$};
  \node[anchor=south] at (4.1+1.8,.25+.858) {$\lambda$};
  \node[anchor=south] at (4.4+1.8,.25+.806) {$\mu$};
  \node[anchor=south] at (4.7+1.8,.25+.858) {$\nu$};
 \node at (-.15,-.15) {$x$};
 \node at (7,1.8) {$\cdots$};
 \node at (3.2,3.83) {$\vdots$};
 \node[gray!80] at (2.3,3.9) {$\vdots$};
 \node[gray!80] at (2.3,-.4) {$\vdots$};
\end{tikzpicture}
\end{center}
\caption{A sample $x\in\textup{Per}_H(\OT)$, where $H=\big\langle (3,1)\big\rangle$ and $\c^x=(6,3)$.}
\label{fig:Case2}
\end{figure}

\begin{definition}
For each $\m=(m_1,m_2)\in\Z_+^2$ and $r\in m_1$, define the Borel function $$\sigma_{r,\m}:(2^\omega)^{m_1\times\Z}\to (2^\omega)^{m_1\times\Z}$$ as follows: for each $x\in (2^\omega)^{m_1\times\Z}$ and $(i,j)\in m_1\times\Z$, put 
\[
\sigma_{r,\m}(x)(i,j) \ := \ 
\begin{cases}
    x(i+r,j) & \text{if $i+r<m_1$;} \\
    x(i+r-m_1, \, j-m_2) & \text{if $i+r\geq m_1$.}
\end{cases}
\]
\end{definition}

Intuitively, $\sigma_{r,\m}(x)$ is the bi-infinite vertical strip in $\Omega^{m_1\times\Z}$ obtained by extending $x$ to an element $\hat{x}$ in $\Omega^{\omega\times\Z}$ using $\m$-periodicity and then shifting the strip $x$ within $\hat{x}$ to the right by $r$. We use $2^\omega$ in the definition of $\sigma_{r,\m}$ since we can encode elements of $\Omega$ in $2^\omega$ no matter what $\Omega$ is. Now we define an equivalence relation $\E_{\m}$ on the space $(2^\omega)^{m_1\times\Z}$ of bi-infinite vertical $m_1$-strips as follows.

\begin{definition}
For each $\m=(m_1,m_2)\in\Z_+^2$ and $x,y\in (2^\omega)^{m_1\times\Z}$, let
\[
x\mathrel{\mathcal{E}_{\m}}y \quad :\Longleftrightarrow \quad (\exists r\in m_1)(\exists a,b\in\N)(\forall i\in m_1)(\forall j\in\Z) \ x(i, a+j) = \sigma_{r,\m}(y)(i, b+j).
\]
\end{definition}

Intuitively, if $H=\langle\m\rangle\leq\Z^2$ with $\m\in\Z_+^2$, then for all $x,y\in\textup{Per}_H\big((2^\omega)^{\omega^2}\big)$ we have that $x\mathrel{E_2^{2^\omega}}y$ if and only if the bi-infinite vertical strips determined by $x$ and $y$ as in Figure \ref{fig:Case2} are $\E_{\m}$-equivalent to each other. We can therefore show $\EOT\res\textup{Per}_H(\OT)$ is hypersmooth by coding $\Omega$ into $2^\omega$, verifying the preceding statement, and checking that $\E_{\m}$ is a hypersmooth equivalence relation. Proving all this  in detail is unfortunately rather technical, so the reader is encouraged to keep this intuitive picture in mind throughout the argument.

Before turning to the details, let us also consider the intuition behind case (iii), since cases (ii) and (iii) are similar and we will be able to combine some aspects of their proofs.  In case (iii), we have $H=\langle\m\rangle$ where $\m=(m_1,m_2)\in\Z_+\times-\N$. This time, the $\EOT$-class of $x\in\text{Per}_H(\OT)$ is determined by an upward-infinite vertical strip of values of width $m_1$, again considered up to a horizontal shift.  A visual example is shown in Figure \ref{fig:Case3}, where now $H=\langle\m\rangle$ with $\m=(3,-1)$ and once again $\c^x\cdot x$ is strongly $H$-periodic with $\c^x=(6,3)$. Note that an $\omega$-ordered vertical strip suffices in this case, in place of the $\Z$-ordered vertical strip we require in case (ii). Because of this, we will define $\sigma_{r,\m}$ and $\E_{\m}$ slightly differently in this case, so we add a prime in the notation to avoid confusion, but the basic idea behind their definitions is the same.

\begin{figure}[h]
\begin{center}
\begin{tikzpicture}[x=1.2cm, y=1.22cm]
     \draw[->,gray!120] (4.55,1.45) -- (2.79,2.18);
     \fill [fill=gray!25] (1.85,1.12) rectangle (2.75,3.6); 
     \draw(4.55,1.15) rectangle (5.45,1.45);
 \draw [dashed,gray!80] (0,0) -- (1.85,1.12);
 \draw [thick] (0,0) -- (6.7,0);
 \draw [thick] (0,0) -- (0,3.6);
  \node at (0,0) [circle,fill,inner sep=1.2pt]{};
  \draw [gray!80] (1.85,1.12) -- (6.7,1.12);
  \draw [gray!80] (1.85,1.12) -- (1.85,3.6);
    \node at (1.85,1.12) [circle,fill=gray!80,inner sep=1pt]{};
 \foreach \x in {0,.3,...,6.6} \foreach \y in {0,.35,.7} \node at (.2+\x,.25+\y) {$0$};
 \foreach \x in {0,.3,.6,.9,1.2,1.5} \foreach \y in {1.05,1.4,...,3.15} \node at (.2+\x,.25+\y) {$0$};
  \foreach \x in {0} \node[anchor=south] at (.2+1.8+3*\x,.25+.858+21*\x/18) {$a$};
  \foreach \x in {0} \node[anchor=south] at (.5+1.8+3*\x,.25+.858+21*\x/18) {$b$};
  \foreach \x in {0} \node[anchor=south] at (.8+1.8+3*\x,.25+.858+21*\x/18) {$c$};
  \foreach \x in {0,.3} \node[anchor=south] at (.2+1.8+3*\x,.25+1.208-21*\x/18) {$d$};
  \foreach \x in {0,.3} \node[anchor=south] at (.5+1.8+3*\x,.25+1.208-21*\x/18) {$e$};
  \foreach \x in {0,.3} \node[anchor=south] at (.8+1.8+3*\x,.25+1.15-21*\x/18) {$f$};
  \foreach \x in {0,.3,.6} \node[anchor=south] at (.2+1.8+3*\x,.25+1.51-21*\x/18) {$g$};
  \foreach \x in {0,.3,.6} \node[anchor=south] at (.5+1.8+3*\x,.25+1.558-21*\x/18) {$h$};
  \foreach \x in {0,.3,.6} \node[anchor=south] at (.8+1.8+3*\x,.25+1.558-21*\x/18) {$i$};
  \foreach \x in {0,.3,.6,.9} \node[anchor=south] at (.2+1.8+3*\x,.25+1.86-21*\x/18) {$j$};
  \foreach \x in {0,.3,.6,.9} \node[anchor=south] at (.5+1.8+3*\x,.25+1.908-21*\x/18) {$k$};
  \foreach \x in {0,.3,.6,.9} \node[anchor=south] at (.8+1.8+3*\x,.25+1.908-21*\x/18) {$l$};
  \foreach \x in {0,.3,.6,.9,1.2} \node[anchor=south] at (.2+1.8+3*\x,.25+2.258-21*\x/18) {$m$};
  \foreach \x in {0,.3,.6,.9,1.2} \node[anchor=south] at (.5+1.8+3*\x,.25+2.258-21*\x/18) {$n$};
  \foreach \x in {0,.3,.6,.9,1.2} \node[anchor=south] at (.8+1.8+3*\x,.25+2.258-21*\x/18) {$o$};
  \foreach \x in {0,.3,.6,.9,1.2,1.5} \node[anchor=south] at (.2+1.8+3*\x,.25+2.56-21*\x/18) {$p$};
  \foreach \x in {0,.3,.6,.9,1.2} \node[anchor=south] at (.5+1.8+3*\x,.25+2.56-21*\x/18) {$q$};
  \foreach \x in {0,.3,.6,.9,1.2} \node[anchor=south] at (.8+1.8+3*\x,.25+2.608-21*\x/18) {$r$};
  \foreach \x in {0,.3,.6,.9,1.2,1.5} \node[anchor=south] at (.2+1.8+3*\x,.25+2.958-21*\x/18) {$s$};
  \foreach \x in {0,.3,.6,.9,1.2} \node[anchor=south] at (.5+1.8+3*\x,.25+2.958-21*\x/18) {$t$};
  \foreach \x in {0,.3,.6,.9,1.2} \node[anchor=south] at (.8+1.8+3*\x,.25+2.958-21*\x/18) {$u$};
  \foreach \x in {0,.3,.6,.9,1.2} \node[anchor=south] at (1.1+1.8+3*\x,.25+2.958-21*\x/18) {$\alpha$};
  \foreach \x in {0,.3,.6,.9} \node[anchor=south] at (1.4+1.8+3*\x,.25+2.9-21*\x/18) {$\beta$};
  \foreach \x in {0,.3,.6,.9} \node[anchor=south] at (1.7+1.8+3*\x,.25+2.9-21*\x/18) {$\gamma$};
  \foreach \x in {0,.3,.6,.9} \node[anchor=south] at (2+1.8+3*\x,.25+2.955-21*\x/18) {$\delta$};
  \foreach \x in {0,.3,.6} \node[anchor=south] at (2.3+1.8+3*\x,.25+2.958-21*\x/18) {$\epsilon$};
  \foreach \x in {0,.3,.6} \node[anchor=south] at (2.6+1.8+3*\x,.25+2.9-21*\x/18) {$\zeta$};
  \foreach \x in {0,.3,.6} \node[anchor=south] at (2.9+1.8+3*\x,.25+2.92-21*\x/18) {$\eta$};
  \foreach \x in {0,.3} \node[anchor=south] at (3.2+1.8+3*\x,.25+2.955-21*\x/18) {$\theta$};
  \foreach \x in {0,.3} \node[anchor=south] at (3.5+1.8+3*\x,.25+2.958-21*\x/18) {$\iota$};
  \foreach \x in {0,.3} \node[anchor=south] at (3.8+1.8+3*\x,.25+2.958-21*\x/18) {$\kappa$};
  \foreach \x in {0} \node[anchor=south] at (4.1+1.8+3*\x,.25+2.958-21*\x/18) {$\lambda$};
  \foreach \x in {0} \node[anchor=south] at (4.4+1.8+3*\x,.25+2.906-21*\x/18) {$\mu$};
  \node[anchor=south] at (4.7+1.8,.25+2.958) {$\nu$};
  %
 \node at (-.15,-.15) {$x$};
 \node at (7,1.8) {$\cdots$};
 \node at (3.2,3.9) {$\vdots$};
 \node[gray!80] at (2.3,3.9) {$\vdots$};
 \draw(1.85,2.2) rectangle (2.75,2.5);
\end{tikzpicture}
\end{center}
\caption{A sample $x\in\textup{Per}_H(\OT)$, where $H=\big\langle (3,-1)\big\rangle$ and $\c^x=(6,3)$.}
\label{fig:Case3}
\end{figure}

\begin{definition}
For each $\m=(m_1,m_2)\in\Z_+\times-\N$ and $r\in m_1$, define the Borel function $$\sigma_{r,\m}':(2^\omega)^{m_1\times\omega}\to (2^\omega)^{m_1\times\omega}$$ by the same rule as $\sigma_{r,\m}$, that is,
\[
\sigma_{r,\m}'(x)(i,j) \ := \ 
\begin{cases}
    x(i+r,j) & \text{if $i+r<m_1$;} \\
    x(i+r-m_1, \, j-m_2) & \text{if $i+r\geq m_1$}
\end{cases}
\]
for each $x\in (2^\omega)^{m_1\times\omega}$ and $(i,j)\in m_1\times\omega$.
\end{definition}

\begin{definition}
For each $\m=(m_1,m_2)\in\Z_+\times-\N$ and $x,y\in (2^\omega)^{m_1\times\omega}$, let
\[
x\mathrel{\mathcal{E}_{\m}'}y \quad :\Longleftrightarrow \quad (\exists r\in m_1)(\exists a,b\in\N)(\forall i\in m_1)(\forall j\in\omega) \ x(i, a+j) = \sigma_{r,\m}'(y)(i, b+j).
\]
\end{definition}

Note that $\E_{\m}$ lives on $(2^\omega)^{m_1\times\Z}$ and is defined when $\m\in\Z_+^2$, while $\E_{\m}'$ lives on $(2^\omega)^{m_1\times\omega}$ and is defined when $\m\in\Z_+\times-\N$, but the only differences in their defining conditions are whether $j$ ranges over $\omega$ or $\Z$ and whether $\sigma_{r,\m}$ or $\sigma_{r,\m}'$ is used for the ``horizontal shifts." Once again, the intuition is that if $H=\langle\m\rangle\leq\Z^2$ with $\m\in\Z_+\times-\N$, then for all $x,y\in\textup{Per}_H\big((2^\omega)^{\omega^2}\big)$ we have that $x\mathrel{E_2^{2^\omega}}y$ if and only if the upward-infinite vertical strips determined by $x$ and $y$ as in Figure \ref{fig:Case3} are $\E_{\m}'$-equivalent to each other. With these intuitions in place, we now turn to the details of the proofs.  There are three main things we need to prove: that $\E_{\m}$ and $\E_{\m}'$ are equivalence relations, that they are hypersmooth, and that $\EOT\res\textup{Per}_H(\OT)$ does indeed reduce to the appropriate one, depending on the case.  We start by establishing some basic facts about the maps $\sigma_{r,\m}$ and $\sigma_{r,\m}'$, often omitting $\m$ from the notation when it is clear from context. 

\begin{lemma}\label{lem:sigma}
Let $\m\in\Z_+^2$ and $x\in(2^\omega)^{m_1\times\Z}$ and $r,s\in m_1$ with $t:=(r+s)\mod m_1$.
\begin{enumerate}
 \item[(i)] If $r+s<m_1$, then $\sigma_r\big(\sigma_s(x)\big)=\sigma_t(x)$.
 \item[(ii)] If $r+s\geq m_1$, then $\sigma_r\big(\sigma_s(x)\big)\big(i,\, j\big)=\sigma_t(x)\big(i, \, j-m_2\big)$ for all $i\in m_1$ and $j\in\Z$.
\end{enumerate}
\end{lemma}

\begin{proof}
For (i), assume $r+s<m_1$, so $r+s=t$. Let $i\in m_1$ and $j\in\Z$ be arbitrary. We consider several cases. If $i+t<m_1$, so also $i+r<m_1$, then
\[
\sigma_r\big(\sigma_s(x)\big)(i,j) \ = \ \sigma_s(x)(i+r,\, j) \ = \ x(i+t,\, j) \ = \ \sigma_t(x)(i,j).
\]
If $i+t\geq m_1$ and $i+r<m_1$, then
\[
\sigma_r\big(\sigma_s(x)\big)(i,j) \ = \ \sigma_s(x)(i+r,\, j) \ = \ x(i+t-m_1,\, j-m_2) \ = \ \sigma_t(x)(i,j).
\]
Finally, if $i+t\geq m_1$ and $i+r\geq m_1$, then $i+t-m_1<m_1$ and 
\begin{align*}
\sigma_r\big(\sigma_s(x)\big)(i,j) \ = \ \sigma_s(x)(i+r-m_1, \, j-m_2) \ = \ x(i+t-m_1, \, j-m_2) \ = \ \sigma_t(x)(i, j).
\end{align*}
So $\sigma_r\big(\sigma_s(x)\big)(i,j)=\sigma_t(x)(i, j)$ in all cases, and we have shown $\sigma_r\big(\sigma_s(x)\big)=\sigma_t(x)$.

For (ii), assume $r+s\geq m_1$, so $t=r+s-m_1$. Let $i\in m_1$ and $j\in\Z$ be arbitrary. Once again we consider cases. If $i+r<m_1$, then also $i+t<m_1$ and 
\[
\sigma_r\big(\sigma_s(x)\big)(i,\, j) \ = \ \sigma_s(x)(i+r,\, j) \ = \ x(i+t, \, j-m_2) \ = \ \sigma_t(x)(i,\, j-m_2).
\]
If $i+r\geq m_1$ and $i+t<m_1$, then
\[
\sigma_r\big(\sigma_s(x)\big)(i,j) \ = \ \sigma_s(x)(i+r-m_1,\, j-m_2) \ = \ x(i+t,\,j-m_2) \ = \ \sigma_t(x)(i,\, j-m_2). 
\]
Finally, if $i+r\geq m_1$ and $i+t\geq m_1$, then
\[
\sigma_r\big(\sigma_s(x)\big)(i,\, j) \ = \ \sigma_s(x)(i+r-m_1,\, j-m_2) \ = \ x(i+t-m_1, \, j-2m_2) \ = \ \sigma_t(x)(i, \, j-m_2). 
\]
Thus $\sigma_r\big(\sigma_s(x)\big)(i,j)=\sigma_t(x)(i,\, j-m_2)$ in all cases, completing the proof of (ii) and of the lemma.
\end{proof}

\begin{lemma}\label{lem:sigma'}
Let $\m\in\Z_+\times-\N$ and $x\in(2^\omega)^{m_1\times\omega}$ and $r,s\in m_1$ with $t:=(r+s)\mod m_1$.
\begin{enumerate}
 \item[(i)] If $r+s<m_1$, then $\sigma_r'\big(\sigma_s'(x)\big)=\sigma_t'(x)$.
 \item[(ii)] If $r+s\geq m_1$, then $\sigma_r'\big(\sigma_s'(x)\big)\big(i,\, j\big)=\sigma_t'(x)\big(i, \, j-m_2\big)$ for all $i\in m_1$ and $j\in\omega$.
\end{enumerate}
\end{lemma}

\begin{proof}
The proof is identical to that of \ref{lem:sigma} once we let $j$ vary over $\omega$ instead of $\Z$, and replace $\sigma$ with $\sigma'$ throughout. Note that while we work on $(2^\omega)^{m_1\times\omega}$ here instead of $(2^\omega)^{m_1\times\Z}$, we need not worry about subtracting $m_2$ in the second component because now $m_2\leq 0$.
\end{proof}

\begin{lemma}\label{lem:Eequiv}
For all $\m\in\Z_+^2$, $\mathcal{E}_{\m}$ is an equivalence relation on $(2^\omega)^{m_1\times\Z}$.
\end{lemma}

\begin{proof}
Clearly $\mathcal{E}_{\m}$ is reflexive, since we can take $r=a=b=0$. For symmetry, fix $x,y\in (2^\omega)^{m_1\times\Z}$ such that  $x\mathrel{\mathcal{E}_{\m}}y$, and fix $r\in m_1$ and $a,b\in\N$ such that
\[
x(i,\, a+j) = \sigma_{r}(y)(i,\, b+j) \quad \text{for all $i\in m_1$ and $j\in\Z$.}
\]
If $r=0$, then $\sigma_r(y)=y$ and it is obvious that $y\mathrel{\mathcal{E}_{\m}}x$; otherwise, put $\tilde{r}:=m_1-r$, so $(r+\tilde{r})\mod m_1=0$, and let $i\in m_1$ and $j\in\Z$ be arbitrary. If $\tilde{r}+i<m_1$, then 
\[
\sigma_{\tilde{r}}\big(\sigma_r(y)\big)(i, \, b+j) \ = \ \sigma_r(y)\big(\tilde{r}+i, \, b+j\big) \ = \ x\big(\tilde{r}+i, \, a+j\big) \ = \ \sigma_{\tilde{r}}(x)(i,a+j),
\]
and if $\tilde{r}+i\geq m_1$ then again
\[
\sigma_{\tilde{r}}\big(\sigma_r(y)\big)(i, \, b+j) \ = \ \sigma_r(y)\big(\tilde{r}+i-m_1,\, b+j-m_1\big) \ = \ x\big(\tilde{r}+i-m_1,\, a+j-m_1\big) \ = \ \sigma_{\tilde{r}}(x)(i,a+j).
\]
Thus by Lemma \ref{lem:sigma}(ii) we have
\[
y(i, \, b-m_2+j) \ = \ \sigma_{\tilde{r}}\big(\sigma_r(y)\big)(i, \, b+j) \ = \ \sigma_{\tilde{r}}(x)(i, \, a+j),
\]
showing $y\mathrel{\E_{\m}}x$ as required.

For transitivity, fix also $z\in (2^\omega)^{m_1\times\Z}$ such that $y\mathrel{\mathcal{E}_{\m}}z$, and suppose $s\in m_1$ and $c,d\in\N$ satisfy $y(i,\, c+j)=\sigma_{s}(z)(i,\, d+j)$ for all $i\in m_1$ and $j\in\Z$. We must show $x\mathrel{\mathcal{E}_{\m}}z$. Put $t:=(r+s)\mod m_1$, and let $i\in m_1$ and $j\in\Z$ be arbitrary. We consider two cases. 

\smallskip

\underline{Case 1}: $r+s<m_1$, so $r+s=t$. Using Lemma \ref{lem:sigma}(i), if $i+r<m_1$ then
\[
x(i, \, a+c+j) \ = \ \sigma_r(y)(i, \, b+c+j) \ = \ y(i+r, \, b+c+j) \ = \ \sigma_s(z)(i+r, \, b+d+j) 
\ = \ \sigma_t(z)(i, \, b+d+j),
\]
and if $i+r\geq m_1$ then again
\begin{align*}
x(i, \, a+c+j) \ &= \ \sigma_r(y)(i, \, b+c+j) \ = \ y(i+r-m_1, \, b+c-m_2+j) \\ &= \ \sigma_s(z)(i+r-m_1, \, b+d-m_2+j) \ = \ \sigma_r\big(\sigma_s(z)\big)(i, \, b+d+j) 
\ = \ \sigma_t(z)(i, \, b+d+j).
\end{align*}
So $x(i, \, a+c+j)=\sigma_t(z)(i, \, b+d+j)$ for all $i\in m_1$ and $j\in\Z$, which shows $x\mathrel{\E_{\m}}z$ in this case.

\smallskip

\underline{Case 2}: $r+s\geq m_1$. Using Lemma \ref{lem:sigma}(ii), if $i+r<m_1$ then
\begin{align*}
x(i, \, a+c+j) \ &= \ \sigma_r(y)(i, \, b+c+j) \ = \ y(i+r, \, b+c+j) \ = \  \sigma_s(z)(i+r, \, b+d+j) \\ &= \ \sigma_r\big(\sigma_s(z)\big)(i, \, b+d+j) \ = \ \sigma_t(z)(i, \, b+d-m_2+j),
\end{align*}
and if $i+r\geq m_1$ then again
\begin{align*}
x(i, \, a+c+j) \ &= \ \sigma_r(y)(i, \, b+c+j) \ = \ y(i+r-m_1, \, b+c-m_2+j) \\
&= \ \sigma_s(z)(i+r-m_1, \, b+d-m_2+j) \ = \ \sigma_r\big(\sigma_s(z)\big)(i, \, b+d+j) \\ &= \ \sigma_t(z)(i, \, b+d-m_2+j).
\end{align*}
So $x(i, \, a+c+j)=\sigma_t(z)(i, \, b+d-m_2+j)$ for all $i\in m_1$ and $j\in\Z$, showing $x\mathrel{\E_{\m}}z$ in this case as well to complete the proof.
\end{proof}

\begin{lemma}\label{lem:Eequiv'}
For all $\m\in\Z_+\times-\N$, $\mathcal{E}_{\m}'$ is an equivalence relation on $(2^\omega)^{m_1\times\omega}$.
\end{lemma}

\begin{proof}
Again, the proof is identical to that of \ref{lem:Eequiv} once we let $j$ vary over $\omega$ instead of $\Z$, replace $\sigma$ with $\sigma'$ throughout, and use Lemma \ref{lem:sigma'} in place of Lemma \ref{lem:sigma}.
\end{proof}

The following facts, which we will not need, are easy to prove from Lemmas \ref{lem:sigma} and \ref{lem:Eequiv}, and might help to clarify the meaning of $\sigma_r$ and $\mathcal{E}_{\m}$. Analogous statements hold for $\sigma_r'$ and $\mathcal{E}_{\m}'$.

\begin{corollary}
For all $\m\in\Z_+^2$, $r\in m_1$, and $x\in (2^\omega)^{m_1\times\Z}$, we have $x\mathrel{\mathcal{E}}_{\m}\sigma_{r,\m}(x)$. In particular, the map $\sigma_{r,\m}$ is a Borel reduction from $\mathcal{E}_{\m}$ to $\mathcal{E}_{\m}$. \hfill $\square$
\end{corollary}


\begin{corollary}
For all $\m\in\Z_+^2$ and $x,y\in (2^\omega)^{m_1\times\Z}$, we have $x\mathrel{\mathcal{E}_{\m}}y$ if and only if there exist $r,s\in m_1$ and $a,b\in\N$ such that for all $i\in m_1$ and $j\in\Z$, $\sigma_r(x)(i, a+j) = \sigma_s(y)(i, b+j)$. \hfill $\square$
\end{corollary}

Now we will show that for each $\m\in\Z_+^2$ and $\m'\in\Z_+\times-\N$, the equivalence relations $\E_{\m}$ and $\E_{\m'}'$ are hypersmooth. In order to do this, we first define functions that code the ``vertical strips'' from Figures \ref{fig:Case2} and \ref{fig:Case3} into single columns, row by row. 

\begin{definition}
For each $m\in\Z_+$, let $\alpha_m:(2^\omega)^{m}\to 2^\omega$ be a fixed Borel bijection, and define the Borel bijections $\hat{\alpha}_m:(2^\omega)^{m\times\Z}\to (2^\omega)^{\Z}$ and $\hat{\alpha}_m':(2^\omega)^{m\times\omega}\to (2^\omega)^{\omega}$ by the (same) rules 
\begin{align*}
\hat{\alpha}_m(x)(j) \ &:= \ \alpha_m\Big(x(0,j),\ldots, x(m-1, j)\Big) \quad \text{for all $x\in(2^\omega)^{m\times\Z}$ and $j\in\Z$,} \quad \text{and} \\
\hat{\alpha}_m'(x)(j) \ &:= \ \alpha_m\Big(x(0,j),\ldots, x(m-1, j)\Big) \quad \text{for all $x\in(2^\omega)^{m\times\omega}$ and $j\in\omega$.}
\end{align*}
\end{definition}

Next, given $\m\in\Z_+^2$ and $\m'\in\Z_+\times-\N$, let $\F_{\m}$ and $\F_{\m'}'$ be the Borel equivalence relations on $(2^\omega)^\Z$ and $(2^\omega)^\omega$, respectively, defined by
\begin{equation}
x\mathrel{\F_{\m}}y \ :\Longleftrightarrow \ \hat{\alpha}_{m_1}^{-1}(x)\mathrel{\mathcal{E}_{\m}}\hat{\alpha}_{m_1}^{-1}(y) \qquad \text{and} \qquad 
x\mathrel{\F_{\m'}'}y \ :\Longleftrightarrow \ \hat{\alpha}_{m_1'}'^{-1}(x)\mathrel{\mathcal{E}_{\m'}'}\hat{\alpha}_{m_1'}'^{-1}(y),
\end{equation}
so $\mathcal{E}_{\m}\cong_B\F_{\m}$ via $\hat{\alpha}_{m_1}$ and $\mathcal{E}_{\m'}'\cong_B\F_{\m'}'$ via $\hat{\alpha}_{m_1}'$.

\begin{lemma}\label{lem:EmHF}
For all $\m\in\Z_+^2$, the equivalence relation $\E_{\m}$ is hyperfinite.
\end{lemma}

\begin{proof}
Fix $\m\in\Z_+^2$, and write $E_{\Z}$ for the hyperfinite orbit equivalence relation arising from the shift action of $\Z$ on $(2^\omega)^\Z$. We show first that $E_{\Z}\subseteq\F_{\m}$. To see this, fix $x,y\in (2^\omega)^{m_1\times\Z}$ such that $\hat{\alpha}_{m_1}(x)\mathrel{E_{\Z}}\hat{\alpha}_{m_1}(y)$, along with $k\in\Z$ such that $k\cdot\hat{\alpha}_{m_1}(x)=\hat{\alpha}_{m_1}(y)$. Then for all $j\in\Z$ we have
\[
\alpha_{m_1}\Big(x(0,k+j),\ldots,x(m_1-1, k+j)\Big) \ = \ \hat{\alpha}_{m_1}(x)(k+j) \ = \ \hat{\alpha}_{m_1}(y)(j) \ = \ \alpha_{m_1}\Big(y(0,j),\ldots,y(m_1-1, j)\Big),
\]
which means $x(i,k+j)=y(i,j)$ for all $i\in m_1$ and $j\in\Z$. This shows $x\mathrel{\mathcal{E}_{\m}}y$ (with $r=0$), so $\hat{\alpha}_{m_1}(x)\mathrel{\F_{\m}}\hat{\alpha}_{m_1}(y)$, and we conclude $E_{\Z}\subseteq\F_{\m}$.

Now we show $[\F_{\m}:E_{\Z}]\leq m_1$, so $\F_{\m}$ is hyperfinite by \cite[1.3(vii)]{JKL}, which means $\mathcal{E}_{\m}$ is hyperfinite as well. Recycling notation, fix $x,y\in (2^\omega)^{m_1\times\Z}$ such that $x\mathrel{\E_{\m}}y$, along with $r\in m_1$ and $a,b\in\N$ such that $x(i,a+j)=\sigma_r(y)(i,b+j)$ for all $i\in m_1$ and $j\in\Z$. Note that if $r=0$, then $a\cdot\hat{\alpha}_{m_1}(x)=b\cdot\hat{\alpha}_{m_1}(y)$, so $\hat{\alpha}_{m_1}\big(\sigma_0(x)\big)\mathrel{E_{\Z}}\hat{\alpha}_{m_1}(y)$. Assuming $r\ne 0$, we put $\tilde{r}:=m_1-r$ and $\tilde{x}:=\sigma_{\tilde{r}}(x)$, and show $\hat{\alpha}_{m_1}(\tilde{x})\mathrel{E_{\Z}}\hat{\alpha}_{m_1}(y)$. To show this, let $i\in m_1$ and $j\in\Z$ be arbitrary. If $\tilde{r}+i<m_1$, then
\[
\tilde{x}(i, \, a+j+m_2) \ = \ x(\tilde{r}+i, \, a+j+m_2) \ = \ \sigma_r(y)(\tilde{r}+i, \, b+j+m_2) \ = \ y(i, \, b+j),
\]
and if $\tilde{r}+i\geq m_1$, then again
\[
\tilde{x}(i, \, a+j+m_2) \ = \ x(\tilde{r}+i-m_1, \, a+j) \ = \ \sigma_r(y)(\tilde{r}+i-m_1, \, b+j) \ = \ y(i, \, b+j).
\]
This implies $(a+m_2)\cdot\hat{\alpha}_{m_1}(\tilde{x})=b\cdot\hat{\alpha}_{m_1}(y)$, so $\hat{\alpha}_{m_1}(\tilde{x})\mathrel{E_{\Z}}\hat{\alpha}_{m_1}(y)$ as claimed. It follows that for each $x\in (2^\omega)^{m_1\times\Z}$, 
\[
[\hat{\alpha}_{m_1}(x)]_{\F_{\m}}\subseteq\bigcup_{r\mkern1mu\in\mkern1mu m_1}\big[ \hat{\alpha}_{m_1}\big(\sigma_r(x)\big) \big]_{E_{\Z}},
\]
showing $[\F_{\m}:E_{\Z}]\leq m_1$ to complete the proof.
\end{proof}

\begin{lemma}\label{lem:Em'HS}
For all $\m\in\Z_+\times-\N$, the equivalence relation $\E_{\m}'$ is hypersmooth.
\end{lemma}

\begin{proof}
This proof is essentially identical to that of Lemma \ref{lem:EmHF}, with $E_{1,t}$, $\F_{\m}'$, and $\hat{\alpha}_{m_1}'$ taking the place of $E_{\Z}$, $F_{\m}$, and $\hat{\alpha}_{m_1}$, respectively, and $j$ ranging over $\omega$ instead of $\Z$. Thus $E_{1,t}$ is a sub-equivalence relation of $\F_{\m}'$ of index at most $m_1$, and we use Corollary \ref{cor:fi1} in place of \cite[1.3(vii)]{JKL} to conclude that $\F_{\m}'$, and consequently $\E_{\m}'$, is hypersmooth.
\end{proof}

Now we are ready to finish showing that $\EOT\res\POT$ is hypersmooth. We handle case (iii) first, since its proof avoids a slight technical complication that shows up in case (ii), although otherwise the proofs are similar.

\begin{proposition}\label{prop:Case3}
For every $\m\in\Z_+\times-\N$, $\EOT\res\textup{Per}_{\langle\m\rangle}(\OT)$ is hypersmooth.
\end{proposition}

\begin{proof}
Fix $\m\in\Z_+\times-\N$ and write $H=\langle\m\rangle\leq\Z^2$. For each $x\in\textup{Per}_H(\OT)$, let $\bar{c}^x=(c^x_1,c^x_2)\in\N^2$ be least such that $\bar{c}^x\cdot x$ is strongly $\m$-periodic, so $x\mapsto\bar{c}^x$ is Borel. Fix a Borel injection $\mu:\Omega\to 2^\omega$, and let $f:\textup{Per}_H(\OT)\to (2^\omega)^{m_1\times\omega}$ be the Borel map defined by
\[
f(x)(i,j) \ = \ \mu\Big((\bar{c}^x\cdot x)(i,j)\Big) \qquad \text{for all $i\in m_1$ and $j\in\omega$.}
\]
We claim that $f$ is a Borel reduction from $\EOT\res\textup{Per}_H(\OT)$ to $\mathcal{E}_{\m}'$, showing by Lemma \ref{lem:Em'HS} that $\EOT\res\textup{Per}_H(\OT)$ is hypersmooth.

To prove this, fix $x,y\in\textup{Per}_H(\OT)$, and suppose first that $x\mathrel{\EOT}y$. Fix $\bar{a},\bar{b}\in\N^2$ such that $\bar{a}\cdot x=\bar{b}\cdot y$, along with $k_0\in\N$ and $r_0\in m_1$ such that $c_1^y+a_1=k_0m_1-r_0$, and $k\in\N$ and $r\in m_1$ such that $c_1^x+b_1+r_0=km_1+r$. Then for all $i\in m_1$ and $j\in\omega$,
\begin{align*}
(\bar{c}^x\cdot x)\big(i, \, c_2^y+a_2-k_0m_2+j\big) \ &= \ (\bar{c}^x\cdot x)\big(i+r_0+k_0m_1-r_0, \, c_2^y+a_2+j\big) \\
&= \ \Big((\bar{c}^x+\bar{c}^y+\bar{a})\cdot x\Big)(i+r_0,\,j) \\
&= \ \Big((\bar{c}^x+\bar{c}^y+\bar{b})\cdot y\Big)(i+r_0,\,j) \\
&= \  (\bar{c}^y\cdot y)\big(c_1^x+b_1+i+r_0, \, c^x_2+b_2+j\big) \\
&= \ (\bar{c}^y\cdot y)\big(km_1+i+r, \, c_2^x+b_2+j\big) \\ 
&= \ (\bar{c}^y\cdot y)\big(i+r, \, c^x_2+b_2-km_2+j\big).
\end{align*}
Now, let $i\in m_1$ and $j\in\omega$ be arbitrary. If $i+r<m_1$, then 
\begin{align*}
f(x)(i,\, c_2^y+a_2-k_0m_2+j) \ &= \ \mu\Big( (\bar{c}^x\cdot x)\big(i, \, c_2^y+a_2-k_0m_2+j\big)\Big) \\
&= \ \mu\Big( (\bar{c}^y\cdot y)\big(i+r, \, c_2^x+b_2-km_2+j\big)\Big) \\
&= \ f(y)(i+r, \, c_2^x+b_2-km_2+j) \\ 
&= \ \sigma_r'\big(f(y)\big)(i, \, c_2^x+b_2-km_2+j),
\end{align*}
and if $i+r\geq m_1$, then again
\begin{align*}
f(x)(i,\, c_2^y+a_2-k_0m_2+j) \ &= \ \mu\Big( (\bar{c}^y\cdot y)\big(i+r, \, c_2^x+b_2-km_2+j\big)\Big) \\
&= \ \mu\Big( (\bar{c}^y\cdot y)\big(i+r-m_1, \, c_2^x+b_2-(k+1)m_2+j\big)\Big) \\
&= \ f(y)\big(i+r-m_1, \, c_2^x+b_2-(k+1)m_2+j\big) \\ 
&= \ \sigma_r'\big(f(y)\big)(i, \, c_2^x+b_2-km_2+j).
\end{align*}
This shows $f(x)\mathrel{\mathcal{E}_{\m}'}f(y)$.

Conversely, suppose $f(x)\mathrel{\mathcal{E}_{\m}'}f(y)$, and fix $r\in m_1$ and $a,b\in\N$ such that for all $i\in m_1$ and $j\in\omega$, 
\[
f(x)(i, \,a+j) \, = \, \sigma_r'\big(f(y)\big)(i, \, b+j) \, = \, \begin{cases} f(y)(i+r, \, b+j) & \text{if $i+r<m_1$;} \\ f(y)(i+r-m_1, \, b+j-m_2) & \text{if $i+r\geq m_1$.} \end{cases}
\]
Let $\bar{i}=(i_1,i_2)\in\omega^2$ be arbitary, and fix $i\in m_1$ and $k\in\N$ such that $i_1=i+km_1$. Let $j:=i_2-km_2\in\omega$. If $i+r<m_1$, then
\begin{align*}
(\c^x\cdot x)(i_1, \, a+i_2) \ &= \ (\c^x\cdot x)(i+km_1, \, a+j+km_2) \ = \ (\c^x\cdot x)(i, \, a+j) \ = \ (\bar{c}^y\cdot y)(i+r, \, b+j) \\
&= \ (\c^y\cdot y)(i_1+r-km_1, \, b+i_2-km_2) \ = \ (\c^y\cdot y)(i_1+r,\,b+i_2),
\end{align*}
and if $i+r\geq m_1$, then again
\begin{align*}
(\c^x\cdot x)(i_1, \, a+i_2) \ &= \ (\c^x\cdot x)(i, \, a+j) \ = \ (\bar{c}^y\cdot y)(i+r-m_1, \, b+j-m_2) \ = \ (\bar{c}^y\cdot y)(i+r, \, b+j) \\ 
&= \ (\c^y\cdot y)(i_1+r,\,b+i_2).
\end{align*}
Since $\bar{i}\in\omega^2$ was arbitrary, it follows that $(c^x_1,\, c^x_2+a)\cdot x = (c^y_1+r,\, c^y_2+b)\cdot y$, that is, $x\mathrel{\EOT}y$ as desired.
\end{proof}

Recall that the Borel equivalence relation $E$ is \emph{essentially hyperfinite} if $E$ Borel reduces to a hyperfinite equivalence relation.

\begin{proposition}\label{prop:Case2}
For every $\m\in\Z_+^2$, $\EOT\res\textup{Per}_{\langle\m\rangle}(\OT)$ is essentially hyperfinite.
\end{proposition}

\begin{proof}
Fix $\m\in\Z_+^2$ and write $H=\langle\m\rangle\leq\Z^2$. Our proof will mirror that of Proposition \ref{prop:Case3}, although we will need to tweak our definition of $f$, since now $f(x)(i,j)$ must be defined for $j\in\Z$ rather than just $j\in\omega$. Again fix a Borel injection $\mu:\Omega\to 2^\omega$, and for each $x\in\textup{Per}_H(\OT)$ let $\bar{c}^x=(c^x_1,c^x_2)\in\N^2$ be least such that $\bar{c}^x\cdot x$ is strongly $\m$-periodic. Now define the Borel map $f:\textup{Per}_H(\OT)\to (2^\omega)^{m_1\times\Z}$ by 
\[
f(x)(i,j) \ = \ \mu\Big((\bar{c}^x\cdot x)(nm_1+i, \, nm_2+j)\Big) \quad \text{where $n\in\omega$ is least such that $nm_2+j\in\omega$.}
\]
We claim that $f$ is a Borel reduction from $\EOT\res\textup{Per}_H(\OT)$ to $\mathcal{E}_{\m}$, showing by Lemma \ref{lem:EmHF} that $\EOT\res\textup{Per}_H(\OT)$ is essentially hyperfinite, in particular hypersmooth.

To see this, fix $x,y\in\textup{Per}_H(\OT)$, and suppose first that $x\mathrel{\EOT}y$. Fix $\bar{a},\bar{b}\in\N^2$ such that $\bar{a}\cdot x=\bar{b}\cdot y$, along with $k_0\in\N$ and $r_0\in m_1$ such that $c_1^y+a_1=k_0m_1-r_0$, and $k\in\N$ and $r\in m_1$ such that $c_1^x+b_1+r_0=km_1+r$. Let $i\in m_1$ and $j\in\Z$ be arbitrary. Fix $n_0,n_1,n_1'\in\omega$ least such that
\begin{align*}
n_0m_2+c^y_2+a_2-k_0m_2+j \ &\geq \ 0, \\
n_1m_2+c^x_2+b_2-km_2+j \ &\geq \ 0, \\
\text{and} \qquad n_1'm_2+c^x_2+b_2-(k+1)m_2+j \ &\geq \ 0.
\end{align*}
Also fix an integer $n\geq\max\{n_0,n_1,n_1'\}$ large enough so that $nm_2+j\in\omega$. Then
\begin{align*}
(\bar{c}^x\cdot x)\big(nm_1+i, \, nm_2+c_2^y+a_2-k_0m_2+j\big) \ &= \ (\bar{c}^x\cdot x)\big(nm_1+i+r_0+k_0m_1-r_0, \, nm_2+c_2^y+a_2+j\big) \\
&= \ \Big((\bar{c}^x+\bar{c}^y+\bar{a})\cdot x\Big)(nm_1+i+r_0,\,nm_2+j) \\
&= \ \Big((\bar{c}^x+\bar{c}^y+\bar{b})\cdot y\Big)(nm_1+i+r_0,\,nm_2+j) \\
&= \  (\bar{c}^y\cdot y)\big(nm_1+c_1^x+b_1+i+r_0, \, nm_2+c^x_2+b_2+j\big) \\
&= \ (\bar{c}^y\cdot y)\big(nm_1+km_1+i+r, \, nm_2+c_2^x+b_2+j\big) \\ 
&= \ (\bar{c}^y\cdot y)\big(nm_1+i+r, \, nm_2+c^x_2+b_2-km_2+j\big).
\end{align*}
Now, if $i+r<m_1$ then 
\begin{align*}
f(x)(i,\, c_2^y+a_2-k_0m_2+j) \ &= \ \mu\Big( (\bar{c}^x\cdot x)\big(n_0m_1+i, \, n_0m_2+c_2^y+a_2-k_0m_2+j\big)\Big) \\
&= \ \mu\Big( (\bar{c}^x\cdot x)\big(nm_1+i, \, nm_2+c_2^y+a_2-k_0m_2+j\big)\Big) \\
&= \ \mu\Big( (\bar{c}^y\cdot y)\big(nm_1+i+r, \, nm_2+c_2^x+b_2-km_2+j\big)\Big) \\
&= \ \mu\Big( (\bar{c}^y\cdot y)\big(n_1m_1+i+r, \, n_1m_2+c_2^x+b_2-km_2+j\big)\Big) \\
&= \ f(y)(i+r, \, c_2^x+b_2-km_2+j) \\ 
&= \ \sigma_r\big(f(y)\big)(i, \, c_2^x+b_2-km_2+j),
\end{align*}
and if $i+r\geq m_1$ then again
\begin{align*}
f(x)(i,\, c_2^y+a_2-k_0m_2+j) \ &= \ \mu\Big( (\bar{c}^y\cdot y)\big(nm_1+i+r, \, nm_2+c_2^x+b_2-km_2+j\big)\Big) \\
&= \ \mu\Big( (\bar{c}^y\cdot y)\big(n_1'm_1+i+r-m_1, \, n_1'm_2+c_2^x+b_2-(k+1)m_2+j\big)\Big) \\
&= \ f(y)(i+r-m_1, \, c_2^x+b_2-(k+1)m_2+j) \\ 
&= \ \sigma_r\big(f(y)\big)(i, \, c_2^x+b_2-km_2+j).
\end{align*}
Since $i\in m_1$ and $j\in\Z$ were arbitrary, this shows $f(x)\mathrel{\mathcal{E}_{\m}}f(y)$.

Conversely, suppose $f(x)\mathrel{\mathcal{E}_{\m}}f(y)$, and fix $r\in m_1$ and $a,b\in\N$ such that for all $i\in m_1$ and $j\in\Z$,
\begin{equation}\label{eq:ugly1}
f(x)(i, \,a+j) \, = \, \sigma_r\big(f(y)\big)(i, \, b+j) \, = \, \begin{cases} f(y)(i+r, \, b+j) & \text{if $i+r<m_1$;} \\ f(y)(i+r-m_1, \, b+j-m_2) & \text{if $i+r\geq m_1$.} \end{cases}
\end{equation}
Let $\bar{i}=(i_1,i_2)\in\omega^2$ be arbitrary, and fix $i\in m_1$ and $k\in\N$ such that $i_1=i+km_1$. Let $j:=i_2-km_2\in\Z$. Fix $n_0,n_1\in\omega$ least such that
\[
n_0m_2+a+j \geq 0 \qquad \text{and} \qquad n_1m_2+b+j\geq 0,
\]
and let $n:=\max\{k,n_0,n_1\}$ and $n':=n-k\in\N$. Then
\begin{align}\label{eq:ugly2}
f(x)(i, \, a+j) \ &= \ \mu\Big((\bar{c}^x\cdot x)(n_0m_1+i, \, n_0m_2+a+j)\Big) \ = \ \mu\Big((\bar{c}^x\cdot x)(nm_1+i, \, nm_2+a+j)\Big) \\\notag
&= \ \mu\Big((\bar{c}^x\cdot x)(n'm_1+i_1, \, n'm_2+a+i_2)\Big) \ = \ \mu\Big((\bar{c}^x\cdot x)(i_1, \, a+i_2)\Big).
\end{align}
Furthermore, if $i+r<m_1$ then
\begin{align*}
f(y)(i+r, \, b+j) \ &= \ \mu\Big((\bar{c}^y\cdot y)(n_1m_1+i+r, \, n_1m_2+b+j)\Big) \ = \ \mu\Big((\bar{c}^y\cdot y)(nm_1+i+r, \, nm_2+b+j)\Big) \\
&= \ \mu\Big((\bar{c}^y\cdot y)(n'm_1+i_1+r, \, n'm_2+b+i_2)\Big) \ = \ \mu\Big((\bar{c}^y\cdot y)(i_1+r, \, b+i_2)\Big),
\end{align*}
and similarly if $i+r\geq m$ then
\begin{align*}
f(y)(i+r-m_1,\, b+j-m_2) \ &= \ \mu\Big((\bar{c}^y\cdot y)\big((n_1+1)m_1+i+r-m_1, \, (n_1+1)m_2+b+j-m_2\big) \Big) \\
&= \ \mu\Big((\bar{c}^y\cdot y)(n_1m_1+i+r, \, n_1m_2+b+j)\Big) \ = \ \mu\Big((\bar{c}^y\cdot y)(i_1+r, \, b+i_2)\Big).
\end{align*}
Combining the previous two calculations with Equation \eqref{eq:ugly1} gives
\[
f(x)(i,\, a+j) \ = \ \mu\Big((\bar{c}^y\cdot y)(i_1+r, \, b+i_2)\Big),
\]
which together with Equation \eqref{eq:ugly2} implies 
\[
(\c^x\cdot x)(i_1, \, a+i_2) \ = \ (\c^y\cdot y)(i_1+r, \, b+i_2).
\]
Since $\bar{i}\in\omega^2$ was arbitrary, we conclude that $(c^x_1,\, c^x_2+a)\cdot x = (c^y_1+r,\, c^y_2+b)\cdot y$, that is, $x\mathrel{\EOT}y$ as desired.
\end{proof}

\begin{proposition}\label{prop:FinalCase}
For all $m\in\Z_+$, $\EOT\res\textup{Per}_{\langle m\bar{e}_2\rangle}(\OT)$ is hypersmooth.
\end{proposition}

\begin{proof}
Here we can use symmetry to reduce to case (iii).  Formally, for each $i\in\{1,2\}$ let $X_i$ be the set of all $x\in\POT$ such that $H_x=\langle m\bar{e}_i\rangle$ for some $m\in\Z_+$. We have already shown $\EOT\res X_1$ is hypersmooth in Proposition \ref{prop:Case3}, and $\delta:X_2\to X_1$ is a Borel reduction from $\EOT\res X_2$ to $\EOT\res X_1$.
\end{proof}

\begin{theorem}\label{thm:Periodic}
For every standard Borel space $\Omega$, $\EOT\res\POT$ is hypersmooth.
\end{theorem}

\begin{proof}
By Lemmas \ref{lem:HBorel} and \ref{lem:piecewiseHS}, it suffices to show $\EOT\res\textup{Per}_H(\OT)$ is hypersmooth for each nontrivial subgroup $H\leq\Z^2$. This is obvious when $H=\Z^2$, and all other cases have been covered in Propositions \ref{prop:Case1}, \ref{prop:Case3}, \ref{prop:Case2}, and \ref{prop:FinalCase}.
\end{proof}


\section{The graphs $\Ge$, $\Gee$, and $\Gp$}\label{sec:ep}

In this section, we consider the graphs $\Ge$, $\Gee$, and $\Gp$ on two-dimensional product spaces defined over a finite base space, which we always take to be an integer $J\geq 2$. These graphs are the symmetrizations of the directed graphs $\vec{\G}_{\bar{e}_1}$, $\vec{\G}_{\vec{e}_2}$, and $\vec{\G}_{\p}$, which are defined over any base space by
\[
x\mathrel{\vec{\G}_{\bar{e}_1}}y \ :\Longleftrightarrow \ \bar{e}_1\cdot x=y, \qquad 
x\mathrel{\vec{\G}_{\bar{e}_2}}y \ :\Longleftrightarrow \ \bar{e}_2\cdot x=y, \qquad \text{and} \qquad
x\mathrel{\vec{\G}_{\p}}y \ :\Longleftrightarrow \ \p\cdot x=y, 
\]
where
\[
\bar{e}_1 \ := \ \begin{bmatrix} 1 \\ 0 \end{bmatrix}, \qquad 
\bar{e}_2 \ := \ \begin{bmatrix} 0 \\ 1 \end{bmatrix}, \qquad \text{and} \qquad
\p \ := \ \begin{bmatrix} 1 \\ 1 \end{bmatrix}.
\]
Should we ever need to emphasize the base space $J$, we will write, e.g., $\G_{\bar{e}_1}(J)$ for the graph $\G_{\bar{e}_1}$ defined on $\JT$, but usually the value of $J$ is unimportant and we suppress it from the notation when convenient.

We already know $\chi_B(\Ge)=\aleph_0$ by the proof of Observation \ref{claim:NoColoring}, and similar arguments show $\chi_B(\Gee)=\chi_B(\Gp)=\aleph_0$ as well. Nevertheless, we will attempt to find large $\EJT$-invariant Borel subsets of $\JT$ on which these graphs admits finite Borel colorings. Specifically, we prove the following:

\begin{theorem}\label{thm:e1}
For each $J\in\Z_{\geq 2}$ and $\bar{v}\in\{\bar{e}_1,\bar{e}_2,\p\}$, there is a partition of $\neg\textup{Per}_{\bar{v}}(\JT)$ into a pair of $\EJT$-invariant Borel sets $X_{\bar{v}}$ and $W_{\bar{v}}$ such that $\chi_B\big(\G_{\bar{v}}(J)\res X_{\bar{v}}\big)\leq 3$ and $\EJT\res W_{\bar{v}}$ is hypersmooth.
\end{theorem}

Recall that $\delta:\JT\to\JT$ is the coordinate-swapping involution defined by $\delta(x)(i,j):=x(j,i)$ for all $x\in\JT$. By $\delta$-symmetry, the arguments we apply to $\Ge$ also apply to $\Gee$, and we will see that similar arguments can be adapted to treat $\Gp$.  Most of our attention, therefore, will be devoted to $\Ge$. The basic idea of the proof of Theorem \ref{thm:e1} in this case is that if $x\in\JT$ has infinitely many rows that are not eventually constant (which is an $\EJT$-invariant condition), then we can color $\Ge$ by applying a Borel $3$-coloring of $\G_1(J)$ to the least such row; otherwise, if almost all rows of $x$ are eventually constant, then as we move up through the rows in $x$, the sequence of first points beyond which that row is constant must be unbounded to the right since $x$ is not $\bar{e}_1$-periodic. Thus $x$ determines a sequence of points in $[x]_{\EJT}$, darkly shaded in Figure \ref{fig:Ge1} below, and we shall argue that if $x$ and $y$ are $\EJT$-equivalent then they produce the same sequence through their $\EJT$-class, up to tail equivalence.

\begin{figure}[h]
\begin{center}
\begin{tikzpicture}[x=1.2cm, y=1.22cm]
     \fill [fill=gray!40] (.05+1*.3,.07+1*.35) rectangle (.35+1*.3,.4+1*.35); 
     \fill [fill=gray!40] (.05+4*.3,.07+4*.35) rectangle (.35+4*.3,.4+4*.35); 
     \fill [fill=gray!40] (.05+5*.3,.07+6*.35) rectangle (.35+5*.3,.4+6*.35); 
     \fill [fill=gray!40] (.05+6*.3,.07+7*.35) rectangle (.35+6*.3,.4+7*.35); 
     \fill [fill=gray!40] (.05+9*.3,.07+9*.35) rectangle (.35+9*.3,.4+9*.35); 
     \fill [fill=gray!15] (.05+2*.3,.07+1*.35) rectangle (.35+21*.3,.4+1*.35); 
     \fill [fill=gray!15] (.05+1*.3,.07+2*.35) rectangle (.35+21*.3,.4+2*.35); 
     \fill [fill=gray!15] (.05+0*.3,.07+3*.35) rectangle (.35+21*.3,.4+3*.35); 
     \fill [fill=gray!15] (.05+5*.3,.07+4*.35) rectangle (.35+21*.3,.4+4*.35); 
     \fill [fill=gray!15] (.05+2*.3,.07+5*.35) rectangle (.35+21*.3,.4+5*.35); 
     \fill [fill=gray!15] (.05+6*.3,.07+6*.35) rectangle (.35+21*.3,.4+6*.35); 
     \fill [fill=gray!15] (.05+7*.3,.07+7*.35) rectangle (.35+21*.3,.4+7*.35); 
     \fill [fill=gray!15] (.05+2*.3,.07+8*.35) rectangle (.35+21*.3,.4+8*.35); 
     \fill [fill=gray!15] (.05+10*.3,.07+9*.35) rectangle (.35+21*.3,.4+9*.35);      
     \draw (.05+1*.3,.07+1*.35) rectangle (.35+1*.3,.4+1*.35); 
     \draw (.05+4*.3,.07+4*.35) rectangle (.35+4*.3,.4+4*.35); 
     \draw (.05+5*.3,.07+6*.35) rectangle (.35+5*.3,.4+6*.35); 
     \draw (.05+6*.3,.07+7*.35) rectangle (.35+6*.3,.4+7*.35); 
     \draw (.05+9*.3,.07+9*.35) rectangle (.35+9*.3,.4+9*.35); 
 \draw [thick] (0,0) -- (6.7,0);
 \draw [thick] (0,0) -- (0,3.6);
  \node at (0,0) [circle,fill,inner sep=1.2pt]{};
 \foreach \x in {0,.6,...,6.6} \node at (.2+\x,.25) {$0$}; \foreach \x in {0,.6,...,6.6} \node at (.5+\x,.25) {$1$};
 \node at (.2,.6) {$1$}; \foreach \x in {.3,.6,...,6.6} \node at (.2+\x,.6) {$0$};
 \node at (.2,.95) {$0$}; \foreach \x in {.3,.6,...,6.6} \node at (.2+\x,.95) {$1$};
 \foreach \x in {0,.3,...,6.6} \node at (.2+\x,1.3) {$0$};
 \node at (.2,1.65) {$0$}; \node at (.5,1.65) {$0$}; \node at (.8,1.65) {$1$}; \node at (1.1,1.65) {$0$}; \foreach \x in {1.2,1.5,...,6.6} \node at (.2+\x,1.65) {$1$};
 \node at (.2,2) {$1$}; \node at (.5,2) {$0$}; \foreach \x in {.6,.9,...,6.3} \node at (.2+\x,2) {$1$};
 \node at (.2,2.35) {$1$}; \node at (.5,2.35) {$0$}; \node at (.8,2.35) {$0$}; \node at (1.1,2.35) {$1$}; \node at (1.4,2.35) {$1$}; \foreach \x in {1.5,1.8,...,6.6} \node at (.2+\x,2.35) {$0$};
 \node at (.2,2.7) {$0$}; \node at (.5,2.7) {$1$}; \node at (.8,2.7) {$0$}; \node at (1.1,2.7) {$1$}; \node at (1.4,2.7) {$0$}; \node at (1.7,2.7) {$0$}; \foreach \x in {1.8,2.1,...,6.6} \node at (.2+\x,2.7) {$1$};
 \node at (.2,3.05) {$1$}; \node at (.5,3.05) {$1$}; \foreach \x in {.6,.9,...,6.3} \node at (.2+\x,3.05) {$0$};
  \node at (.2,3.4) {$0$}; \node at (.5,3.4) {$1$}; \node at (.8,3.4) {$1$}; \node at (1.1,3.4) {$0$}; \node at (1.4,3.4) {$0$}; \node at (1.7,3.4) {$1$}; \node at (2,3.4) {$0$}; \node at (2.3,3.4) {$0$}; \node at (2.6,3.4) {$0$}; \foreach \x in {2.7,3,...,6.6} \node at (.2+\x,3.4) {$1$};
 \node at (-.15,-.15) {$x$};
 \node at (7,1.8) {$\cdots$};
 \node at (3.2,3.9) {$\vdots$};
\end{tikzpicture}
\end{center}
\caption{A sample $x\in\We(2)$, with constant tails of rows shaded.}
\label{fig:Ge1}
\end{figure}
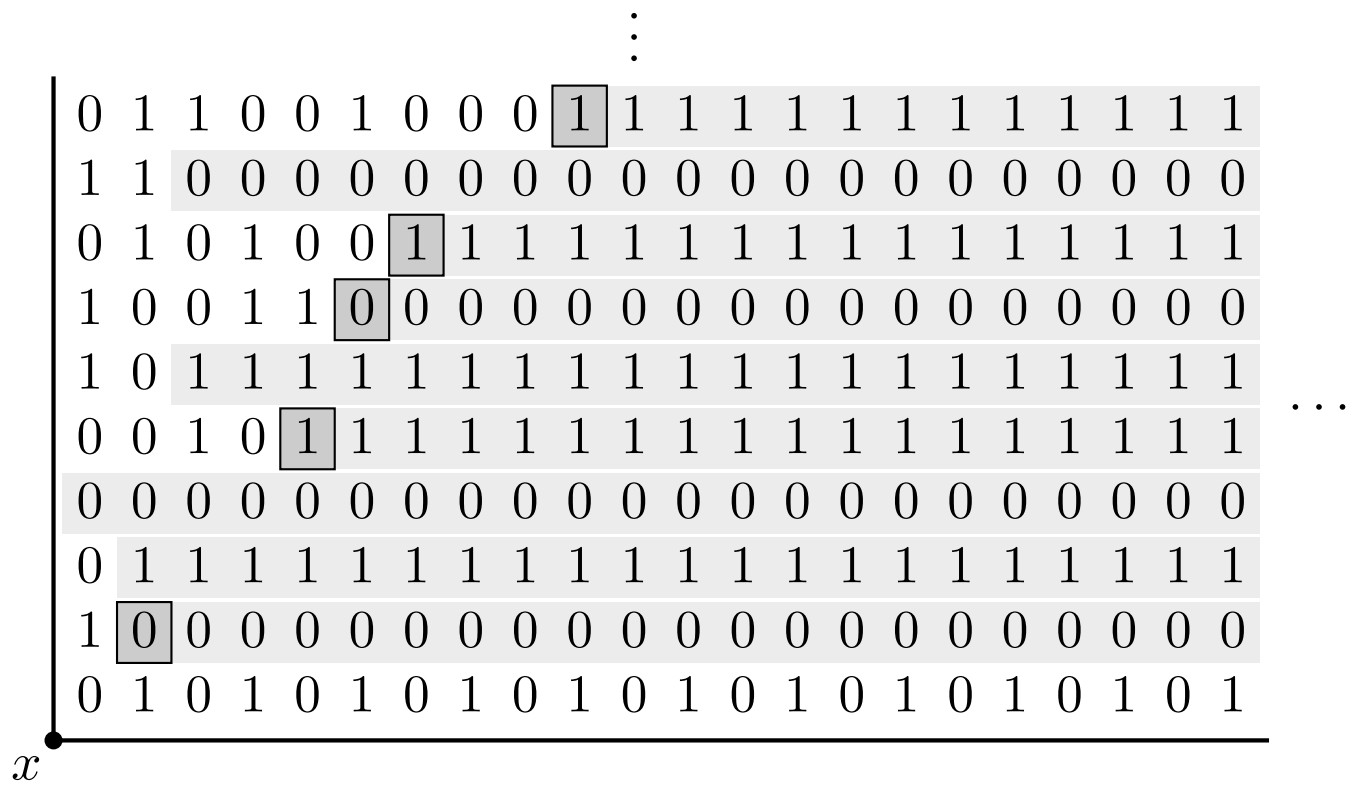

We turn now to the details of the argument, beginning with definitions of the sets $\Xe$ and $\We$.

\begin{definition}\label{df:XeWe}
For each $J\in\Z_{\geq 2}$, let
\[
\Xe \ = \ \Xe(J) \ := \ \Big\{ x\in\neg\textup{Per}_{\bar{e}_1}(J^{\omega^2})\st (\exists^{\infty}j\in\omega) \text{  $\langle x(i,j)\st i\in\omega\rangle$ is not eventually constant}\Big\},
\]
and let
\begin{align*}
\We \ = \ \We(J) \ :\!&= \ \neg\textup{Per}_{\bar{e}_1}(J^{\omega^2})\setminus\Xe(J) \\
&= \ \Big\{ x\in\neg\textup{Per}_{\bar{e}_1}(J^{\omega^2})\st (\forall^{\infty}j\in\omega) \text{  $\langle x(i,j)\st i\in\omega\rangle$ is eventually constant}\Big\}.
\end{align*}
\end{definition}

Note that $\Xe$ and $\We$ are Borel and $\EJT$-invariant.

\begin{proposition}\label{prop:ColoringGe1}
For every $J\in\Z_{\geq 2}$, we have $\chi_B\big(\Ge\res\Xe\big)=3$.
\end{proposition}

\begin{proof}
Fix $J\in\Z_{\geq 2}$. 
Using Fact \ref{fact:KST5.2}, let $c_J:J^\omega\to 3$ be a Borel coloring of $\G_1(J)$. For each $x\in\Xe$, let $j(x)$ be the least $j\in\omega$ such that the sequence $x_j:=\langle x(i,j)\st i\in\omega\rangle$ is not eventually constant.  Noting that the map $x\mapsto j(x)$ is Borel, define the Borel function $c_{\bar{e}_1}:\Xe\to 3$ by $c_{\bar{e}_1}(x):=c_J\big(x_{j(x)}\big)$. Let $x,y\in\Xe$ be arbitrary such that $x\mathrel{\Ge}y$. Then $j(x)=j(y)$ and $x_{j(x)}\mathrel{\G_1}y_{j(y)}$, which means $c_J(x_{j(x)})\ne c_J(y_{j(y)})$ and therefore $c_{\bar{e}_1}(x)\ne c_{\bar{e}_1}(y)$. This shows $\chi_B\big(\Ge\res\Xe\big)\leq 3$, and since $\G_1(J)$ embeds into $\Ge\res\Xe$ and $\chi_B\big(\G_1(J)\big)=3$ (Fact \ref{fact:KST5.2}), we also have $\chi_B\big(\Ge\res\Xe\big)\geq 3$, so $\chi_B\big(\Ge\res\Xe\big)=3$.
\end{proof}

By tweaking the proof of Fact \ref{fact:KST5.2} that is given in \cite{KST}, one can show that the restriction of $\G_1(J)$ to the set of non-constant sequences in $J^{\omega}$ admits a continuous $3$-coloring. On the other hand, the map $x\mapsto j(x)$ in the proof of Proposition \ref{prop:ColoringGe1} is Borel, but not continuous. Finite Borel colorings suffice for our purposes by Theorem \ref{thm:connection} and Corollary \ref{cor:main}, but since continuous colorings lead to clopen marker sets, we ask:

\begin{question}
Does there exist a continuous finite coloring of $\Ge\res\Xe$?
\end{question}

Towards showing $\EJT\res\We$ is hypersmooth, we distinguish certain subsets of $\We$. Define
\begin{equation}
\We' \ = \ \We'(J) \ := \ \Big\{x\in\We(J) \st (\forall j\in\omega) \text{  $\langle x(i,j)\st i\in\omega\rangle$ is eventually constant}\Big\}
\end{equation}
and
\begin{equation}
\We'' \ = \ \We''(J) \ := \ \Big\{x\in\We'(J) \st \text{  $\langle x(i,0)\st i\in\omega\rangle$ is constant}\Big\}.
\end{equation}
Note that $\We'$ and $\We''$ are Borel but not $\EJT$-invariant. In Figure \ref{fig:Ge1}, assuming every row above those pictured is eventually constant, we have $x\in\We\setminus\We'$ and $\bar{e}_2\cdot x\in\We'\setminus\We''$ and $\d{1}\cdot x\in\We''$. 

\begin{lemma}\label{lem:e1Easy}
For all $x\in\We$ there exists $n\in\omega$ such that: 
\begin{enumerate}[(i)]
    \item for all $j\geq n$, the sequence $\langle x(i,j)\st i\in\omega\rangle$ is eventually constant; but
    \item the sequence $\langle x(i,n)\st i\in\omega\rangle$ is not constant.
\end{enumerate}
\end{lemma}

\begin{proof}
Fix $x\in\We$. By definition of $\We$, we may fix $n_0\in\omega$ such that $\langle x(i,j)\st i\in\omega\rangle$ is eventually constant for all $j\geq n_0$. But $\We\subseteq\neg\textup{Per}_{\bar{e}_1}(\JT)$ by definition, so there must be some $n\geq n_0$ such that $\langle x(i,n)\st i\in\omega\rangle$ is non-constant or else $x$ would be $\bar{e}_1$-periodic.
\end{proof}

Using Lemma \ref{lem:e1Easy}, for each $x\in\We$ let $\nu(x)$ be the least $n\in\omega$ such that for all $j\geq n$ the sequence $\langle x(i,j)\st i\in\omega\rangle$ is eventually constant but $\langle x(i,n)\st i\in\omega\rangle$ is non-constant, and let $\mu(x)$ be the least $m\in\omega$ such that the sequence $\big\langle x\big(m+i,\nu(x)\big)\st i\in\omega\big\rangle$ is constant. Note that the maps $\nu$ and $\mu$ are Borel, that $\mu(x)\geq 1$ for all $x\in\We$, and that $\nu(x)\geq 1$ for all $x\in\We''$. Define the Borel function $\sigma:\We\to\We''$ by
\begin{equation}\label{eq:sigma}
    \sigma(x) \ := \ \big(\mu(x), \, \nu(x)\big)\cdot x.
\end{equation}
It is easy to check that $\sigma$ really does take values in $\We''$, and also that $\sigma$ is a Borel reduction from $\EJT\res\We$ to $\EJT\res\We''$. Intuitively, iterating $\sigma$ takes $x$ through the sequence of darkly shaded points in Figure \ref{fig:Ge1}. With these definitions in place, we now state and prove the key technical lemma we will use to show that both $\EJT\res\We$ and $\EJT\res\Wp$ are hypersmooth.

\begin{lemma}\label{lem:Ge1HSTechLem}
There is a Borel map $s:\We''\to (\N^2)^\omega$ that assigns to each $x\in\We''$ a sequence $\big(s_k(x)\big)_{k\in\omega}$ in $\N^2$ with the following property: for all $x,y\in\We''$ and $\bar{u},\bar{v}\in\N^2$ satisfying $\bar{u}\cdot x=\bar{v}\cdot y$, there exist $m,n\in\omega$ such that for all $k\in\omega$, $s_{m+k}(x)-\bar{u}=s_{n+k}(y)-\bar{v}\in\N^2$.
\end{lemma}

\begin{proof}
For each $x\in\We$ and $j\in\omega$, write $x_j:=\sigma^j(x)$. We define $s:\We''\to (\N^2)^\omega$ as follows. Given $x\in\We''$, first put $s_0(x)=s(x)(0):=\bar{0}$; then, for each $k\in\Z_+$ define
\[
s_k(x) \ = \ s(x)(k) \ := \ \left( \sum_{j=0}^{k-1}\mu(x_j), \; \sum_{j=0}^{k-1}\nu(x_j)\right) \ \in \ \N^2.
\]
Since the maps $\mu$, $\nu$, and $\sigma$ are Borel, so is $s$. Note that $s_k(x)\cdot x=x_k$ for all $k\in\omega$ by a straightforward induction.

To see that $s$ has the desired property, fix $x,y\in\We''$ and $\bar{u},\bar{v}\in\N^2$ such that $\bar{u}\cdot x=\bar{v}\cdot y$. Note that if we can find $m,n\in\omega$ such that $s_m(x)-\bar{u}=s_n(y)-\bar{v}\in\N^2$, then 
\[
x_m \ = \ s_m(x)\cdot x \ = \ \big(s_m(x)-\bar{u}\big)\cdot (\bar{u}\cdot x) \ = \ \big(s_n(y)-\bar{v}\big)\cdot (\bar{v}\cdot y) \ = \ s_n(y)\cdot y \ = \ y_n,
\]
which means $x_{m+k}=y_{n+k}$ and $s_{m+k}(x)-s_m(x)=s_{n+k}(y)-s_n(y)$ for all $k\in\omega$, so $s$ is as desired. Thus it will suffice to find such $m,n$. To improve notation for the remainder of the proof, write $(a_k,b_k):=s_k(x)$ and $(c_k,d_k):=s_k(y)$ for all $k\in\omega$, and additionally write $x_m<_1y_n$ and $x_m<_2y_n$ as a geometrically suggestive shorthand for $(a_m-u_1)<(c_n-v_1)$ and $(b_m-u_2)<(d_n-v_2)$, respectively, for any $m,n\in\omega$. Intuitively, if we imagine locating $x_m$ and $y_n$ in Figure \ref{fig:equivalence}, we have $x_m<_1y_n$ if $x_m$ lies to the left of $y_n$, and $x_m<_2y_n$ if $x_m$ lies below $y_n$. Likewise define, e.g., $\leq_1$ and $>_2$ and $=_2$ and other variations; this will make our arguments below more intelligible.

Now, noting that all four sequences $(a_k)$, $(b_k)$, $(c_k)$, and $(d_k)$ strictly increase, let $m_0\in\omega$ be least such that $a_{m_0}>u_1$ and $b_{m_0}>u_2$, and $n_0\in\omega$ least such that $c_{n_0}>v_1$ and $d_{n_0}>v_2$. We prove the following two claims.

\smallskip

\noindent \underline{Claim (i)}: For all integers $m\geq m_0$ and $n\geq n_0$, if $x_m=_2y_n$ then also $x_m=_1y_n$; that is, if $b_{m}-u_2=d_{n}-v_2$ then also $a_{m}-u_1=c_{n}-v_1$.

\smallskip

\noindent To prove this, fix $m\geq m_0$ and $n\geq n_0$ and assume $h:=b_{m}-u_2 = d_{n}-v_2$. Also write $t:=a_{m}-u_1-1\in\omega$, and let $x':=(a_{m}-1, \, b_{m})\cdot x$. Then $x'(0,0)\ne x'(1,0)$ but $x'(i,0)=x'(1,0)$ for all $i\geq 1$. Additionally, for each $i\in\omega$ we have
\begin{equation}\label{eq:characterization}
    x'(i,0) \ = \ \Big( (u_1+t, \, u_2+h)\cdot x\Big)(i,0) \ = \ \Big( (v_1+t, \, v_2+h)\cdot y\Big)(i,0) \ = \ \Big((v_1+t, \, d_{n})\cdot y\Big)(i,0).
\end{equation}
Since $c_{n}$ is the unique $c\in\omega$ such that
\[
\big((c-1,d_{n})\cdot y\big)(0,0)\ne \big((c-1,d_{n})\cdot y\big)(1,0) \quad \text{but}  \quad (\forall i\in\Z_+) \ \big((c-1,d_{n})\cdot y\big)(i,0) = \big((c-1,d_{n})\cdot y\big)(1,0),
\]
Equation \eqref{eq:characterization} shows that $v_1+t=c_{n}-1$, so $c_{n}-v_1=t+1=a_{m}-u_1$, establishing Claim (i).

\smallskip

\noindent \underline{Claim (ii)}: For all integers $m\geq m_0$ and $n\geq n_0$, if $x_m<_2y_n$ but $x_{m+1}>_2y_n$, then $y_n\leq_1x_m$; symmetrically, if $y_n<_2x_m$ but $y_{n+1}>_2x_m$, then $x_m\leq_1 y_n$.

\smallskip

\noindent To prove the first statement, fix $m\geq m_0$ and $n\geq n_0$, and suppose $x_m<_2y_n$; we show that $x_m<_1 y_n$ implies $x_{m+1}\leq_2 y_n$. Indeed, assuming $x_m<_1 y_n$, let
\[
\bar{r} \ := \ \Big((c_n-v_1)-(a_m-u_1)-1,\;(d_n-v_2)-(b_m-u_2)\Big) \ \in \ \N\times\Z_+,
\]
so $\bar{e}_1\cdot (\bar{r}\cdot x_m)=y_n$, and write $x':=\bar{r}\cdot x_m$. Then $x'(0,0)\ne x'(1,0)$ but $x'(i,0)=x'(1,0)$ for all $i\geq 1$ since $y_n\in\We''$. Therefore $\nu(x_m)\leq r_2$, so $x_{m+1}\leq_2 y_n$ as desired. 
The second statement follows by symmetry, concluding the proof of Claim (ii).

\smallskip

Now we are ready to finish proving the lemma. We must find $m,n\in\omega$ such that $x_m=_1y_n$ and $x_m=_2y_n$.  Assume without loss of generality that $x_{m_0}\geq_2y_{n_0}$. If $x_{m_0}=_2y_{n_0}$, then we are done by Claim (i), so assume $x_{m_0}>_2y_{n_0}$. Let $n_1\in\omega$ be least such that $x_{m_0}\leq_2y_{n_1}$, so $n_0<n_1$. 
If $x_{m_0}=_2y_{n_1}$, then again we are done by Claim (i), so assume $x_{m_0}<_2y_{n_1}$. Then Claim (ii) implies $x_{m_0}\leq_1 y_{n_1-1}$. Now let $m_1>m_0$ be least such that $x_{m_1}\geq_2 y_{n_1}$. Then $x_{m_1-1}<_2y_{n_1}$, and also $y_{n_1-1}<_2 x_{m_0}\leq_2 x_{m_1-1}$, so $x_{m_1-1}\leq_1 y_{n_1-1}<_1y_{n_1}$ by Claim (ii) and the fact that $\mu(y_{n_1-1})>0$. 
But now $x_{m_1}>_2 y_{n_1}$ would contradict Claim (ii), so in fact $x_{m_1}=_2 y_{n_1}$ and we are done by Claim (i).
\end{proof}

\begin{proposition}\label{prop:WeHS}
For every $J\in\Z_{\geq 2}$, the restriction of $\EJT$ to $\We$ is hypersmooth.
\end{proposition}

\begin{proof}
Fix $J\geq 2$, and let $s:\We''\to(\N^2)^\omega$ be as in the statement of Lemma \ref{lem:Ge1HSTechLem}. Let $\rho:\JT\to 2^\omega$ be a fixed Borel bijection, and define the Borel function $g:\We''\to (2^\omega)^\omega$ by
\[
g(x)(n) \ := \ \rho\Big(s_n(x)\cdot x\Big).
\]
We claim that $g$ is a Borel reduction from $\EJT\res\We''$ to $E_{1,t}$.

To see this, fix $x,y\in\We''$, and suppose first that $g(x)\mathrel{E_{1,t}}g(y)$. Then there exist $m,n\in\omega$ such that $s_m(x)\cdot x=s_n(y)\cdot y$. Since $x\mathrel{\EJT}s_m(x)\cdot x$ and $s_n(y)\cdot y\mathrel{\EJT}y$, it follows that $x\mathrel{\EJT}y$. Conversely, suppose $x\mathrel{\EJT}y$, and fix $\bar{u},\bar{v}\in\N^2$ such that $\bar{u}\cdot x=\bar{v}\cdot y$. Using Lemma \ref{lem:Ge1HSTechLem}, fix $m,n\in\omega$ such that $s_{m+k}(x)-\bar{u}=s_{n+k}(y)-\bar{v}\in\N^{\ell}$ for all $k\in\omega$. Then for all $k\in\omega$ we have
\begin{align*}
\rho^{-1}\Big(g(x)(m+k)\Big) \ &= \ s_{m+k}(x)\cdot x \ = \ \big(s_{m+k}(x)-\bar{u}\big)\cdot (\bar{u}\cdot x) \\ &= \ \big(s_{n+k}(y)-\bar{v}\big)\cdot (\bar{v}\cdot y) \ = \ s_{n+k}(y)\cdot y \ = \ \rho^{-1}\Big(g(y)(n+k)\Big),
\end{align*}
which implies $g(x)\mathrel{E_{1,t}}g(y)$. Thus $g$ is indeed a Borel reduction from $\EJT\res\We''$ to $E_{1,t}$, so $\EJT\res\We''$ is hypersmooth. Since $\sigma:\We\to\We''$ is a Borel reduction from $\EJT\res\We$ to $\EJT\res\We''$, it follows that $\EJT\res\We$ is hypersmooth as well.
\end{proof}

This completes the proof of Theorem \ref{thm:e1} for the case $\bar{v}=\bar{e}_1$. The case $\bar{v}=\bar{e}_2$ can be proven similarly, but we can also use $\delta$-symmetry to deduce this case directly from what we have already shown. Indeed, note that $\delta[\textup{Per}_{\bar{e}_1}(J^{\omega^2})]=\textup{Per}_{\bar{e}_2}(J^{\omega^2})$, and define
\begin{equation}
X_{\bar{e}_2} := \delta\big[\Xe\big] \qquad \text{and} \qquad W_{\bar{e}_2} := \delta\big[\We\big],
\end{equation}
so $\neg\textup{Per}_{\bar{e}_2}(J^{\omega^2})=X_{\bar{e}_2}\sqcup W_{\bar{e}_2}$. 

\begin{corollary}\label{cor:e2}
For every $J\in\Z_{\geq 2}$, $\chi_B\big(\G_{\bar{e}_2}\res X_{\bar{e}_2}\big)\leq 3$ and $\EJT\res W_{\bar{e}_2}$ is hypersmooth.
\end{corollary}

\begin{proof}
Let $c_{\bar{e}_1}:X_{\bar{e}_1}\to 3$ be the Borel $3$-coloring of $\G_{\bar{e}_1}\res X_{\bar{e}_1}$ given in Proposition \ref{prop:ColoringGe1}. It is easy to check using Lemma \ref{lem:delta} that $c_{\bar{e}_2}:=c_{\bar{e}_1}\circ\delta$ is a Borel $3$-coloring of $\G_{\bar{e}_2}\res X_{\bar{e}_2}$, and also that $\delta$ is a Borel reduction from $\EJT\res W_{\bar{e}_2}$ to $\EJT\res\We$.
\end{proof}

\begin{corollary}
Let $J\in\Z_{\geq 2}$ and write $X_2=X_2(J):=\Xe(J)\cap\Xee(J)$. Then $\chi_B\big(\vec{\G}_{2}\res X_2\big)\leq 9$.
\end{corollary}

\begin{proof}
The Borel map $c_2:X_2\to 3\times 3$ defined by $c_2(x)=\big(c_{\bar{e}_1}(x), \, c_{\bar{e}_2}(x)\big)$ is a Borel $9$-coloring of $\vec{\G}_{2}\res X_2$.
\end{proof}

To complete the proof of Theorem \ref{thm:e1}, we turn our attention to the graph $\Gp$ on $\JT$. Recall that we write $\p:=(1,1)\in\N^2$. We will be able to apply to $\Gp$ essentially the same methods we used to treat $\Ge$, including in particular Lemma \ref{lem:Ge1HSTechLem}, which was written so as to be applicable in both cases. Note the similarity of Definition \ref{df:XpWp} below to that of Definition \ref{df:XeWe}.

\begin{definition}\label{df:XpWp}
For each $J\in\Z_{\geq 2}$, let
\[
\Xp \ = \ \Xp(J) \ := \ \Big\{ x\in\neg\textup{Per}_{\p}(J^{\omega^2})\st (\exists\bar{a}\in\N^2) \text{  $\langle (\bar{a}\cdot x)(i,i)\st i\in\omega\rangle$ is not eventually constant}\Big\},
\]
and let
\begin{align*}
\Wp \ = \ \Wp(J) \ :\!&= \ \neg\textup{Per}_{\p}(J^{\omega^2})\setminus\Xp(J) \\
&= \ \Big\{ x\in\neg\textup{Per}_{\p}(J^{\omega^2})\st (\forall\bar{a}\in\N^2) \text{ $\langle (\bar{a}\cdot x)(i,i)\st i\in\omega\rangle$ is eventually constant}\Big\}.
\end{align*}
\end{definition}

Note that $\Xp$ and $\Wp$ are Borel and $\EJT$-invariant.

\begin{proposition}\label{prop:GpColoring}
For every $J\in\Z_{\geq 2}$, we have $\chi_B\big(\Gp\res\Xp\big)\leq 3$.
\end{proposition}

\begin{proof}
Fix $J\in\Z_{\geq 2}$. 
Using Fact \ref{fact:KST5.2}, let $c_J:J^\omega\to 3$ be a fixed Borel coloring of $\G_1(J)$. For each $x\in\Xp$, let $\bar{a}(x)$ be the least $\bar{a}\in\N^2$ such that the sequence $s(x):=\langle (\bar{a}\cdot x)(i,i)\st i\in\omega\rangle$ is not eventually constant. Noting that the map $x\mapsto\bar{a}(x)$ is Borel, define the Borel function $c_{\p}:\Xp\to 3$ by $c_{\p}(x)=c_J\big(s(x)\big)$. Let $x,y\in\Xp$ be arbitrary such that $x\mathrel{\Gp}y$. Then $\bar{a}(x)=\bar{a}(y)$ and $s(x)\mathrel{\G_1}s(y)$, so $c_J\big(s(x)\big)\ne c_J\big(s(y)\big)$ and therefore $c_{\p}(x)\ne c_{\p}(y)$. 
\end{proof}

In order to prove $\EJT\res\Wp$ is hypersmooth, it will be helpful to have some additional tools at our disposal. First, define the map $T:\N^2\to\N^2$ by $T(a_1,a_2):=(a_1, a_1+a_2)$, and, analogously, let $h:\JT\to\JT$ be the continuous function defined by $h(x)(i,j) := x(i, \, i+j)$. Note that
\begin{equation}\label{eq:Th}
\bar{a}\cdot h(x) \ = \ h\big(T(\bar{a})\cdot x\big) \qquad \text{for all $x\in\JT$ and $\bar{a}\in\N^2$.}
\end{equation}
Also note that for each $x\in\JT$ and $j\in\omega$ we have
\begin{equation}\label{eq:stupid0}
\big\langle h(x)(i,j)\st i\in\omega\big\rangle \ = \ \big\langle (j\bar{e}_2\cdot x)(i,i)\st i\in\omega\big\rangle.
\end{equation}
Next, for each $x\in\We'$ and $n\in\omega$, let $\alpha(x)(n)$ be the least $m\in\omega$ such that $x(i,n)=x(m,n)$ for all $i\geq m$. Then $\limsup\alpha(x)(n)=\infty$ for all $x\in\We'$, since $\We'\subseteq\neg\textup{Per}_{\bar{e}_1}(\JT)$. Similarly, for each $x\in\Wp$ and $n\in\omega$, let $\beta(x)(n)$ be the least $m\in\omega$ such that $x(i,n+i)=x(m,n+m)$ for all $i\geq m$, and let $\beta(x)(-n)$ be the least $m\in\omega$ such that $x(n+i,i)=x(n+m,m)$ for all $i\geq m$. Note that $\delta[\Wp]=\Wp$ and for all $x\in\Wp$ and $n\in\omega$ we have $\beta(x)(-n)=\beta\big(\delta(x)\big)(n)$. Furthermore, using Equation \eqref{eq:stupid0} we see that for all $x\in\Wp$ and $n\in\omega$,
\[
\beta(x)(n) = \alpha\big(h(x)\big)(n) \qquad \text{and} \qquad \beta(x)(-n) = \alpha\big(h(\delta(x))\big)(n).
\]
If $x\in\Wp$, then $\limsup\beta(x)(n)=\infty$ \emph{or} $\limsup\beta(x)(-n)=\infty$, and we have 
\[
\limsup\beta(x)(n)<\infty \quad \Longleftrightarrow \quad h(x)\in\textup{Per}_{\bar{e}_1}(\JT).
\]
Now we are ready to complete the proof of Theorem \ref{thm:e1}.

\begin{proposition}\label{prop:WpHS}
For every $J\in\Z_{\geq 2}$, the restriction of $\EJT$ to $\Wp$ is hypersmooth.
\end{proposition}

\begin{proof}
Define the Borel subsets $\Wp''\subseteq\Wp'\subseteq\Wp$ by
\[
\Wp' \ := \ \big\{x\in\Wp\st h(x)\in\We'\big\} \qquad \text{and} \qquad \Wp'' \ := \ \big\{x\in\Wp\st h(x)\in\We''\big\}.
\]
Note that by \eqref{eq:stupid0}, if $x\in\Wp$ then $h(x)$ has the property that for all $j\in\omega$ the sequence $\langle h(x)(i,j)\st i\in\omega\rangle$ is eventually constant. Thus
\[
\Wp' \ = \ \Big\{x\in\Wp\st h(x)\in\neg\textup{Per}_{\bar{e}_1}(\JT)\Big\} \ = \ \Big\{x\in\Wp\st \limsup \beta(x)(n)=\infty\Big\}.
\]
Given $x,y\in\Wp$ such that $x\mathrel{\EJT}y$, we have $\limsup\beta(x)(n)=\infty$ if and only if $\limsup\beta(y)(n)=\infty$, which shows $\Wp'$ is $\EJT$-invariant. Next, note that
\[
x\in \Wp\setminus\Wp' \ \Longrightarrow \ \limsup\beta(x)(n)<\infty \
\Longrightarrow \ \limsup\beta(x)(-n)=\infty \ \Longrightarrow \ \limsup\beta\big(\delta(x)\big)(n)=\infty,
\]
so $\delta\big[\Wp\setminus\Wp']\subseteq\Wp'$, which means $\delta$ is a Borel reduction from $\EJT\res(\Wp\setminus\Wp')$ to $\EJT\res\Wp'$. Furthermore, if for $x\in\Wp'$ we let $k_x:=\beta(x)(0)=\mu\big(h(x)\big)$, then the map $x\mapsto \dot{k}_x\cdot x$ is a Borel reduction from $\EJT\res\Wp'$ to $\EJT\res\Wp''$. It will suffice, therefore, to show $\EJT\res\Wp''$ is hypersmooth.

We define a Borel reduction $g:\Wp''\to (2^\omega)^\omega$ from $\EJT\res\Wp''$ to $E_{1,t}$ as follows. 
Given $x\in\Wp''$, so $h(x)\in\We''$, again let $\rho:J^{\omega^2}\to 2^\omega$ be a fixed Borel bijection, and define
\[
g(x)(n) \ := \ \rho\Big(\big(T\circ s_n\circ h\big)(x)\cdot x\Big),
\]
where $s:\We''\to (\N^2)^\omega$ is the Borel function given in Lemma \ref{lem:Ge1HSTechLem}. To check that $g$ is indeed a Borel reduction from $\EJT\res\Wp''$ to $E_{1,t}$, fix $x,y\in\Wp''$ and suppose first that $g(x)\mathrel{E_{1,t}}g(y)$. Then there exist $m,n\in\N$ such that $g(x)(m)=g(y)(n)$, which means
\[
\big(T\circ s_m\circ h\big)(x)\cdot x \ = \ \big(T\circ s_n\circ h\big)(y)\cdot y,
\]
and therefore $x\mathrel{\EJT}y$. Conversely, suppose $x\mathrel{\EJT}y$, and fix $\bar{a},\bar{b}\in\N^2$ such that $\bar{a}\cdot x=\bar{b}\cdot y$. Then also
\[
\big(a_1, \, a_1+b_1+a_2\big)\cdot x \ = \ \big(b_1, \, a_1+b_1+b_2\big)\cdot y,
\]
which by Equation \eqref{eq:Th} implies
\[
(a_1, \, b_1+a_2)\cdot h(x) \ = \ h\Big((a_1, \, a_1+b_1+a_2)\cdot x\Big) \ = \ h\Big((b_1, \, a_1+b_1+b_2)\cdot y\Big) \ = \ (b_1, \, a_1+b_2)\cdot h(y).
\]
Write $\bar{u}:=(a_1,\, b_1+a_2)$ and $\bar{v}:=(b_1, \, a_1+b_2)$, and also $x':=h(x)$ and $y':=h(y)$, so $\bar{u}\cdot x'=\bar{v}\cdot y'$ and $T(\bar{u})\cdot x=T(\bar{v})\cdot y$. Applying Lemma \ref{lem:Ge1HSTechLem} to the elements $x',y'\in\We''$, fix $m,n\in\omega$ such that $s_{m+k}(x')-\bar{u}=s_{n+k}(y')-\bar{v}\in\N^2$ for all $k\in\omega$. Then for all $k\in\omega$,
\begin{align*}
\big(T\circ s_{m+k}\circ h\big)(x)\cdot x \ &= \ \Big(T\big(s_{m+k}(x')-\bar{u}\big)+T(\bar{u})\Big)\cdot x \\
&= \ \Big(T\big(s_{m+k}(x')-\bar{u}\big)\Big) \cdot \big(T(\bar{u})\cdot x\big) \\
&= \ \Big(T\big(s_{n+k}(y')-\bar{v}\big)\Big) \cdot \big(T(\bar{v})\cdot y\big) \\
&= \ \Big(T\big(s_{n+k}(y')-\bar{v}\big)+T(\bar{v})\Big)\cdot y \ = \ \big(T\circ s_{n+k}\circ h\big)(y)\cdot y.
\end{align*}
It follows that $g(x)(m+k)=g(y)(n+k)$ for all $k\in\omega$, that is, $g(x)\mathrel{E_{1,t}}g(y)$, as desired.
\end{proof}

We conclude this section by briefly considering whether Theorem \ref{thm:e1} can be extended to $2$-dimensional product spaces with infinite base space. Given $x\in (2^\omega)^{\omega^2}$ and $n\in\omega$, define the ``$n$th slice'' $x_n\in 2^{\omega^2}$ of $x$ by $x_n(i,j) := x(i,j)(n)$. One can attempt to apply the arguments of this section slice-by-slice to elements of $(2^\omega)^{\omega^2}$ and show, for instance, that a version of Theorem \ref{thm:e1} holds on the $E^{2^\omega}_2$-invariant Borel set of all $x\in (2^\omega)^{\omega^2}$ for which there exists $n\in\omega$ such that $x_n$ is not $\bar{e}_1$-periodic. By a series of additional arguments along these lines, we are able to isolate an apparent obstruction to extending Theorem \ref{thm:e1} to $E^{2^\omega}_2$ or even $E^\omega_2$ which provides  what is essentially the simplest example we do not know how to handle.

\begin{example}\label{ex:stuck}
With $[\omega]^{\omega}$ as in Example \ref{ex:BadAssGraph}, let $\eta:\omega^2\to\omega$ be the computable bijection with inverse
\[
\big((0,0), (1,0), (0,1), (2,0), (1,1), (0,2), (3,0),(2,1),(1,2),(0,3),.\ldots\big),
\]
and define the continuous injection $f:[\omega]^\omega\to \omega^{\omega^2}$ by $f(x)(i,j) := x\big(\eta(i,j)\big)$. Let $X:=\ran(f)\subseteq \omega^{\omega^2}$, so $X$ is Borel and $\N^2$-invariant. It is not too difficult to show that $\chi_B\big(\G_{\bar{v}}(\omega)\res X\big)=\aleph_0$ for each $\bar{v}\in\{\bar{e}_1,\bar{e}_2,\p\}$. Is $E^\omega_2\res X$ hypersmooth? Note that we can code this example into $(2^\omega)^{\omega^2}$ by composing $f$ with the map $g:\omega^{\omega^2}\to(2^\omega)^{\omega^2}$ defined by $g(x)(i,j)(n)=1 \; :\Leftrightarrow \; x(i,j)=n$.
\end{example}


%

\section{The graph $\Gq$}\label{sec:se}

In this section we consider the graph $\Gq$, which is defined over any base space as in Equation \eqref{eq:mq} by 
\begin{equation}\label{eq:Gse}
x\mathrel{\Gq}y \quad :\Longleftrightarrow \quad \bar{e}_1\cdot x=\bar{e}_2\cdot y \ \vee \ \bar{e}_1\cdot y=\bar{e}_2\cdot x.
\end{equation}
Recall that $\q:=(1,-1)$, and that we write $x\mathrel{\Gsea}y$ to mean $\bar{e}_1\cdot x=\bar{e}_2\cdot y$, so $\Gq$ is the symmetrization of $\Gsea$. Unlike $\vec{\G}_{\bar{e}_1}$, $\vec{\G}_{\bar{e}_2}$, and $\vec{\G}_{\p}$, the directed graph $\Gsea$ has not just uncountable in-degree, but uncountable out-degree as well, regardless of the base space.

As in Section \ref{sec:ep}, we work over a finite base space which we take to be an integer $J\geq 2$. The particular value of $J$ will occasionally be important in this section, and when we wish to emphasize it we will write $\Gq(J)$ or $\vec{\G}_{\q}(J)$, as in Sections \ref{sec:graphs} and \ref{sec:ep}. Our goal will be to prove the following analogue of Theorem \ref{thm:e1}.

\begin{theorem}\label{thm:MainSE}
For each $J\in\Z_{\geq 2}$, there is a partition of $\neg\textup{Per}_{\q}(\JT)$ into a pair of $\EJT$-invariant Borel sets $\Xq(J)$ and $\Wq(J)$ such that $\chi_B\big(\Gq(J)\res\Xq(J)\big)<\infty$ and $\EJT\res\Wq(J)$ is hypersmooth; more precisely, we have $\chi_C\big(\Gq(J)\res\Xq(J)\big)\leq J^2+3J+19$.
\end{theorem}

Here, in contrast to Theorem \ref{thm:e1}, we obtain a continuous coloring but require more colors. We make no attempt to find $\chi_B\big(\Gq(J)\res\Xq(J)\big)$ exactly; any finite bound would suffice for our purposes by Corollary \ref{cor:main}. 
Proving Theorem \ref{thm:MainSE} will require different methods, and somewhat more effort, than we used in proving Theorem \ref{thm:e1}, and this entire section is devoted to its proof.  Throughout, even when we do not say so explicitly, assume $J\in\Z_{\geq 2}$. Also recall that $\delta:J^{\omega^2}\to J^{\omega^2}$ is the coordinate-swapping involution defined in Equation \ref{eq:delta} by $\delta(x)(i,j):=x(j,i)$.
 
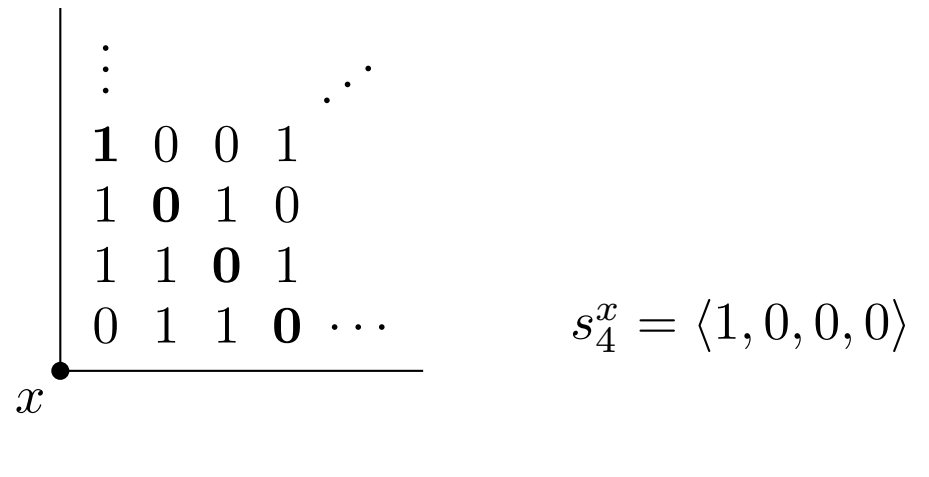
\begin{figure}[h]
\begin{center}
\begin{tikzpicture}
\draw (0,0) -- (2.4,0);
\draw (0,0) -- (0,2.4);
\node at (-.2,-.2) {$x$};
\node at (0,0) [circle,fill,inner sep=1.2pt]{};
\node at (.3,.3) {$0$};
\node at (.7,.3) {$1$};
\node at (1.1,.3) {$1$};
\node at (1.5,.3) {$\bf{0}$};
\node at (.3,.7) {$1$};
\node at (.7,.7) {$1$};
\node at (1.1,.7) {$\bf{0}$};
\node at (1.5,.7) {$1$};
\node at (.3,1.1) {$1$};
\node at (.7,1.1) {$\bf{0}$};
\node at (1.1,1.1) {$1$};
\node at (1.5,1.1) {$0$};
\node at (.3,1.5) {$\bf{1}$};
\node at (.7,1.5) {$0$};
\node at (1.1,1.5) {$0$};
\node at (1.5,1.5) {$1$};
\node at (.3, 2.1) {$\vdots$};
\node at (2, .28) {$\cdots$};
\node at (1.9,2) {\reflectbox{$\ddots$}};
\node at (4.5,.3) {$s^x_4=\langle 1, 0, 0, 0\rangle$};
\end{tikzpicture}
\end{center}
\caption{The sequence $s^x_4$ for a given $x\in\TT$.}
\label{fig:sxk}
\end{figure}

We begin with the following definition: for each $x\in\JTT$ and $k\in\Z_+$, let $s^x_k$ be the sequence of length $k$ defined by 
\begin{equation}\label{df:sxk}
s^x_k(i) \ := \ x(i,\,k-1-i) \quad \text{for all $0\leq i<k$}.
\end{equation}
The sequences $s^x_k$ may be pictured as in Figure \ref{fig:sxk}. Note that $s_k^x$ is constant for all $k$ if and only if $\bar{e}_1\cdot x=\bar{e}_2\cdot x$.

\begin{claim}\label{claim:SEequivalent}
For all $x\in\JT$, the following are equivalent:
\begin{enumerate}[(i)]
 \item $x$ is not $\q$-periodic;
 \item for all $\bar{c}\in\N^2$ there exist infinitely many $k\in\Z_+$ for which $s_k^{\bar{c}\cdot x}$ is non-constant.
\end{enumerate}
\end{claim}

\begin{proof}
We prove the contrapositives. 
Fix $x\in\JT$. For (i) $\Rightarrow$ (ii), assume $\bar{c}\in\N^2$ is such that $s_k^{\bar{c}\cdot x}$ is constant for all but finitely many $k$. Fix $n$ such that $s_k^{\bar{c}\cdot x}$ is constant for all $k\geq n$, and let $y:=(\bar{c}+\d{n})\cdot x$, so $s_k^y$ is constant for all $k$. Then $\bar{e}_1\cdot y=\bar{e}_2\cdot y$, so $(\bar{c}+\d{n}+\bar{e}_1)\cdot x=(\bar{c}+\d{n}+\bar{e}_2)\cdot x$, showing that $x$ is $\q$-periodic.

Conversely, for (ii) $\Rightarrow$ (i), suppose $x$ is $\q$-periodic, and fix $\bar{a},\bar{b}\in\N^2$ with $\bar{a}-\bar{b}=\q$ such that $\bar{a}\cdot x=\bar{b}\cdot x$. Let $n:=\max\{|\bar{a}|,|\bar{b}|\}$. Then $\d{n}\cdot x=(\d{n}-\bar{b})\cdot (\bar{a}\cdot x)=\big(\d{n}+\q\big)\cdot x$, so $(\d{n}+\bar{e}_2)\cdot x=(\d{n}+\bar{e}_1)\cdot x$, which means $s^{\d{n}\cdot x}_k$ is constant for all $k$.
\end{proof}

Using Claim \ref{claim:SEequivalent}, we make the following definitions.

\begin{definition}\label{df:XqWq}
For each $x\in\npseJ$, let $r(x)$ be the least $k\in\Z_+$ for which $s^x_k$ is non-constant, and write $s^x:=s^x_{r(x)}$. Then define
\[
\Xse \ = \ \Xse(J) \ := \ \left\{ x\in\npseJ\st\liminf_{n\to\infty} r(\d{n}\cdot x)<\infty\right\},
\]
and let
\begin{align*}
\Wse \ = \ \Wse(J) \ :\!&= \ \npseJ\setminus\Xse(J) \\
&= \ \left\{ x\in\npseJ\st\lim_{n\to\infty} r(\d{n}\cdot x)=\infty\right\}.
\end{align*}
\end{definition}

To improve notation below, we write $\d{r}(x)$ for the constant tuple with value $r(x)$, that is, $\d{r}(x):=\d{n}$ where $n=r(x)$. 
Note that $x\mapsto r(x)$ is continuous, and that $r(x)\geq 2$ for all $x\in\npseJ$. Since
\[
\Wse \ = \ \Big\{ x\in\npseJ\st (\forall N)(\exists M)(\forall k\geq M) \ r(\d{k}\cdot x)>N \Big\}
\]
and ``$r(\d{k}\cdot x)>N$'' is a clopen condition, the set $\Wse$ is Borel, hence so is $\Xse$. Since $s^{\delta(x)}_k$ is just $s^x_k$ in reverse order, we have $r(\delta(x))=r(x)$ for all $x\in\npseJ$, so the sets $\Xq$ and $\Wq$ are $\delta$-invariant. These sets are also $\EJT$-invariant, although we postpone a formal proof of this until the next section (see Corollary \ref{cor:WSEinvar}), when we establish another way of thinking about $\Wse$ which we shall need anyway and which makes proving invariance a bit simpler.

\begin{claim}\label{claim:postponed}
The sets $\Xse$ and $\Wse$ are Borel, $\delta$-invariant, and $\EJT$-invariant. \hfill $\square$
\end{claim}

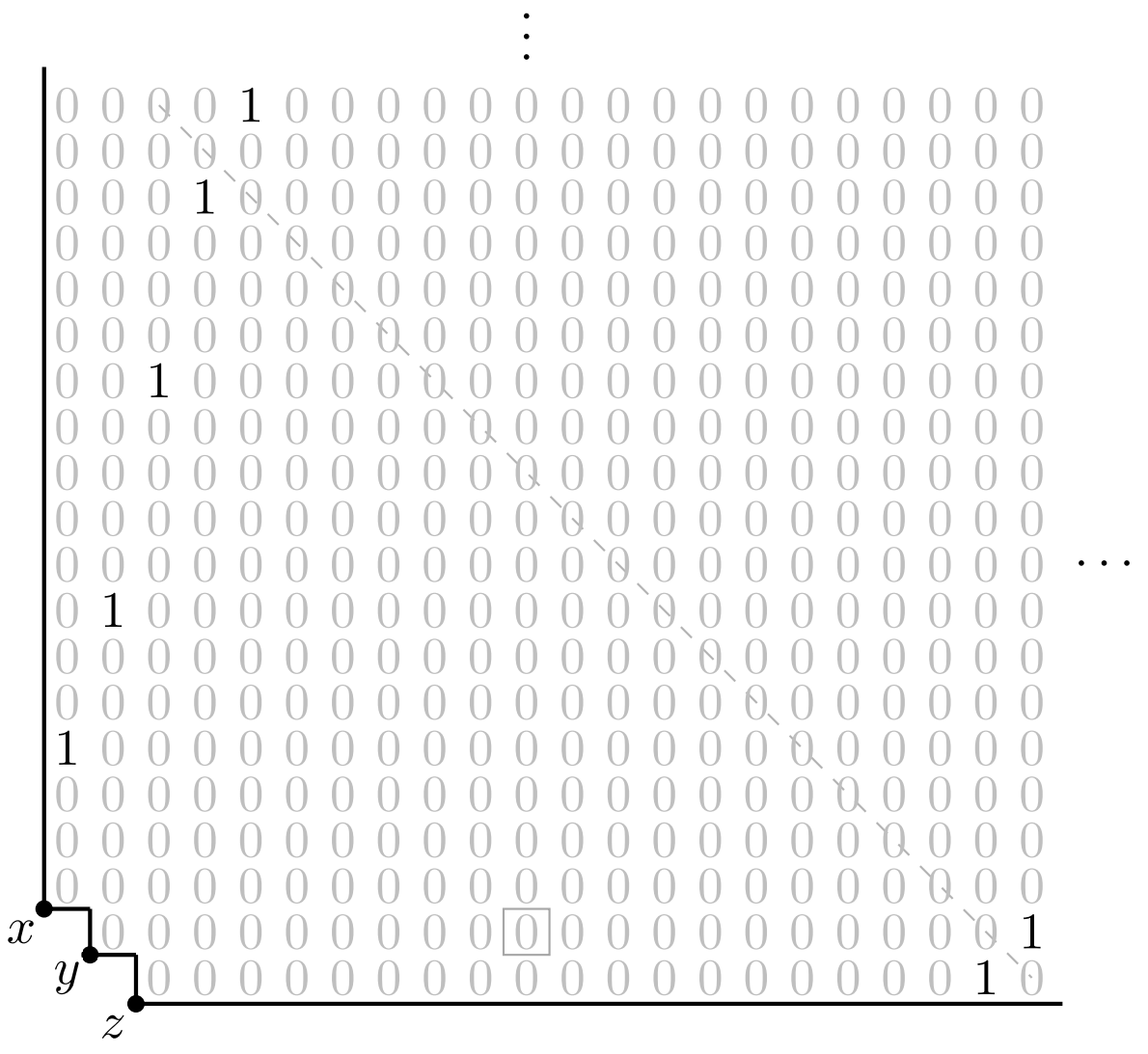
\begin{figure}[h]
\begin{center}
\begin{tikzpicture}[x=1.05cm, y=1.05cm]
     %
     \draw [gray!60,dashed] (.2+.3*2,.2+.3*17) -- (.2+.3*21,.2-.3*2);
     \draw [gray!70] (.05+.3*10,.05-.3*1) rectangle (.35+.3*10,.35-.3*1); 
 \draw [thick] (.05,.05) -- (.05,5.55);
 \draw [thick] (.05,.05) -- (0.35,.05);
 \draw [thick] (0.35,.05) -- (0.35,-0.25);
 \draw [thick] (0.35,-0.25) -- (0.65,-0.25);
 \draw [thick] (0.65,-0.25) -- (0.65,-0.57);
 \draw [thick] (0.65,-0.57) -- (6.7,-0.57);
  \node at (.05,.05) [circle,fill,inner sep=1.2pt]{};
  \node at (0.35,-0.25) [circle,fill,inner sep=1.2pt]{};
  \node at (0.65,-0.57) [circle,fill,inner sep=1.2pt]{};
 \foreach \x in {0,1,...,21} \foreach \y in {-2,-1,0,1,...,17} \node at (.2+.3*\x,.2+.3*\y) {{\color{lightgray}$0$}};
 \fill [fill=white] (.1+.3*0,.07-.3*2) rectangle (.3+.3*0,.33-.3*2);
 \fill [fill=white] (.1+.3*0,.07-.3) rectangle (.3+.3*0,.33-.3);
 \fill [fill=white] (.1+.3*1,.07-.3*2) rectangle (.3+.3*1,.33-.3*2);
\fill [fill=white] (.1+.3*0,.07+.3*3) rectangle (.3+.3*0,.33+.3*3); 
 \node at (.2+.3*0,.2+.3*3) {$1$};
\fill [fill=white] (.1+.3*1,.07+.3*6) rectangle (.3+.3*1,.33+.3*6); 
 \node at (.2+.3*1,.2+.3*6) {$1$};
\fill [fill=white] (.1+.3*2,.07+.3*11) rectangle (.3+.3*2,.33+.3*11); 
 \node at (.2+.3*2,.2+.3*11) {$1$};
 \fill [fill=white] (.1+.3*3,.07+.3*15) rectangle (.3+.3*3,.33+.3*15); 
 \node at (.2+.3*3,.2+.3*15) {$1$};
  \fill [fill=white] (.1+.3*4,.07+.3*17) rectangle (.3+.3*4,.33+.3*17); 
 \node at (.2+.3*4,.2+.3*17) {$1$};
\fill [fill=white] (.1+.3*21,.07-.3) rectangle (.3+.3*21,.33-.3); 
 \node at (.2+.3*21,.2-.3) {$1$};
\fill [fill=white] (.1+.3*20,.07-.3*2) rectangle (.3+.3*20,.33-.3*2); 
 \node at (.2+.3*20,.2-.3*2) {$1$};
 \node at (-.1,-.1) {$x$};
 \node at (.2,-.4) {$y$};
 \node at (.5,-.72) {$z$};
 \node at (7,2.3) {$\cdots$};
 \node at (3.2,5.85) {$\vdots$};
 \node at (9.6,2.5)
{
\begin{tabular}{|c|c|c|}
\hline
 & $r$ & $\textup{tp}$ \\ \hline
$x$ & $4$ & $\L$ \\ \hline
$\bar{e}_1\cdot x$ & $7$ & $\L$ \\ \hline
$2\bar{e}_1\cdot x$ & $12$ & $\L$ \\ \hline
$3\bar{e}_1\cdot x$ & $16$ & $\L$ \\ \hline
$y$ & $8$ & $\L$ \\ \hline
$\bar{e}_1\cdot y$ & $13$ & $\L$ \\ \hline
$2\bar{e}_1\cdot y$ & $17$ & $\L$ \\ \hline
$3\bar{e}_1\cdot y$ & $18$ & $\R$ \\ \hline
$z$ & $14$ & $\L$ \\ \hline
$\bar{e}_1\cdot z$ & $18$ & $\mathcal{B}$ \\ \hline
$2\bar{e}_1\cdot z$ & $17$ & $\R$ \\ \hline
$3\bar{e}_1\cdot z$ & $16$ & $\R$ \\ \hline
\end{tabular}
};
\end{tikzpicture}
\end{center}
\caption{Sample points $x,y,z\in\L_1$.}
\label{fig:L1}
\end{figure}

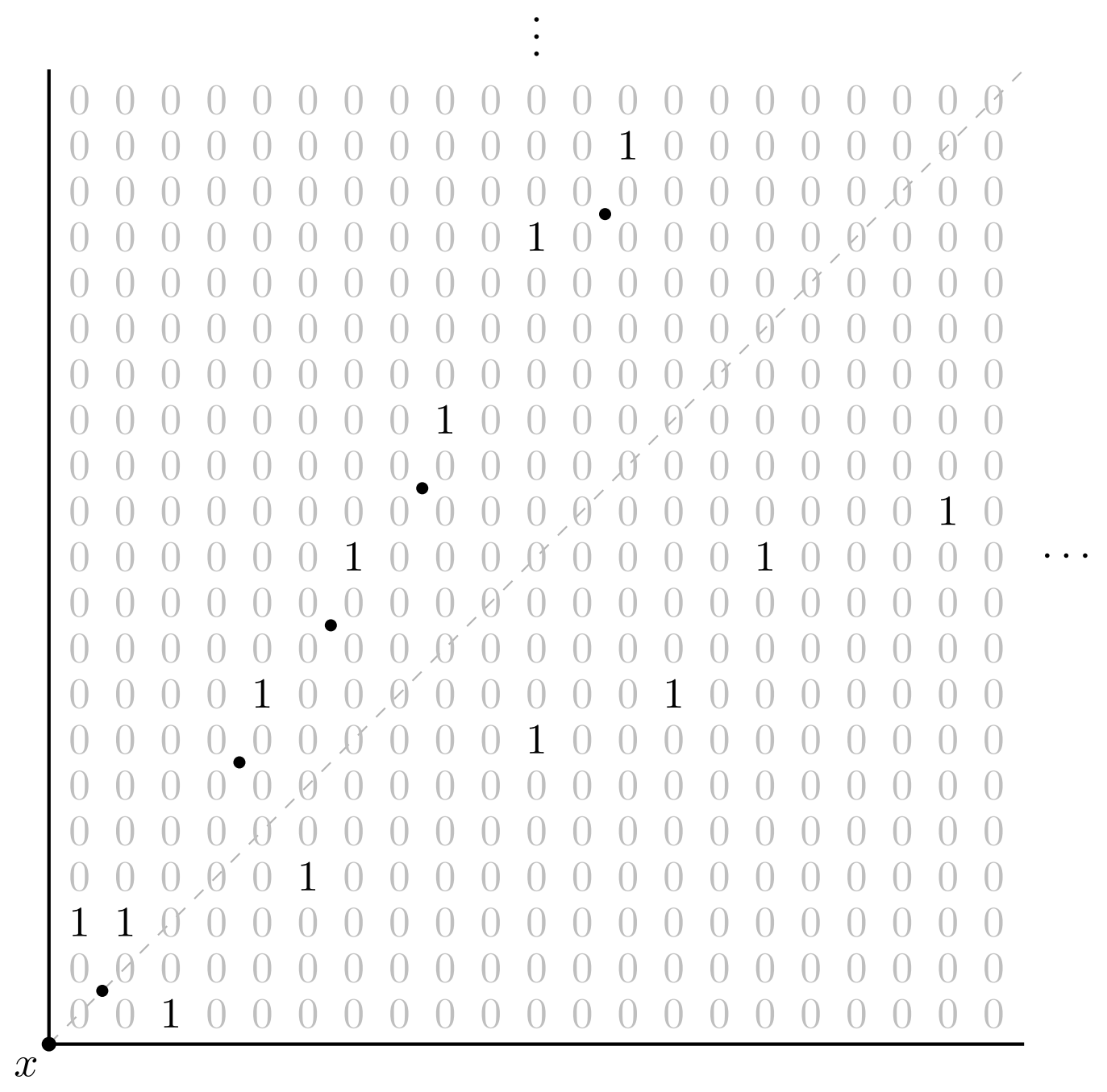
\begin{figure}[h]
\begin{center}
\begin{tikzpicture}[x=1.25cm, y=1.25cm]
     %
     %
     \draw [gray!60,dashed] (0,0) -- (6.4,6.4);
 \draw [thick] (0,0) -- (6.4,0);
 \draw [thick] (0,0) -- (0,6.4);
  \node at (0,0) [circle,fill,inner sep=1.2pt]{};
  \node at (.05+.3*1,.05+.3*1) [circle,fill,inner sep=1pt]{};
  \node at (.05+.3*4,.05+.3*6) [circle,fill,inner sep=1pt]{};
  \node at (.05+.3*6,.05+.3*9) [circle,fill,inner sep=1pt]{};
  \node at (.05+.3*8,.05+.3*12) [circle,fill,inner sep=1pt]{};
   \node at (.05+.3*12,.05+.3*18) [circle,fill,inner sep=1pt]{};
     %
  %
 \foreach \x in {0,1,...,20} \foreach \y in {0,1,...,20} \node at (.2+.3*\x,.2+.3*\y) {{\color{lightgray}$0$}};
\fill [fill=white] (.1+.3*0,.07+.3*2) rectangle (.3+.3*0,.33+.3*2); 
 \node at (.2+.3*0,.2+.3*2) {$1$};
\fill [fill=white] (.1+.3*2,.07+.3*0) rectangle (.3+.3*2,.33+.3*0); 
 \node at (.2+.3*2,.2+.3*0) {$1$};
\fill [fill=white] (.1+.3*1,.07+.3*2) rectangle (.3+.3*1,.33+.3*2); 
 \node at (.2+.3*1,.2+.3*2) {$1$};
\fill [fill=white] (.1+.3*5,.07+.3*3) rectangle (.3+.3*5,.33+.3*3); 
 \node at (.2+.3*5,.2+.3*3) {$1$};
\fill [fill=white] (.1+.3*4,.07+.3*7) rectangle (.3+.3*4,.33+.3*7); 
 \node at (.2+.3*4,.2+.3*7) {$1$};
\fill [fill=white] (.1+.3*6,.07+.3*10) rectangle (.3+.3*6,.33+.3*10); 
 \node at (.2+.3*6,.2+.3*10) {$1$};
\fill [fill=white] (.1+.3*10,.07+.3*6) rectangle (.3+.3*10,.33+.3*6); 
 \node at (.2+.3*10,.2+.3*6) {$1$};
\fill [fill=white] (.1+.3*13,.07+.3*7) rectangle (.3+.3*13,.33+.3*7); 
 \node at (.2+.3*13,.2+.3*7) {$1$}; 
\fill [fill=white] (.1+.3*8,.07+.3*13) rectangle (.3+.3*8,.33+.3*13); 
 \node at (.2+.3*8,.2+.3*13) {$1$};
\fill [fill=white] (.1+.3*15,.07+.3*10) rectangle (.3+.3*15,.33+.3*10); 
 \node at (.2+.3*15,.2+.3*10) {$1$};
\fill [fill=white] (.1+.3*19,.07+.3*11) rectangle (.3+.3*19,.33+.3*11); 
 \node at (.2+.3*19,.2+.3*11) {$1$}; 
\fill [fill=white] (.1+.3*10,.07+.3*17) rectangle (.3+.3*10,.33+.3*17); 
 \node at (.2+.3*10,.2+.3*17) {$1$};
\fill [fill=white] (.1+.3*12,.07+.3*19) rectangle (.3+.3*12,.33+.3*19); 
 \node at (.2+.3*12,.2+.3*19) {$1$};
 \node at (-.15,-.15) {$x$};
 \node at (6.7,3.2) {$\cdots$};
 \node at (3.2,6.7) {$\vdots$};
%
\node at (9,4.4)
{
\begin{tabular}{|c|c|c|}
\hline
$n$ & $r(\d{n}\cdot x)$ & $\textup{tp}(\d{n}\cdot x)$ \\ \hline\hline
$0$ & $3$ & $\mathcal{A}_3$ \\ \hline
$1$ & $2$ & $\mathcal{A}_2$ \\ \hline
$2$ & $5$ & $\mathcal{M}$ \\ \hline
$3$ & $3$ & $\mathcal{A}_3$ \\ \hline
\rowcolor{gray!30} $4$ & $4$ & $\mathcal{L}$ \\ \hline
$5$ & $7$ & $\mathcal{M}$ \\ \hline
\rowcolor{gray!30} $6$ & $5$ & $\mathcal{B}$ \\ \hline
$7$ & $7$ & $\mathcal{R}$ \\ \hline
$8$ & $6$ & $\mathcal{L}$ \\ \hline
$9$ & $8$ & $\mathcal{M}$ \\ \hline
\rowcolor{gray!30} $10$ & $6$ & $\mathcal{R}$ \\ \hline
$11$ & $9$ & $\mathcal{R}$ \\ \hline
\rowcolor{gray!30} $12$ & $8$ & $\mathcal{L}$ \\ \hline
\end{tabular}
};
\node at (9,0.5)
{
\begin{tabular}{|c|c|c|c|c|c|c|}
\hline
$n$ & $0$ & $1$ & $2$ & $3$ & $4$ & $5$ \\
\hline\hline
$j^x_n$ & $0$ & $0$ & $2$ & $1$ & $1$ & $2$ \\
\hline
$k^x_n$ & $0$ & $0$ & $2$ & $1$ & $1$ & $3$ \\
\hline
$J^x_n$ & $0$ & $0$ & $2$ & $3$ & $4$ & $6$ \\
\hline
$K^x_n$ & $0$ & $1$ & $4$ & $6$ & $8$ & $12$ \\
\hline
\end{tabular}
};
\end{tikzpicture}
\end{center}
\caption{A sample $x\in\Wqnw(2)$, with dots representing the first few terms of $s_n(x)\cdot x$.}
\label{fig:Wq}
\end{figure}

The remainder of this section is devoted to proving that $\Xse$ and $\Wse$ satisfy the conclusion of Theorem \ref{thm:MainSE}. We concentrate first on coloring $\Gse\res\Xse$. Our initial strategy will be to color elements $x$ based on where in $s^x$ we can find witnesses that $s^x$ is non-constant. As our proof that $\EJT\res\Wq$ is hypersmooth also uses this idea, and this portion of our coloring argument does not actually require working in $\Xq$, we make the following definitions in $\npseJ$ rather than just in $\Xq$. Let us define
\[
\L_0 \ := \ \big\{x\in\npseJ\st \forall i,j\in\{1,\ldots,r(x)-1\} \ s^x(i)=s^x(j)\big\} \qquad \text{and} \qquad \L \ := \ \big\{x\in\L_0\st r(x)\geq 4\big\},
\]
\[
\R_0 \ := \ \big\{x\in\npseJ\st \forall i,j\in\{0,\ldots,r(x)-2\} \ s^x(i)=s^x(j)\big\} \qquad \text{and} \qquad \R \ := \ \big\{x\in\R_0\st r(x)\geq 4\big\},
\]
and
\[
\B \ := \ \big\{x\in\npseJ\st r(x)\geq 4 \, \wedge \, \forall i,j\in\{1,\ldots,r(x)-2\} \ s^x(i)=s^x(j)\big\} \setminus\big(\mathcal{L}\cup\mathcal{R}\big).
\]
Additionally define
\begin{align*}
\M \ :\!&= \ \big\{x\in\npseJ\st r(x)\geq 4\big\} \setminus \big(\mathcal{L}\cup\mathcal{R}\cup\mathcal{B}\big) \\
&= \ \big\{x\in\npseJ\st r(x)\geq 4 \, \wedge \, \big\langle s^x(1),\ldots,s^x\big(r(x)-2\big)\big\rangle \text{ is non-constant}\big\}.
\end{align*}
Intuitively, $\L$, $\R$, $\B$, and $\M$ mean \emph{left}, \emph{right}, \emph{both}, and \emph{middle}, respectively. As an example, the element $x$ in Figure \ref{fig:sxk} belongs to $\L$; further examples can be found in Figures \ref{fig:L1} and \ref{fig:Wq}. Note that $\L\cap\R=\emp$, so the sets $\L$, $\R$, $\B$, and $\M$ are pairwise disjoint. Let us also write $\A_2:=\{x\in\npseJ\st r(x)=2\}$ and $\A_3:=\{x\in\npseJ\st r(x)=3\}$, so we have 
\begin{equation}
\npseJ \ = \ \A_2\sqcup\A_3\sqcup\L\sqcup\R\sqcup\B\sqcup\M.
\end{equation}
For each $x\in\npseJ$, let $\textup{tp}(x)$, the \emph{type} of $x$, be the unique set in this partition to which $x$ belongs.

\begin{claim}\label{claim:clopenSE}
The sets $\A_2$, $\A_3$, $\L_0$, $\R_0$, $\L$, $\R$, $\B$, and $\M$ are relatively clopen subsets of $\npseJ$, and, consequently, their intersections with $\Xse$ are relatively clopen in $\Xse$.
\end{claim}

\begin{proof}
In each case this can be checked by verifying that membership of $x\in\npseJ$ in the given set is determined by its values on the finite subset $\R\big(\d{r}(x)\big)$ of $\omega^2$. For example, fixing $x\in\npseJ$ and writing $s:=x\res\R\big(\d{r}(x)\big)$, if $x\in\L_0$ then $x\in [s]\subseteq\L_0$, and if $x\not\in\L_0$ then $x\in [s]\subseteq\npseJ\setminus\L_0$, showing $\L_0$ to be both open and closed in $\npseJ$. The other cases are identical.
\end{proof}

We will incorporate $\textup{tp}(x)$ into our coloring scheme, and indeed type together with a bit more information will enable us to distinguish $x$ from its $\Gse$-neighbors when $x\in\A_2\sqcup\A_3\sqcup\B\sqcup\M$, but we will have to work harder when $x\mathrel{\Gse}y$ and both $x$ and $y$ belong either to $\L$ or to $\R$. The symmetry between $\L$ and $\R$ will allow us to focus our attention on just one of them, so we will concentrate on $\L$.  Indeed, we have $\delta[\L]=\R$ and $\delta[\R]=\L$, while $\delta[\B]=\B$ and $\delta[\M]=\M$, and for any $\delta$-invariant Borel set $X\subseteq\npseJ$, the map $c_{\L}$ colors $\Gse\res(\L\cap X)$ if and only if $c_{\R}:=c_{\L}\circ\delta$ colors $\Gse\res(\R\cap X)$. For this reason, most of our work below will be devoted to ensuring that $\Gse$-neighbors $x$ and $y$ in $\Xq$ have different colors when $x,y\in\L$. This will take some effort, and we begin with a series of lemmas. The reader is encouraged to draw pictures in order to visualize them and see why they are true. 

\begin{lemma}\label{lem:M}
For all $x,y\in\M$, if $x\mathrel{\Gse}y$ then $r(x)=r(y)$.
\end{lemma}

\begin{proof}
Fix $x,y\in\M$ such that $x\mathrel{\Gsea}y$, along with $0<i<r(x)-2$ such that $s^x(i)\ne s^x(i+1)$. Then
\[
s^y_{r(x)}(i-1) \ = \ s^x(i) \ \ne \ s^x(i+1) \ = \ s^y_{r(x)}(i),
\]
which shows $r(y)\leq r(x)$.  Similarly, fixing $0<j<r(y)-2$ such that $s^y(j)\ne s^y(j+1)$, we have
\[
s^x_{r(y)}(j+1) \ = \ s^y(j) \ \ne \ s^y(j+1) \ = \ s^x_{r(y)}(j+2),
\]
which shows $r(x)\leq r(y)$.
\end{proof}

\begin{lemma}\label{lem:move1}
For all $x\in\npseJ$, we have the following:
\begin{enumerate}[(i)]
 \item $r(x)\leq r(\bar{e}_1\cdot x)+1$ and $r(x)\leq r(\bar{e}_2\cdot x)+1$;
 \item $r(x)\leq r(\bar{e}_1\cdot x)$ if and only if $x\in\L_0$, and $r(x)\leq r(\bar{e}_2\cdot x)$ if and only if $x\in\R_0$.
\end{enumerate}
\end{lemma}

\begin{proof}
Fix $x\in\npseJ$. For (i), if $s^{\bar{e}_1\cdot x}_k$ is non-constant, then $s^x_{k+1}$ is non-constant as well since $1\cdot s^x_{k+1}=s^{\bar{e}_1\cdot x}_k$. Thus $r(x)\leq r(\bar{e}_1\cdot x)+1$, and the proof that $r(x)\leq r(\bar{e}_2\cdot x)+1$ is similar. For (ii), assume first that $x\in\L_0$. We already know $2\leq r(x)\leq r(\bar{e}_1\cdot x)+1$ from (i), but $1\cdot s^x=s^{\bar{e}_1\cdot x}_{r(x)-1}$ is constant since $x\in\L_0$, which means $r(\bar{e}_1\cdot x)\ne r(x)-1$ and thus $r(x)\leq r(\bar{e}_1\cdot x)$. Conversely, suppose $k:=r(x)\leq r(\bar{e}_1\cdot x)$. Then $s^x_k$ is non-constant, but $s^{\bar{e}_1\cdot x}_{k-1}=1\cdot s^x_k$ must be constant since $k-1<r(\bar{e}_1\cdot x)$, which means $x\in\L_0$. That $r(x)\leq r(\bar{e}_2\cdot x)$ if and only if $x\in\R_0$ can be shown similarly, or, alternatively, follows by $\delta$-symmetry.
\end{proof}

\begin{lemma}\label{lem:up2}
Let $x\in\npseJ$, and write $x':=\bar{e}_2\cdot x$.
\begin{enumerate}[(i)]
 \item Suppose $x\in\L$.  Then $x'\in\L_0$ and $r(x')=r(x)-1$; in particular, $x'\in\L$ if and only if $r(x)\geq 5$.
 \item If $x'\in\L$, then exactly one of the following holds:
 \begin{enumerate}[(a)]
  \item $r(x)\leq r(x')$ and $x\in\R_0$; 
  \item $r(x)=r(x')+1$ and $x\in\L$;  
  \item $r(x)=r(x')+1$ and $x\in\B$.
 \end{enumerate}
\end{enumerate}
\end{lemma}

\begin{proof}
For (i), assume $x\in\L$, so $1\cdot s^x$ is constant. Then $x\big(0,r(x)-1\big)\ne x\big(1,r(x)-2\big)$, so $x'\big(0,r(x)-2\big)\ne x'\big(1,r(x)-3\big)$, which shows $r(x')\leq r(x)-1$. Thus $r(x')=r(x)-1$ by Lemma \ref{lem:move1}. 
It follows that $s^{x'}=s^{x}\res r(x')$, so $1\cdot s^{x'}$ is constant since $1\cdot s^x$ is constant. Hence $x'\in\L_0$, and therefore $x'\in\L \ \Leftrightarrow \ r(x')\geq 4 \ \Leftrightarrow \ r(x)\geq 5$.

For (ii), suppose $x'\in\L$. By Lemma \ref{lem:move1} and the fact that $\R_0$, $\L$, and $\B$ are pairwise disjoint, it will suffice to show that if $r(x)=r(x')+1$ then $x\in\L\cup\B$. So assume $r(x)=r(x')+1$, which means $r(x)\geq 5$ since $x'\in\L$. Then $s^x\res r(x')=s^{x'}$, so $s^x(0)\ne s^x(1)$ but $1\cdot\big(s^x\res r(x')\big)$ is constant since $x'\in\L$. Now if $s^x\big(r(x)-2\big)=s^x\big(r(x)-1\big)$ then $x\in\L$, and if $s^x\big(r(x)-2\big)\ne s^x\big(r(x)-1\big)$ then $x\in\B$.
\end{proof}

\begin{lemma}\label{lem:upn}
For all $x\in\L\cup\B$ and $0<n\leq r(x)-2$, we have $n\bar{e}_2\cdot x\in\L_0$ and $r(n\bar{e}_2\cdot x)=r(x)-n$; in particular, $n\bar{e}_2\cdot x\in\L$ if and only if $n\leq r(x)-4$.
\end{lemma}

\begin{proof}
Fix $x\in\L\cup\B$ and $0<n\leq r(x)-2$, and write $x':=n\bar{e}_2\cdot x$. Then for all $k\in\Z_+$ we have $s^{x'}_k=s^{x}_{n+k}\res k$, so $s^x_{n+k}$ is non-constant if $s^{x'}_k$ is non-constant, which shows $r(x)\leq r(x')+n$. On the other hand, we have
\[
x'\big(0,\,r(x)-n-1\big) \ = \ x\big(0,\,r(x)-1\big) \ \ne \ x\big(1,\,r(x)-2\big) \ = \ x'\big(1,\,r(x)-n-2\big),
\]
which shows $r(x')\leq r(x)-n$. Thus $r(x')=r(x)-n$. This implies $s^{x'}=s^x\res r(x')$, which since $x\in\L\cup\B$ and $2\leq r(x')<r(x)$ implies $1\cdot s^{x'}$ is constant and $s^{x'}(0)\ne s^{x'}(1)$, that is, $x'\in\L_0$. It now follows from the fact that $r(x')=r(x)-n$ that $x'\in\L$ if and only if $n\leq r(x)-4$.
\end{proof}

The points $x,y\in\L$ satisfying $x\mathrel{\Gsea}y$ that will give us the most trouble are those for which $r(x)\ll r(y)$. To specify these ``bad'' points precisely, we make the following definition which will be central to our argument. Define 
\begin{equation}
\L_1 \ := \ \Big\{x\in\L\st \d{1}\cdot x\in\L \, \wedge \ r(x)+2\leq r(\d{1}\cdot x)\Big\}.
\end{equation}
Sample $\Gsea$-related points $x$ and $y$ in $\L_1$ are pictured in Figure \ref{fig:L1}. As mentioned above, the case $x,y\in\L$ with $x\mathrel{\Gq}y$ will have a ``dual'' version where $x,y\in\R$, but we need only consider one of these two cases by $\delta$-symmetry. The dual version of $\L_1$ consisting of the ``bad'' points $x\in\R$ is given by
\begin{equation}\label{eq:dual}
\R_1' \ := \ \big\{x\in\R\st\delta(x)\in\L_1\big\} \ 
= \ \Big\{x\in\R\st \d{1}\cdot x\in\R \, \wedge \ r(x)+2\leq r(\d{1}\cdot x)\Big\}.
\end{equation}
Below, we will focus on the (arbitrarily chosen) primary case $x,y\in\L$, and indicate ``dual'' definitions by including a prime in the notation. These dual definitions can always be made either indirectly using $\delta$ or directly by dualizing the primary definition, as in \eqref{eq:dual} above. Note that for each $x\in\npseJ$, membership of $x$ in $\L_1$ is determined by $x\res\R(\d{k})$ where $k=\max\{r(x),\, 1+r(\d{1}\cdot x)\}$, so $\L_1$ is clopen in $\npseJ$, as is $\R_1'$, and both $\L_1\cap\Xse$ and $\R_1'\cap\Xse$ are clopen in $\Xse$.

\begin{lemma}\label{lem:LL}
For all $x,y\in\L$ such that $x\mathrel{\Gsea}y$, if $x\not\in\L_1$ then $r(x)+1\leq r(y)\leq r(x)+3$. 
\end{lemma}

\begin{proof}
Suppose $x,y$ are as stated, with $x\not\in\L_1$. For the first inequality, note that since $y\in\L$ we have 
\[
x\big(1,r(y)-2\big) \ = \ y\big(0,r(y)-1\big) \ = \ s^y(0) \ \ne \ s^y(1) \ = \ y\big(1,r(y)-2\big) \ = \ x\big(2, r(y)-3\big),
\]
which shows $r(x)<r(y)$ since $x\in\L$. For the second, using $x\not\in\L_1$, suppose first that $2\bar{e}_2\cdot y=\d{1}\cdot x\not\in\L$. Then $r(y)<6$ by Lemma \ref{lem:upn}, but also $4\leq r(x)<r(y)$, so in fact $r(y)=5=r(x)+1$. If instead $\d{1}\cdot x\in\L$, then $r(\d{1}\cdot x)\leq r(x)+1$ since $x\not\in\L_1$, so by Lemma \ref{lem:upn} we have $r(y)=r(\d{1}\cdot x)+2\leq r(x)+3$.
\end{proof}

At this point we are ready to color the restriction of $\Gse$ to $\npseJ\setminus (\L_1\cup\R_1')$, and we pause to do this before considering $\Gq$ on $\Xq\cap(\L_1\cup\R_1')$, which will take considerably more effort. Our eventual goal is to define a finite set $\Sse$ and a continuous map $\cse:\Xse\to \Sse$ such that $\cse(x)\ne \cse(y)$ whenever $x\mathrel{\Gse}y$. For each $x\in\Xse\setminus (\L_1\cup\R_1')$, $\cse(x)$ will be a finite sequence with first term $\textup{tp}(x)$. The set $\Sse$ and function $\cse$ will depend on $J$ here, so we occasionally include $J$ in the notation and write, e.g., $\Sse(J)$ and $\cse^J$ in place of the abbreviations $\Sse$ and $\cse$ that we wrote above. Let us define 
\[
S_2(J):=\big\{(\A_2,j)\st j\in J\big\} \qquad \text{and}\qquad S_3(J):=\big\{ (\A_3, i, j)\st i,j\in J\big\},
\]
so $|S_2(J)|=J$ and $|S_3(J)|=J^2$, along with the functions
\[
c_2^J:\A_2\to S_2(J) \qquad \text{and} \qquad c_3^J:\A_3\to S_3(J)
\]
defined by
\[
c_2^J(x) \ := \ \big(\A_2, \, s^x(0)\big) \qquad \text{and} \qquad c_3^J(x) \ := \ \big(\A_3, \, s^x(0), \, s^x(1)\big).
\]
Next, for each $x\in\M$ let $\epsilon(x)\in\{0,1\}$ be the parity of the least $i<r(x)$ such that $s^x(i)\ne s^x(0)$, and define the function $c_4^J$ on $\M\cup\B$ by
\[
c_4^J(x) \ := \ 
\begin{cases}
\Big(\M, \, s^x(0),\, \epsilon(x)\Big) & \text{if $x\in\M$;} \\
(\B) & \text{if $x\in\B$.} 
\end{cases}
\]
Write $S_4(J):=\ran(c_4^J)$, so $|S_4(J)|=2J+1$. Finally, define the function $c_0$ on $(\L\setminus\L_1)\cup(\R\setminus\R_1')$ by 
\[
c_0(x) \ := \ 
\begin{cases}
\big(\L,\, r(x)\!\mod 4\big) & \text{if $x\in\L\setminus\L_1$;} \\
\big(\R,\, r(x)\!\mod 4\big) & \text{if $x\in\R\setminus\R_1'$.}
\end{cases}
\]
Put $S_0:=\ran(c_0)$, so $|S_0|=8$. Note that $S_0$ and $c_0$ do not depend on $J$, and that
\begin{equation}\label{eq:Jcount}
|S_0\sqcup S_2(J)\sqcup S_3(J)\sqcup S_4(J)| \ = \ J^2+3J+9.
\end{equation}

\begin{claim}\label{claim:continuousSE}
The maps $c_2$, $c_3$, $c_4$, and $c_0$ defined above are continuous.
\end{claim}

\begin{proof}
This is routine to check, using Claim \ref{claim:clopenSE} and arguments similar to the proof of Claim \ref{claim:clopenSE}. Note first that for all $k\in\Z_+$, the set $\{x\in\npseJ\st r(x)=k\}$ is relatively clopen in $\npseJ$, since membership of $x$ in this set is determined by $x\res \R\big(\d{r}(x)\big)$. It follows that the function $x\mapsto r(x)$ from $\npseJ$ to $\omega$ is continuous, as are the functions $x\mapsto r(x)\!\mod k$ for $k\in\Z_+$. Next, note that for all $k\in\Z_+$ and $i<k$, the map $x\mapsto s^x(i)$ is continuous on the clopen set $\{x\in\npseJ\st r(x)\geq k\}$. These observations suffice to show that $c_2$, $c_3$, and $c_0$ are continuous on their domains. Finally, note that $\epsilon^{-1}[\{0\}]$ and $\epsilon^{-1}[\{1\}]$ are clopen subsets of $\M$, and $\M$ and $\B$ are clopen subsets of $\npseJ$, so $c_4$ is continuous as well.
\end{proof}

Eventually we will define a finite set $S_1$ and a continuous map $c_1:\Xse(J)\cap(\L_1\cup\R_1')\to S_1$, and, given this, put $\Sse=\Sse(J):=S_0\sqcup S_1\sqcup S_2^J\sqcup S_3^J\sqcup S_4^J$ and define our coloring $\cse^J:\Xse(J)\to\Sse(J)$ by 
\begin{equation}\label{eq:cse}
\cse(x) \ = \ \cse^J(x) \ := \ 
\begin{cases}
c_0(x) & \text{if $x\in(\L\setminus\L_1)\cup(\R\setminus\R_1')$;} \\
c_1(x) & \text{if $x\in\L_1\cup\R_1'$;} \\
c_2^J(x) & \text{if $x\in\A_2$;} \\
c_3^J(x) & \text{if $x\in\A_3$;} \\
c_4^J(x) & \text{if $x\in\B\cup\M$.}
\end{cases}
\end{equation}
That $\cse^J$ is continuous on $\Xse(J)$ will then follow from Claims \ref{claim:clopenSE} and \ref{claim:continuousSE} and continuity of $c_1$. 
The definitions of $S_1$ and $c_1$ will be given later, and for now we write $\tilde{c}_{\q}=\tilde{c}_{\q}^J:=c_0\sqcup c_2^J\sqcup c_3^J\sqcup c_4^J$ for the coloring we have already defined, and prove the following:

\begin{proposition}\label{prop:InitialColoring}
The continuous map $\tilde{c}_{\q}^J$ is a coloring of $\Gse(J)\res\big(\npseJ\setminus(\L_1\cup\R_1')\big)$.
\end{proposition}

\begin{proof}
Fix arbitrary $x,y\in\npseJ\setminus(\L_1\cup\R_1')$ such that $x\mathrel{\Gsea}y$. We will show $\tilde{c}_{\q}(x)\ne \tilde{c}_{\q}(y)$, arguing by cases depending on type. We may, of course, assume $\textup{tp}(x)=\textup{tp}(y)$, since otherwise $\tilde{c}_{\q}(x)\ne \tilde{c}_{\q}(y)$ trivially.

\smallskip

\noindent \underline{Case 1}: $\textup{tp}(x)=\textup{tp}(y)=\A_2$. Then $\tilde{c}_{\q}(x)=\big(\A_2,s^x(0)\big)$ and $\tilde{c}_{\q}(y)=\big(\A_2,s^y(0)\big)$, where 
\[
x(0,1)=s^x(0)\ne s^x(1)=x(1,0) \qquad \text{and} \qquad y(0,1)=s^y(0)\ne s^y(1)=y(1,0).
\]
Since $x(1,0)=y(0,1)$, it follows that $s^x(0)\ne s^y(0)$, so $\tilde{c}_{\q}(x)\ne \tilde{c}_{\q}(y)$.

\smallskip

\noindent \underline{Case 2}: $\textup{tp}(x)=\textup{tp}(y)=\A_3$. Then $\tilde{c}_{\q}(x)=\big( \A_3,\, s^x(0),\, s^x(1)\big)$ and $\tilde{c}_{\q}(y)=\big( \A_3,\, s^y(0),\, s^y(1)\big)$. If $s^x(0)\ne s^x(1)$, then $s^y(0)=s^x(1)\ne s^x(0)$, so $\tilde{c}_{\q}(x)$ and $\tilde{c}_{\q}(y)$ differ in their second term. If instead $s^x(0)=s^x(1)$, then since $s^x$ is non-constant we must have $s^y(1)=s^x(2)\ne s^x(1)$, so $\tilde{c}_{\q}(x)$ and $\tilde{c}_{\q}(y)$ differ in their third term.

\smallskip

\noindent \underline{Case 3}: $\textup{tp}(x)=\textup{tp}(y)=\B$. This case cannot arise, since $x\in\B$ implies $r(y)\leq r(x)$ and $y\in\B$ implies $r(x)\leq r(y)$, but if $x\in\B$ and $r(x)=r(y)$ then $y\in\M$. In other words, if $x\in\B$, then this fact alone distinguishes $x$ from its $\Gse$-neighbors.

\smallskip

\noindent \underline{Case 4}: $\textup{tp}(x)=\textup{tp}(y)=\M$. Using Lemma \ref{lem:M}, write $k:=r(x)=r(y)$. By the definition of $c_4$, if $s^x(0)\ne s^y(0)$ we are done, so assume otherwise. Then $s^x(0)=s^y(0)=s^x(1)$, so if $i_x$ and $i_y$ are least such that $s^x(i_x)\ne s^x(0)$ and $s^y(i_y)\ne s^y(0)$, respectively, then $i_x=i_y+1$. Thus $\epsilon(x)\ne\epsilon(y)$, so $\tilde{c}_{\q}(x)\ne \tilde{c}_{\q}(y)$ as desired.

\smallskip

\noindent \underline{Case 5}: $\textup{tp}(x)=\textup{tp}(y)=\L$. Since we are assuming $x,y\not\in\L_1$, the conclusion $\tilde{c}_{\q}(x)\ne \tilde{c}_{\q}(y)$ follows immediately from Lemma \ref{lem:LL} and the definition of $c_0$.

\smallskip

\noindent \underline{Case 6}: $\textup{tp}(x)=\textup{tp}(y)=\R$. Again, since we are assuming $x,y\not\in\R_1'$, the conclusion $\tilde{c}_{\q}(x)\ne \tilde{c}_{\q}(y)$ follows from Lemma \ref{lem:LL}, the definition of $c_0$, and the following obvious facts: $\delta[\L]=\R$ and $\delta[\R]=\L$; $\delta[\L_1]=\R_1'$ and $\delta[\R_1']=\L_1$; and $r(\delta(x))=r(x)$ for all $x\in\npseJ$.
\end{proof}

\begin{remark}
Our proof of Proposition \ref{prop:InitialColoring} and the lemmas leading up to it never used the definition of $\Xq$, but this will change when we consider $\L_1\cup\R_1'$. Indeed, by arguments similar to the proof of Observation \ref{claim:NoColoring}, one can show that $\chi_B\big(\Gq\res(\Wq\cap\L_1)\big)=\chi_B\big(\Gq\res(\Wq\cap\R_1')\big)=\aleph_0$.
\end{remark}

Now we turn our attention to defining $\cse$ on $\Xq\cap(\L_1\cup\R_1')$, where, again, we can focus on $\L_1$ by $\delta$-symmetry. Our approach is motivated by the following question: given $x,y\in\Xq\cap\L_1$ with $x\mathrel{\Gsea}y$, what happens if we repeatedly apply $\bar{e}_1$ to $x$ and $y$? The following lemmas provide some answers to this question.


\begin{lemma}\label{lem:L1bounded}
For all $x\in\Xq\cap\L_1$ and $m\in\omega$, if $k\bar{e}_1\cdot x\in\L_1$ for all $k\leq m$ then $r(m\bar{e}_1\cdot x)\geq 4+3m$; consequently, there is $m\in\omega$ such that $m\bar{e}_1\cdot x\not\in\L_1$.
\end{lemma}

\begin{proof}
Fix $x\in\Xse\cap\L_1$, and for each $k\in\omega$ write $x_k:=k\bar{e}_1\cdot x$. We argue by induction on $m$, the base case $m=0$ holding because $x\in\L$, so $r(x)\geq 4$ by definition. For the inductive step, fix $m\geq 0$ and assume that if $x_k\in\L_1$ for all $k\leq m$, then $r(x_m)\geq 4+3m$. Now suppose $x_k\in\L_1$ for all $k\leq m+1$, so in particular $r(x_m)\geq 4+3m$ by inductive hypothesis. Since $x_m\in\L_1$ we have $\d{1}\cdot x_m\in\L$ and $4+3m\leq r(x_m)\leq r(\d{1}\cdot x_m)-2$, so $6+3m\leq r(\d{1}\cdot x_m)$. As $x_{m+1}\in\L$, it then follows from Lemma \ref{lem:up2} that $r(x_{m+1})=r(\d{1}\cdot x_m)+1\geq 7+3m=4+3(m+1)$, completing the induction. But now if $x_m\in\L_1$ for all $m\in\omega$, then $r(\d{m}\cdot x)=r(m\bar{e}_2\cdot x_m)=r(x_m)-m\geq 4+2m$ for all $m$ by Lemma \ref{lem:upn}, contradicting $x\in\Xse$; hence there is some $m\in\omega$ for which $x_m\not\in\L_1$.
\end{proof}

\begin{DefinitionNotation}\label{df:mx*}
Using Lemma \ref{lem:L1bounded}, for each $x\in\Xse\cap\L_1$ let $m_x$ be the least $m\in\omega$ such that $m\bar{e}_1\cdot x\not\in\L_1$, and write $x^* := m_x\bar{e}_1\cdot x$.
\end{DefinitionNotation}

For example, with $x,y,z$ as in Figure \ref{fig:L1}, we have $m_x=3$, $m_y=2$, and $m_z=1$. 
Our coloring of the point $x\in\Xq\cap\L_1$ will depend on what we can say about $x^*$, and on how $x^*$ and $y^*$ are related when $x\mathrel{\Gsea}y$. Let us gather some additional information about $m_x$ and $x^*$ before outlining our strategy further.

\begin{lemma}\label{lem:xstar} 
For all $x\in\Xse\cap\L_1$, we have $x^*\in\L\cup\R\cup\B$. 
\end{lemma}

\begin{proof}
Fix $x\in\Xse\cap\L_1$, and again write $x_k:=k\bar{e}_1\cdot x$ for each $k\in\omega$, so $x^*=x_{m_x}$. Note that $m_x\geq 1$. Since $x_{m_x-1}\in\L_1$, we have $\bar{e}_2\cdot x^*=\d{1}\cdot x_{m_x-1}\in\L$, so  $x^*\in\L\cup\R_0\cup\B$ by Lemma \ref{lem:up2}(ii). And if $x^*\in\R_0$, then since $x_{m_x-1}\in\L$ we have $4\leq r(x_{m_x-1})\leq r(x^*)$ by Lemma \ref{lem:move1}(ii), so in fact $x^*\in\R$.
\end{proof}

\begin{lemma}\label{lem:L1tail}
For all $x,y\in\Xse\cap\L_1$ such that $x\mathrel{\Gsea}y$, we have $k\bar{e}_1\cdot x\in\L_1$ for all $k\leq m_y$; consequently, $m_y<m_x$.
\end{lemma}

\begin{proof}
Fix $x,y\in\Xq\cap\L_1$ with $x\mathrel{\Gsea}y$, and once again write $x_k:=k\bar{e}_1\cdot x$ for each $k\in\omega$. We show $x_k\in\L_1$ for each $k\leq m_y$ by induction on $k$, the base case $k=0$ holding by assumption. Fix $0\leq k<m_y$ and assume $x_k\in\L_1$. Since $k<m_y$, also $y_k:=k\bar{e}_1\cdot y\in\L_1$. Since $y_k\in\L$ and $2\bar{e}_2\cdot y_k=\d{1}\cdot x_k\in\L$, we have $2\leq r(y_k)-4$ by Lemma \ref{lem:upn}, so $x_{k+1}=\bar{e}_2\cdot y_k\in\L$ again by \ref{lem:upn}. Furthermore, since $y_k\in\L_1$, we have $\d{1}\cdot y_k\in\L$ with $4\leq r(y_k)\leq r(\d{1}\cdot y_k)-2$, so $\d{1}\cdot x_{k+1}=\bar{e}_2\cdot (\d{1}\cdot y_k)\in\L$ by Lemma \ref{lem:up2}(i). Finally, also by Lemma \ref{lem:up2} and the fact that $y_k\in\L_1$, we have $r(x_{k+1})=r(y_k)-1\leq r(\d{1}\cdot y_k)-3=r(\d{1}\cdot x_{k+1})-2$, so in fact $x_{k+1}\in\L_1$, completing the induction.
\end{proof}

\begin{lemma}\label{lem:L1first}
For all $x,y\in\Xse\cap\L_1$ such that $x\mathrel{\Gsea}y$, if $y^*\in\L$ then $m_y<m_x\leq m_y+2$.
\end{lemma}

\begin{proof}
Fix $x,y\in\Xse\cap\L_1$ with $x\mathrel{\Gsea}y$, and write $x':=m_y\bar{e}_1\cdot x$ and $x_1':=\bar{e}_1\cdot x'=\bar{e}_2\cdot y^*$, so $x'\mathrel{\Gsea}y^*$. Assume $y^*\in\L\setminus\L_1$, so $\d{1}\cdot y^*\not\in\L$ or $r(\d{1}\cdot y^*)\leq r(y^*)+1$. If $\d{1}\cdot y^*\not\in\L$, then $m_x\leq m_y+2$ since $(m_y+2)\bar{e}_1\cdot x=\d{1}\cdot y^*$. So assume $\d{1}\cdot y^*\in\L$ but $r(\d{1}\cdot y^*)\leq r(y^*)+1$. Then $r(\d{1}\cdot x_1')=r(\d{1}\cdot y^*)-1$ and $r(x_1')=r(y^*)-1$ by Lemma \ref{lem:up2}(i), so $r(\d{1}\cdot x_1') \leq r(y^*)=r(x_1')+1$ and thus $x_1'\not\in\L_1$, showing in this case that $m_x\leq m_y+1$. Since $m_y<m_x$ by Lemma \ref{lem:L1tail}, this completes the proof.
\end{proof}

We take a moment now to motivate what comes next. Suppose $x,y\in\Xq\cap\L_1$ with $x\mathrel{\Gsea}y$. We know $m_y<m_x$ by Lemma \ref{lem:L1tail}, and if $m_x$ is not too much larger than $m_y$, say $m_x\leq m_y+2$ as in the conclusion of Lemma \ref{lem:L1first}, then we can incorporate $m_x\!\mod 3$ into our coloring scheme in order to distinguish $x$ from $y$. But if $m_y\ll m_x$, which is possible if $y^*\in\R\cup\B$, then we cannot rely on ``local'' information about $x^*$ and $y^*$ to distinguish $x$ from $y$ since they may be very far apart.  To address this we will define, for each $x\in\Xq\cap\L_1$, a finite sequence $(x_i)$ with the property that when $m_y\ll m_x$, we have $y_i=x_{i+1}$ for all $i>m_y$. The lengths of the sequences determined by $x$ and $y$ will then have different parity, enabling us to distinguish $x$ from $y$ in this case.

The terms of the sequence $(x_i)$ just described will be iterates of $x$ under a function $\sigma$ which sends $x_i$ to $\bar{e}_1\cdot x_i$ if $x_i\in\L_1$, and to either $\d{1}\cdot x_i$ or $\bot$ (``halt'') otherwise. Given $x,y\in\Xq\cap\L_1$ with $x\mathrel{\Gsea}y$, we want $\sigma(y^*)=\d{1}\cdot y^*$ if $m_y\ll m_x$, but we also want the sequence $\sigma^n(y)$ to halt eventually, and this is where the definition of $\Xq$ will come into play. Let us turn to the details.

To define the conditions under which $\sigma(x)=\d{1}\cdot x$, let $\Phi\subseteq\npseJ$ be the clopen condition 
\begin{equation}
\Phi(x) \quad :\Longleftrightarrow \quad \bar{e}_2\cdot x\in\L_1 \ \wedge \ \d{1}\cdot x\in\L\cup\R\cup\B \ \wedge \ r(x)+2\leq r(\d{1}\cdot x),
\end{equation}
and define 
\begin{equation}
\R_2 \ := \ \big\{x\in\R\st \Phi(x)\big\} \qquad \text{and} \qquad \B_2 \ := \ \big\{x\in\B\st \Phi(x)\big\},
\end{equation}
so $\R_2$ and $\B_2$ are clopen in $\Xq$. Then define the continuous map $\sigma:\npseJ\cup\{\bot\}\to\npseJ\cup\{\bot\}$ by
\[
\sigma(x) \ := \ 
\begin{cases}
\bar{e}_1\cdot x & \text{if $x\in\L_1$;} \\
\d{1}\cdot x & \text{if $x\in\R_2\cup\B_2$;} \\
\bot & \text{otherwise.}
\end{cases}
\]

\begin{example}
Suppose the boxed ``$0$" in Figure \ref{fig:L1} were a ``$1$" instead. Then we would have
\[
\sigma^0(y)=y\in\L_1, \qquad \sigma^1(y)=\bar{e}_1\cdot y\in\R_2, \qquad \sigma^2(y)=(2,1)\cdot y\in\L_1,
\]
and
\[
\sigma^0(x)=x\in\L_1, \qquad \sigma^1(x)=\bar{e}_1\cdot x\in\L_1, \qquad \sigma^2(x)=2\bar{e}_1\cdot x\in\L_1, \qquad \sigma^3(x)=3\bar{e}_1\cdot x=\sigma^2(y)\in\L_1,
\]
so $\sigma^n(y)=\sigma^{n+1}(x)$ for all $n\geq 2$.
\end{example}

Symmetrically, we could define the condition $\Phi'(x)$, the sets $\L_2'$ and $\B_2'$, and the function $\sigma'$, but by $\delta$-symmetry we need only carry out the argument below once, and will not need to use $\Phi'$, $\L_2'$, $\B_2'$, or $\sigma'$. We do need to verify that for every $x\in\Xq\cap\L_1$, there is $n\in\Z_+$ such that $\sigma^n(x)=\bot$. Note that for all $x\in\Xse\cap\L_1$ and $n\in\Z_+$ such that $\sigma^n(x)\ne\bot$, we have $\sigma^k(x)\in\L_1\cup\R_2\cup\B_2$ for all $k<n$. Let us define, for each $x\in\npseJ$ and $n\in\omega$, the integers $K_n(x)$ and $M_n(x)$ by
\[
K_n(x) \ := \ \left|\big\{0<k\leq n\st \sigma^k(x)\in\L_1\big\}\right|
\]
and
\[
M_n(x) \ := \ \left|\big\{0<k\leq n\st \sigma^{k-1}(x)\in\R_2\cup\B_2\wedge\sigma^k(x)\in\R_2\cup\B_2\big\}\right|.
\]

Note that $K_0(x)=M_0(x)=0$ for all $x$. Note also by an easy induction that for all $x\in\L_1$ and $n\in\omega$, if $\sigma^{n+1}(x)\ne\bot$ then $\sigma^{n+1}(x)=\big(n+1, \, n-K_n(x)\big)\cdot x$.

\begin{lemma}\label{lem:rsigma}
For all $x\in\L_1$ and $n\in\omega$, if $\sigma^{n+1}(x)\ne\bot$ then $r\big(\sigma^n(x)\big)\geq 4+2K_n(x)+M_n(x)$. 
\end{lemma}

\begin{proof}
Fix $x\in\L_1$, and write $K_n:=K_n(x)$ and $M_n:=M_n(x)$. We argue by induction on $n$, the base case $n=0$ holding because $x\in\L_1$, so $r(x)\geq 4$. Fix $n\in\omega$ and assume for inductive hypothesis that $\sigma^{n+1}(x)\ne\bot$ implies $r\big(\sigma^n(x)\big)\geq 4+2K_n+M_n$. Suppose $\sigma^{n+2}(x)\ne\bot$. Then $\sigma^{n+1}(x)\ne\bot$, so $r\big(\sigma^n(x)\big)\geq 4+2K_n+M_n$. We must show $r\big(\sigma^{n+1}(x)\big)\geq 4+2K_{n+1}+M_{n+1}$. Write $y:=\sigma^n(x)$, so both $y$ and $\sigma(y)$ belong to $\L_1\cup\R_2\cup\B_2$. We consider the various cases.

First, suppose $y\in\L_1$, so $\sigma(y)=\bar{e}_1\cdot y$. If $\sigma(y)\in\L_1$, so $K_{n+1}=K_n+1$ and $M_{n+1}=M_n$, then since $\d{1}\cdot y\in\L$ with $r(\d{1}\cdot y)\geq r(y)+2$, by Lemma \ref{lem:up2}(i) we have
\[
r\big(\sigma^{n+1}(x)\big) \ = \ r(\bar{e}_1\cdot y) \ = \ r(\d{1}\cdot y)+1 \ \geq \ r(y)+3 \ \geq \ 4+2K_n+M_n+3 \ > \ 4+2K_{n+1}+M_{n+1}.
\]
On the other hand if $\sigma(y)\in\R_2\cup\B_2$, so $K_{n+1}=K_n$ and $M_{n+1}=M_n$, then by Lemma \ref{lem:move1} we have
\[
r\big(\sigma^{n+1}(x)\big) \ = \ r(\bar{e}_1\cdot y) \ \geq \ r(y) \ \geq \ 4+2K_n+M_n \ = \ 4+2K_{n+1}+M_{n+1}.
\]

Next, suppose $y\in\R_2\cup\B_2$, so $\sigma(y)=\d{1}\cdot y$ and $r(y)+2\leq r(\d{1}\cdot y)$. If $\sigma(y)\in\L_1$, then $K_{n+1}=K_n+1$ and $M_{n+1}=M_n$, while if $\sigma(y)\in\R_2\cup\B_2$, then $K_{n+1}=K_n$ and $M_{n+1}=M_n+1$; either way, we have
\[
r\big(\sigma^{n+1}(x)\big) \ = \ r(\d{1}\cdot y) \ \geq \ r(y)+2 \ \geq \ 4+2K_n+M_n+2 \ \geq \ 4+2K_{n+1}+M_{n+1}.
\]
This completes the induction, and the proof of the lemma.
\end{proof}

\begin{lemma}\label{lem:unbounded}
For all $x\in\L_1$ and $n\in\omega$, if $\sigma^{n+1}(x)\ne\bot$ then $r(\d{n}\cdot x)\geq 4+K_n(x)+M_n(x)$. 
\end{lemma}

\begin{proof}
Our proof is similar to that of the previous lemma. 
Again we fix $x\in\L_1$, write $K_n:=K_n(x)$ and $M_n:=M_n(x)$, and argue by induction on $n$, the base case $n=0$ holding because $x\in\L_1$. For the inductive step, fix $n\in\omega$ and assume $r(\d{n}\cdot x)\geq 4+K_n+M_n$ and $\sigma^{n+2}(x)\ne\bot$. Write $y:=\sigma^n(x)$, so both $y$ and $\sigma(y)$ belong to $\L_1\cup\R_2\cup\B_2$. Write $m:=n-K_n$, and note that $\sigma(y)=(n+1,m)\cdot x$. Suppose first that $y\in\L_1$, so $\d{1}\cdot y\in\L$ and $r(y)+2\leq r(\d{1}\cdot y)$, and additionally $\sigma(y)=\bar{e}_1\cdot y$ and $M_{n+1}=M_n$. By Lemma \ref{lem:rsigma}, we know $4+2K_n+M_n\leq r(y)$, so $6+2K_n+M_n\leq r(\d{1}\cdot y)$. Thus by Lemma \ref{lem:upn} and the fact that $K_{n+1}\leq K_n+1$, we have
\begin{align}\label{eq:gleeful}
r\big((\d{n}+\d{1})\cdot x\big) \ &= \ r\big((0,n+1-m)\cdot\sigma(y)\big) \ = \ r\big((0,n-m)\cdot (\d{1}\cdot y)\big) \ = \ r(\d{1}\cdot y)-(n-m) \\ &= \ r(\d{1}\cdot y)-K_n \ \geq \ 6+K_n+M_n \ > \ 4+K_{n+1}+M_{n+1}, \nonumber
\end{align}
regardless of whether $\sigma(y)$ belongs to $\L_1$ or $\R_2\cup\B_2$.

Now suppose $y\in\R_2\cup\B_2$, so $\sigma(y)=\d{1}\cdot y$ and $r(y)+2\leq r(\d{1}\cdot y)$. Then, using Lemma \ref{lem:rsigma}, we have
\[
4+2K_n+M_n \ \leq \ r(y) \ \leq \ r(\d{1}\cdot y)-2.
\]
If $\d{1}\cdot y\in\L_1$ then $K_{n+1}=K_n+1$ and $M_{n+1}=M_n$, and if $\d{1}\cdot y\in\R_2\cup\B_2$ then $K_{n+1}=K_n$ and $M_{n+1}=M_n+1$; either way, arguing as in \eqref{eq:gleeful} above, we have
\[
r\big((\d{n}+\d{1})\cdot x\big) \ = \ r(\d{1}\cdot y)-K_n-1 \ \geq \ 5+K_n+M_n \ = \ 4+K_{n+1}+M_{n+1}.
\]
This completes the induction, and the proof of the claim.
\end{proof}

\begin{corollary}\label{cor:terminates}
For all $x\in\Xse\cap\L_1$, there is $n\in\Z_+$ such that $\sigma^n(x)=\bot$.
\end{corollary}

\begin{proof}
Fix $x\in\L_1$, and suppose there is no $n\in\Z_+$ such that $\sigma^n(x)=\bot$. Then $\sigma^n(x)\in\L_1\cup\R_2\cup\B_2$ for all $n$. If there exist infinitely many $n$ for which $\sigma^n(x)\in\L_1$, then $\lim K_n(x)=\infty$, and if there are only finitely many $n$ for which $\sigma^n(x)\in\L_1$, then $\lim M_n(x)=\infty$. Either way, we have $\lim r(\d{n}\cdot x)=\infty$ by Lemma \ref{lem:unbounded}, so $x\not\in\Xse$.
\end{proof}

Using Corollary \ref{cor:terminates}, we introduce the following notation.

\begin{DefinitionNotation}\label{df:nx}
For each $x\in\Xse\cap\L_1$, let $n_x$ be the least $n\in\Z_+$ such that $\sigma^n(x)=\bot$.
\end{DefinitionNotation}

\begin{claim}\label{claim:c1continuity}
The maps $x\mapsto n_x$, $x\mapsto m_x$, and $x\mapsto x^*$ are continuous. 
\end{claim}

\begin{proof}
Continuity of $x\mapsto n_x$ follows from continuity of $\sigma$, and whether $m_x=m$ is determined by $x\res\R\big(\d{m}+\d{r}(m\bar{e}_1\cdot x)\big)$, so $x\mapsto m_x$ is continuous as well. Finally, $x\mapsto x^*$ is the composition of the continuous functions $x\mapsto (x,m_x)$ and $(x,m_x)\mapsto m_x\bar{e}_1\cdot x$.
\end{proof}

We are just about ready now to define our coloring of $\Gse$ on $\L_1\cup\R_1'$.  The next lemma may be viewed as a continuation of Lemma \ref{lem:L1first} now that we have defined $\R_2$, $\B_2$, and $\sigma$; together, these two lemmas will ensure that our coloring actually works.

\begin{lemma}\label{lem:L1combined}
Let $x,y\in\Xse\cap\L_1$ with $x\mathrel{\Gsea}y$.
\begin{enumerate}[(i)]
 \item If $y^*\in\big(\R\setminus\R_2\big)\cup\big(\B\setminus\B_2\big)$, then $m_y<m_x\leq m_y+2$.
 \item If $y^*\in\R_2\cup\B_2$, then $n_x=n_y+1$.
\end{enumerate}
\end{lemma}

\begin{proof}
Fix $x,y\in\Xse\cap\L_1$ with $x\mathrel{\Gsea}y$. For (i), assume $y^*\in\big(\R\setminus\R_2\big)\cup\big(\B\setminus\B_2\big)$. We already know $m_y<m_x$ by Lemma \ref{lem:L1tail}. If $\bar{e}_2\cdot y^*\in\L_1$ and $\d{1}\cdot y^*\in\L_1$, then $(1,2)\cdot y^*\in\L$ and by Lemmas \ref{lem:move1}(i) and \ref{lem:up2}(i) and the definition of $\L_1$ we have
\[
r(y^*)+2 \ \leq \ r(\bar{e}_2\cdot y^*)+3 \ \leq \ r\big((1,2)\cdot y^*\big)+1 \ = \ r(\d{1}\cdot y^*).
\]
Therefore, since $y^*\not\in\R_2\cup\B_2$, we must have either $\bar{e}_2\cdot y^*\not\in\L_1$ or $\d{1}\cdot y^*\not\in\L_1$; since $\bar{e}_2\cdot y^*=(m_y+1)\bar{e}_1\cdot x$ and $\d{1}\cdot y^*=(m_y+2)\bar{e}_1\cdot x$, the former implies $m_x\leq m_y+1$ and the latter implies $m_x\leq m_y+2$.

For (ii), assume $y^*\in\R_2\cup\B_2$ and write $x':=m_y\bar{e}_1\cdot x$. Then $x'=\sigma^{m_y}(x)\in\L_1$ by Lemma \ref{lem:L1tail}, and $\bar{e}_2\cdot y^*=\bar{e}_1\cdot x'\in\L_1$ since $y^*\in\R_2\cup\B_2$. 
Thus $\sigma^{m_y+2}(x)=2\bar{e}_1\cdot x'=\d{1}\cdot y^*=\sigma^{m_y+1}(y)\ne\bot$, so $n_x>m_y+2$ and $n_y>m_y+1$, but also $\sigma^{k+2}(x)=\sigma^{k+1}(y)$ for all $k\geq m_y$. It follows that $n_x=n_y+1$.
\end{proof}

Now we complete our coloring of $\Gse\res\Xse$ by defining the set $S_1$ and the continuous function $c_1:\Xq\cap (\L_1\cup\R_1')\to S_1$, which will not depend on $J$. The cases $\L_1$ and $\R_1'$ are symmetric, and we start with $\L_1$. For each $x\in\Xse(J)\cap\L_1$, define 
\begin{equation}
\hat{c}_1(x) \ := \ 
\begin{cases}
m_x\!\mod 3 & \text{if $x^*\in\L\cup(\R\setminus\R_2)\cup(\B\setminus\B_2)$;} \\
3+(n_x\!\mod 2) & \text{if $x^*\in\R_2\cup\B_2$.}
\end{cases}
\end{equation}
Here we are using the fact (Lemma \ref{lem:xstar}) that $x^*\in\L\cup\R\cup\B$ for all $x\in\Xse\cap\L_1$. 
Then, given arbitrary $x\in\Xse(J)\cap (\L_1\cup\R_1')$, let
\[
c_1(x) \ := \ 
\begin{cases}
 \big(\L_1, \, \hat{c}_1(x)\big) & \text{if $x\in\L_1$;} \\
 \big(\R_1', \, \hat{c}_1(\delta(x))\big) & \text{if $x\in\R_1'$.}
\end{cases}
\]
Note that $c_1$ is indeed continuous by Claim \ref{claim:c1continuity} and the fact that all the sets involved in the piecewise definitions of $c_1$ and $\hat{c}_1$ are clopen. Letting $S_1:=\ran(c_1)$, we see that $|S_1|=10$, so by \eqref{eq:Jcount} we have
\[
|\Sse(J)| \ = \ (J^2+3J+9)+10 \ = \ J^2+3J+19.
\]
This completes the defintion of $\cse=\cse^J$ that was given in Equation \eqref{eq:cse}. 

\begin{proposition}\label{prop:colorGSE}
For all $J\in\Z_{\geq 2}$, the continuous map $\cse^J$ is a coloring of $\Gse(J)\res\Xse(J)$. In particular,
\[
\chi_C\big(\Gse(J)\res\Xse(J)\big) \ \leq \ J^2+3J+19.
\]
\end{proposition}

\begin{proof}
This is, of course, just a continuation of the proof of Proposition \ref{prop:InitialColoring}. Fixing $J\in\Z_{\geq 2}$ and $x,y\in\Xse(J)$ such that $x\mathrel{\Gsea}y$, it will suffice by that result to consider the following two cases.

\smallskip

\noindent \underline{Case 7}: $x,y\in\L_1$. It will suffice to show $\hat{c}_1(x)\ne\hat{c}_1(y)$, so we may assume $x^*$ and $y^*$ either both belong to $\L\cup(\R\setminus\R_2)\cup(\B\setminus\B_2)$ or both belong to $\R_2\cup\B_2$. In the former case, $m_y<m_x\leq m_y+2$ by Lemmas \ref{lem:L1first} and \ref{lem:L1combined}(i), so $m_x\not\equiv m_y\!\mod 3$ and thus $\hat{c}_1(x)\ne\hat{c}_1(y)$; in the latter, $n_x\not\equiv n_y\!\mod 2$ by Lemma \ref{lem:L1combined}(ii), so again $\hat{c}_1(x)\ne\hat{c}_1(y)$. 

\smallskip

\noindent \underline{Case 8}: $x,y\in\R_1'$. Then $\delta(x),\delta(y)\in\L_1$ and $\delta(y)\mathrel{\Gsea}\delta(x)$ by Lemma \ref{lem:delta}, so $\hat{c}_1(\delta(x))\ne \hat{c}_1(\delta(y))$ by Case 7, which implies $c_1(x)\ne c_1(y)$ and therefore $\cse(x)\ne\cse(y)$.
\end{proof}

\begin{question}
Proposition \ref{prop:colorGSE} implies $\chi_B\big(\Gq(2)\res\Xq(2)\big)\leq 29$, but what is the exact value of $\chi_B\big(\Gq(2)\res\Xq(2)\big)$, or more generally $\chi_B\big(\Gq(J)\res\Xq(J)\big)$? Does the value depend on $J$?
\end{question}


We have now finished proving the first half of Theorem \ref{thm:MainSE}, and are ready to turn our attention to the second half and show that $\EJT\res\Wse$ is hypersmooth. Before describing our proof strategy, we state a few more definitions and prove some facts about them.

For each $x\in\Wse$, let
\[
A_x \ := \ \big\{n\in\omega\st (\forall m>n)\; 4\leq r(\d{n}\cdot x)<r(\d{m}\cdot x)\big\},
\]
so $A_x$ is infinite and the sequence $\big(r(\d{n}\cdot x)\st n\in A_x\big)$ strictly increases. Then let
\[
\Wqnw \ := \ \big\{x\in\Wse\st (\exists^\infty n\in A_x) \, \d{n}\cdot x\in\L\cup\B\big\},
\]
and put
\[
\tWqnw \ := \ \big[\Wqnw\big]_{\EJT} \qquad \text{and} \qquad \tWqse \ := \ \Wse\setminus\tWqnw.
\]

\begin{lemma}\label{lem:nxLRB}
For all $x\in\Wse$ and $n\in A_x$, we have $\d{n}\cdot x\in\L\cup\R\cup\B$.
\end{lemma}

\begin{proof}
Fix $x\in\Wse$ and $n\in A_x$, and write $k:=r(\d{n}\cdot x)\geq 4$ and $x':=(\d{n}+\d{1})\cdot x$. The sequence
\[
s^{x'}_{k-2} \ = \ \big\langle s^{\d{n}\cdot x}(1)\,,\, \ldots\,,\, s^{\d{n}\cdot x}(k-2)\big\rangle
\]
is constant since $r(x')>k$, which means $\d{n}\cdot x\in\L\cup\R\cup\B$ as claimed.
\end{proof}

\begin{lemma}\label{lem:SEtoNW}
For all $x\in\tWqse$ we have $\delta(x)\in\Wqnw$.
\end{lemma}

\begin{proof}
Let $x\in\tWqse$. Then $x\in\Wse$, so $\delta(x)\in\Wse$ as well by Claim \ref{claim:postponed}. Since $x\not\in\tWqnw$, by Lemma \ref{lem:nxLRB} we have $\d{n}\cdot x\in\R$ for all but finitely many $n\in A_x$, and therefore $\delta(\d{n}\cdot x)\in\L$ for all but finitely many $n\in A_x$. But $\delta(\d{n}\cdot x)=\d{n}\cdot\delta(x)$ and $A_x=A_{\delta(x)}$, so $\d{n}\cdot\delta(x)\in\L$ for all but finitely many $n\in A_{\delta(x)}$, which means $\delta(x)\in\Wqnw$ as desired.
\end{proof}

\begin{lemma}\label{lem:NWclosure}
The sets $\tWqnw$ and $\tWqse$ are Borel and $\EJT$-invariant. Furthermore, for all $x\in\Wqnw$, all $\bar{c}\in\N^2$ satisfying $c_1\leq c_2$, and all $n\geq c_1$ in $A_x$ such that $\d{n}\cdot x\in\L\cup\B$ and $r(\d{n}\cdot x)\geq 4+c_2-c_1$, we have $n-c_1\in A_{\bar{c}\cdot x}$ and $(\d{n}-\d{c}_1)\cdot (\bar{c}\cdot x)\in\L\cup\B$; in particular, $\bar{c}\cdot x\in\Wqnw$.
\end{lemma}

\begin{proof}
Clearly $\Wqnw$ is Borel. To see that $\tWqnw$ is Borel, note that $\Wqnw$ is closed under the diagonal $\N$-action, so if $x\in\Wqnw$ and $\bar{a}\cdot x=\bar{b}\cdot y$, then $\d{n}\cdot x=(\d{n}-\bar{a}+\bar{b})\cdot y\in\Wqnw$ where $n:=|\bar{a}|$; this shows $\tWqnw=\bigcup_{\bar{c}\mkern1mu\in\mkern1mu\N^2}\bar{c}^{-1}\big[\Wqnw\big]$, so $\tWqnw$ is indeed Borel. Hence $\tWqse$ is Borel as well, and both are $\EJT$-invariant by definition.

For the second claim, fix $x\in\Wqnw$ and $\bar{c}=(c_1,c_2)\in\N^2$ with $c_1\leq c_2$, and fix $n\geq c_1$ in $A_x$ such that $\d{n}\cdot x\in\L\cup\B$ and $r(\d{n}\cdot x)\geq 4+c_2-c_1$. Write $x':=(\d{n}-\d{c}_1)\cdot (\bar{c}\cdot x)=(c_2-c_1)\bar{e}_2\cdot (\d{n}\cdot x)$. Since $c_2-c_1\leq r(\d{n}\cdot x)-4$, we have $x'\in\L\cup\B$ with $r(x')=r(\d{n}\cdot x)-c_2+c_1$ by Lemma \ref{lem:upn}. Now let $m>n-c_1$ be arbitrary. Then $r(\d{n}\cdot x)<r\big((\d{m}+\d{c}_1)\cdot x\big)$ since $n\in A_x$, and $r\big((\d{m}+\d{c}_1)\cdot x\big) \leq r\big(\d{m}\cdot (\bar{c}\cdot x)\big)+c_2-c_1$ by repeated use of Lemma \ref{lem:move1}. Thus 
\[
4 \ \leq \ r(x') \ = \ r(\d{n}\cdot x)-c_2+c_1 \ < \ r\big((\d{m}+\d{c}_1)\cdot x\big)-c_2+c_1 \ \leq \ r\big(\d{m}\cdot (\bar{c}\cdot x)\big),
\]
which since $m>n-c_1$ was arbitrary shows $n-c_1\in A_{\bar{c}\cdot x}$. 
Since there are infinitely many $n\geq c_1$ in $A_x$ for which $\d{n}\cdot x\in\L\cup\B$ and $r(\d{n}\cdot x)\geq 4+c_2-c_1$, it follows that $\bar{c}\cdot x\in\Wqnw$.
\end{proof}

A sample element $x\in\Wqnw(2)$ is pictured in Figure \ref{fig:Wq}, along with the values $r(\d{n}\cdot x)$ and $\textup{tp}(\d{n}\cdot x)$ for small $n$ in the first table in that figure, where rows with $n\in A_x$ are shaded. We fill $x$ almost entirely with zeros here to make the picture more intelligible; while this is not representative of a generic element of $\Wq$, it will make our proof strategy easier to visualize. 
We will show $\EJT\res\Wqnw$ is hypersmooth by associating to each $x\in\Wqnw$, in a manner similar to our proof of Lemma \ref{lem:Ge1HSTechLem} and Proposition \ref{prop:WeHS}, a sequence $s(x)\in (\N^2)^\omega$ such that the map $x\mapsto \big(s_n(x)\cdot x\big)_{n\in\omega}$ is (after coding) a Borel reduction from $\EJT\res\Wqnw$ to $E_{1,t}$; then we reduce the rest of $\EJT\res\Wq$ to this case. The solid dots in Figure \ref{fig:Wq} represent the first few points in the sequence $s_n(x)$, and are determined by the positions of the 1s above the main diagonal. In order to describe this sequence for an arbitrary element of $\Wqnw$, we make some more definitions and prove a couple more lemmas, after which the second table in Figure \ref{fig:Wq} will make sense.

Given $x\in\JT$, let us say that $x$ \emph{splits at $k$} if $x(k,k+1)\ne x(k+1,k)$, and that $x$ \emph{splits} if there is $k\in\omega$ such that $x$ splits at $k$. Note that $x$ splits at $k$ if and only if $r(\d{k}\cdot x)=2$.

\begin{lemma}
Let $x\in\Wqnw$.
\begin{enumerate}[(i)]
 \item There is $j\in\omega$ such that $j\bar{e}_2\cdot x$ splits.
 \item There exist at most finitely many $k\in\omega$ such that $x$ splits at $k$.
\end{enumerate}
\end{lemma}

\begin{proof}
Fix $x\in\Wqnw$. 
For (i), if $n\in A_x$ is such that $\d{n}\cdot x\in\L\cup\B$ and $j:=r(\d{n}\cdot x)-2$, then $j\geq 2$ and
\[
(j\bar{e}_2\cdot x)(n, n+1) \ = \ (\d{n}\cdot x)\big(0, \, r(\d{n}\cdot x)-1\big) \ \ne \ (\d{n}\cdot x)\big(1, \, r(\d{n}\cdot x)-2\big) \ = \
(j\bar{e}_2\cdot x)(n+1,n),
\]
so $j\bar{e}_2\cdot x$ splits at $n$. For (ii), using $\lim r(\d{n}\cdot x)=\infty$, if $N\in\omega$ is such that $r(\d{n}\cdot x)>2$ for all $n\geq N$, then there is no $k\geq N$ for which $x$ splits at $k$.
\end{proof}

For each $x\in\Wqnw$, let $j^x$ be the least $j\in\omega$ such that $j\bar{e}_2\cdot x$ splits, and let $k^x$ be the largest $k\in\omega$ such that $j^x\bar{e}_2\cdot x$ splits at $k$, so $(k^x, \, j^x+k^x)\cdot x$ splits at $0$ but does not split at any $n>0$. Now let $f,g:\Wqnw\to\Wqnw$ be the Borel maps defined by
\[
f(x) \ := \ \big(k^x, \, j^x+k^x\big)\cdot x \qquad \text{and} \qquad g(x) \ := \ f(\d{1}\cdot x).
\]

Note that $f$ and $g$ really do take values in $\Wqnw$ by Lemma \ref{lem:NWclosure}, and that $\ran(g)\subseteq\ran(f)$. For each $x\in\Wqnw$, define sequences $j^x_n$ and $k^x_n$ inductively for $n\in\omega$ by $j^x_0:=0$ and $k^x_0:=0$ and
\begin{equation}
j^x_{n+1} \ := \ j^{\d{1}\cdot g^n(x)} \qquad \text{and} \qquad k^x_{n+1} \ := \ k^{\d{1}\cdot g^n(x)}.
\end{equation}
Finally, for each $x\in\Wqnw$ and $n\in\omega$ let 
\[
J^x_n:=\sum_{i=0}^nj^x_i \qquad \text{and} \qquad K^x_n:=n+\sum_{i=0}^nk^x_i,
\]
so $\lim K^x_n=\infty$.

The sequences $(K_n^x)$ and $(J^x_n)$ are defined so that $g^n(x)=(K^x_n, \, J^x_n+K^x_n)\cdot x$ for all $n\in\omega$, and we verify this along with some other basic properties of these sequences in Lemmas \ref{lem:jk} and \ref{lem:NoSplitting} below. For a visual example the reader may again consult Figure \ref{fig:Wq}, where the first few values of $j^x_n$, $k^x_n$, $J^x_n$, and $K^x_n$ are included in the second table. Our goal will be to show that $s_n(x):=(K^x_n, \, J^x_n+K^x_n)$ is the desired sequence described above, so that the map $x\mapsto \big(g^n(x)\big)_{n\in\omega}$ is (after coding) a Borel reduction from $\EJT$ to $E_{1,t}$. Notice the similarity between this approach and our proof of Proposition \ref{prop:WeHS}, and compare Figure \ref{fig:Wq} with Figure \ref{fig:Ge1}.

\begin{lemma}\label{lem:jk}
Let $x\in\Wqnw$.
\begin{enumerate}[(i)]
 \item For all $n\in\omega$, $g^n(x) = \big(K^x_n,\, J^x_n+K^x_n\big)\cdot x$.
 \item We have $j^x_n>0$ for all $n\in\Z_{\geq 2}$, and $J^x_{n+1}>J^x_n\geq n-1$ for all $n\in\Z_{+}$; in particular, $\lim J^x_n=\infty$.
 \item For all $n\in\Z_{+}$, $J^x_n\bar{e}_2\cdot x$ splits at $K^x_n$ but does not split at any $k>K^x_n$.
\end{enumerate}
\end{lemma}

\begin{proof}
Fix $x\in\Wqnw$. We prove (i) by induction on $n$, with the base case $n=0$ holding because $K^x_n=J^x_n=0$. For the induction, fix $n\in\omega$ and assume $g^n(x)=(K^x_n, \, J_n^x+K_n^x)\cdot x$. Then 
\begin{align*}
g^{n+1}(x) \ 
&= \ g\big(g^n(x)\big) \ = \ f\big(\d{1}\cdot g^n(x)\big) 
\ = \ \Big(k^x_{n+1}, \, j^x_{n+1}+k^x_{n+1}\Big)\cdot \big(\d{1}\cdot g^n(x)\big) \\
&= \ \Big(K^x_n+k^x_{n+1}, \, J^x_n+K^x_n+j^x_{n+1}+k^x_{n+1}\Big)\cdot (\d{1}\cdot x) \ = \ \big(K^x_{n+1}, \, J^x_{n+1}+K^x_{n+1}\big)\cdot x.
\end{align*}

For (ii), given arbitrary $n\in\Z_+$ we have $g^n(x)\in\ran(f)$, so $g^n(x)$ does not split at any $k>0$, so $\d{1}\cdot g^n(x)$ does not split at all, which means $j^x_{n+1}=j^{\d{1}\cdot g^n(x)}>0$. This shows $j^x_n>0$ for all $n\in\Z_{\geq 2}$. To show $J^x_{n+1}>J^x_n\geq n-1$ for all $n\in\Z_{+}$, we argue by induction, the base case $n=1$ holding because $J^x_2=J^x_1+j^x_2>J^x_1\geq 0$ by what we just showed. For the induction, if $J^x_{n+1}>J_n^x\geq n-1$, so $J_{n+1}^x\geq n$, then $J^x_{n+2}=J^x_{n+1}+j^x_{n+2}>J^x_{n+1}\geq n$ as required.

Finally, for (iii),  let $n\in\omega$ be arbitrary. Then $y:=\dot{K}^x_{n+1}\cdot (J^x_{n+1}\bar{e}_2\cdot x)=g^{n+1}(x)\in\ran(f)$, so $y$ splits at $0$ but not at any $k>0$, which means $J^x_{n+1}\bar{e}_2\cdot x$ splits at $K^x_{n+1}$ but not at any $k>K^x_{n+1}$. 
\end{proof}

\begin{lemma}\label{lem:NoSplitting}
For all $x\in\Wqnw$, $n\in\Z_+$, and $j<J^x_n$, there is no $k\geq K^x_n$ such that $j\bar{e}_2\cdot x$ splits at $k$.
\end{lemma}

\begin{proof}
Fixing $x\in\Wqnw$, we argue by induction on $n$. For the base case $n=1$, we have $J^x_1=j^{\d{1}\cdot x}$ and $K^x_1=1+k^{\d{1}\cdot x}\geq 1$. By definition of $j^{\d{1}\cdot x}$, there is no $j<j^{\d{1}\cdot x}$ such that $j\bar{e}_2\cdot (\d{1}\cdot x)$ splits; therefore, given any $j<J^x_1$, there is no $k\in\omega$ such that $j\bar{e}_2\cdot (\d{1}\cdot x)$ splits at $k$, and hence no $k\geq 1$ such that $j\bar{e}_2\cdot x$ splits at $k$. Since $K^x_1\geq 1$, this establishes the base case. For the induction, fix $n\in\Z_+$ and assume for inductive hypothesis that for all $j<J^x_n$ there is no $k\geq K^x_n$ such that $j\bar{e}_2\cdot x$ splits at $k$. We must show that for all $j<J^x_{n+1}$,
\begin{equation}\tag{$\ast$}
\text{there is no $k\geq K^x_{n+1}$ such that $j\bar{e}_2\cdot x$ splits at $k$.}
\end{equation}
Since $K^x_{n+1}>K^x_n$, $(\ast)$ holds for $j<J^x_n$ by inductive hypothesis, and Lemma \ref{lem:jk}(iii) implies that $(\ast)$ holds for $j=J^x_n$. So assume $J^x_n<j<J^x_{n+1}=J^x_n+j^x_{n+1}$. Write $j':=j-J^x_n$, and also $x':=g^n(x)=\big(K^x_n,\,J^x_n+K^x_n\big)\cdot x$, so $j^x_{n+1}=j^{\d{1}\cdot x'}$. Then $0<j'<j^{\d{1}\cdot x'}$, so $j'\bar{e}_2\cdot (\d{1}\cdot x')$ does not split. Since $(\d{K}^x_n+\d{1})\cdot (j\bar{e}_2\cdot x)=j'\bar{e}_2\cdot (\d{1}\cdot x')$, this means there is no $k\geq K_n^x+1$ such that $j\bar{e}_2\cdot x$ splits at $k$, and therefore no $k\geq K^x_{n+1}$ such that $j\bar{e}_2\cdot x$ splits at $k$. This completes the induction, and the proof.
\end{proof}

\begin{figure}[h]
\begin{center}
\begin{tikzpicture}[x=1.25cm, y=1.25cm]
     \fill [fill=gray!10] (2.1,.3) rectangle (5.9,7.6);
  \foreach \x in {.3,.6} \draw [lightgray] (\x,0) -- (\x,7.6);
  \foreach \x in {.9,1.2,...,6} \draw [lightgray] (\x,-.6) -- (\x,7.6);
  \draw [lightgray] (.9,-.3) -- (5.9,-.3);
  \foreach \y in {.3,.6,...,7.8} \draw [lightgray] (0,\y) -- (5.9,\y);
 \draw [gray!60,dashed] (0,0) -- (5.9,6);
 \draw [gray!60,dashed] (.9,-.6) -- (5.9,4.5);
 \draw [gray!60,dashed] (0,.3*10) -- (.3*15.3,.3*25.3);
 \draw [dotted] (0,0) -- (.9,-.6);
 \draw [thick] (0,0) -- (5.9,0);
 \draw [thick] (0,0) -- (0,7.6);
 \draw [thick] (.9,-.6) -- (5.9,-.6);
 \draw [thick] (.9,-.6) -- (.9,7.6);
   %
  \draw (0,0) -- (2.1,.3);
  \draw (.9,-.6) -- (2.1,.3);
  \draw (0,0) -- (.3*13,.3*23);
  \draw (.9,-.6) -- (.3*13,.3*23);
  \draw [gray] (.3*4,.3*5) -- (.3*4,.3*12);
  \draw [gray] (.3*4,.3*12) -- (.3*6,.3*14);
  \draw [gray] (0,0) -- (.3*6,.3*14);
  \draw [gray] (0,0) -- (.3*3,.3*4);
  \node at (0,0) [circle,fill,inner sep=1.2pt]{};
  \node at (.9,-.6) [circle,fill,inner sep=1pt]{};
(.05+.3*3,.05-.3*2) [circle,fill,inner sep=1pt]{};
  \node at (.3*13,.3*23) [circle,fill,inner sep=1pt]{};
  \node at (.3*4,.3*5) [circle,fill,inner sep=1pt]{};
  \node at (.3*6,.3*14) [circle,fill=gray!80,inner sep=1pt]{};
  \node at (.17+.3*4,.17+.3*5) {$x'$};
    \node at (.155+.3*13,.155+.3*24) {{\color{gray}$0$}};
    \node at (.155+.3*14,.155+.3*23) {{\color{gray}$1$}};
    \node at (.155+.3*6,.155+.3*15) {{\color{gray}$0$}};
    \node at (.155+.3*7,.155+.3*14) {{\color{gray}$1$}};
   \node at (1.6,-.3) {$\bar{b}$};
   \node at (.74,.23) {$\bar{a}$};
   \node at (3.2,4.66) {$\bar{b}'$};
   \node at (2.65,4.95) {$\bar{a}'$};
   \node at (.741,1.03) {$\bar{c}$};
   \node at (.704,1.96) {$\bar{d}$};
   \node at (1.24,3.3) {$j^{x'}$};
   \node at (1.47,4.05) {$\d{k}^{x'}$};
   \node at (3.52,7.05) {$g^n(y)$};
 \node at (.33,-.41) {$\bar{m}$};
 \node at (-.15,-.15) {$x$};
 \node at (.75,-.75) {$y$};
%
\end{tikzpicture}
\end{center}
\caption{A picture for the proofs of Lemma \ref{lem:technical} and Proposition \ref{prop:NWhs}.}
\label{fig:technical}
\end{figure}
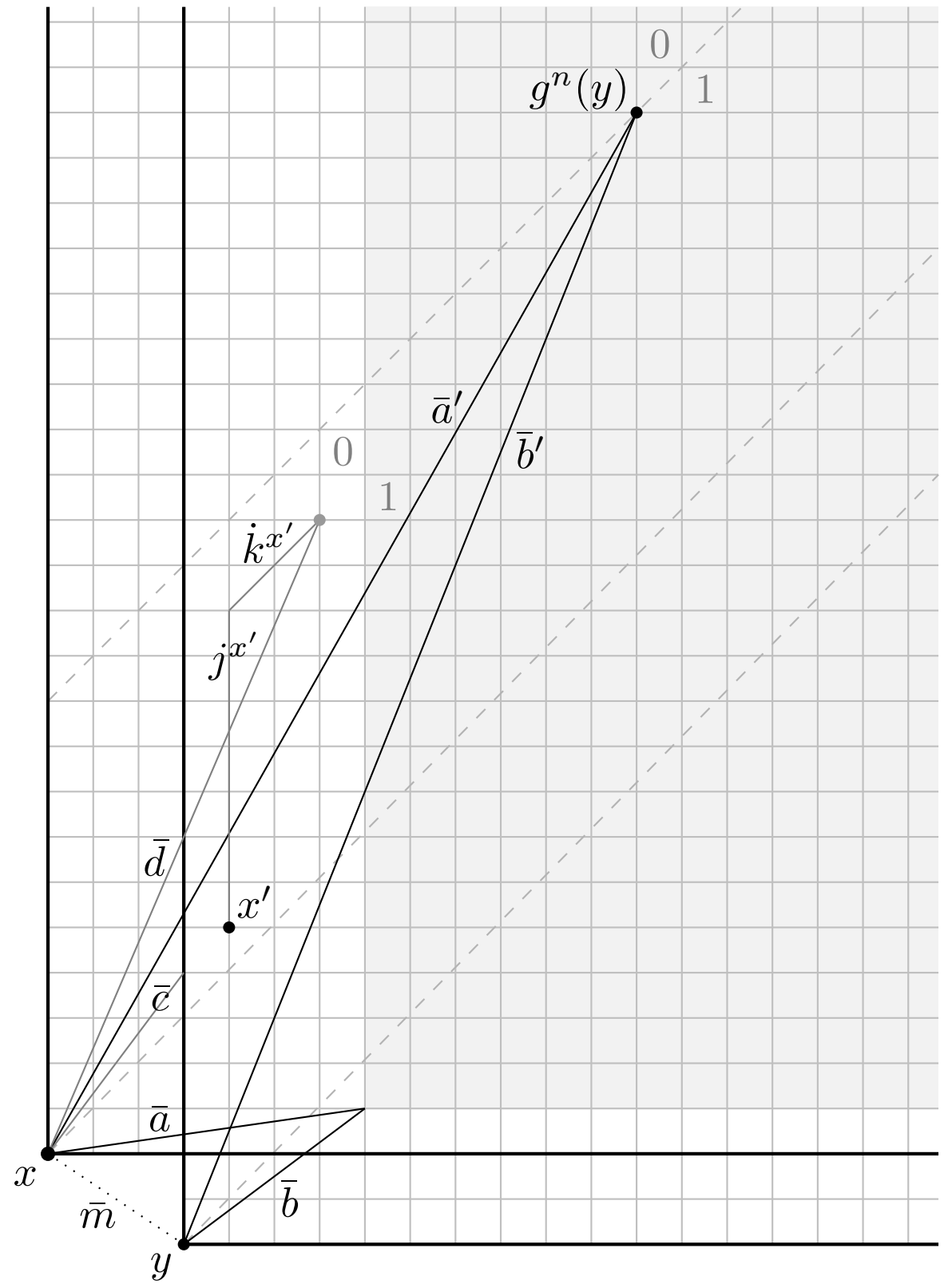

We are now ready to show that $\EJT\res\Wqnw$ is hypersmooth. We split the proof into two pieces, consisting of a technical lemma and the main argument, both of which are best understood by reference to Figure \ref{fig:technical}. First, the technical lemma:

\begin{lemma}\label{lem:technical}
Let $x,y\in\Wqnw$ and $\bar{a},\bar{b}\in\N^2$ with $\bar{m}:=\bar{a}-\bar{b}$, and suppose $\bar{a}\cdot x=\bar{b}\cdot y$ and $m_1\geq m_2$. Fix $n\in\Z_+$ such that $K^y_n>\max\{|\bar{a}|,|\bar{b}|\}$, and let $\bar{b}'=(b_1',b_2'):=(K_n^y, J_n^y+K_n^y)$ and $\bar{a}'=(a_1',a_2'):=\bar{a}+\bar{b}'-\bar{b}$, so $\bar{a}'\cdot x=\bar{b}'\cdot y=g^n(y)$. Suppose $\bar{c}=(c_1,c_2)\in\N^2$ satisfies $c_1\leq c_2$ and $0<a_1'-c_1\leq a_2'-c_2$, and write $x':=(\d{1}+\bar{c})\cdot x$ and $\bar{d}:=\d{1}+\bar{c}+(k^{x'},\, j^{x'}+k^{x'})$. Then either $\bar{d}=\bar{a}'$ or else $0<a_1'-d_1<a_2'-d_2$.
\end{lemma}

\begin{proof}
Assume the notation and hypotheses in the statement of the lemma. Write $u:=a_1'-c_1-1$ and $v:=a_2'-c_2-1-u$,  so $u,v\in\N$ and $(u,u+v)\cdot x'=\bar{a}'\cdot x=g^n(y)$. Since $g^n(y)$ splits at $0$, we have that $v\bar{e}_2\cdot x'$ splits at $u$, so $j^{x'}\leq v$. If $j^{x'}=v$, then we must have $k^{x'}=u$ since $g^n(y)$ splits at $0$ but 
does not split at any $k>0$. Thus if $j^{x'}=v$, then $\bar{d}=\d{1}+\bar{c}+(u,u+v)=\bar{a}'$ and we are done. Assume, therefore, that $j^{x'}<v$. We will show $k^{x'}<u$. This will suffice, since then
\[
0 \ < \ u-k^{x'} \ = \ a_1'-d_1 \ < \ u-k^{x'}+v-j^{x'} \ = \ a_2'-d_2.
\]

Assume for contradiction that $k^{x'}\geq u$. Since $j^{x'}\bar{e}_2\cdot x'$ splits at $k^{x'}$, we have that $(j^{x'}+c_2-c_1)\bar{e}_2\cdot x$ splits at $k^{x'}+c_1+1$. Let us write $j_0:=j^{x'}+c_2-c_1$ and $k_0:=k^{x'}+c_1+1$, so that $j_0\bar{e}_2\cdot x$ splits at $k_0$. Note that
\[
k_0 \ = \ k^{x'}+c_1+1 \ \geq \ u+c_1+1 \ = \ a_1',
\]
so let us write $k_1:=k_0-a_1'\in\N$, so $(a_1',\, a_1'+j_0)\cdot x$ splits at $k_1$. Let us also note that 
\[
a_1'+v+c_2-c_1 \ = \ a_1'+(a_2'-c_2-1-a_1'+c_1+1)+c_2-c_1 \ = \ a_2',
\]
so 
\[
a_1'+j_0 \ < \ a_1'+v+c_2-c_1 \ = \ a_2'.
\]
Now recall that we are assuming $m_1\geq m_2$, that is, $a_1-b_1\geq a_2-b_2$. Therefore,
\[
a_1' \ = \ a_1+b_1'-b_1 \ = \ (a_1-b_1)+K_n^y \ > \ a_2-b_2+b_2 \ = \ a_2,
\]
but also $a_1'=a_1+(b_1'-b_1)>a_1$. From $\bar{a}\cdot x=\bar{b}\cdot y$ and $a_1'>a_1$, we get
\[
(a_1', a_2)\cdot x \ = \ (b_1+a_1'-a_1, \, b_2)\cdot y, \qquad \text{where} \qquad a_2<a_1'+j_0<a_2'.
\]
Therefore,
\[
(a_1', \, a_1'+j_0)\cdot x \ = \ \big(b_1+a_1'-a_1, \, b_2-a_2+a_1'+j_0\big)\cdot y \ = \ (b_1', t)\cdot y,
\]
where $t:=b_2-a_2+a_1'+j_0=b_1'+j_0+m_1-m_2$, so $(b_1',t)\cdot y$ splits at $k_1$ and $t$ satisfies 
\[
b_1' \ \leq \ b_1'+j_0+m_1-m_2 \ = \ t \qquad \text{and} \qquad t \ = \ b_2-a_2+a_1'+j_0 \ < \ b_2-a_2+a_2' \ = \ b_2' \ = \ J^y_n+b_1'.
\]
It follows that $(t-b_1')\bar{e}_2\cdot y$ splits at $b_1'+k_1$, with $0\leq t-b_1'<J^y_n$ and $b_1'+k_1=K^y_n+k_1\geq K^y_n$, contradicting Lemma \ref{lem:NoSplitting}.
\end{proof}

\begin{proposition}\label{prop:NWhs}
$\EJT\res\Wqnw$ is hypersmooth.
\end{proposition}

\begin{proof}
We will prove $\EJT\res\Wqnw$ is hypersmooth by Borel reducing it to $E_{1,t}$. Let $\rho:\JT\to 2^\omega$ be a fixed Borel bijection, and define the Borel map $\phi:\Wqnw\to (2^\omega)^\omega$ by $\phi(x)(n) := \rho\big( g^n(x)\big)$. We claim that $\phi$ is a Borel reduction from $\EJT\res\Wqnw$ to $E_{1,t}$.

To prove this, let $x,y\in\Wqnw$ be arbitrary.  If $\phi(x)\mathrel{E_{1,t}}\phi(y)$, then there exist $m,n\in\omega$ such that $g^m(x)=g^n(y)$, but since $z\mathrel{\EJT}g(z)$ for all $z\in\Wqnw$, this easily implies that $x\mathrel{\EJT}y$. Conversely, assume $x\mathrel{\EJT}y$, and fix $\bar{a},\bar{b}\in\N^2$ such that $\bar{a}\cdot x=\bar{b}\cdot y$. Write $\bar{m}:=\bar{a}-\bar{b}\in\Z^2$, and without loss of generality suppose $m_1\geq m_2$. Using the fact that $\lim K^y_n=\infty=\lim J^y_n$, fix $n\in\Z_+$ such that
\[
K^y_n > \max\big\{|\bar{a}|,|\bar{b}|\big\} \qquad \text{and} \qquad J^y_n>m_1-m_2.
\]
We claim that there is $m\in\omega$ such that $g^m(x)=g^n(y)$, from which $\phi(x)\mathrel{E_{1,t}}\phi(y)$ easily follows.

To show this, let us write $\bar{b}'=(b_1',b_2'):=(K_n^y, J_n^y+K_n^y)$ and $\bar{a}'=(a_1',a_2'):=\bar{a}+\bar{b}'-\bar{b}$ as in Lemma \ref{lem:technical}, so $\bar{a}\prec\bar{a}'$, $\bar{b}\prec\bar{b}'$, and $\bar{a}'\cdot x=\bar{b}'\cdot y=g^n(y)$. (Recall that $\preceq$, $\precneq$, and $\prec$ were defined in Equation \ref{eq:prec}). 
Note that $a_2'-a_1'=J^y_n-(m_1-m_2)>0$. We now show by induction that for all $m\in\omega$, 
\begin{equation}\tag{$\ast$}
\text{if $(K^x_m, \, J^x_m+K^x_m)\precneq\bar{a}'$, then $0<a_1'-K^x_m<a_2'-J^x_m-K^x_m$.}
\end{equation}
The base case $m=0$ holds trivially since $K^x_0=J^x_0=0<a_1'< a_2'$. For the induction, fix $m\in\omega$ and assume ($\ast$) holds for $m$. Assume $(K^x_{m+1}, \, J^x_{m+1}+K^x_{m+1})\precneq\bar{a}'$, which implies also $(K^x_m, \, J^x_m+K^x_m)\precneq\bar{a}'$, so $0<a_1'-K^x_m<a_2'-J^x_m-K^x_m$ by inductive hypothesis. As in Lemma \ref{lem:technical}, write $\c:=(K^x_m, \, J^x_m+K^x_m)$, $x':=(\d{1}+\bar{c})\cdot x$ and $\bar{d}:=\d{1}+\bar{c}+(k^{x'},\, j^{x'}+k^{x'})$, so $\bar{d}=(K^x_{m+1}, \, J^x_{m+1}+K^x_{m+1})\preceq\bar{a}'$. Since we are assuming $\bar{d}\ne\bar{a}'$, it now follows from Lemma \ref{lem:technical} that $0<a_1'-d_1<a_2'-d_2$, completing the induction.

Now let $m\in\omega$ be largest such that $(K^x_m, \, J^x_m+K^x_m)\preceq\bar{a}'$. Such $m$ must exist because $(K_0^x, \, J_0^x+K_0^x)=\bar{0}\preceq\bar{a}'$ 
but $\lim K^x_m=\infty$. We claim that in fact $\big(K^x_m, \, J^x_m+K^x_m\big)=\bar{a}'$, so that $g^m(x)=g^n(y)$. To see this, suppose otherwise for contradiction. Again let $\c:=(K^x_m, \, J^x_m+K^x_m)$, so $c_1\leq c_2$ and $0<a_1'-c_1< a_2'-c_2$ by ($\ast$). But then $\bar{c}$ satisfies the hypotheses of Lemma \ref{lem:technical}, so, again writing $x':=(\d{1}+\bar{c})\cdot x$ and $\bar{d}:=\d{1}+\bar{c}+(k^{x'}, j^{x'}+k^{x'})$, we have $\bar{d}=(K^x_{m+1},\, J^x_{m+1}+K^x_{m+1})\preceq\bar{a}'$ by Lemma \ref{lem:technical}, contradicting maximality of $m$. Thus $g^m(x)=g^n(y)$ as desired, and $\phi(x)(m+k)=\phi(y)(n+k)$ for all $k$ to complete the proof.
\end{proof}

\begin{corollary}\label{prop:SENWhs}
$\EJT\res\Wse$ is hypersmooth.
\end{corollary}

\begin{proof}
Fix $y\in\Wqnw$ and $x\in\JT$ such that $x\mathrel{\EJT}y$, along with $\bar{a},\bar{b}\in\N^2$ such that $\bar{a}\cdot x=\bar{b}\cdot y$. Let $\bar{a}':=(a_1,a_2+b_1)$ and $\bar{b}':=(b_1,b_2+b_1)$. Then $\bar{a}'\cdot x=\bar{b}'\cdot y\in\Wqnw$ by Lemma \ref{lem:NWclosure}. This shows that for all $x\in\tWqnw$ there is $\bar{c}\in\N^2$ such that $\bar{c}\cdot x\in\Wqnw$. For each $x\in\tWqnw$, let $h(x):=\bar{c}\cdot x$ where $\bar{c}\in\N^2$ is least such that $\bar{c}\cdot x\in\Wqnw$. Then $h:\tWqnw\to\Wqnw$ is a Borel reduction from $\EJT\res\tWqnw$ to $\EJT\res\Wqnw$, which shows $\EJT\res\tWqnw$ is hypersmooth by Proposition \ref{prop:NWhs}. Furthermore, Lemma \ref{lem:SEtoNW} implies that $\delta\res\tWqse$ is a Borel reduction from $\EJT\res\tWqse$ to $\EJT\res\Wqnw$, so $\EJT\res\tWqse$ is also hypersmooth. Since $\Wse=\tWqnw\sqcup\tWqse$, we are done by Lemma \ref{lem:piecewiseHS}.
\end{proof}

\section{The graphs $\G_{\m}$ and $\G_{2,n}$}\label{sec:Gm}

In this section, we complete the proof of Theorem \ref{thm:main} by extending our results from Sections \ref{sec:ep} and \ref{sec:se} first to the graphs $\G_{\m}$ where $\bar{0}\ne\m\in\N^2\cup (\Z_+\times\Z_-)$, and then ultimately to the graphs $\G_{2,n}$ for $n\in\Z_+$. 
Once again we work over an aribitrary finite base space $J\in\Z_{\geq 2}$.  Recall from Equations \ref{eq:Gl} and \ref{eq:mq} that for each $\bar{0}\ne\m\in\N^2$, the directed graph $\vec{\G}_{\m}(J)$ is defined on $\JT$ by 
\[
x\mathrel{\vec{\G}_{\m}(J)}y \quad :\Longleftrightarrow \quad \m\cdot x=y, 
\]
while for each $\m\in\Z^+\times\Z^-$ the directed graph $\vec{\G}_{\m}(J)$ is defined on $\JT$ by
\[
x\mathrel{\vec{\G}_{\m}(J)}y \quad :\Longleftrightarrow \quad m_1\bar{e}_1\cdot x=-m_2\bar{e}_2\cdot y. 
\]
In both cases, the undirected graph $\G_{\m}(J)$ is the symmetrization of $\vec{\G}_{\m}(J)$.  As in previous sections, we omit $J$ from the notation when this is convenient.  Most of our effort in this section will be directed at proving the following:

\begin{theorem}\label{thm:m}
For each nonzero $\m\in\N^2\cup (\Z_+\times\Z_-)$, there is a partition of $\neg\textup{Per}_{\m}(\TT)$ into a pair of $\ETT$-invariant Borel sets $X_{\m}$ and $W_{\m}$ such that $\chi_B\big(\G_{\m}\res\Xm\big)<\infty$ and $\ETT\res\Wm$ is hypersmooth.
\end{theorem}

We have, of course, already proven this for $\m\in\{\bar{e}_1,\bar{e}_2,\p,\q\}$, and our proofs in this section will make use of those results. Roughly speaking, our argument will split into three cases: first we consider $\m\in\Z_+^2$, then $\m\in\Z_+\times\Z_-$, and finally the case where $\m$ is a multiple of $\bar{e}_1$ or $\bar{e}_2$. There are strong similarities between the three cases, and for each our strategy will be to refer back to one of our results from Sections \ref{sec:ep} or \ref{sec:se}, via an appropriate coding function.

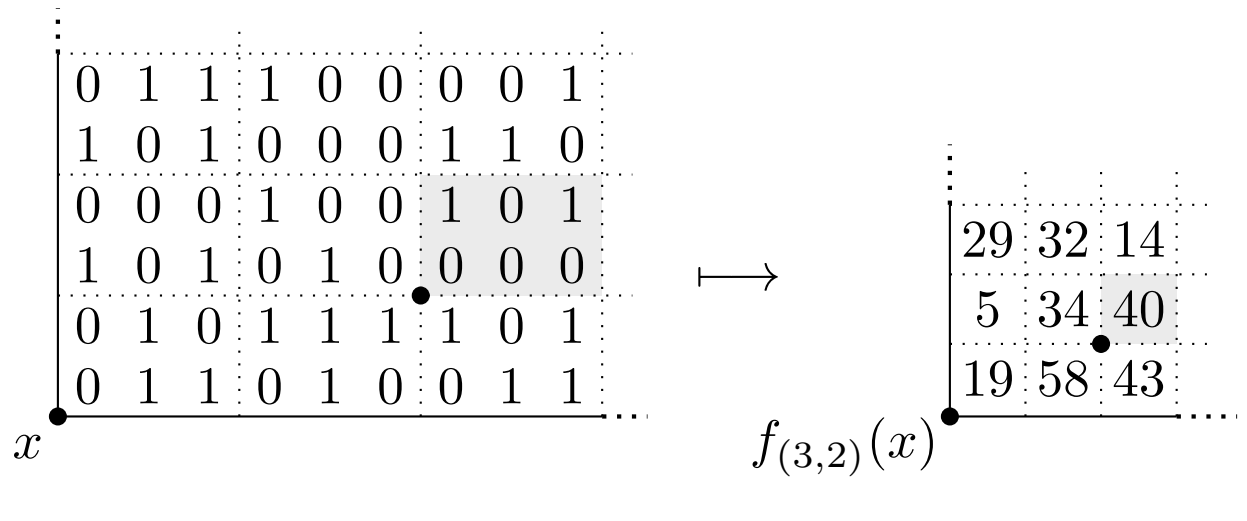
\begin{figure}[h]
\begin{center}
\begin{tikzpicture}
\fill [fill=gray!16] (2.5,.9) rectangle (3.7,1.7);
\draw (.1,.1) -- (3.7,.1);
\draw [dotted, thick] (3.7,.1) -- (4,.1);
\draw (.1,.1) -- (.1,2.5);
\draw [dotted, thick] (.1,2.5) -- (.1,2.8);
\draw [dotted] (.1,.9) -- (3.9,.9);
\draw [dotted] (.1,1.7) -- (3.9,1.7);
\draw [dotted] (.1,2.5) -- (3.9,2.5);
\draw [dotted] (1.3,.1) -- (1.3,2.7);
\draw [dotted] (2.5,.1) -- (2.5,2.7);
\draw [dotted] (3.7,.1) -- (3.7,2.7);
\node at (-.1,-.1) {$x$};
\node at (.1,.1) [circle,fill,inner sep=1.2pt]{};
\node at (2.5,.9) [circle,fill,inner sep=1.2pt]{};
\node at (.3,.3) {$0$};
\node at (.7,.3) {$1$};
\node at (1.1,.3) {$1$};
\node at (1.5,.3) {$0$};
\node at (1.9,.3) {$1$};
\node at (2.3,.3) {$0$};
\node at (2.7,.3) {$0$};
\node at (3.1,.3) {$1$};
\node at (3.5,.3) {$1$};
\node at (.3,.7) {$0$};
\node at (.7,.7) {$1$};
\node at (1.1,.7) {$0$};
\node at (1.5,.7) {$1$};
\node at (1.9,.7) {$1$};
\node at (2.3,.7) {$1$};
\node at (2.7,.7) {$1$};
\node at (3.1,.7) {$0$};
\node at (3.5,.7) {$1$};
\node at (.3,1.1) {$1$};
\node at (.7,1.1) {$0$};
\node at (1.1,1.1) {$1$};
\node at (1.5,1.1) {$0$};
\node at (1.9,1.1) {$1$};
\node at (2.3,1.1) {$0$};
\node at (2.7,1.1) {$0$};
\node at (3.1,1.1) {$0$};
\node at (3.5,1.1) {$0$};
\node at (.3,1.5) {$0$};
\node at (.7,1.5) {$0$};
\node at (1.1,1.5) {$0$};
\node at (1.5,1.5) {$1$};
\node at (1.9,1.5) {$0$};
\node at (2.3,1.5) {$0$};
\node at (2.7,1.5) {$1$};
\node at (3.1,1.5) {$0$};
\node at (3.5,1.5) {$1$};
\node at (.3,1.9) {$1$};
\node at (.7,1.9) {$0$};
\node at (1.1,1.9) {$1$};
\node at (1.5,1.9) {$0$};
\node at (1.9,1.9) {$0$};
\node at (2.3,1.9) {$0$};
\node at (2.7,1.9) {$1$};
\node at (3.1,1.9) {$1$};
\node at (3.5,1.9) {$0$};
\node at (.3,2.3) {$0$};
\node at (.7,2.3) {$1$};
\node at (1.1,2.3) {$1$};
\node at (1.5,2.3) {$1$};
\node at (1.9,2.3) {$0$};
\node at (2.3,2.3) {$0$};
\node at (2.7,2.3) {$0$};
\node at (3.1,2.3) {$0$};
\node at (3.5,2.3) {$1$};
\node at (4.6,1) {$\longmapsto$};
\fill [fill=gray!16] (7,.58) rectangle (7.5,1.04);
\draw (6,.1) -- (7.5,.1);
\draw [dotted, thick] (7.5,.1) -- (7.9,.1);
\draw (6,.1) -- (6,1.5);
\draw [dotted, thick] (6,1.5) -- (6,1.9);
\node at (5.3,-.1) {$f_{(3,2)}(x)$};
\node at (6,.1) [circle,fill,inner sep=1.2pt]{};
\node at (7,.58) [circle,fill,inner sep=1.2pt]{};
\node at (6.25,.35) {$19$};
\node at (6.75,.35) {$58$};
\node at (7.25,.35) {$43$};
\node at (6.25,.81) {$5$};
\node at (6.75,.81) {$34$};
\node at (7.25,.81) {$40$};
\node at (6.25,1.27) {$29$};
\node at (6.75,1.27) {$32$};
\node at (7.25,1.27) {$14$};
\draw [dotted] (6,.58) -- (7.7,.58);
\draw [dotted] (6,1.04) -- (7.7,1.04);
\draw [dotted] (6,1.5) -- (7.7,1.5);
\draw [dotted] (6.5,.1) -- (6.5,1.75);
\draw [dotted] (7,.1) -- (7,1.75);
\draw [dotted] (7.5,.1) -- (7.5,1.75);
\end{tikzpicture}
\end{center}
\caption{The coding function $f_{(3,2)}:\TT\to 64^{\omega^2}$.}
\label{fig:compression}
\end{figure}

Let us consider the general strategy for $\m\in\Z_+^2$, before turning to the details in all three cases. For the sake of illustration, fix some $\m\in\Z_+^2$, say $\m=(3,2)$, and consider the continuous bijective coding function $f_{(3,2)}:\TT\to 64^{\omega^2}$ defined by
\begin{equation}\label{eq:(3,2)-CodingFunction}
f_{(3,2)}(x)(i,j) \ := \ \beta\bigg( \Big((3i, 2j)\cdot x\Big)\res\R\big((3,2)\big)\bigg),
\end{equation}
where $\beta:2^{3\times 2}\to 64$ is a fixed bijection. We can understand this map by partitioning $\omega^2$ into $3\times 2$ rectangles that we collapse into single entries by coding their bits as natural numbers less than 64.  If $\beta$ converts from binary, reading row by row, this map can be pictured as in Figure \ref{fig:compression}, where, for instance, $f_{(3,2)}(x)(2,1)=40$.

We can define $f_{\m}$ analogously for any $\m\in\Z_+^2$, and our motivation for introducing these maps is that the following hold, with $J_{\m}:=J^{m_1m_2}$: 
\begin{enumerate}[(i)]
 \item for all $x,y\in\TT$, we have $x\mathrel{\Gm}y$ if and only if $f_{\m}(x)\mathrel{\Gp}f_{\m}(y)$;
 \item for all $x\in\TT$, $x$ is $\m$-periodic if and only if $f_{\m}(x)$ is $\p$-periodic;
 \item $c_{\p}\circ f_{\m}$ is a Borel $3$-coloring of $\Gm\res f_{\m}^{-1}\big[\Xp\big]$;
 \item $f_{\m}$ Borel reduces $\ETT\res f_{\m}^{-1}\big[\Wp\big]$ to a finite-index Borel super-equivalence relation of $E_2^{J}\res\Wp$ where $J=2^{m_1m_2}$, so $\ETT\res f_{\m}^{-1}\big[\Wp\big]$ itself is hypersmooth by Proposition \ref{prop:WpHS} and the results of Section \ref{sec:hs}.
\end{enumerate}

Of course, we will need to state all this precisely and prove it, and (iv) in particular will take a bit of work.  We will need to pay attention to the base space $J$, since the domain and codomain of $f_{\m}$ have different base spaces.  Once we finish the case $\m\in\Z_+^2$, which relies on our results from Section \ref{sec:ep} for $\Gp$, we will be able to make similar arguments for $\m\in\Z_+\times\Z_-$ that rely on our results from Section \ref{sec:se} for $\Gq$. Finally, we can use our results for $\Ge$ and $\G_{\bar{e}_2}$ to treat the symmetric cases $\m\in\Z_+\times\{0\}$ and $\m\in\{0\}\times\Z_+$. With Figure \ref{fig:compression} and this preliminary discussion as a guide, we now turn to the details. For clarity of presentation we split what follows into five subsections, with the first containing some uniform definitions and results, the next three treating our three cases in the proof of Theorem \ref{thm:m} in parallel, and a short final one treating the graphs $\G_{2,n}$ and completing the proof of Theorem \ref{thm:main}.

\subsection{The coding functions $f_{J,\m}$} 

For $\m\in\Z_+^2$, we can define $f_{J,\m}$ by generalizing Equation \eqref{eq:(3,2)-CodingFunction} above.  For the other cases, it turns out that when $\m\in\Z_+\times\Z_-$, we will want to use the same coding function that we use for $(m_1,-m_2)\in\Z_+^2$, and when $\m\in\Z_+\times\{0\}$, we want the same coding function that we use for $(m_1,1)\in\Z_+^2$. 
In order to avoid having to state everything three times, let us employ the following notational convention: for each $m\in\Z$ and $\m\in\Z^2$, let
\begin{equation}\label{eq:star2}
m^{\star} := \begin{cases} m & \text{if $m>0$} \\ -m & \text{if $m<0$} \\ 1 &\text{if $m=0$} \end{cases} \qquad \text{and} \qquad \m^{\star} := (m_1^\star,m_2^\star) \ \in \ \Z_+^2.
\end{equation}
Using this, we can make the following uniform definitions.

\begin{DefinitionNotation}
For each $J\in\Z_{\geq 2}$ and $\m\in\N^2\cup(\Z_+\times\Z_-)$, let
\[
J_{\m} \ := \ J^{m_1^{\star}m_2^{\star}} \ \in \ \Z_+,
\]
let
\[
\beta_{J,\m}:J^{\R(\m^{\star})}\to J_{\m}
\]
be a fixed bijection, and define the continuous bijection $f_{J,\m}:\JT\to J_{\m}^{\omega^2}$ by
\begin{equation}
f_{J,\m}(x)(i,j) \ := \ \beta_{J,\m}\bigg( \Big(\big(m_1^\star i, \, m_2^\star j\big)\cdot x\Big)\res\R(\m^{\star})\bigg).
\end{equation}
Let us agree that in the degenerate case $\m^{\star}=\p=(1,1)$, $\beta_{J,\m}$ is the natural bijection from $J^{\{\bar{0}\}}$ to $J_{\p}=J^1=J$, so that in this case $f_{J,\m}$ is the identity map on $\JT$.  Also define the Borel equivalence relation $E_{J,\m}$ on $J_{\m}^{\omega^2}$ by
\begin{equation}
x\mathrel{E_{J,\m}}y \quad :\Longleftrightarrow \quad f_{J,\m}^{-1}(x)\mathrel{\EJT}f_{J,\m}^{-1}(y),
\end{equation}
so $f_{J,\m}$ is a continuous isomorphism from $\EJT$ to $E_{J,\m}$. 
\end{DefinitionNotation}

We now prove two lemmas about the coding functions $f_{J,\m}$ and equivalence relations $E_{J,\m}$ before specializing to our three cases in Subsections \ref{sec:SubSec2}, \ref{sec:SubSec3}, and \ref{sec:SubSec4}. 
Note that $\m^\star=\m$ for all $\m\in\Z_+^2$; furthermore, since the star operation is idempotent, for all $J\in\Z_{\geq 2}$ and $\m\in\N^2\cup(\Z_+\times\Z_-)$ we have
\begin{equation}
J_{\m}=J_{\m^{\star}}, \qquad \beta_{J,\m}=\beta_{J,\m^{\star}}, \qquad f_{J,\m}=f_{J,\m^{\star}}, \qquad \text{and} \qquad E_{J,\m}=E_{J,\m^{\star}},
\end{equation}
so in the proofs of Lemmas \ref{lem:boring1} and \ref{claim:compressionFI} below, we could assume without loss of generality that $\m\in\Z_+^2$ (or, equivalently, ignore every occurrence of $\star$).

\begin{lemma}\label{lem:boring1}
For all $J\in\Z_{\geq 2}$, $\m\in\N^2\cup(\Z_+\times\Z_-)$, $\bar{a}=(a_1,a_2)\in\N^2$, and $x\in\JT$, we have
\[
\bar{a}\cdot f_{J,\m}(x) \ = \ f_{J,\m}\Big(\big(m_1^\star a_1, \, m_2^\star a_2\big)\cdot x\Big).
\]
\end{lemma}

\begin{proof}
Let $J$, $\m$, $\bar{a}$, and $x$ be given. Then for all $i,j\in\omega$ we have
\begin{align*}
\Big(\bar{a}\cdot f_{J,\m}(x)\Big)(i,j) \ &= \ \Big(f_{J,\m}(x)\Big)(a_1+i,\, a_2+j) \\ 
&= \ \beta_{J,\m}\bigg( \Big( \big( m_1^{\star}(a_1+i), \, m_2^{\star}(a_2+j) \big)\cdot x\Big) \res \R(\m^{\star})\bigg) \\
&= \ \beta_{J,\m}\bigg(\Big(\big(m_1^{\star}i, \, m_2^{\star}j\big)\cdot \big( (m_1^{\star}a_1, \, m_2^{\star}a_2) \cdot x\big)\Big)\res\R(\m^{\star})\bigg) \\ 
&= \ f_{J,\m}\Big(\big(m_1^{\star}a_1, \, m_2^{\star}a_2\big)\cdot x\Big)(i,j). \qedhere
\end{align*}
\end{proof}

\begin{lemma}\label{claim:compressionFI}
For all $J\in\Z_{\geq 2}$ and $\m\in\N^2\times(\Z_+\times\Z_-)$, we have $\EJmtwo\subseteq E_{J,\m}$; furthermore, every $E_{J,\m}$-equivalence class contains at most $|m_1^{\star}m_2^{\star}|$ $\EJmtwo$-equivalence classes.
\end{lemma}

\begin{proof}
Let $J$ and $\m$ be given, and fix $x,y\in J_{\m}^{\omega^2}$ such that $x\mathrel{\EJmtwo}y$. Fix $\bar{a},\bar{b}\in\N^2$ such that $\bar{a}\cdot x=\bar{b}\cdot y$, and write $x_0:=f_{J,\bar{m}}^{-1}(x)$ and $y_0:=f_{J,\bar{m}}^{-1}(y)$. Let $\bar{a}':=(m_1^{\star}a_1, \, m_2^{\star}a_2)$ and $\bar{b}':=(m_1^{\star}b_1,\,m_2^{\star}b_2)$.  Then by Lemma \ref{lem:boring1} we have 
\[
f_{J,\m}(\bar{a}'\cdot x_0) \ = \ \bar{a}\cdot f_{J,\m}(x_0) \ = \ \bar{b}\cdot f_{J,\m}(y_0) \ = \ f_{J,\m}(\bar{b}'\cdot y_0).
\]
Thus $\bar{a}'\cdot x_0=\bar{b}'\cdot y_0$, so $x_0\mathrel{\EJT}y_0$, which means $x\mathrel{E_{J,\m}}y$ as desired.

For the second claim, we show that for all $x,y\in J_{\m}^{\omega^2}$ such that $x\mathrel{E_{J,\m}}y$, there is $\bar{c}\in\R(\m^{\star})$ such that $f_{J,\m}\big(\bar{c}\cdot f_{J,\m}^{-1}(x)\big)\mathrel{\EJmtwo}y$. This will suffice to complete the proof, since it shows the $E_{J,\m}$-class of $x$ is covered by the $|m_1^{\star}m_2^{\star}|$ $\EJmtwo$-classes containing the points $f_{J,\m}\big(\bar{c}\cdot f^{-1}_{J,\m}(x)\big)$ as $\bar{c}$ ranges over $\R(\m^{\star})$.  To prove the claim, let $x,y\in J_{\m}^{\omega^2}$ satisfy $x\mathrel{E_{J,\m}}y$, again write $x_0:=f_{J,\m}^{-1}(x)$ and $y_0:=f_{J,\m}^{-1}(y)$, and this time fix $\bar{a},\bar{b}\in\N^2$ such that $\bar{a}\cdot x_0=\bar{b}\cdot y_0$.  Let $c_1:=(a_1-b_1)\mod m_1^{\star}$ and $c_2:=(a_2-b_2)\mod m_2^{\star}$, and let $r_1:=-a_1\mod m_1^{\star}$ and $r_2:=-a_2\mod m_2^{\star}$.  Then $m_1^{\star}$ divides $b_1+c_1+r_1$ and $m_2^{\star}$ divides $b_2+c_2+r_2$. Fix $\bar{u}=(u_1,u_2)$ and $\bar{v}=(v_1,v_2)$ in $\N^2$ such that $\bar{a}+\bar{r}=(m_1^{\star}u_1,\, m_2^{\star}u_2)$ and $\bar{b}+\bar{c}+\bar{r}=(m_1^{\star}v_1,\, m_2^{\star}v_2)$.  Then
\[
(m_1^{\star}u_1, \, m_2^{\star}u_2)\cdot (\bar{c}\cdot x_0) \ = \ (\bar{a}+\bar{r})\cdot (\bar{c}\cdot x_0) \ = \ (\bar{b}+\bar{c}+\bar{r})\cdot y_0 \ = \ (m_1^{\star}v_1,\, m_2^{\star}v_2)\cdot y_0,
\]
so by Lemma \ref{lem:boring1} we have
\[
\bar{u}\cdot f_{J,\m}\big(\bar{c}\cdot x_0\big) \ = \ f_{J,\m}\Big((m_1^{\star}u_1, \, m_2^{\star}u_2)\cdot (\bar{c}\cdot x_0)\Big) \ = \ f_{J,\m}\Big((m_1^{\star}v_1, \, m_2^{\star}v_2)\cdot y_0\Big) \ = \ \bar{v}\cdot f_{J,\m}(y_0),
\]
showing $\c\in\R(\m^{\star})$ satisfies $f_{J,\m}\big(\bar{c}\cdot f_{J,\m}^{-1}(x)\big)\mathrel{\EJmtwo}y$ to complete the proof.
\end{proof}

\subsection{The graphs $\Gm$ for $\m\in\Z_+^2$}\label{sec:SubSec2}

In this subsection, we prove Theorem \ref{thm:m} in the case $\m\in\Z_+^2$. 
We start by establishing (i) and (ii) from the discussion around Figure \ref{fig:compression}. 

\begin{lemma}\label{claim:fJmHom+}
For all $J\in\Z_{\geq 2}$, $\m\in\Z_+^2$, and $x,y\in\JT$, we have 
\begin{align*}
x\mathrel{\vec{\G}_{\m}(J)}y \quad &\Longleftrightarrow \quad f_{J,\m}(x)\mathrel{\vec{\G}_{\p}(J_{\m})}f_{J,\m}(y) \\
\shortintertext{and}
x\mathrel{\G_{\m}(J)}y \quad &\Longleftrightarrow \quad f_{J,\m}(x)\mathrel{\G_{\p}(J_{\m})}f_{J,\m}(y).
\end{align*}
\end{lemma}

\begin{proof}
Let $J$, $\m$, and $x,y$ be given. Since $\p\cdot f_{J,\m}(x) = f_{J,\m}\big(\m\cdot x\big)$ by Lemma \ref{lem:boring1}, we have, using bijectivity of $f_{J,\m}$,
\[
\m\cdot x=y \quad 
\Longleftrightarrow \quad f_{J,\m}(\m\cdot x) = f_{J,\m}(y) \quad \Longleftrightarrow \quad \p\cdot f_{J,\m}(x)=f_{J,\m}(y). 
\]
This proves the first claim, and the second follows easily from the first.
\end{proof}

\begin{lemma}\label{claim:fJmPeriodicity+}
For all $J\in\Z_{\geq 2}$, $\m\in\Z_+^2$, and $x\in\JT$, $x$ is $\m$-periodic if and only if $f_{J,\m}(x)$ is $\p$-periodic.
\end{lemma}

\begin{proof}
Let $J$, $\m$, and $x$ be given. Using Lemmas \ref{lem:boring1} and \ref{lem:per}, we have 
\begin{align*}
\text{$x$ is $\m$-periodic} \quad &\Longleftrightarrow \quad (\exists\bar{c}\in\N^2) \ \bar{c}\cdot x = (\bar{c}+\m)\cdot x \\
&\Longleftrightarrow \quad (\exists\bar{a}\in\N^2) \ (a_1m_1,\, a_2m_2)\cdot x = \big((a_1+1)m_1, \, (a_2+1)m_2\big)\cdot x \\
&\Longleftrightarrow \quad (\exists\bar{a}\in\N^2) \ f_{J,\m}\Big((a_1m_1,\, a_2m_2)\cdot x\Big) \ = \ f_{J,\m}\Big( \big((a_1+1)m_1, \, (a_2+1)m_2\big)\cdot x\Big) \\
&\Longleftrightarrow \quad (\exists\bar{a}\in\N^2) \ (a_1,a_2)\cdot f_{J,\m}(x) = (a_1+1, \, a_2+1)\cdot f_{J,\m}(x) \\
&\Longleftrightarrow \quad \text{$f_{J,\m}(x)$ is $\p$-periodic.} \qedhere
\end{align*}
\end{proof}

Now we generalize to $\m\in\Z_+^2$ the definitions of $\Xp$ and $\Wp$ that were given in Definition \ref{df:XpWp}.

\begin{definition}\label{df:WmXm+}
Given $J\in\Z_{\geq 2}$ and $\m=(m_1,m_2)\in\Z_+^2$, let
\[
\Xm(J) \ := \ \Big\{ x\in\neg\textup{Per}_{\m}(\JT)\st (\exists\bar{a}\in\N^2) \text{  $\langle (\bar{a}\cdot x)(im_1,\, im_2)\st i\in\omega\rangle$ is not eventually constant}\Big\},
\]
and let
\begin{align*}
\Wm(J) \ :\!&= \ \neg\textup{Per}_{\m}(\JT)\setminus\Xm(J) \\
&= \ \Big\{ x\in\neg\textup{Per}_{\m}(\JT)\st (\forall\bar{a}\in\N^2) \text{ $\big\langle (\bar{a}\cdot x)(im_1,\, im_2)\st i\in\omega\big\rangle$ is eventually constant}\Big\}.
\end{align*}
\end{definition}

Note that Definition \ref{df:WmXm+} agrees with Definiton \ref{df:XpWp} when $\m=\p$, and also that the sets $\Xm(J)$ and $\Wm(J)$ are Borel and $\EJT$-invariant.

\begin{lemma}\label{claim:annoying+}
For all $J\in\Z_{\geq 2}$ and $\m\in\Z_+^2$, we have
\[
\Xm(J)=\big\{x\in\JT\st f_{J,\m}(x)\in\Xp(J_{\m})\big\} \qquad \text{and} \qquad \Wm(J)=\big\{x\in\JT\st f_{J,\m}(x)\in\Wp(J_{\m})\big\}.
\]
\end{lemma}

\begin{proof} Let $J$ and $\m$ be given. By definition we have
\[
\neg\textup{Per}_{\m}(\JT) = \Xm(J)\sqcup\Wm(J) \qquad \text{and} \qquad \neg\textup{Per}_{\p}(J_{\m}^{\omega^2}) = \Xp(J_{\m})\sqcup\Wp(J_{\m}).
\]
Therefore, it will suffice by Lemma \ref{claim:fJmPeriodicity+} to show $f_{J,\m}[\Xm(J)]\subseteq\Xp(J_{\m})$ and $f_{J,\m}[\Wm(J)]\subseteq\Wp(J_{\m})$.

To show these, first fix $x\in\Xm(J)$, along with $\bar{a}\in\N^2$ such that the sequence $\langle (\bar{a}\cdot x)(im_1,im_2)\st i\in\omega\rangle$ is not eventually constant. Write $a_1=b_1m_1+r_1$ and $a_2=b_2m_2+r_2$ where $b_1,b_2\in\N$, $r_1\in m_1$, and $r_2\in m_2$. 
Since for each $i\in\omega$ the rectangle $\big(m_1(b_1+i), m_2(b_2+i)\big)+\R(\m)$ contains $(a_1+im_1, a_2+im_2)$, the sequence $\big\langle \big(\bar{b}\cdot f_{J,\m}(x)\big)(i,i)\st i\in\omega\big\rangle$ is not eventually constant, so $f_{J,\m}(x)\in\Xp(J_{\m})$.  For the second inclusion, fix $y\in\Wm(J)$, so for all $\bar{a}\in\N^2$ the sequence $\big\langle (\bar{a}\cdot y)(im_1, im_2)\st i\in\omega\big\rangle$ is eventually constant. Let $\bar{b}\in\N^2$ be arbitrary, and consider the rectangle $R:=(b_1m_1, b_2m_2)+\R(\m)$.  For every $\bar{a}\in R$, there is $N_{\bar{a}}\in\N$ such that the sequence $\big\langle (\bar{a}\cdot y)(im_1, im_2)\st i\geq N_{\bar{a}}\big\rangle$ is constant. Let $N=\max\{N_{\bar{a}}\st\bar{a}\in R\}$. Then the sequence $\big\langle \big(\bar{b}\cdot f_{J,\m}(y)\big)(i,i)\st i\geq N\big\rangle$ is constant, so $f_{J,\m}(y)\in\Wp(J_{\m})$. We have shown $f_{J,\m}[\Xm(J)]\subseteq\Xp(J_{\m})$ and $f_{J,\m}[\Wm(J)]\subseteq\Wp(J_{\m})$, completing the proof.
\end{proof}

\begin{corollary}\label{cor:EJmInvar+}
For all $J\in\Z_{\geq 2}$ and $\m\in\Z_+^2$, the set $\Wp(J_{\m})$ is $E_{J,\m}$-invariant.
\end{corollary}

\begin{proof}
Let $J$ and $\m$ be given, fix arbitrary $x,y\in J_{\m}^{\omega^2}$, and write $x_0:=f_{J,\m}^{-1}(x)$ and $y_0:=f_{J,\m}^{-1}(y)$. Suppose $x\in\Wp(J_{\m})$ and $x\mathrel{E_{J,\m}}y$. Then $x_0\in\Wm(J)$ by Lemma \ref{claim:annoying+}, and $x_0\mathrel{\EJT}y_0$. Since $\Wm(J)$ is $\EJT$-invariant, we also have $y_0\in\Wm(J)$, so $y\in\Wp(J_{\m})$ again by \ref{claim:annoying+}.
\end{proof}

\begin{proposition}\label{prop:XmColorable+}
For all $\bar{m}\in\Z_+^2$, we have $\chi_B\big(\G_{\m}(2)\res\Xm(2)\big)\leq 3$.
\end{proposition}

\begin{proof}
Fix $\bar{m}\in\Z_+^2$. Using Proposition \ref{prop:GpColoring}, let $c_{\p}:\Xp(2_{\m})\to 3$ be a Borel $3$-coloring of $\Gp(2_{\m})\res\Xp(2_{\m})$, and define the Borel function $c_{\m}:\Xm(2)\to 3$ by 
\[
c_{\m}(x) \ := \ c_{\p}\big(f_{2,\m}(x)\big).
\]
Let $x,y\in\Xm(2)$ be arbitrary such that $x\mathrel{\G_{\m}(2)}y$. Then $f_{2,\m}(x)\in\Xp(2_{\m})$ and $f_{2,\m}(y)\in\Xp(2_{\m})$ by Lemma \ref{claim:annoying+}, and $f_{2,\m}(x)\mathrel{\Gp(2_{\m})}f_{2,\m}(y)$ by Lemma \ref{claim:fJmHom+}. This means $c_{\p}\big(f_{2,\m}(x)\big)\ne c_{\p}\big(f_{2,\m}(y)\big)$, so $c_{\m}(x)\ne c_{\m}(y)$. We have shown that $c_{\m}$ is a Borel $3$-coloring of $\G_{\m}(2)\res\Xm(2)$.
\end{proof}

\begin{proposition}\label{prop:WmHS+}
For all $\bar{m}\in\Z_+^2$, $\ETT\res\Wm(2)$ is hyperselective.
\end{proposition}

\begin{proof}
Fix $\bar{m}\in\Z_+^2$, and let $J:=2_{\m}=2^{m_1m_2}\in\Z_+$. By Corollary \ref{cor:EJmInvar+}, $\Wp(J)$ is $E_{2,\m}$-invariant.  By Proposition \ref{prop:WpHS} and Observation \ref{ob:UniformlyContinuous}, $\EJT\res\Wp(J)$ is hypersmooth, uniformly continuous, and has $K_{\sigma}$ classes. Therefore, since by Lemma \ref{claim:compressionFI} $E_{2,\m}$ contains $\EJT$ as a finite-index sub-equivalence relation on $\Wp(J)$, the equivalence relation $E_{2,\m}\res\Wp(J)$ is hypersmooth by Corollary \ref{cor:fi1}. But $f_{2,\m}\res\Wm(2)$ is a Borel reduction from $\ETT\res\Wm(2)$ to $E_{2,\m}\res\Wp(J)$, so $\ETT\res\Wm(2)$ is hypersmooth as well, and is therefore hyperselective by Observations \ref{ob:UniformlyContinuous} and \ref{ob:UniformlyContinuousAndKSigma}.
\end{proof}

\subsection{The graphs $\Gm$ for $\m\in\Z_+\times\Z_-$}\label{sec:SubSec3}

We proceed as in the previous subsection, starting with analogues of Lemmas \ref{claim:fJmHom+} and \ref{claim:fJmPeriodicity+}. 

\begin{lemma}\label{claim:fJmHom-}
For all $J\in\Z_{\geq 2}$, $\m\in\Z_+\times\Z_-$, and $x,y\in\JT$, we have 
\begin{align*}
x\mathrel{\vec{\G}_{\m}(J)}y \quad &\Longleftrightarrow \quad f_{J,\m}(x)\mathrel{\vec{\G}_{\q}(J_{\m})}f_{J,\m}(y) \\
\shortintertext{and}
x\mathrel{\G_{\m}(J)}y \quad &\Longleftrightarrow \quad f_{J,\m}(x)\mathrel{\G_{\q}(J_{\m})}f_{J,\m}(y).
\end{align*}
\end{lemma}

%

\begin{proof}
Let $J$, $\m$, and $x,y$ be given. Since $\bar{e}_1\cdot f_{J,\m}(x)=f_{J,\m}(m_1\bar{e}_1\cdot x)$ and $\bar{e}_2\cdot f_{J,\m}(y)=f_{J,\m}(-m_2\bar{e}_2\cdot y)$ by Lemma \ref{lem:boring1}, we have, using bijectivity of $f_{J,\m}$,
\begin{align*}
x\mathrel{\vec{\G}_{\m}(J)}y \quad &\Longleftrightarrow \quad m_1\bar{e}_1\cdot x \ = \ -m_2\bar{e}_2\cdot y \\
&\Longleftrightarrow \quad f_{J,\m}(m_1\bar{e}_1\cdot x) \ = \ f_{J,\m}(-m_2\bar{e}_2\cdot y) \\
&\Longleftrightarrow \quad \bar{e}_1\cdot f_{J,\m}(x) \ = \ \bar{e}_2\cdot f_{J,\m}(y) \quad \Longleftrightarrow \quad
f_{J,\m}(x)\mathrel{\vec{\G}_{\q}(J_{\m})}f_{J,\m}(y).
\end{align*}
This proves the first claim, and the second follows easily from the first.
\end{proof}

\begin{lemma}\label{claim:fJmPeriodicity-}
For all $J\in\Z_{\geq 2}$, $\m\in\Z_+\times\Z_-$, and $x\in\JT$, $x$ is $\m$-periodic if and only if $f_{J,\m}(x)$ is $\q$-periodic.
\end{lemma}

\begin{proof}
Let $J$ and $\m$ be given. Note first that for any $x\in\JT$,
\begin{align*}
\text{$x$ is strongly $\m$-periodic} \quad &\Longleftrightarrow \quad m_1\bar{e}_1\cdot x = -m_2\bar{e}_2\cdot x \\
&\Longleftrightarrow \quad f_{J,\m}(m_1\bar{e}_1\cdot x) = f_{J,\m}(-m_2\bar{e}_2\cdot x) 
\\
&\Longleftrightarrow \quad \bar{e}_1\cdot f_{J,\m}(x) = \bar{e}_2\cdot f_{J,\m}(x) \qquad &\text{(by Lemma \ref{lem:boring1})} \\
&\Longleftrightarrow \quad \text{$f_{J,\m}(x)$ is strongly $\q$-periodic.}
\end{align*}
Next, note by Lemma \ref{lem:per} that for each $x\in\JT$,
\begin{align*}
\text{$x$ is $\m$-periodic} \quad &\Longleftrightarrow \quad \text{$\exists\bar{c}\in\N^2$ such that $\bar{c}\cdot x$ is strongly $\m$-periodic} \\
&\Longleftrightarrow \quad \text{$\exists\bar{a}\in\N^2$ such that $\big(m_1a_1, \, -m_2a_2\big)\cdot x$ is strongly $\m$-periodic.}
\end{align*}
Combining these and using Lemma \ref{lem:boring1} again, we see that for all $x\in\JT$,
\begin{align*}
\text{$x$ is $\m$-periodic} \quad &\Longleftrightarrow \quad \text{$\exists\bar{a}\in\N^2$ such that $\big(m_1a_1, \, -m_2a_2\big)\cdot x$ is strongly $\m$-periodic} \\
&\Longleftrightarrow \quad \text{$\exists\bar{a}\in\N^2$ such that $f_{J,\m}\Big((m_1a_1, \, -m_2a_2)\cdot x\Big)$ is strongly $\q$-periodic} \\
&\Longleftrightarrow \quad \text{$\exists\bar{a}\in\N^2$ such that $\bar{a}\cdot f_{J,\m}(x)$ is strongly $\q$-periodic} \\
&\Longleftrightarrow \quad \text{$f_{J,\m}(x)$ is $\q$-periodic.} \qedhere
\end{align*}
\end{proof}

Now we generalize to $\m\in\Z_+\times\Z_-$ the definitions of $\Wse$ and $\Xse$ that were given in Definition \ref{df:XqWq}.  We use here the notion of a rectangular patch being \emph{$\m$-symmetric} that was introduced in Section \ref{sec:periodicity}, in the discussion between Lemmas \ref{lem:HBorel} and \ref{lem:delta}.

\begin{definition}\label{df:WmXm}
Given $J\in\Z_{\geq 2}$ and $\m\in\Z_+\times\Z_-$, let
\[
W_{\m}(J) \ := \ \big\{ x\in\neg\textup{Per}_{\m}(\JT)\st (\forall M)(\exists N)(\forall n\geq N) \; (n\m^{\star}\cdot x)\res\R(\d{M}) \text{ is $\m$-symmetric}\big\}
\]
and $X_{\m}(J):=\neg\textup{Per}_{\m}(\JT)\setminus W_{\m}(J)$.
\end{definition}

That Definitions \ref{df:WmXm} and \ref{df:XqWq} agree when $\m=\q$ is not immediately obvious, so we will verify it below.  In the meantime, to avoid confusion, we take $\Xm$ and $\Wm$ to be defined as in Definition \ref{df:WmXm} even when $\m=\q$, and temporarily write 
\[
\hat{X}_{\q} \ := \ \left\{ x\in\npseJ\st\liminf_{n\to\infty} r(n\cdot x)<\infty\right\} \quad \text{and} \quad \hat{W}_{\q} \ := \ \npseJ\setminus\hat{X}_{\q},
\]
for the sets in Definition \ref{df:XqWq}.  But first we check that $X_{\m}$ and $W_{\m}$ are Borel and $\EJT$-invariant.

\begin{lemma}\label{claim:WmInvar}
For all $J\in\Z_{\geq 2}$ and $\m\in\Z_+\times\Z_-$, the sets $\Xm(J)$ and $\Wm(J)$ are Borel and $\EJT$-invariant.
\end{lemma}

\begin{proof}
Let $J$ and $\m$ be given.  It is clear from the definition that $W_{\m}$, hence also $X_{\m}$, is Borel.  For invariance of $\Wm$, fix $x\in\Wm$, $y\in\neg\textup{Per}_{\m}(\JT)$, and $\bar{a},\bar{b}\in\N^2$ such that $\bar{a}\cdot x=\bar{b}\cdot y$. Write $k:=\max\{|\bar{a}|,|\bar{b}|\}$, and let $M\in\omega$ be arbitrary. Put $M':=M+2k|\m|$. Fix $N'\geq k$ such that $(n\m^{\star}\cdot x)\res\R(\d{M}')$ is $\m$-symmetric for all $n\geq N'$. Let $N:=N'+k$, and suppose $n\geq N$. Then
\begin{align*}
n\m^{\star}\cdot y \ &= \ \big(m_1n-b_1+a_1, \, -m_2n-b_2+a_2)\cdot x \\ &= \ \big(m_1k-b_1+a_1, \, -m_2k-b_2+a_2\big)\cdot \big((n-k)\m^{\star}\cdot x\big).
\end{align*}
Since $\big((n-k)\m^{\star}\cdot x\big)\res\R(\d{M}')$ is $\m$-symmetric and $k\m^{\star}-\bar{b}+\bar{a}+\d{M}\preceq\d{M}'$, it follows that $(n\m^{\star}\cdot y)\res\R(\d{M})$ is $\m$-symmetric. Since $M$ was arbitrary, we conclude that $y\in\Wm$. This shows $\Wm$ is $\EJT$-invariant, from which invariance of $\Xm$ easily follows.
\end{proof}

The following analogue of Lemma \ref{claim:annoying+} is a bit harder to show when $\m\in\Z_+\times\Z_-$.

\begin{lemma}\label{claim:annoying-}
For all $J\in\Z_{\geq 2}$ and $\m\in\Z_+\times\Z_-$, we have 
\[
W_{\m}(J)=\big\{x\in\JT\st f_{J,\m}(x)\in\hat{W}_{\q}(J_{\m})\big\} \qquad \text{and} \qquad X_{\m}(J)=\big\{x\in\JT\st f_{J,\m}(x)\in\hat{X}_{\q}(J_{\m})\big\}.
\]
\end{lemma}

\begin{proof}
Let $J$ and $\m$ be given, and suppose first that $x\in\W_{\m}(J)$. Then $x$ is not $\m$-periodic, so $f_{J,\m}(x)\in\npseJm$ by Lemma \ref{claim:fJmPeriodicity-}. With $r$ as in Definition \ref{df:XqWq}, we show $\displaystyle{\lim_{n\to\infty}r\big(\d{n}\cdot f_{J,\m}(x)\big)=\infty}$. Let $M\in\omega$, write $M':=M|\m|$, and fix $N\in\omega$ such that $(n\m^{\star}\cdot x)\res\R(\d{M}')$ is $\m$-symmetric for all $n\geq N$. Let $n\geq N$, $k<M$, and $0\leq i<k-1$ be arbitrary. Then for all $\bar{a}\in\R(\m^{\star})$, we have
\begin{align*}
x\Big(m_1(n+i)+a_1, \; -m_2(n+k-i-1)+a_2\Big) \ &= \ (n\m^\star\cdot x)\Big(im_1+a_1, \, (k-i-1)m_2^{\star}+a_2\Big) \\
= \ (n\m^\star\cdot x)\Big((i+1)m_1+a_1, \, (k-i-2)m_2^{\star}+a_2\Big) \ &= \ x\Big(m_1(n+i+1)+a_1, \; -m_2(n+k-i-2)+a_2\Big)
\end{align*}
since $(n\m^{\star}\cdot x)\res\R(\d{M}')$ is $\m$-symmetric and
\[
|(i+1)m_1+a_1| \ \leq \ km_1 \ < \ M' \qquad \text{and} \qquad |(k-i-1)m_2^{\star}+a_2| \ \leq \ km_2^{\star} \ < \ M'.
\]
Therefore, 
\begin{align*}
\big(\d{n}\cdot f_{J,\m}(x)\big)(i, \, k-i-1) \ &= \ f_{J,\m}(x)(n+i, \, n+k-i-1) \\
&= \ \beta_{J,\m}\bigg( \Big(\big(m_1(n+i),\, -m_2(n+k-i-1)\big)\cdot x\Big)\res\R(\m^{\star})\bigg) \\
&= \ \beta_{J,\m}\bigg( \Big(\big(m_1(n+i+1),\, -m_2(n+k-i-2)\big)\cdot x\Big)\res\R(\m^{\star})\bigg) \\
&= \ f_{J,\m}(x)(n+i+1, \, n+k-i-2) \ = \ \big(\d{n}\cdot f_{J,\m}(x)\big)(i+1, \, k-i-2).
\end{align*}
Since $0\leq i<k-1$ was arbitrary, this shows $s^{y_n}_k$ is constant where $y_n:=\d{n}\cdot f_{J,\m}(x)$. Since $k<M$ was arbitrary, this shows $r\big(\d{n}\cdot f_{J,\m}(x)\big)\geq M$. But $n\geq N$ was arbitrary as well, so $r\big(\d{n}\cdot f_{J,\m}(x)\big)\geq M$ for all $n\geq N$, and since $M$ was arbitary, we conclude $\displaystyle{\lim_{n\to\infty}r\big(\d{n}\cdot f_{J,\m}(x)\big)=\infty}$ as desired. 

Conversely, fix $x\in\JT$ and suppose $f_{J,\m}(x)\in\hat{W}_{\q}(J_{\m})$, so in particular $f_{J,\m}(x)\in\neg\textup{Per}_{\q}(J_{\m}^{\omega^2})$ and thus $x\in\neg\textup{Per}_{\m}(J^{\omega^2})$ by Lemma \ref{claim:fJmPeriodicity-}. Let $M\in\omega$, fix $N\in\omega$ such that $r\big(\d{n}\cdot f_{J,\m}(x)\big)>2M$ for all $n\geq N$, and suppose $n\geq N$. We show $(n\m^{\star}\cdot x)\res\R(\d{M})$ is $\m$-symmetric. Let $\bar{a},\bar{b}\in\R(\d{M})$ be arbitrary with $\bar{b}-\bar{a}=\m$. Let $r_1:=a_1\mod m_1$ and $r_2:=a_2\mod m_2^{\star}$, and write $t_1:=(a_1-r_1)/m_1$ and $t_2:=(a_2-r_2)/m_2^{\star}$, so $(r_1,r_2)\in\R(\m^{\star})$ and $t_1,t_2\in\omega$ with $t_1+t_2\leq a_1+a_2<2M$. Note that $t_2\geq 1$, since $(1-t_2)m_2+r_2=b_2\in\omega$. Since $r\big(\d{n}\cdot f_{J,\m}(x)\big)>2M$, we have 
\begin{align*}
&\beta_{J,\m}\bigg( \Big(\big(m_1(n+t_1),\, -m_2(n+t_2)\big)\cdot x\Big)\res\R(\m^{\star})\bigg) \ = \ \big(\d{n}\cdot f_{J,\m}(x)\big)(t_1,t_2) \\
= \ &\big(\d{n}\cdot f_{J,\m}(x)\big)(t_1+1, t_2-1) \ = \ \beta_{J,\m}\bigg( \Big(\big(m_1(n+t_1
+1),\, -m_2(n+t_2-1)\big)\cdot x\Big)\res\R(\m^{\star})\bigg).
\end{align*}
Therefore,
\begin{align*}
(n\m^{\star}\cdot x)(\bar{a}) \ &= \ x\Big(m_1(n+t_1)+r_1, \, -m_2(n+t_2)+r_2\Big) \\
&= \ x\Big(m_1(n+t_1+1)+r_1, \, -m_2(n+t_2-1)+r_2\Big) \ = \ (n\m^{\star}\cdot x)(\bar{b}).
\end{align*}
This shows $(n\m^{\star}\cdot x)\res\R(\d{M})$ is indeed $\m$-symmetric, completing the proof that 
\begin{equation}\label{eq:half}
W_{\m}(J) \ = \ \big\{x\in\JT\st f_{J,\m}(x)\in\hat{W}_{\q}(J_{\m})\big\}.
\end{equation}
It now follows from Lemma \ref{claim:fJmPeriodicity-}, Equation \eqref{eq:half}, and the decompositions
\[
\neg\textup{Per}_{\m}(\JT) = X_{\m}(J)\sqcup W_{\m}(J) \qquad \text{and} \qquad \neg\textup{Per}_{\q}(J_{\m}^{\omega^2}) = \hat{X}_{\q}(J_{\m})\sqcup\hat{W}_{\q}(J_{\m})
\]
that $X_{\m}(J)=\big\{x\in\JT\st f_{J,\m}(x)\in\hat{X}_{\q}(J_{\m})\big\}$ as well.
\end{proof}

\begin{corollary}\label{cor:WSEinvar}
For all $J\in\Z_{\geq 2}$, $X_{\q}(J)=\hat{X}_{\q}(J)$ and $W_{\q}(J)=\hat{W}_{\q}(J)$; that is, Definitions \ref{df:XqWq} and \ref{df:WmXm} coincide when $\m=\q$. Consequently, the sets $\hat{X}_{\q}(J)$ and $\hat{W}_{\q}(J)$ are Borel and $\EJT$-invariant.
\end{corollary}

\begin{proof}
That the two definitions coincide when $\m=\q$ follows immediately from Lemma \ref{claim:annoying-}, since $f_{J,\q}$ is the identity map. Thus $\hat{X}_{\q}(J)$ and $\hat{W}_{\q}(J)$ are Borel and $\EJT$-invariant by Lemma \ref{claim:WmInvar}.
\end{proof}

In fact, $\Wse(J_{\m})$ is not just $E_2^{J_{\m}}$-invariant, but actually $E_{J,\m}$-invariant, as we now verify in the analogue for $\m\in\Z_+\times\Z_-$ of Corollary \ref{cor:EJmInvar+}.

\begin{corollary}\label{cor:EJmInvar}
For all $J\in\Z_{\geq 2}$ and $\m\in\Z_+\times\Z_-$, the set $\Wse(J_{\m})$ is $E_{J,\m}$-invariant.
\end{corollary}

\begin{proof}
Let $J$ and $\m$ be given, fix arbitrary $x,y\in J_{\m}^{\omega^2}$, and write $x_0:=f_{J,\m}^{-1}(x)$ and $y_0:=f_{J,\m}^{-1}(y)$. Suppose $x\in\Wq(J_{\m})$ and $x\mathrel{E_{J,\m}}y$. Then $x_0\in\Wm(J)$ by Lemma \ref{claim:annoying-}, and $x_0\mathrel{\EJT}y_0$. Since $\Wm(J)$ is $\EJT$-invariant by Lemma \ref{claim:WmInvar}, we also have $y_0\in\Wm(J)$, so $y\in\Wq(J_{\m})$ again by \ref{claim:annoying-}.
\end{proof}

And now we have the analogues of Propositions \ref{prop:XmColorable+} and \ref{prop:WmHS+} for the case $\m\in\Z_+\times\Z_-$.

\begin{proposition}\label{prop:XmColorable-}
For all $\bar{m}\in\Z_+\times\Z_-$, we have $\chi_C\big(\G_{\m}(2)\res\Xm(2)\big)\leq 4^{|m_1m_2|}+3\cdot 2^{|m_1m_2|}+19$. 
\end{proposition}

\begin{proof}
Fix $\bar{m}\in\Z_+\times\Z_-$, and write $J:=2_{\m}=2^{|m_1m_2|}\in\omega$. Using Proposition \ref{prop:colorGSE}, let $c_{\q}^{J}:\Xq(J)\to S_{\q}(J)$ be a continuous coloring of $\Gq(J)\res\Xq(J)$, and define the continuous function $c_{\m}:\Xm(2)\to\Sse(J)$ by 
\[
c_{\m}(x) \ := \ \cse^J\big(f_{2,\m}(x)\big).
\]
Recall that $f_{2,\m}$ maps $\TT$ to $J^{\omega^2}$, and that
\[
|\Sse(J)| \ = \ J^2+3J+19 \ = \ 4^{|m_1m_2|}+3\cdot 2^{|m_1m_2|}+19
\]
by the discussion preceding Proposition \ref{prop:colorGSE}. 
Let $x,y\in\Xm(2)$ be arbitrary such that $x\mathrel{\Gm(2)}y$. Then $f_{2,\m}(x)$ and $f_{2,\m}(y)$ belong to $\Xq(J)$ by Lemma \ref{claim:annoying-}, and $f_{2,\m}(x)\mathrel{\Gq(J)}f_{2,\m}(y)$ by Lemma \ref{claim:fJmHom-}. This means $\cse^J\big(f_{2,\m}(x)\big)\ne\cse^J\big(f_{2,\m}(y)\big)$, so $c_{\m}(x)\ne c_{\m}(y)$, and we conclude that $c_{\m}$ is a coloring of $\G_{\m}(2)\res\Xm(2)$.
\end{proof}

\begin{proposition}\label{prop:WmHS-}
For all $\bar{m}\in\Z_+\times\Z_-$, $\ETT\res\Wm(2)$ is hyperselective.
\end{proposition}

\begin{proof}
Fix $\bar{m}\in\Z_+\times\Z_-$, and let $J:=2_{\m}=2^{|m_1m_2|}\in\omega$. By Corollary \ref{cor:EJmInvar}, $\Wse(J)$ is $E_{2,\m}$-invariant.  By Proposition \ref{prop:SENWhs} and Observation \ref{ob:UniformlyContinuous}, $\EJT\res\Wse(J)$ is hypersmooth, uniformly continuous, and has $K_{\sigma}$ classes.  Therefore, since by Lemma \ref{claim:compressionFI} $E_{2,\m}$ contains $\EJT$ as a finite-index sub-equivalence relation on $\Wse(J)$, the equivalence relation $E_{2,\m}\res\Wse(J)$ is hypersmooth by Corollary \ref{cor:fi1}. But $f_{2,\m}\res\Wm(2)$ is a Borel reduction from $\ETT\res\Wm(2)$ to $E_{2,\m}\res\Wse(J)$, so $\ETT\res\Wm(2)$ is hypersmooth as well, and is therefore hyperselective by Observations \ref{ob:UniformlyContinuous} and \ref{ob:UniformlyContinuousAndKSigma}.
\end{proof}

\subsection{The graphs $\Gm$ for $\m$ a nonzero multiple of $\bar{e}_1$ or $\bar{e}_2$}\label{sec:SubSec4}

In what should now be a familiar pattern of argument, we present the analogues for $\m\in\Z_+\times\{0\}$ of Lemmas \ref{claim:fJmHom+} and \ref{claim:fJmPeriodicity+}, Definition \ref{df:WmXm+}, Lemma \ref{claim:annoying+} and Corollary \ref{cor:EJmInvar+}, and finally Propositions \ref{prop:XmColorable+} and \ref{prop:WmHS+}. The proofs of Lemmas \ref{claim:fJmHom0} and \ref{claim:TooManyClaims0} are nearly identical to those of Lemmas \ref{claim:fJmHom+} and \ref{claim:fJmPeriodicity+}, respectively.

\begin{lemma}\label{claim:fJmHom0}
For all $J\in\Z_{\geq 2}$, $\m=(m_1,0)\in\Z_+\times\{0\}$, and $x,y\in\JT$, we have 
\begin{flalign*}
& & x\mathrel{\vec{\G}_{\m}(J)}y \quad &\Longleftrightarrow \quad f_{J,\m}(x)\mathrel{\vec{\G}_{\bar{e}_1}(J_{\m})}f_{J,\m}(y) & &  \\
\shortintertext{and}
& & x\mathrel{\G_{\m}(J)}y \quad &\Longleftrightarrow \quad f_{J,\m}(x)\mathrel{\G_{\bar{e}_1}(J_{\m})}f_{J,\m}(y). & & \makebox[0pt][r]{$\square$}
\end{flalign*}
\end{lemma}

\begin{lemma}\label{claim:TooManyClaims0}
For all $J\in\Z_{\geq 2}$, $\m=(m_1,0)\in\Z_+\times\{0\}$, and $x\in\JT$, $x$ is $\m$-periodic if and only if $f_{J,\m}(x)$ is $\bar{e}_1$-periodic. \hfill $\square$
\end{lemma}

\begin{definition}\label{df:XWm0}
Given $J\in\Z_{\geq 2}$ and $\m=(m_1,0)\in\Z_+\times\{0\}$, define
\begin{align*}
\Xm(J) \ := \ \Big\{ &x\in\neg\textup{Per}_{\m}(J^{\omega^2})\st \big(\exists l\in m_1\big)\,\big(\exists^{\infty}j\in\omega\big) \\& \text{$\big\langle x(l+im_1,j)\st i\in\omega\big\rangle$ is not eventually constant}\Big\}
\end{align*}
and
\[
\Wm(J) \ := \ \neg\textup{Per}_{\m}(\JT)\setminus\Xm(J).
\]
\end{definition}

Note that Definition \ref{df:XWm0} agrees with Definition \ref{df:XeWe} when $\m=\bar{e}_1$, and also that the sets $\Xm$ and $\Wm$ are Borel and $\EJT$-invariant.

\begin{lemma}\label{claim:ugh}
For all $J\in\Z_{\geq 2}$ and $\m=(m_1,m_2)\in\Z_+\times\{0\}$, we have 
\[
\Xm(J) \ = \ \big\{x\in\JT\st f_{J,\m}(x)\in\Xe(J_{\m})\big\} \qquad \text{and} \qquad 
\Wm(J) \ = \ \big\{x\in\JT\st f_{J,\m}(x)\in\We(J_{\m})\big\}.
\]
\end{lemma}

\begin{proof}
Let $J$ and $\m$ be given. By definition we have
\[
\neg\textup{Per}_{\m}(\JT) = \Xm(J)\sqcup\Wm(J) \qquad \text{and} \qquad \neg\textup{Per}_{\bar{e}_1}(J_{\m}^{\omega^2}) = \Xe(J_{\m})\sqcup\We(J_{\m}).
\]
Therefore, it will suffice by Lemma \ref{claim:TooManyClaims0} to show $f_{J,\m}[\Xm(J)]\subseteq\Xe(J_{\m})$ and $f_{J,\m}[\Wm(J)]\subseteq\We(J_{\m})$.

To show these, first note that for each $x\in\neg\textup{Per}_{\m}(\JT)$ and $j\in\omega$, the sequence $\big\langle f_{J,\m}(x)(i,j)\st i\in\omega\big\rangle$ is eventually constant if and only if for all $l\in m_1$ the sequence $\big\langle x(l+im_1, j)\st i\in\omega\big\rangle$ is eventually constant. In particular, if for all $l\in m_1$ and all but finitely many $j\in\omega$ the sequence $\big\langle x(l+im_1, j)\st i\in\omega\big\rangle$ is eventually constant, then likewise the sequence $\big\langle f_{J,\m}(x)(i,j)\st i\in\omega\big\rangle$ is eventually constant for all but finitely many $j$; that is, we have $f_{J,\m}[\Wm(J)]\subseteq\We(J_{\m})$. Similarly, if for some $l\in m_1$ there are infinitely many $j$ for which $\big\langle x(l+im_1, j)\st i\in\omega\big\rangle$ is not eventually constant, then there are infinitely many $j$ for which $\big\langle f_{J,\m}(x)(i,j)\st i\in\omega\big\rangle$ is not eventually constant; that is, we have $f_{J,\m}[\Xm(J)]\subseteq\Xe(J_{\m})$. 
\end{proof}

\begin{corollary}\label{claim:EJmInvar+}
For all $J\in\Z_{\geq 2}$ and $\m=(m_1,0)\in\Z_+\times\{0\}$, the set $\We(J_{\m})$ is $E_{J,\m}$-invariant.
\end{corollary}

\begin{proof}
Let $J$ and $\m$ be given, fix arbitrary $x,y\in J_{\m}^{\omega^2}$, and write $x_0:=f_{J,\m}^{-1}(x)$ and $y_0:=f_{J,\m}^{-1}(y)$. Suppose $x\in\We(J_{\m})$ and $x\mathrel{E_{J,\m}}y$. Then $x_0\in\Wm(J)$ by Lemma \ref{claim:ugh}, and $x_0\mathrel{\EJT}y_0$. Since $\Wm(J)$ is $\EJT$-invariant, we also have $y_0\in\Wm(J)$, so $y\in\We(J_{\m})$ again by \ref{claim:ugh}.
\end{proof}

\begin{proposition}\label{prop:XmColorable0}
For all $\m=(m_1,0)\in\Z_+\times\{0\}$, we have $\chi_B\big(\Gm(2)\res\Xm(2)\big)\leq 3$.
\end{proposition}

\begin{proof}
Fix $\m\in\Z_+\times\{0\}$. Using Proposition \ref{prop:ColoringGe1}, let $c_{\bar{e}_1}:\Xe(2_{\m})\to 3$ be a Borel $3$-coloring of $\Ge(2_{\m})\res\Xe(2_{\m})$, and define the Borel function $c_{\m}:\Xm(2)\to 3$ by 
\[
c_{\m}(x) \ := \ c_{\bar{e}_1}\big(f_{2,\m}(x)\big).
\]
Let $x,y\in\Xm(2)$ be arbitrary such that $x\mathrel{\G_{\m}(2)}y$. Then $f_{2,\m}(x)\in\Xe(2_{\m})$ and $f_{2,\m}(y)\in\Xe(2_{\m})$ by Lemma \ref{claim:ugh}, and $f_{2,\m}(x)\mathrel{\Ge(2_{\m})}f_{2,\m}(y)$ by Lemma \ref{claim:fJmHom0}. This means $c_{\bar{e}_1}\big(f_{2,\m}(x)\big)\ne c_{\bar{e}_1}\big(f_{2,\m}(y)\big)$, so $c_{\m}(x)\ne c_{\m}(y)$, and we have shown that $c_{\m}$ is a Borel $3$-coloring of $\G_{\m}(2)\res\Xm(2)$.
\end{proof}

\begin{proposition}\label{prop:WmHS0}
For all $\m=(m_1,0)\in\Z_+\times\{0\}$, $\ETT\res\Wm(2)$ is hyperselective.
\end{proposition}

\begin{proof}
Fix $\m=(m_1,0)\in\Z_+\times\{0\}$, and let $J:=2_{\m}=2^{m_1}\in\omega$. By Corollary \ref{claim:EJmInvar+}, $\We(J)$ is $E_{2,\m}$-invariant. By Proposition \ref{prop:WeHS} and Observation \ref{ob:UniformlyContinuous}, $\EJT\res\We(J)$ is hypersmooth, uniformly continuous, and has $K_{\sigma}$ classes. Therefore, since by Lemma \ref{claim:compressionFI} $E_{2,\m}$ contains $\EJT$ as a finite-index sub-equivalence relation on $\We(J)$, the equivalence relation $E_{2,\m}\res\We(J)$ is hypersmooth by Corollary \ref{cor:fi1}. But $f_{2,\m}\res\Wm(2)$ is a Borel reduction from $\ETT\res\Wm(2)$ to $E_{2,\m}\res\We(J)$, so $\ETT\res\Wm(2)$ is hypersmooth as well, and is therefore hyperselective by Observations \ref{ob:UniformlyContinuous} and \ref{ob:UniformlyContinuousAndKSigma}.
\end{proof}

Finally, we consider the case $\m\in\{0\}\times\Z_+$. For each $J\in\Z_{\geq 2}$ and $\m=(0,m_2)\in\{0\}\times\Z_+$, let
\[
\Xm(J) := \delta\big[X_{\m^\delta}(J)\big] \qquad \text{and} \qquad \Wm(J) := \delta\big[W_{\m^\delta}(J)\big],
\]
so $\Xm(J)$ and $\Wm(J)$ are Borel and $\ETT$-invariant, and $\neg\textup{Per}_{\m}(\JT)=\Xm(J)\sqcup\Wm(J)$.

\begin{corollary}\label{cor:XmColorable0'}
For all $\bar{m}=(0,m_2)\in\{0\}\times\Z_+$, $\chi_B\big(\Gm(2)\res\Xm(2)\big)\leq 3$ and $\ETT\res\Wm(2)$ is hyperselective.
\end{corollary}

\begin{proof}
Fix $\m=(0,m_2)\in\{0\}\times\Z_+$, so $\m^\delta=(m_2,0)$. Using Proposition \ref{prop:XmColorable0}, let $c_{\m^{\delta}}$ be a Borel $3$-coloring of $\G_{\m^{\delta}}(2)\res X_{\m^{\delta}}(2)$.  Then $c_{\m}:=c_{\m^{\delta}}\circ\delta$ is a Borel $3$-coloring of $\Gm(2)\res\Xm(2)$, and $\delta\res\Wm(2)$ is a Borel reduction from $\ETT\res\Wm(2)$ to $\ETT\res W_{\m^{\delta}}(2)$, so $\ETT\res\Wm(2)$ is hypersmooth by Proposition \ref{prop:WmHS0} and hyperselective by Observations \ref{ob:UniformlyContinuous} and \ref{ob:UniformlyContinuousAndKSigma}.
\end{proof}

\subsection{The graphs $\G_{2,n}$ and the proof of Theorem \ref{thm:main}}

Now we tie together the results obtained so far in this section in order to complete our proof of Theorem \ref{thm:main}. 
Recall from Section \ref{sec:graphs} that for each $n\in\Z_+$ and nonzero $\m\in\Z^2$, $\G_{\m,n}$ is the (undirected) graph on $\TT$ defined by
\[
x\mathrel{\G_{\m,n}}y \quad :\Longleftrightarrow \quad x\ne y \, \wedge \, 
\big(\exists \bar{a},\bar{b}\in\N^{2}_{<n}\big) \, \Big(\bar{a}\cdot x=\bar{b}\cdot y \, \wedge \, \bar{b}-\bar{a}=\pm\m\Big),
\]
and $\G_{2,n}=\G_{2,n,+}\cup\G_{2,n,-}$ where
\[
\G_{2,n,+} \ := \ \bigcup_{\m \mkern1mu \in \mkern1mu \N^{2}_{<n}\setminus\{\bar{0}\}}\G_{\m,n} \qquad \text{and} \qquad 
\G_{2,n,-} \ := \ \bigcup_{\m \mkern1mu \in \mkern1mu (\Z_+\times\Z_-)_{<n}}\G_{\m,n}.
\]
We can color the graphs $\G_{2,n,+}$ and $\G_{2,n,-}$ by combining our colorings of the graphs $\G_{\m}$ that we obtained in Propositions \ref{prop:XmColorable+}, \ref{prop:XmColorable-}, and \ref{prop:XmColorable0} and Corollary \ref{cor:XmColorable0'}.

\begin{corollary}\label{cor:+good}
For all $n\in\Z_+$, there is a partition of $\FTT$ into a pair of $\ETT$-invariant Borel sets $\FTT=X_n^+\sqcup W_n^+$ such that $\ETT\res W_n^+$ is hypersmooth and $\G_{2,n,+}\res X_n^+$ admits a Borel $3^{n^2}$-coloring.
\end{corollary}

\begin{proof}
Fix $n\in\Z_+$, and write $\mathbb{M}_n^+=\N^2_{<n}\setminus\{\bar{0}\}$, noting that $|\mathbb{M}_n^+|=n^2-1$. Using Propositions \ref{prop:XmColorable+} and \ref{prop:XmColorable0} and Corollary \ref{cor:XmColorable0'}, for each $\m\in\mathbb{M}_n^+$ let $c_{\m}:\Xm(2)\to 3$ be a Borel coloring of $\Gm(2)\res\Xm(2)$. Let
\[
X_n^+ \, := \, \bigcap_{\m\mkern1mu\in\mkern1mu\mathbb{M}_n^+}X_{\m}^{\textup{Fr}}(2) \qquad\quad \text{and} \qquad\quad 
W_n^+ \, := \, \bigcup_{\m\mkern1mu\in\mkern1mu\mathbb{M}_n^+}W_{\m}^{\textup{Fr}}(2),
\]
and define the Borel function $c_n^+:X^+_n\to 3^{\mathbb{M}_n^+}$ by $c_n^+(x)(\m)=c_{\m}(\d{n}\cdot x)$, so $|\ran(c_n^+)|\leq 3^{n^2}$. Then $\ETT\res W_n^+$ is hypersmooth by Propositions \ref{prop:WmHS+} and \ref{prop:WmHS0}, Corollary \ref{cor:XmColorable0'}, and Lemma \ref{lem:piecewiseHS}, and we claim that $c_n^+$ colors $\G_{2,n,+}\res X_n^+$. To see this, fix $x,y\in X_n^+$ such that $x\mathrel{\G_{2,n,+}}y$, say without loss of generality that $\m:=y-x\in\mathbb{M}_n^+$. Fix $\bar{a},\bar{b}\in\N^2_{<n}$ such that $\bar{a}\cdot x=\bar{b}\cdot y$, so $\m=\bar{a}-\bar{b}$. Then 
\[
\d{n}\cdot y \ = \ (\d{n}-\bar{b})\cdot (\bar{b}\cdot y) \ = \ (\d{n}-\bar{b})\cdot (\bar{a}\cdot x) \ = \ \m\cdot (\d{n}\cdot x),
\]
so $(\d{n}\cdot x)\mathrel{\Gm}(\d{n}\cdot y)$. Thus $c_{\m}(\d{n}\cdot x)\ne c_{\m}(\d{n}\cdot y)$, which means $c_n^+(x)\ne c_n^+(y)$. 
\end{proof}

\begin{corollary}\label{cor:-good}
 For all $n\in\Z_+$, there is a partition of $\FTT$ into a pair of $\ETT$-invariant Borel sets $\FTT=X_n^-\sqcup W_n^-$ such that $\ETT\res W_n^-$ is hypersmooth and $\G_{2,n,-}\res X_n^-$ admits a continuous $4^{n^4}$-coloring.
\end{corollary}

\begin{proof}
Fix $n\in\Z_+$, and write $\mathbb{M}_n^-=(\Z_+\times\Z_-)_{<n}$, noting that $|\mathbb{M}_n^-|=(n-1)^2$. Using Proposition \ref{prop:XmColorable-}, for each $\m\in\mathbb{M}_n^-$ let $c_{\m}:\Xm(2)\to |S_{\q}(2_{\m})|$ be a continuous coloring of $\Gm(2)\res\Xm(2)$ using natural numbers as colors. Put
\[
S_{\q,n} := \max\big\{|S_{\q}(2_{\m})|\st \m\in\mathbb{M}_n^-\big\}, 
 \qquad X_n^- := \bigcap_{\m\mkern1mu\in\mkern1mu\mathbb{M}_n^-}X_{\m}^{\textup{Fr}}(2), \qquad \text{and} \qquad 
W_n^- := \bigcup_{\m\mkern1mu\in\mkern1mu\mathbb{M}_n^-}W_{\m}^{\textup{Fr}}(2),
\]
and define the continuous function $c_n^-:X^-_n\to\left(S_{\q,n}\right)^{\mathbb{M}_n^-}$ by $c_n^-(x)(\m)=c_{\m}(\d{n}\cdot x)$. 
Recall that for all $\m\in\mathbb{M}_n^-$ we have $|S_{\q}(2_{\m})|=4^{|m_1m_2|}+3\cdot 2^{|m_1m_2|}+19\leq 4^{n^2+1}$, so $|\ran(c_n^-)|\leq \big(4^{n^2+1}\big)^{(n-1)^2}\leq 4^{n^4}$. 

Now, $\ETT\res W_n^-$ is hypersmooth by Proposition \ref{prop:WmHS-} and Lemma \ref{lem:piecewiseHS}, and we claim that $c_n^-$ colors $\G_{2,n,-}\res X_n^-$. To see this, fix $x,y\in X_n^-$ such that  $x\mathrel{\G_{2,n,-}}y$, say without loss of generality that $\m:=y-x\in\mathbb{M}_n^-$. Fix $\bar{a},\bar{b}\in\N^2_{<n}$ such that $\bar{a}\cdot x=\bar{b}\cdot y$, so $\m=\bar{a}-\bar{b}$. 
Then
\[
\d{n}\cdot y \ = \ (\d{n}-\bar{b})\cdot (\bar{b}\cdot y) \ = \ (\d{n}-\bar{b})\cdot (\bar{a}\cdot x) \ = \ (\d{n}+\bar{m})\cdot x,
\]
so $m_1\bar{e}_1\cdot (\d{n}\cdot x)=-m_2\bar{e}_2\cdot (\d{n}\cdot y)$. This shows $(\d{n}\cdot x)\mathrel{\Gm}(\d{n}\cdot y)$, which means $c_{\m}(\d{n}\cdot x)\ne c_{\m}(\d{n}\cdot y)$. Thus $c_n^-(x)(\m)\ne c_n^-(y)(\m)$, so $c_n^-(x)\ne c_n^-(y)$, completing the proof.
\end{proof}

\begin{proof}[Proof of Theorem \ref{thm:main}]
Fix $n\in\Z_+$. Let $\mathbb{M}_n^+$ and $\mathbb{M}_n^-$ be as above, along with $c_n^+$ and $c_n^-$ and the sets $X_n^+$, $X_n^-$, $W_n^+$, $W_n^-$, and $S_{\q,n}$. Let
\[
X_n \ := \ X_n^+\cap X_n^- \qquad \text{and} \qquad W_n \ := \ W_n^+\cup W_n^-,
\]
so $\FTT=X_n\sqcup W_n$. Then $\ETT\res W_n$ is hypersmooth by Corollaries \ref{cor:+good} and \ref{cor:-good} and Lemma \ref{lem:piecewiseHS}. Now define the Borel function $c_n:X_n\to 3^{\mathbb{M}_n^+}\times \left(S_{\q,n}\right)^{\mathbb{M}_n^-}$ by $c_n(x)=\big(c_n^+(x), \, c_n^-(x)\big)$. Then $c_n$ colors $\G_{2,n}\res X_n$, and $|\ran(c_n)|\leq 3^{n^2}\cdot 4^{n^4}$.
\end{proof}

\section{Hypersmoothness in measure of $\ETl$ for $\ell\geq 3$}\label{sec:measure}

In this section, we prove that with respect to the natural $(\frac{1}{2},\frac{1}{2})$-probability measures $\mu_{\ell}$ on the spaces $\Tl$, Conjecture \ref{con:HigherDimensions} holds on an invariant conull set for all $\ell\geq 3$. After showing this, we also provide a proof that Conjecture \ref{con:HigherDimensions} implies Conjecture \ref{con:final}.

\begin{theorem}\label{thm:measure}
For all $\ell\in\Z_{\geq 3}$, there is a $\mu_{\ell}$-conull, $\ETl$-invariant Borel set $X\subseteq\Tl$ such that $\ETl\res X$ is hyperselective.
\end{theorem}

The same conclusion then holds for any Borel probability measure $\mu_{\ell}'$ on $\Tl$ having the same null sets as $\mu_{\ell}$. We could use any finite base space here, but for simplicity in this section we stick to $2=\{0,1\}$ and omit the base space from our notation. Our argument ultimately relies on Corollary \ref{cor:main}. Intuitively, the ``wild'' sets where we could not color our graphs in Sections \ref{sec:ep}--\ref{sec:Gm} all have measure zero, so we can ignore them, obtain finite Borel graph colorings on their complements, and then appeal to Corollary \ref{cor:main}. We begin with a simple lemma.

\begin{lemma}\label{lem:unlikely}
Let $\ell\in\Z_{\geq 3}$. Let $(\bar{a}_k)_{k\in\omega}$ and $(\bar{b}_k)_{k\in\omega}$ be disjoint sequences of distinct points in $\N^{\ell}$, so as maps $k\mapsto\bar{a}_k$ and $k\mapsto\bar{b}_k$ are injective with disjoint ranges. Then
\[
\mu_{\ell}\Big(\big\{ x\in\Tl\st (\forall k\in\omega) \, x(\bar{a}_k)=x(\bar{b}_k)\big\}\Big)=0.
\]
\end{lemma}

\begin{proof}
Fix $\ell\in\Z_{\geq 3}$, and let $(\bar{a}_k)_{k\in\omega}$ and $(\bar{b}_k)_{k\in\omega}$ be as stated. For each $k\in\omega$, write
\[
A_k := \big\{x\in\Tl\st x(\bar{a}_k)=x(\bar{b}_k)\big\} \qquad \text{and} \qquad
A_{<k} := \big\{x\in\Tl\st (\forall i<k) \; x(\bar{a}_i)=x(\bar{b}_i)\big\}.
\]
For each $k\in\omega$, we have $\mu_{\ell}(A_k)=\frac{1}{2}$ since $\bar{a}_k\ne\bar{b}_k$, and therefore
\[
\mu_{\ell}(A_{<k}) \ = \ \prod_{i<k}\mu_{\ell}(A_i) \ = \ 2^{-k}
\]
since $\bar{a}_0,\ldots,\bar{a}_{k-1},\bar{b}_0,\ldots,\bar{b}_{k-1}$ are all distinct. 
As $\big\{ x\in\Tl\st (\forall k\in\omega) \, x(\bar{a}_k)=x(\bar{b}_k)\big\}$ is the intersection of the decreasing sequence of sets $(A_{<k})$, it follows that
\[
\mu_{\ell}\Big(\big\{ x\in\Tl\st (\forall k\in\omega) \, x(\bar{a}_k)=x(\bar{b}_k)\big\}\Big) \ = \ \lim_{k\to\infty}\mu_{\ell}\big(A_{<k}\big) \ = \ \lim_{k\to\infty}2^{-k} \ = \ 0. \qedhere
\]
\end{proof}

Lemma \ref{lem:unlikely} enables us to ignore the periodic parts in our proof of Theorem \ref{thm:measure}:

\begin{corollary}
For all $\ell\in\Z_+$ and nonzero $\m\in\Z^{\ell}$, we have $\mu_{\ell}\big(\textup{Per}_{\m}(\Tl)\big)=0$. Consequently, $\mu_{\ell}\Big(\textup{Per}(\Tl)\Big)=0$ and $\mu_{\ell}\Big(\FTl\Big)=1$. \hfill $\square$
\end{corollary}

We consider now the graphs $\G_{\m}$ for nonzero $\m\in\Z^{\ell}$, splitting our argument into the two cases $\m\in\pm\N^{\ell}$ and otherwise. For the first case, we may assume $\m\in\N^{\ell}$ since $\G_{\m}=\G_{-\m}$.

\begin{proposition}\label{prop:forward}
For each $\ell\in\Z_{\geq 3}$ and nonzero $\m\in\N^{\ell}$, there is an $\ETl$-invariant Borel set $\tilde{X}_{\m}\subseteq\neg\textup{Per}_{\m}(\Tl)$ such that $\mu_{\ell}(\tilde{X}_{\m})=1$ and $\chi_B\big(\G_{\m}\res\tilde{X}_{\m}\big)\leq 3$.
\end{proposition}

\begin{proof}
Given $\ell\in\Z_{\geq 3}$ and nonzero $\m\in\N^{\ell}$, let 
\[
\tilde{X}_{\m} \ := \ \Big\{x\in\neg\textup{Per}_{\m}(\Tl)\st (\forall\bar{a}\in\N^{\ell})(\exists\bar{b}\in\N^{\ell})\, \langle (\bar{a}\cdot x)(\bar{b}+i\m)\st i\in\omega\rangle \text{ is not eventually constant}\Big\}.
\]
Clearly $\tilde{X}_{\m}$ is Borel and $\ETl$-invariant. 
For each fixed $\bar{a},\bar{b}\in\N^{\ell}$, let $A_{\bar{a},\bar{b}}$ be the set of all $x\in\Tl$ such that $\langle (\bar{a}\cdot x)(\bar{b}+i\m)\st i\in\omega\rangle$ is not eventually constant. Then for each $\bar{a},\bar{b}\in\N^{\ell}$,  $\Tl\setminus A_{\bar{a},\bar{b}}$ has measure zero by Lemma \ref{lem:unlikely}, so $\mu(A_{\bar{a},\bar{b}})=1$, which means
\[
A \ := \ \bigcap_{\bar{a} \mkern1mu \in \mkern1mu \N^{\ell}}\bigcup_{\bar{b} \mkern1mu \in \mkern1mu \N^{\ell}}A_{\bar{a},\bar{b}}
\]
has measure one as well. Since
\[
\tilde{X}_{\m} \ = \ A\cap\neg\textup{Per}_{\m}(\Tl)
\]
and $\mu\big(\neg\textup{Per}_{\m}(\Tl)\big)=1$, we have $\mu(\tilde{X}_{\m})=1$.

To complete the proof, we define a Borel $3$-coloring $\tilde{c}_{\m}:\tilde{X}_{\m}\to 3$ of $\G_{\m}\res\tilde{X}_{\m}$ as follows. Using Fact \ref{fact:KST5.2}, let $c:2^\omega\to 3$ be a Borel coloring of $\G_1$. Given $x\in\tilde{X}_{\m}$, let $\bar{b}_x\in\N^{\ell}$ be least such that $\langle x(\bar{b}_x+i\m)\st i\in\omega\rangle$ is not eventually constant, and put
\[
\tilde{c}_{\m}(x) \ := \ c\Big(\langle x(\bar{b}_x+i\m)\st i\in\omega\rangle\Big).
\]
To see that $\tilde{c}_{\m}$ colors $\G_{\m}$ on $\tilde{X}_{\m}$, let $x,y\in\tilde{X}_{\m}$ be arbitrary such that $x\mathrel{\G_{\m}}y$, say without loss of generality that $\m\cdot x=y$. Note that for all $\bar{b}\in\N^{\ell}$, we have 
\[
1\cdot \big\langle x(\bar{b}+i\m)\st i\in\omega\big\rangle \ = \ \big\langle y(\bar{b}+i\m)\st i\in\omega\big\rangle.
\]
Thus $\bar{b}_x=\bar{b}_y$ and
\[
\big\langle x(\bar{b}_x+i\m)\st i\in\omega\big\rangle \mathrel{\G_1} 
\big\langle y(\bar{b}_y+i\m)\st i\in\omega\big\rangle,
\]
which implies $\tilde{c}_{\m}(x)\ne\tilde{c}_{\m}(y)$, completing the proof.
\end{proof}

Next we consider the case $\m\in\Z^{\ell}\setminus\pm\N^{\ell}$, so both $\m^+$ and $\m^-$ are nonzero (recall that this notation was introduced in Equation \eqref{eq:mplusminus}). For each such $\m$, define the linear maps $T_{\m}:\mathbb{R}^2\to\mathbb{R}^{\ell}$ and $\gamma_{\m}:\mathbb{R}^2\to\mathbb{R}^{\ell}$ and the continuous function $g_{\m}:2^{\omega^\ell}\to 2^{\omega^2}$ by
\[
T_{\m}(i,j) \ := \ i\m^++j\m^-+(i+j)\d{1}, \qquad \gamma_{\m}(a_1,a_2) \ := \ a_1(\m^++\d{1})+a_2(\m^-+\d{1}),
\]
and 
\[
g_{\m}(x)(i,j) \ := \ x\Big(T_{\m}(i,j)\Big) \ = \ x\Big(i\m^++j\m^-+(i+j)\d{1}\Big).
\]
By a straightforward calculation, for all $\ell\in\Z_{\geq 3}$, $x\in\Tl$, $\m\in\Z^{\ell}\setminus\pm\N^{\ell}$, and $\bar{a}\in\N^2$, we have 
\begin{equation}\label{eq:gmNewUseful}
\bar{a}\cdot g_{\m}(x) \ = \ g_{\m}\Big(\gamma_{\m}(\bar{a})\cdot x\Big).
\end{equation}

From this we can easily show the following:

\begin{lemma}\label{lem:SatBorel}
Let $\ell\in\Z_{\geq 3}$ and $\m\in\Z^{\ell}\setminus\pm\N^{\ell}$, and suppose the set $X\subseteq\TT$ is $\ETT$-dense in its $\ETT$-saturation. Then the set $g_{\m}^{-1}[X]$ is $\ETl$-dense in its $\ETl$-saturation; therefore,
\[
\Big[g_{\m}^{-1}[X]\Big]_{\ETl} \ = \ \bigcup_{\bar{a} \mkern1mu \in \mkern1mu \N^{\ell}}\bar{a}^{-1}\big[g_{\m}^{-1}[X]\big],
\]
so the $\ETl$-saturation of $g_{\m}^{-1}[X]$ is Borel if $X$ is Borel.
\end{lemma}

\begin{proof}
Let $\ell$ and $\m$ be given, and suppose $X\subseteq\TT$ is $\ETT$-dense. Fix $x\in g_{\m}^{-1}[X]$ and $\bar{a}\in\N^{\ell|}$. Choose $c_1,c_2\in\N$ large enough so that $\bar{a}\preceq\gamma_{\m}(\c)$ and $\c\cdot g_{\m}(x)\in X$. Then $g_{\m}(\gamma_{\m}(\c)\cdot x)=\c\cdot g_{\m}(x)\in X$, showing $g_{\m}^{-1}[X]$ is $\ETl$-dense as claimed.
\end{proof}

Now, one can show that for all $\m\in\Z^{\ell}\setminus\pm\N^{\ell}$ and $x\in\textup{Per}_{\m}(\Tl)$, we have $g_{\m}(x)\in\textup{Per}_{\q}(\TT)$. If $x\in\neg\textup{Per}_{\m}(\Tl)$, however, we cannot necessarily expect $g_{\m}(x)\in\neg\textup{Per}_{\q}(\TT)$; indeed, we cannot even expect that there is $\bar{a}\in\N^{\ell}$ such that $g_{\m}(\bar{a}\cdot x)\in\neg\textup{Per}_{\q}(\TT)$. Nevertheless, we have the following:

\begin{claim}\label{claim:PerMeasureZero}
For all $\ell\in\Z_{\geq 3}$ and $\m\in\Z^{\ell}\setminus\pm\N^{\ell}$, the $\ETl$-saturation of the set $g_{\m}^{-1}\big[\textup{Per}_{\q}(\TT)\big]$ is Borel and has measure zero in $\Tl$.
\end{claim}

\begin{proof}
Fix $\ell\in\Z_{\geq 3}$ and $\m\in\Z^{\ell}\setminus\pm\N^{\ell}$, and write $Y_{\m}:=g_{\m}^{-1}\big[\textup{Per}_{\q}(\TT)\big]$ and $\tilde{Y}_{\m}:=[Y_{\m}]_{\ETl}$. That $\tilde{Y}_{\m}$ is Borel follows from Lemma \ref{lem:SatBorel}. We will show $\mu_{\ell}\big(\bar{a}^{-1}[Y_{\m}]\big)=0$ for each $\bar{a}\in\N^{\ell}$, implying $\mu_{\ell}(\tilde{Y}_{\m})=0$ again by Lemma \ref{lem:SatBorel}. 

Let $\bar{a}\in\N^{\ell}$ and $x\in \bar{a}^{-1}[Y_{\m}]$ be arbitrary. Since $g_{\m}(\bar{a}\cdot x)$ is $\q$-periodic, we may fix $\c\in\N^{2}$ such that $(\c+\bar{e}_1)\cdot g_{\m}(\bar{a}\cdot x)=(\c+\bar{e}_2)\cdot g_{\m}(\bar{a}\cdot x)$. For each $n\in\omega$, write 
\[
\bar{b}_n \ := \ (c_1+n)\m^++(c_2+n)m^-+(c_1+c_2+2n+1)\d{1} \ \in \ \N^{\ell},
\]
so $n\mapsto \bar{b}_n$ is injective. Then for all $n\in\omega$ we have
\begin{align*}
x(\bar{a}+\bar{b}_n+\m^+) \ &= \ g_{\m}(\bar{a}\cdot x)\big(c_1+n+1, \, c_2+n\big) \\ &= \ g_{\m}(\bar{a}\cdot x)\big(c_1+n, \, c_2+n+1\big) \ = \ x(\bar{a}+\bar{b}_n+\m^-).
\end{align*}
The set $Z_{\m,\bar{a},\bar{c}}$ of all $x\in\Tl$ such that $x(\bar{a}+\bar{b}_n+\m^+)=x(\bar{a}+\bar{b}_n+\m^-)$ for all $n\in\omega$ has measure zero by Lemma \ref{lem:unlikely}; therefore, since the argument above shows 
\[
\bar{a}^{-1}\big[Y_{\m}\big] \ \subseteq \ \bigcup_{\bar{c} \mkern1mu \in \mkern1mu \N^2}Z_{\m,\bar{a},\bar{c}},
\]
we conclude that $\bar{a}^{-1}[Y_{\m}]$ has measure zero as well to complete the proof.
\end{proof}

\begin{claim}\label{claim:WqMeasureZero}
For all $\ell\in\Z_{\geq 3}$ and $\m\in\Z^{\ell}\setminus\pm\N^{\ell}$, the $\ETl$-saturation of the set $g_{\m}^{-1}[\Wq]$ is Borel and has measure zero in $\Tl$.
\end{claim}

\begin{proof}
Our proof is similar to that of Claim \ref{claim:PerMeasureZero}, using Lemma \ref{lem:SatBorel}. Again fixing $\ell\in\Z_{\geq 3}$ and $\m\in\Z^{\ell}\setminus\pm\N^{\ell}$, we write $Y_{\m}:=g_{\m}^{-1}[\Wq]$ and $\tilde{Y}_{\m}:=[Y_{\m}]_{\ETl}$ and show $\mu_{\ell}\big(\bar{a}^{-1}[Y_{\m}]\big)=0$ for each $\bar{a}\in\N^{\ell}$.

Let $\bar{a}\in\N^{\ell}$ and $x\in\bar{a}^{-1}[Y_{\m}]$ be arbitrary. By the definition of $\Wq$, we may fix $N\in\N$ such that for all $n\geq N$ we have $g_{\m}(\bar{a}\cdot x)(n+1,n) = g_{\m}(\bar{a}\cdot x)(n,n+1)$. For each $n\geq N$, write
\[
\bar{b}_n \ := \ \bar{a}+n\m^++n\m^-+(2n+1)\d{1},
\]
so $n\mapsto\bar{b}_n$ is injective. Then for all $n\geq N$ we have
\begin{align*}
x(\bar{b}_n+\m^+) \ = \ g_{\m}(\bar{a}\cdot x)(n+1,n) \ = \ g_{\m}(\bar{a}\cdot x)(n,n+1) \ = \ x(\bar{b}_n+\m^-).
\end{align*}
The set $Z_{\m,\bar{a},N}$ of all $x\in\Tl$ satisfying $x(\bar{b}_m+\m^+)=x(\bar{b}_n+\m^-)$ for all $n\geq N$ has measure zero by Lemma \ref{lem:unlikely}; therefore, since the argument above shows
\[
\bar{a}^{-1}[Y_{\m}] \ \subseteq \ \bigcup_{N\in \mkern2mu \omega}Z_{\m,\bar{a},N},
\]
we conclude that $\bar{a}^{-1}[Y_{\m}]$ has measure zero as well. By Lemma \ref{lem:SatBorel}, this completes the proof.
\end{proof}

\begin{proposition}\label{prop:sideways}
For each $\ell\in\Z_{\geq 3}$ and $\m\in\Z^{\ell}\setminus\pm\N^{\ell}$, there is an $\ETl$-invariant Borel set $\tilde{X}_{\m}\subseteq\neg\textup{Per}_{\m}(\Tl)$ such that $\mu_{\ell}(\tilde{X}_{\m})=1$ and $\chi_C\big(\G_{\m}\res\tilde{X}_{\m}\big)\leq 33$.
\end{proposition}

\begin{proof}
Let $\ell$ and $\m$ be given, and let $\tilde{X}_{\m}$ be the complement in $\TT$ of the $\ETl$-saturation of $g_{\m}^{-1}[\Wq]\cup g_{\m}^{-1}[ \textup{Per}_{\q}(\TT)]$. Then $\tilde{X}_{\m}$ is Borel and $\ETl$-invariant and $\mu_{\ell}(\tilde{X}_{\m})=1$ by Claims \ref{claim:PerMeasureZero} and \ref{claim:WqMeasureZero}. Since $\TT=\Xq\sqcup\Wq\sqcup\textup{Per}_{\q}(\TT)$, we have $\tilde{X}_{\m}\subseteq g_{\m}^{-1}[\Xq]$. Now let $c_{\q}$ be the continuous coloring of $\Gq(2)\res\Xq(2)$ given in Proposition \ref{prop:colorGSE}. Since $J=2$ here, we have $|\ran(c_{\q})|=33$. We claim that the continuous map $c_{\q}\circ g_{\m}$ is a coloring of $\G_{\m}\res\tilde{X}_{\m}$. 
To see this, fix $x,y\in\tilde{X}_{\m}$ such that $x\mathrel{\vec{\G}_{\m}}y$, so $\m^+\cdot x=\m^-\cdot y$. Then by Equation \ref{eq:gmNewUseful} we have
\[
\bar{e}_1\cdot g_{\m}(x) \ = \ g_{\m}\big((\m^++\d{1})\cdot x\big) \ = \ g_{\m}\big((\m^-+\d{1})\cdot y\big) \ = \ \bar{e}_2\cdot g_{\m}(y),
\]
so $g_{\m}(x)\mathrel{\vec{\G}_{\q}}g_{\m}(y)$, which means $c_{\q}\big(g_{\m}(x)\big)\ne c_{\q}\big(g_{\m}(y)\big)$ as desired.
\end{proof}

\begin{proof}[Proof of Theorem \ref{thm:measure}]
Fix $\ell\in\Z_{\geq 3}$. Using Propositions \ref{prop:forward} and \ref{prop:sideways}, for each nonzero $\m\in\Z^{\ell}$ fix an $\ETl$-invariant $\mu_{\ell}$-conull Borel set $\tilde{X}_{\m}\subseteq\neg\textup{Per}_{\m}(\Tl)$ 
and a Borel function $c_{\m}:\tilde{X}_{\m}\to 33$ that colors $\G_{\m}\res\tilde{X}_{\m}$. Let
\[
X \ := \ \textup{Fr}(\Tl) \,\cap \bigcap_{\bar{0} \mkern1mu \ne \mkern1mu \m \mkern1mu \in \mkern1mu \Z^{\ell}}\tilde{X}_{\m}.
\]
Then $X$ is Borel and $\ETl$-invariant, and $\mu_{\ell}(X)=1$. Now let $n\in\Z_+$ be arbitrary, and define the Borel map $c_n:X\to 33^{(\Z^{\ell}_{<n}\setminus\{\bar{0}\})}$ by
\[
c_n(x)(\m) \ := \ c_{\m}(\d{n}\cdot x).
\]
An argument identical to the ones in Corollaries \ref{cor:+good} and \ref{cor:-good} shows that $c_n$ colors $\G_{\ell,n}\res X$. Thus $\chi_B(\G_{\ell,n}\res X)$ is finite for every $n$, so $\ETl\res X$ is hyperselective by Corollary \ref{cor:main}.
\end{proof}


Finally, as promised in the Introduction, we show that Conjecture \ref{con:final} follows from Conjecture \ref{con:HigherDimensions}, using the following two lemmas.

\begin{lemma}\label{lem:coding1}
For all $\ell\in\Z_+$, we have $E^{2^\omega}_{\ell}\sqsubseteq_B E^3_{\ell+1}$.
\end{lemma}

\begin{proof}[Proof sketch]
Fix $\ell\in\Z_+$. Writing $\m\in\omega^{\ell}$, define the Borel function $f:(2^\omega)^{\omega^\ell}\to 3^{\omega^{\ell+1}}$ by 
\begin{align*}
f(x)(\bar{m},0) \ &= \ x(\bar{m})(0) \\
f(x)(\bar{m},1) \ &= \ 2 \\
f(x)(\bar{m},2) \ &= \ x(\bar{m})(0) \\
f(x)(\bar{m},3) \ &= \ x(\bar{m})(1) \\
f(x)(\bar{m},4) \ &= \ 2 \\
f(x)(\bar{m},5) \ &= \ x(\bar{m})(0) \\
f(x)(\bar{m},6) \ &= \ x(\bar{m})(1) \\
f(x)(\bar{m},7) \ &= \ x(\bar{m})(2) \\
f(x)(\bar{m},8) \ &= \ 2 \\
&\vdots
\end{align*}
Clearly $f$ is Borel and injective. Fixing $x,y\in (2^\omega)^{\omega^\ell}$ and letting $\bar{a},\bar{b}\in\N^{\ell}$ and $\bar{c},\bar{d}\in\N^{\ell+1}$, it is not too hard to check that if $\bar{a}\cdot x=\bar{b}\cdot y$ then $(\bar{a},0)\cdot f(x)=(\bar{b},0)\cdot f(y)$, and if $\bar{c}\cdot f(x)=\bar{d}\cdot f(y)$ then $c_{\ell+1}=d_{\ell+1}$ and $(c_1,\ldots,c_{\ell})\cdot x=(d_1,\ldots,d_{\ell})\cdot y$.
\end{proof}

\begin{lemma}\label{lem:coding2}
For all $\ell\in\Z_+$, we have $E^3_{\ell}\sqsubseteq_B\ETl$.
\end{lemma}

\begin{proof}
First define the Borel map $f_1:3^\omega\to 2^\omega$ as follows. Intuitively, $f_1$ will send $(x_0,x_1,\ldots)$ to the sequence $(x_0'^\frown x_1'^\frown\ldots)$ where
\[
0' := 1000, \qquad 1':=1100, \qquad \text{and} \qquad 2':=1110.
\]
Formally, this means
\[
f_1(x)(k) \ := \ 
\begin{cases}
0 & \text{if $k\equiv 1\mod 4$ and $x\big((k-1)/4\big)=0$;} \\
0 & \text{if $k\equiv 2\mod 4$ and $x\big((k-2)/4\big)\in\{0,1\}$;} \\
0 & \text{if $k\equiv 3\mod 4$;} \\
1 & \text{otherwise.}
\end{cases}
\] 
We make several observations about $f_1$:
\begin{enumerate}[(i)]
 \item For all $x\in 3^\omega$, $a\in\N$, and $k\in\omega$, we have $f_1(a\cdot x)(k) = f_1(x)(4a+k)$.
 \item For all $x\in 3^\omega$ and $k\in\omega$, we have $f_1(x)(k)=0$ and $f_1(x)(k+1)=1$ if and only if $k\equiv 3\mod 4$.
 \item For all $x,y\in 3^{\omega^{\ell}}$ and $\bar{k}\in\omega^{\ell}$, we have $x(\bar{k})=y(\bar{k})$ if and only if 
 \[
 f_1\Big(\big\langle x(i, k_2,\ldots,k_\ell)\st i\in\omega\big\rangle \Big)(4k_1+j) \ = \ 
f_1\Big(\big\langle y(i, k_2,\ldots,k_\ell)\st i\in\omega\big\rangle \Big)(4k_1+j)
 \]
 for all $j\in\{0,1,2,3\}$.
\end{enumerate}

Next, for each $\ell\in\Z_{\geq 2}$ define the Borel function $f_{\ell}:3^{\omega^\ell}\to 2^{\omega^\ell}$ by
\[
f_{\ell}(x)(k_1,\ldots,k_{\ell}) \ := \ 
\begin{cases}
f_1\Big( \big\langle x\big(i,\frac{k_2}{4},\ldots,\frac{k_{\ell}}{4}\big)\st i\in\omega\big\rangle \Big)(k_1) & \text{if $4\mid k_i$ for each $2\leq i\leq\ell$;} \\
0 & \text{otherwise}.
\end{cases}
\]
The map $f_2$ can be pictured as in Figure \ref{fig:MoreCoding}.

\begin{figure}[h]
\begin{center}
\begin{tikzpicture}
\fill [fill=gray!16] (3.3,1.7) rectangle (4.9,3.3);
\draw (.1,.1) -- (4.9,.1);
\draw [dotted, thick] (4.9,.1) -- (5.2,.1);
\draw (.1,.1) -- (.1,3.3);
\draw [dotted, thick] (.1,3.3) -- (.1,3.6);
\draw [dotted] (.1,1.7) -- (5.15,1.7);
\draw [dotted] (.1,3.3) -- (5.15,3.3);
\draw [dotted] (1.7,.1) -- (1.7,3.55);
\draw [dotted] (3.3,.1) -- (3.3,3.55);
\draw [dotted] (4.9,.1) -- (4.9,3.55);
\node at (-.38,-.1) {$f_2(x)$};
\node at (.1,.1) [circle,fill,inner sep=1.2pt]{};
\node at (.3,.3) {$1$};
\node at (.7,.3) {$0$};
\node at (1.1,.3) {$0$};
\node at (1.5,.3) {$0$};
\node at (1.9,.3) {$1$};
\node at (2.3,.3) {$1$};
\node at (2.7,.3) {$0$};
\node at (3.1,.3) {$0$};
\node at (3.5,.3) {$1$};
\node at (3.9,.3) {$1$};
\node at (4.3,.3) {$1$};
\node at (4.7,.3) {$0$};
\foreach \x in {0,...,11} \node at (.3+.4*\x,.7) {$0$};
\foreach \x in {0,...,11} \node at (.3+.4*\x,1.1) {$0$};
\foreach \x in {0,...,11} \node at (.3+.4*\x,1.5) {$0$};
\node at (.3,1.9) {$1$};
\node at (.7,1.9) {$1$};
\node at (1.1,1.9) {$0$};
\node at (1.5,1.9) {$0$};
\node at (1.9,1.9) {$1$};
\node at (2.3,1.9) {$1$};
\node at (2.7,1.9) {$1$};
\node at (3.1,1.9) {$0$};
\node at (3.5,1.9) {$1$};
\node at (3.9,1.9) {$0$};
\node at (4.3,1.9) {$0$};
\node at (4.7,1.9) {$0$};
\foreach \x in {0,...,11} \node at (.3+.4*\x,2.3) {$0$};
\foreach \x in {0,...,11} \node at (.3+.4*\x,2.7) {$0$};
\foreach \x in {0,...,11} \node at (.3+.4*\x,3.1) {$0$};
\node at (-1.2,1) {$\longmapsto$};
\fill [fill=gray!16] (-3,.58) rectangle (-2.5,1.04);
\draw (-4,.1) -- (-2.5,.1);
\draw [dotted, thick] (-2.5,.1) -- (-2.1,.1);
\draw (-4,.1) -- (-4,1.04);
\draw [dotted, thick] (-4,1.04) -- (-4,1.44);
\node at (-4.2,-.1) {$x$};
\node at (-4,.1) [circle,fill,inner sep=1.1pt]{};
\node at (-3.75,.35) {$0$};
\node at (-3.25,.35) {$1$};
\node at (-2.75,.35) {$2$};
\node at (-3.75,.81) {$1$};
\node at (-3.25,.81) {$2$};
\node at (-2.75,.81) {$0$};
\draw [dotted] (-4,.58) -- (-2.3,.58);
\draw [dotted] (-4,1.04) -- (-2.3,1.04);
\draw [dotted] (-3.5,.1) -- (-3.5,1.25);
\draw [dotted] (-3,.1) -- (-3,1.25);
\draw [dotted] (-2.5,.1) -- (-2.5,1.25);
\end{tikzpicture}
\end{center}
\caption{The coding function $f_2:3^{\omega^2}\to\TT$.}
\label{fig:MoreCoding}
\end{figure}
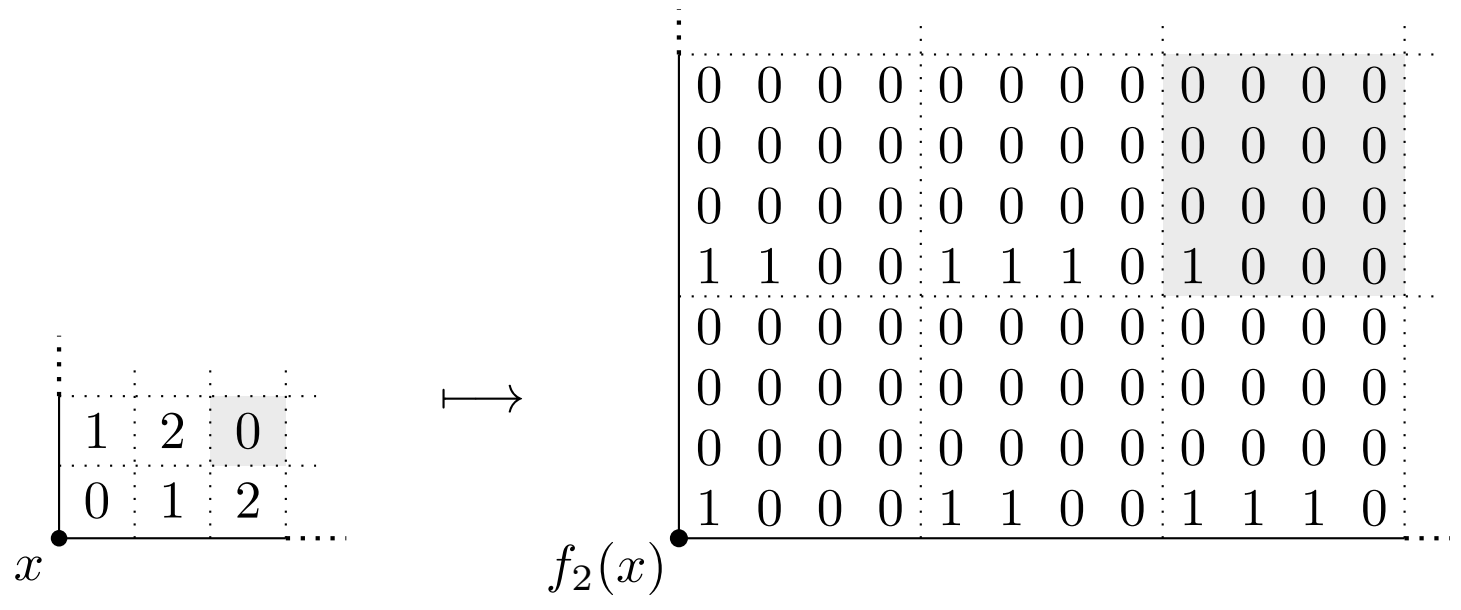

We show now that for each $\ell\in\Z_+$, $f_{\ell}$ is an embedding from $E^3_{\ell}$ to $E^2_{\ell}$. Let $\ell\in\Z_+$ be arbitrary, and for the rest of this proof write $f:=f_{\ell}$ to improve notation. For injectivity of $f$, note that if $x(\bar{k})\ne y(\bar{k})$ then either $f(x)(\bar{e}_1+4\bar{k})\ne f(y)(\bar{e}_1+4\bar{k})$ or $f(x)(2\bar{e}_1+4\bar{k})\ne f(y)(2\bar{e}_1+4\bar{k})$. 

To see that $f$ reduces $E^3_{\ell}$ to $E^2_{\ell}$, fix $x,y\in 3^{\omega^\ell}$ and assume first that $x\mathrel{E^3_{\ell}}y$, say $\bar{a}\cdot x=\bar{b}\cdot y$ where $\bar{a},\bar{b}\in\N^{\ell}$. It should be clear from the picture that $4\bar{a}\cdot f(x)=4\bar{b}\cdot f(y)$. To verify this, let $\bar{k}\in\omega^{\ell}$ be arbitrary. If there is some $i\geq 2$ for which $4\nmid k_i$, then $\big(4\bar{a}\cdot f(x)\big)(\bar{k})=0=\big(4\bar{b}\cdot f(y)\big)(\bar{k})$, so assume otherwise and write $k_i=4t_i$ for each $2\leq i\leq\ell$. Then
\begin{align*}
\big(4\bar{a}\cdot f(x)\big)(\bar{k}) \ &= f(x)\big(4a_1+k_1,4a_2+4t_2,\ldots,4a_{\ell}+4t_{\ell}\big) \\
&= \ f_1\Big(\big\langle x(i, a_2+t_2, \ldots, a_{\ell}+t_{\ell})\st i\in\omega\big\rangle\Big)(4a_1+k_1) \\
&= \ f_1\Big(a_1\cdot \big\langle x(i, a_2+t_2, \ldots, a_{\ell}+t_{\ell})\st i\in\omega\big\rangle\Big)(k_1) \qquad\qquad \text{\big(by Observation (i)\big)} \\
&= \ f_1\Big(\big\langle x(a_1+i, a_2+t_2, \ldots, a_{\ell}+t_{\ell})\st i\in\omega\big\rangle\Big)(k_1) \\
&= \ f_1\Big(\big\langle (\bar{a}\cdot x)(i, t_2, \ldots, t_{\ell})\st i\in\omega\big\rangle\Big)(k_1),
\end{align*}
and likewise
\[
\big(4\bar{b}\cdot f(y)\big)(\bar{k}) \ = \ f_1\Big(\big\langle (\bar{b}\cdot y)(i, t_2, \ldots, t_{\ell})\st i\in\omega\big\rangle\Big)(k_1).
\]
Since $\bar{a}\cdot x=\bar{b}\cdot y$, this shows $\big(4\bar{a}\cdot f(x)\big)(\bar{k})=\big(4\bar{b}\cdot f(y)\big)(\bar{k})$, and we conclude $4\bar{a}\cdot f(x)=4\bar{b}\cdot f(y)$ as desired. 

Conversely, and recycling notation, fix $\bar{a},\bar{b}\in\N^{\ell}$ and suppose $\bar{a}\cdot f(x)=\bar{b}\cdot f(y)$. Note that for each $2\leq i\leq\ell$ and fixed $k_i\in\omega$, the values
\[
f(x)(\bar{a}+\bar{k}) \quad \text{for} \quad k_1,\ldots,k_{i-1}, k_{i+1},\ldots,k_{\ell}\in\omega
\]
are all zero if and only if $a_i+k_i\not\equiv 0\mod 4$. Since the analogous statement holds for $f(y)(\bar{b}+\bar{k})$ and $f(x)(\bar{a}+\bar{k})=f(y)(\bar{b}+\bar{k})$ for all $\bar{k}\in\omega^\ell$, it follows that for each $2\leq i\leq\ell$, we have $a_i\equiv b_i\mod 4$. 

We claim now that $a_1\equiv b_1\mod 4$ as well. Indeed, for each $2\leq i\leq\ell$, let $a_i'\geq a_i$ be least such that $4\mid a_i'$, and likewise let $b_i'\geq b_i$ be least such that $4\mid b_i'$. Then $a_i'-a_i=b_i'-b_i$ for each $2\leq i\leq\ell$, so $(a_1,a_2',\ldots,a_\ell')\cdot f(x)=(b_1,b_2',\ldots,b_{\ell}')\cdot f(y)$. Let us write
\[
x' \ := \ \Big\langle x\big(i,{\textstyle{\frac{a_2'}{4},\ldots, \frac{a_{\ell}'}{4}}}\big)\st i\in\omega\Big\rangle \qquad 
\text{and} \qquad 
y' \ := \ \Big\langle y\big(i,{\textstyle{\frac{b_2'}{4},\ldots, \frac{b_{\ell}'}{4}}}\big)\st i\in\omega\Big\rangle,
\]
so $x',y'\in 3^\omega$. Then for all $k\in\omega$ we have
\[
f_1(x')(a_1+k) \ = \ f(x)(a_1+k, a_2',\ldots,a_{\ell}') \ = \ f(y)(b_1+k, b_2',\ldots,b_{\ell}') \ = \ f_1(y')(b_1+k).
\]
Using Observation (ii), let $k$ be least such that $f_1(x')(a_1+k)=0$ and $f_1(x')(a_1+k+1)=1$, so also $f_1(y')(b_1+k)=0$ and $f_1(y')(b_1+k+1)=1$. Then $a_1+k\equiv 3\mod 4$ and $b_1+k\equiv 3\mod 4$, so $a_1\equiv b_1\mod 4$ as claimed.

Now let $a_1'\geq a_1$ and $b_1'\geq b_1$ be least such that $4\mid a_1'$ and $4\mid b_1'$, and write $\bar{a}':=(a_1',\ldots,a_{\ell}')$ and $\bar{b}':=(b_1',\ldots,b_{\ell}')$, so $\bar{a}'\cdot x=\bar{b}'\cdot y$. For each $1\leq i\leq\ell$, let $c_i:=a_i'/4$ and $d_i:=b_i'/4$. We claim that $\bar{c}\cdot x=\bar{d}\cdot y$. To see this, let $\bar{k}\in\omega^{\ell}$ be arbitrary. Then for each $j\in\{0,1,2,3\}$ we have
\begin{align*}
f_1\Big(\big\langle x(i, c_2+k_2,\ldots,c_{\ell}+k_{\ell})\st i\in\omega\big\rangle\Big)(4c_1+4k_1+j) \ &= \ f(x)\big(4\c+4\bar{k}+j\bar{e}_1\big) \\
= \ \big(\bar{a}'\cdot f(x)\big)(4\bar{k}+j\bar{e}_1) \ &= \ 
\big(\bar{b}'\cdot f(y)\big)(4\bar{k}+j\bar{e}_1) \ = \ f(y)\big(4\bar{d}+4\bar{k}+j\bar{e}_1\big) \\ = \ f_1\Big(\big\langle y(i, d_2+k_2,\mkern1mu\ldots, \,&d_{\ell}+k_{\ell})\st i\in\omega\big\rangle\Big)(4d_1+4k_1+j).
\end{align*}
Therefore $x(\c+\bar{k})=y(\bar{d}+\bar{k})$ by Observation (iii). Since $\bar{k}\in\omega^{\ell}$ was arbitrary, this shows $\bar{c}\cdot x=\bar{d}\cdot y$, completing the proof.
\end{proof}

\begin{proof}[Proof of Conjecture \ref{con:final} from Conjecture \ref{con:HigherDimensions}]
Let $X$, $\ell$, and $U_1,\ldots,U_{\ell}$ be given. For notational convenience, write $\bar{a}\ast x$ for $\big(U_1^{a_1}\circ \cdots\circ U_\ld^{a_\ld}\big)(x)$. Let $\{V_i\st i\in\omega\}$ be a base for a compatible Polish topology on $X$, and define the injective Borel function $f:X\to (2^\omega)^{\omega^\ell}$ by
\[
f(x)(\bar{m})(i) \ = \ \begin{cases} 1 & \text{if $\m\ast x\in V_i$;} \\
0 & \text{otherwise.} \end{cases}
\]
We claim that $f$ is a Borel reduction from $E$ to $E^{2^\omega}_{\ell}$. Indeed, given $x\in X$ and $\bar{a},\bar{m}\in\omega^{\ld}$, observe that since $(\bar{m}+\bar{a})\ast x=\bar{m}\ast(\bar{a}\ast x)$, we have
\[
\big(\bar{a}\cdot f(x)\big)(\bar{m})(i) \ = \ f(x)(\bar{m}+\bar{a})(i) \ = \ f(\bar{a}\ast x)(\bar{m})(i)
\]
for all $i\in\omega$, and therefore $\bar{a}\cdot f(x)=f(\bar{a}\ast x)$. But then since $f$ is injective, for all $x,y\in X$ we have
\[
\begin{array}{rclll}
x\mathrel{E}y & \Longleftrightarrow & (\exists\bar{a},\bar{b}\in\N^\ld) \, \bar{a}\ast x=\bar{b}\ast y & & \\
& \Longleftrightarrow & (\exists\bar{a},\bar{b}\in\N^\ld) \, f(\bar{a}\ast x)=f(\bar{b}\ast y) & & \\
& \Longleftrightarrow & (\exists\bar{a},\bar{b}\in\N^\ld) \, \bar{a}\cdot f(x)=\bar{b}\cdot f(y) & \Longleftrightarrow & f(x)\mathrel{E^{2^\omega}_{\ell}}f(y).
\end{array}
\]
Thus $E\sqsubseteq_B E^{2^\omega}_{\ell}$. But $E^{2^\omega}_{\ell}\sqsubseteq_B E^3_{\ell+1} \sqsubseteq_B E^2_{\ell+1}$ by Lemmas \ref{lem:coding1} and \ref{lem:coding2}, completing the proof.
\end{proof}


\printbibliography

\end{document}